\documentclass[12pt,headings=optiontohead]{book}
\usepackage[english]{babel}
\usepackage[latin1]{inputenc}
\usepackage{amstext}
\usepackage{amsmath}
\usepackage{amsfonts}
\usepackage{amssymb}
\usepackage{amsthm}
\usepackage[all]{xy}
\xyoption{arc}
\CompileMatrices
\usepackage{graphicx}

\newtheorem{thm}{Theorem}[section]
\newtheorem{thmsp}{Teorema}[section]
\newtheorem{cor}[thm]{Corollary}
\newtheorem{lem}[thm]{Lemma}
\newtheorem{prop}[thm]{Proposition}
\newtheorem{dfn}[thm]{Definition}
\newtheorem{rem}[thm]{Remark}
\newtheorem{ex}[thm]{Example}

\newenvironment{pr}[1][Proof]{\noindent\textbf{#1.} }{\ \rule{0.5em}{0.5em}}

\DeclareMathOperator{\kernel}{Ker}
\DeclareMathOperator{\supp}{Supp}

\newcommand{\SO}{{\mathcal{O}}}
\newcommand{\PP}{\mathbb{P}}

\newcommand{\CC}{\mathbb{C}}

\newcommand{\Spec}{\operatorname{Spec}}

\newcommand{\Hilb}{\operatorname{Hilb}}
\newcommand{\Hom}{\operatorname{Hom}}
\newcommand{\Quot}{\operatorname{Quot}}
\newcommand{\Pic}{\operatorname{Pic}}
\newcommand{\Sch}{\operatorname{Sch}}
\newcommand{\Sets}{\operatorname{Sets}}
\newcommand{\Sym}{\operatorname{Sym}}
\newcommand{\id}{\operatorname{id}}
\newcommand{\im}{\operatorname{im}}
\newcommand{\surj}{\twoheadrightarrow}
\newcommand{\inj}{\hookrightarrow}
\newcommand{\too}{\longrightarrow}
\newcommand{\rk}{\operatorname{rk}}
\newcommand{\End}{\operatorname{End}}

\newcommand{\GL}{\operatorname{GL}}
\newcommand{\SL}{\operatorname{SL}}

\newcommand{\Proj}{\operatorname{Proj}}

\selectfont

\DeclareTextFontCommand{\emph}{\bfseries}

 \numberwithin{equation}{section}

\begin{document}

\titlepage

\setlength\topmargin{-1in}
\setlength\textwidth{7in}

\begin{figure}[h]
   \begin{center}
   \includegraphics[width=5cm]{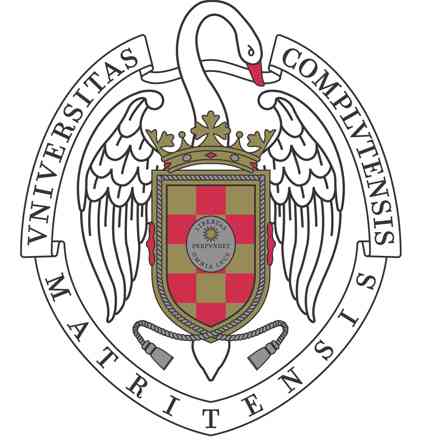}\;\;\;\;\;\;\;\;\;\;\;\;\;\;\;\;\;\;\;\;\;\;\;\;
   \includegraphics[width=5.5cm]{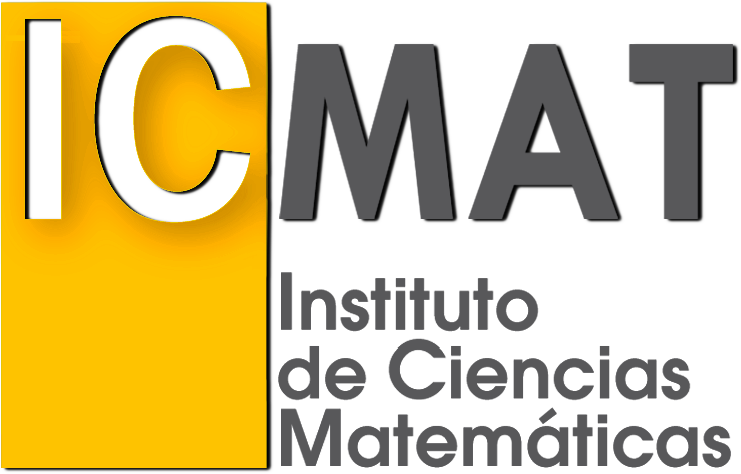}
   \end{center}
\textbf{Universidad Complutense de
Madrid}\;\;\;\;\;\;\textbf{ICMAT(CSIC-UAM-UC3M-UCM)}
\end{figure}

\vspace{1.5cm}

\LARGE{\centerline{\textbf{GIT CHARACTERIZATIONS OF}}} \LARGE{\centerline{\textbf{HARDER-NARASIMHAN
FILTRATIONS}}}
\vspace{0.5cm}
\hrule
\vspace{0.5cm}
\LARGE{\centerline{\textbf{CARACTERIZACIONES GIT DE}}}
\LARGE{\centerline{\textbf{FILTRACIONES DE HARDER-NARASIMHAN}}}

\vspace{2cm}

\LARGE{\centerline{\textbf{Alfonso Zamora Saiz}}}

\vspace{1.5cm}

\Large{\centerline{Dirigida por Tom\'{a}s L. G\'{o}mez}}

\vspace{2cm}

\large{\centerline{\textbf{Tesis presentada para la obtenci\'{o}n
del grado de Doctor}}}
\large{\centerline{\textbf{en
Matem\'{a}ticas por la Universidad Complutense de Madrid}}}


\newpage
\mbox{}
\thispagestyle{empty}

\setlength\topmargin{0in} \setlength\headheight{0in}
\setlength\headsep{1.5cm} \setlength\textheight{21cm}
\setlength\textwidth{6.25in}

\linespread{1.2}

\normalsize

\newpage

\pagenumbering{roman}

\newpage

\pagenumbering{roman}

\setlength\topmargin{0in} \setlength\headheight{0in}
\setlength\headsep{1.5cm} \setlength\textheight{21cm}
\setlength\textwidth{6.25in}

\linespread{1.2}

\normalsize

$$ $$

$$ $$

$$ $$

$$ $$

$$ $$


\begin{flushright}
\LARGE \textit{To my mother}
\end{flushright}

\newpage
\mbox{}

\newpage

\normalsize

$$ $$

$$ $$

\begin{flushright}
TOM\'{A}S.- \textit{Entonces... }\\
ASEL.- \textit{\textexclamdown Entonces hay que salir a la otra c\'{a}rcel!\\
\textexclamdown Y cuando est\'{e}s en ella, salir a otra, y de \'{e}sta, a otra!\\
La verdad te espera en todas, no en la inacci\'{o}n. Te esperaba\\
aqu\'{i}, pero s\'{o}lo si te esforzabas en ver la mentira de la\\
Fundaci\'{o}n que imaginaste. Y te espera en el esfuerzo de ese\\
oscuro t\'{u}nel del s\'{o}tano... En el holograma de esa evasi\'{o}n.}
\end{flushright}

\begin{flushright}
La Fundaci\'{o}n, Antonio Buero Vallejo
\end{flushright}

\vspace{2cm}

\begin{flushright}
THOMAS.- \textit{Then... }\\
ASEL.- \textit{Then you must exit into the other jail!\\
And when you're in it, go into another, and from it, to still another!\\
Truth awaits you in all of them, not in inaction. You found it\\
here, but only by making yourself see the lie of the\\
Foundation you imagined. And it waits for you in the effort of that\\
dark tunnel in the basement... in the hologram of that escape.}
\end{flushright}

\begin{flushright}
The Foundation, Antonio Buero Vallejo, Trans. Marion Peter Holt
\end{flushright}

\newpage
\mbox{}

\chapter*{Acknowledgements}
\addcontentsline{toc}{chapter}{Acknowledgements}

This research has been done thanks to a FPU grant of the Spanish
Ministry of Education (2009-2013) and the participation in the
research grants \textit{Hodge Theory and Moduli of Bundles} (Ref.
MTM2007-63582) and \textit{Moduli spaces, algebraic, arithmetic
and topological questions} (Ref. MTM2010-17389) funded the Spanish
Government.

$ $

First of all, I would like to thank to my advisor, Tom\'as
G\'omez, without whose help this thesis would have not been
possible. Thanks for teaching me and listening to me.

I would also like to thank to Ignacio Sols, for his enthusiastic
proposal, whose expectations I hope to have met.

To the Departmento de \'Algebra of the Facultad de Matem\'aticas,
Universidad Complutense de Madrid, and to the Instituto de
Ciencias Matem\'aticas (CSIC-UAM-UC3M-UCM), for letting me to be
part of them.

I want to thank to the places I have visited scientifically within
these years. Thanks to the Isaac Newton Institute for Mathematical
Sciences in Cambridge, UK, and the organizers of the Moduli Spaces
semester in January-June 2011, for giving me the opportunity to
learn from the best ones. Also, to the Department of Mathematics
of Columbia University in New York City, USA, for hospitality. And
thanks to Professors Michael Thaddeus and Robert Friedman for
listening to me there.

Here I would like to point out some of the Professors from whom I
have learnt or with whom I have discussed about my work, most
directly during these $4$ years. Thanks to Luis
\'Alvarez-C\'onsul, Enrique Arrondo, Marco Castrill\'on, \'Oscar
Garc\'ia-Prada, Jos\'e Fernando, Nitin Nitsure, Ignasi Mundet,
Vicente Mu\~{n}oz, Fran Presas and Markus Reineke. Many thanks to
all the people who gave their lives up during one minute to
listen, think, answer a mail or explaining on the board.

I need a special paragraph to thank to Peter Newstead, for helping
and supporting me so much, in every sense. Also thanks to Peter
and Ann Newstead, for treating me like a son.

I want to recognize special thanks to Jos\'e Manuel Gamboa, for trusting
in me. Many thanks for counseling and consoling me.

I could not forget to acknowledge here to some special people
thanks to whom this story started. To Jos\'e Col\'on, Serapio
Garc\'ia, Bernardino del Campo, Mercedes Fern\'andez, and all the
rest of the organizers of the Albacete, Castilla-La Mancha and Spain Mathematical Olympiads, for letting me think my first problems.

To my fellows from UCM and ICMAT, Alba, Alberto, Alicia, \'Alvaro
A., \'Alvaro S., Ana P., Andrea, Blanca,
 Carlos F., Carlos P., Carlos Q., David, Diego, Emilio, Espe, Giovanni, Javi, Luis, Marina, Mario, Nacho,
Rafa, Roberto, Roger, Silvia, Simone, and many others, for sharing
experiences, conversations, classes, seminars, trips, jokes,
sarcasm, meals, life and mathematics. Thanks for being my family
in Madrid. I have no words to thank, more that specially, to some
of them, for understanding me and lending a shoulder to cry on.
Thanks for putting color and smiles on grey and sadder days. This
thesis also belongs to you.

To my friends from my home village, Tarazona de la Mancha, for
spending holidays with music and good times.

To the Ministry of Education of Spain, for paying $51$ times for
me to think about the relation between the maximal unstability in
the sense of Geometric Invariant Theory and the Harder-Narasimhan
filtration.

To all the people who gave one without expecting two the following
day.

To my father and brother, for loving me and being there.

Finally, to the person this thesis is devoted to. This is for you
and because of you. Without you, this would not have been
possible. Thanks for stay there, day and night, supporting me,
encouraging me and loving me. Thank you.

\newpage

\addcontentsline{toc}{chapter}{Contents}

\tableofcontents

\newpage

\chapter*{Introducci\'{o}n (en espa\~{n}ol)}

\pagenumbering{arabic}

\markboth{\MakeUppercase{Introducci\'{o}n}}{}
\addcontentsline{toc}{chapter}{Introducci\'{o}n (en espa\~{n}ol)}

La presente tesis doctoral est\'{a} dedicada al estudio de la relaci\'{o}n
entre la inestabilidad maximal en el sentido de la Teor\'ia Geom\'etrica de Invariantes
 (que abreviaremos por sus siglas en ingl\'es, GIT) y la filtraci\'{o}n de Harder-Narasimhan
para diferentes problemas de espacios de m\'{o}duli. Muchos de los
problemas de m\'{o}duli en geometr\'{i}a hacen uso de la Teor\'ia Geom\'etrica de Invariantes para la construcci\'{o}n de espacios de m\'{o}duli.
Imponemos, de inicio, una noci\'{o}n de estabilidad en los objetos
para los cuales queremos construir un espacio de m\'{o}duli y,
mediante la Teor\'ia Geom\'etrica de Invariantes, un
objeto estable (resp. semiestable, inestable) se corresponde con un punto GIT
estable (resp. GIT semiestable, inestable) en cierto espacio,
estableciendo una correspondencia entre ambas nociones de
estabilidad. El concepto de m\'{a}xima inestabilidad en el sentido GIT
ha sido estudiado por diferentes autores, y para el presente
prop\'{o}sito consideraremos el tratamiento de Kempf, cuyo art\'{i}culo
\cite{Ke} lo explora. Por otra parte, la filtraci\'{o}n de
Harder-Narasimhan, ampliamente usada en muchos problemas en
geometr\'{i}a algebraica, es el objeto geom\'{e}trico que representa la
idea de inestabilidad maximal para la condici\'{o}n de estabilidad
impuesta de inicio sobre los objetos.

En esta tesis se demuestra que ambas nociones de inestabilidad
maximal coinciden, y se muestra una correspondencia entre ellas en
diferentes casos. El primer cap\'{i}tulo contiene nociones generales
sobre problemas de m\'{o}duli, la Teor\'{i}a Geom\'etrica de Invariantes y
la filtraci\'{o}n de Harder-Narasimhan que usaremos en los cap\'{i}tulos
$2$ y $3$, adem\'{a}s de un ejemplo de construcci\'{o}n de un espacio de
m\'{o}duli para tensores. En el segundo cap\'{i}tulo estudiamos diferentes
problemas de m\'{o}duli relacionados con haces, o con haces con
estructura adicional. Desarrollamos una t\'{e}cnica para probar la
mencionada corres\-pondencia para haces coherentes sin torsi\'{o}n sobre
variedades proyectivas de dimensi\'{o}n arbitraria, pares holomorfos,
haces de Higgs, tensores de rango $2$, y hacemos algunos
comentarios acerca de los tensores de rango $3$, siendo el primer
caso de tensores para el cual la t\'{e}cnica desarrollada no funciona.
En el tercer cap\'{i}tulo estudiamos representaciones de un carcaj, demostrando un resultado similar para
representaciones en la categor\'{i}a de espacios vectoriales y,
nuevamente, haces coherentes vistos como representaciones de un
carcaj de un solo v\'{e}rtice en la categor\'{i}a de haces coherentes.

\section*{Inestabilidad maximal}
Consid\'{e}rese un problema de m\'{o}duli en el que queremos clasificar
una clase de objetos algebro-geom\'{e}tricos m\'{o}dulo una relaci\'{o}n de
equivalencia. Usualmente, tenemos que imponer una condici\'{o}n de
estabilidad en los objetos que clasificamos para obtener un
espacio de m\'{o}duli con buenas propiedades donde cada punto
corresponda a una clase de equivalencia de objetos. Entonces,
a\~{n}adiendo un dato adicional a los objetos (lo que se suele llamar
en la literatura en ingl\'es \textit{to rigidify the data}) los inclu\'{i}mos en un
espacio de par\'{a}metros con el que nos es m\'{a}s sencillo trabajar (como
pueda ser un esquema af\'{i}n o proyectivo). Para deshacernos de este
nuevo dato a\~{n}adido, tenemos que tomar el cociente por la acci\'{o}n de un
grupo que est\'{a} precisamente codificando los cambios en este dato
adicional. Y para tomar este cociente usamos la Teor\'ia Geom\'etrica de Invariantes de Mumford, (v\'{e}ase \cite{Mu} para la primera
edici\'{o}n, \cite{MF} y \cite{MFK} para la segunda y tercera
ediciones) para obtener un espacio de m\'{o}duli proyectivo
clasificando los objetos en el problema de m\'{o}duli. El estudio de
las \'{o}rbitas de la acci\'{o}n de este grupo nos lleva a la noci\'{o}n de
estabilidad GIT definiendo puntos, en el espacio de par\'{a}metros,
que son GIT estables y otros que son GIT inestables.

En cada problema de m\'{o}duli en el que usamos GIT, llegado cierto
momento tenemos que demostrar que ambas nociones de estabilidad
coinciden, con lo que los objetos estables corresponden a los
puntos GIT estables y los objetos inestables corresponden a los
puntos GIT inestables. A tal efecto, Mumford enuncia un criterio
num\'{e}rico (c.f. \cite[Theorem 2.1]{Mu}) basado en ideas de Hilbert
en \cite{Hi}. El Teorema \ref{HMcrit}, conocido como el criterio
de Hilbert-Mumford, caracteriza la estabilidad GIT a trav\'{e}s de
subgrupos uniparam\'{e}tricos, donde una funci\'{o}n num\'{e}rica (que llamamos el \textit{m\'{i}nimo exponente relevante}) toma un valor
positivo o negativo dependiendo de que el subgrupo uniparam\'{e}trico
desestabilice un punto o no, en el sentido GIT. Adem\'{a}s, si un
punto es GIT inestable, podemos hablar de \textit{grados de
inestabilidad}, o de ciertos subgrupos uniparam\'{e}tricos que son m\'{a}s
desestabilizantes que otros.

Basado en el trabajo de Mumford, Tits y otros autores, podemos medir esto
mediante una funci\'{o}n racional en el espacio de subgrupos
uniparam\'{e}tricos, cuyo numerador es la funci\'{o}n num\'{e}rica del
criterio de Hilbert-Mumford y cuyo denominador es una \emph{norma} del
subgrupo uniparam\'{e}trico que escogemos para evitar reescalar la
funci\'{o}n num\'{e}rica (v\'{e}ase secci\'{o}n \ref{kempfsection}). La
\emph{conjetura del centro de Tits} (c.f. \cite[p. 64]{Mu})
establece que existe un \'{u}nico subgrupo uniparam\'{e}trico maximizando
esta funci\'{o}n, que representa la inestabilidad maximal en el
sentido GIT. Kempf explora estas ideas en un art\'{i}culo en $1978$
(c.f. \cite{Ke}), resolviendo lo que \'{e}l llama \emph{la conjetura
de Mumford-Tits} (refiri\'{e}ndose a la conjetura del centro de Tits,
tal como aparece en \cite{Mu}), demostrando que existe un \'{u}nico
subgrupo uniparam\'{e}trico con estas propiedades en \cite[Theorem
2.2, Theorem 3.4]{Ke} (para una correspondencia entre las
definiciones en \cite{Mu} y \cite{Ke} v\'{e}ase \cite[Appendix
2B]{MFK}).

Un objeto inestable proporciona un punto GIT inestable para el
cual existe un \'{u}nico subgrupo uniparam\'{e}trico m\'{a}ximamente
desestabilizante en el sentido GIT. Los di\-fe\-ren\-tes subgrupos
uniparam\'{e}tricos producen, de forma natural, banderas formando el \emph{complejo de banderas} de
un grupo $G$, estudiado por Tits y Mumford. Por tanto, nos gustar\'{i}a
considerar la bandera asociada a ese \'{u}nico subgrupo uniparam\'{e}trico
m\'{a}ximamente desestabilizante en el sentido GIT y construir, a
partir de \'{e}l, una filtraci\'{o}n de subobjetos del objeto inestable
original. Los principales resultados de esta tesis consisten en la
traducci\'{o}n de este subgrupo uniparam\'{e}trico a una
filtraci\'{o}n del objeto, demostrar que esta traducci\'{o}n est\'{a} bien y
un\'{i}vocamente definida (es decir, que no depende de diversas
elecciones hechas durante el proceso) y, finalmente, demostrar que
coincide con la filtraci\'{o}n de Harder-Narasimhan en los casos en
los que esta filtraci\'{o}n es ya conocida, o proporciona una nueva
definici\'{o}n de una tal filtraci\'{o}n en otro caso.

\section*{Teor\'{i}a Geom\'etrica de Invariantes y espacio de m\'{o}duli de tensores}

En el primer cap\'{i}tulo se recogen las nociones preliminares acerca
de espacios de m\'{o}duli y la Teor\'{i}a Geom\'etrica de Invariantes para
presentar el problema estudiado.

La secci\'{o}n \ref{modulis} contiene una descripci\'{o}n de lo que es un
problema de m\'{o}duli, con su formulaci\'{o}n rigurosa. Se proporcionan
ejemplos b\'{a}sicos como el espacio de m\'{o}duli de las c\'{u}bicas
complejas no singulares o la formulaci\'{o}n del problema del espacio
de m\'{o}duli de las curvas algebraicas de g\'{e}nero $g$. Entonces,
recuperamos las nociones de la Teor\'{i}a Geom\'etrica de Invariantes
que necesitaremos en lo sucesivo, los diferentes tipos de
cocientes, linearizaciones de acciones de grupos y el criterio de
Hilbert-Mumford, esencial para tratar con la estabilidad GIT.
Tambi\'{e}n damos un ejemplo, el Ejemplo \ref{hyperboles}, para
resaltar el concepto de $S$-equivalencia para las \'{o}rbitas de la
acci\'{o}n de un grupo, o el ejemplo de
combinaciones de puntos en $\mathbb{P}_{\mathbb{C}}^{1}$ (v\'{e}ase
Ejemplos \ref{exG} y \ref{exG2}) para aplicar el criterio de
Hilbert-Mumford.

La secci\'{o}n \ref{exampletensors} est\'{a} dedicada a presentar un
ejemplo completo de construcci\'{o}n de un espacio de m\'{o}duli usando la
Teor\'{i}a de Invariantes Geom\'{e}tricos: el espacio de m\'{o}duli de tensores
sobre una variedad proyectiva. Esta construcci\'{o}n fue primeramente
estudiada para tensores sobre curvas por Schmitt (c.f. \cite{Sch})
y despu\'{e}s por G\'{o}mez y Sols (c.f. \cite{GS1}). En la secci\'{o}n,
seguimos la referencia \cite{GS1} pero usando el embebimiento de Gieseker (c.f.
subsecci\'{o}n \ref{giesekerembedding}) en vez de el de Simpson, como
se hace en \cite{GS1}. Con este embebimiento, los tensores se
meten en un espacio de par\'{a}metros donde act\'{u}a un grupo. El Teorema
\ref{GIT-delta} establece que los tensores semiestables
corresponden a las \'{o}rbitas GIT semiestables bajo la acci\'{o}n del
grupo, por lo tanto el conciente GIT ser\'{a} el espacio de m\'{o}duli de
los tensores semiestables. En este ejemplo aparecen muchos de los
elementos que suelen encontrarse en las construcciones GIT de
espacios de m\'{o}duli, tales como la dependencia de la estabilidad
con un par\'{a}metro, la necesidad de probar la acotaci\'{o}n del conjunto de objetos semiestables (una prueba que
exige un gran esfuerzo, basada en resultados de \cite{Ma1,Ma2}), o
la identificaci\'{o}n de tensores $S$-equivalentes (v\'{e}ase la
Proposici\'{o}n \ref{Sequivalence}) como puntos semiestables.

Una vez que hemos discutido la correspondencia entre estabilidad y
estabilidad GIT con el ejemplo del m\'{o}duli de tensores, recordamos
la noci\'{o}n de la filtraci\'{o}n de Harder-Narasimhan en la secci\'{o}n
\ref{HNsection}. Explicamos por qu\'{e} esta filtraci\'{o}n captura la
idea de la m\'axima forma de desestabilizar un objeto a trav\'{e}s de casos sencillos (como
en el Ejemplo \ref{HNP1}) y probamos su existencia y unicidad para
el caso de haces coherentes sin torsi\'{o}n. Entonces explicamos la
noci\'{o}n de la filtraci\'{o}n de Harder-Narasimhan en el contexto
abstracto de una categor\'{i}a abeliana, como aparece en \cite{Ru}.

Finalmente, la secci\'{o}n \ref{kempfsection} finaliza el cap\'{i}tulo
explicando las ideas de \cite{Ke}. En ese art\'{i}culo, Kempf prueba
la conjetura de Mumford-Tits, enunciando que un punto GIT
inestable tiene un \'{u}nico subgrupo uniparam\'{e}trico que lo
desestabiliza m\'{a}ximamente, en el sentido de GIT, maximizando una
funci\'{o}n, la \emph{funci\'{o}n de Kempf}, en el Teorema
\ref{kempftheorem0}. Por tanto, teniendo este subgrupo
uniparam\'{e}trico dando la m\'axima forma de desestabilizaci\'on GIT, y la filtraci\'{o}n
de Harder-Narasimhan, podemos conjeturar que ambas corresponden a
la misma noci\'{o}n, y formular la pregunta

$ $

\textbf{\textquestiondown Existe una relaci\'on entre la filtraci\'{o}n de Harder-Narasimhan y el subgrupo uniparam\'{e}trico dado por Kempf?}

$ $

En los cap\'{i}tulos \ref{chaptersheaves} y \ref{chapterquivers}
respondemos positivamente la anterior pregunta en diferentes
casos.

\section*{Correspondencia entre las filtraciones de Kempf y Harder-Narasimhan}

Primero, resumimos c\'{o}mo demostrar la correspondencia para el caso principal,
haces coherentes sin torsi\'{o}n sobre variedades proyectivas. Los
casos restantes ser\'{a}n probados de forma an\'{a}loga a este caso
principal, basado en las mismas ideas y t\'{e}cnicas.

\subsection*{$\bullet$ Haces coherentes sin torsi\'{o}n sobre variedades proyectivas}

Sea $X$ una variedad proyectiva compleja no singular y sea
$\mathcal{O}_X(1)$ un fibrado de l\'{i}nea amplio en $X$. Si $E$ es
un haz coherente en $X$, sea $P_E$ su polinomio de Hilbert con
respecto a $\mathcal{O}_X(1)$, es decir, $P_E(m)=\chi(E\otimes
\mathcal{O}_X(m))$. Si $P$ y $Q$ son polinomios, escribimos $P\leq
Q$ si $P(m)\leq Q(m)$ para $m\gg 0$.

Un haz sin torsi\'{o}n $E$ sobre $X$ se llama \emph{semiestable} si para todo
subhaz propio $F\subset E$, se verifica
$$
\frac{P_F}{\rk F} \leq \frac{P_E}{\rk E} \; .
$$
Si no es semiestable, se llama \emph{inestable}, y posee una filtraci\'{o}n
can\'{o}nica
$$0\subset E_{1} \subset E_{2} \subset \cdots \subset E_{t} \subset
E_{t+1}=E\; ,$$ que satisface las siguientes propiedades, donde
$E^{i}:=E_{i}/E_{i-1}$:
 \begin{enumerate}
   \item Los polinomios de Hilbert verifican
   $$\frac{P_{E^{1}}}{\rk E^{1}}>\frac{P_{E^{2}}}{\rk E^{2}}>\ldots>\frac{P_{E^{t+1}}}{\rk E^{t+1}}$$
   \item Cada $E^{i}$ es semiestable
 \end{enumerate}
que es llamada la \emph{filtraci\'{o}n de Harder-Narasimhan} de $E$
(v\'{e}ase Teorema \ref{HNunique}).

La construcci\'{o}n del espacio de m\'{o}duli para estos objetos es
originalmente debida a Gieseker para superficies, y generalizada a
dimensi\'{o}n superior por Maruyama (c.f. \cite{Gi1,Ma1,Ma2}). Para
construir el espacio de m\'{o}duli de haces sin torsi\'{o}n con polinomio
de Hilbert fijo $P$, elegimos un cierto entero grande $m$ y
consideramos el esquema Quot (siguiendo la nomenclatura de
Grothendieck) que parametriza cocientes
\begin{equation}
\label{quotintsp} V\otimes \SO_X(-m) \too E\; ,
\end{equation}
donde $V$ es un espacio vectorial fijo de dimensi\'{o}n $P(m)$ y $E$
es un haz con $P_{E}=P$. El esquema de cocientes tiene una acci\'{o}n
can\'{o}nica de $\SL(V)$. Gieseker (c.f. \cite{Gi1}) da una
linearizaci\'{o}n de esta acci\'{o}n en un cierto fibrado de l\'{i}nea amplio,
para usar la Teor\'{i}a Geom\'{e}trica de Invariantes para cocientar por
la acci\'{o}n. El espacio de m\'{o}duli de haces semiestables se obtiene
como el cociente GIT.

Sea $E$ un haz sin torsi\'{o}n que es inestable. Eligiendo $m$
suficientemente grande (dependiendo de $E$), y eligiendo un
isomorfismo $V\cong H^0(E(m))$, obtenemos un cociente como en
\eqref{quotintsp}. El correspondiente punto en el esquema Quot ser\'{a}
GIT inestable y, por el criterio de Hilbert-Mumford, habr\'{a} al
menos un subgrupo uniparam\'{e}trico de $\SL(V)$ que lo
\emph{desestabilizar\'{a}} en el sentido de GIT.

Un \emph{subgrupo uniparam\'{e}trico} de $\SL(V)$ es un homomorfismo
no trivial $\CC^*\to \SL(V)$. A un subgrupo uniparam\'{e}trico le
asociamos una filtraci\'{o}n con pesos como sigue. Existe una base de
$V$, $\{e_1,\ldots,e_p\}$, para la cual el subgrupo uniparam\'{e}trico
toma la forma diagonal
$$
t\mapsto \operatorname{diag} \big(
t^{\Gamma_1},\ldots,t^{\Gamma_1},
t^{\Gamma_2},\ldots,t^{\Gamma_2}, \ldots ,
t^{\Gamma_{t+1}},\ldots,t^{\Gamma_{t+1}} \big)
$$
donde $\Gamma_1<\cdots<\Gamma_{t+1}$. Sobre todos los subgrupos
uniparam\'{e}tricos, Kempf muestra que existe una clase de conjugaci\'{o}n
de m\'{a}ximamente desestabilizantes (es decir, que maximicen la
funci\'{o}n de Kempf en la Definici\'{o}n \ref{kempffunction}), todos
ellos dando una \'{u}nica filtraci\'{o}n de $V$ con pesos,
$(V_\bullet,n_\bullet)$,
\begin{equation}
\label{filtVintrosp} 0 \subset V_1 \subset V_2 \subset
\;\cdots\; \subset V_t \subset V_{t+1}=V,
\end{equation}
y n\'{u}meros positivos $n_1,\, n_2,\ldots , \,n_t
> 0$ (c.f. Teorema \ref{kempftheoremsheaves}).

Esta filtraci\'{o}n induce una filtraci\'{o}n de haces de $E$, evaluando
los espacios de secciones globales,
$$0\subseteq E^{m}_{1} \subseteq E^{m}_{2} \subseteq \cdots \subseteq E^{m}_{t} \subseteq
E_{t+1}=E\; ,$$ que llamaremos la \emph{$m$-filtraci\'{o}n de Kempf de
$E$}. Esta filtraci\'{o}n depende del entero $m$ que usamos en la
construcci\'{o}n del espacio de m\'{o}duli y que asegura que los haces
semiestables (resp. inestables) se corresponden con las \'{o}rbitas
GIT semiestables (resp. GIT inestables). Entonces, el punto
principal es probar el siguiente
\begin{thmsp}[c.f. Teorema \ref{kempfstabilizes}]
Existe un entero $m'\gg 0$ tal que la $m$-filtraci\'{o}n de Kempf es
independiente de $m$, para $m'\geq m$.
\end{thmsp}

Dado un entero $m$, la $m$-filtraci\'{o}n de Kempf maximiza una
funci\'{o}n, llamada la \emph{funci\'{o}n de Kempf},
$$
\mu(V_{\bullet},n_{\bullet})=\frac{\sum_{i=1}^{t}  n_{i} (r\dim
V_{i}-r_{i}\dim V)} {\sqrt{\sum_{i=1}^{t+1} {\dim V^{i}_{}}
\Gamma_{i}^{2}}}\; ,$$ la cual identificamos con una funci\'{o}n
geom\'{e}trica (v\'{e}ase Proposici\'{o}n \ref{identification})
$$\mu_{v}(\Gamma)=\frac{(\Gamma,v)}{\|\Gamma\|}\; ,$$
donde $(\; , \;)$ es un producto escalar en $\mathbb{R}^{t+1}$
dado por una matriz diagonal con elementos $\dim V^{i}$ en la
diagonal, y el vector $v$ tiene coordenadas
$$v_{i}=\frac{1}{\dim V^{i}\dim V}\big[r^{i}\dim V-r\dim V^{i}\big]\; .$$

El vector $v$ est\'{a} relacionado con la filtraci\'{o}n
$V_{\bullet}\subset V$ y el vector $\Gamma$ est\'{a} relacionado con
los n\'{u}meros $n_{\bullet}$ (poniendo
$n_{i}=\frac{\Gamma_{i+1}-\Gamma_{i}}{\dim V}$)  en la filtraci\'{o}n
de Kempf (\ref{filtVintrosp}). Entonces, fijando un vector $v$ en un
espacio Eucl\'{i}deo, consideramos la funci\'{o}n $\mu_{v}(\Gamma)$ y nos
preguntamos por el vector $\Gamma$ que proporciona el m\'{a}ximo para
$\mu_{v}$. Ocurre que la respuesta es dada por la envolvente
convexa del grafo determinado por $v$ (v\'{e}ase Teorema
\ref{maxconvexenvelope}).

Las $m$-filtraciones de Kempf, (es decir, las filtraciones de $E$
que obtenemos eva\-luando $(V_{\bullet},n_{\bullet})$ para
diferentes enteros $m$, donde $V\simeq H^{0}(E(m))\;$ ) pueden
diferir para diferentes valores de $m$. Sin embargo, por una parte
son maximales con respecto al valor que la funci\'{o}n de Kempf
alcanza en ellas, y por otra parte verifican propiedades de
convexidad con respecto a la funci\'{o}n $\mu_{v}(\Gamma)$. A partir
de esto, podemos probar diferentes propiedades satisfechas por los
filtros que aparecen en las filtraciones, las cuales caracterizan
la \emph{filtraci\'{o}n de Kempf} y muestran que es independiente de
$m$ (c.f. Teorema \ref{kempfstabilizes}).

La filtraci\'{o}n que obtenemos, la cual de hecho no depende del
entero $m$, se llama la \emph{filtraci\'{o}n de Kempf} de $E$.
Entonces, observamos que las dos propiedades de convexidad que
estaban impl\'{i}citas en los argumentos que condujeron a probar el
Teorema \ref{kempfstabilizes}, propiedades que quedan descritas en
los Lemas \ref{lemmaA} y \ref{lemmaB}, son igualmente sa\-tis\-fechas
por la filtraci\'{o}n de Kempf (c.f Proposiciones
\ref{descendentslopes} y \ref{blocksemistability}). Observamos que
estas propiedades de convexidad, las pendientes descendentes y la
semiestabilidad de los cocientes, son precisamente las propiedades
de la filtraci\'{o}n de Harder-Narasimhan para haces (c.f. Teorema
\ref{HNunique}). Por consiguiente, por unicidad de la filtraci\'{o}n
de Harder-Narasimhan, probamos el siguiente
\begin{thmsp}[c.f. Teorema \ref{kempfisHN}]
La filtraci\'{o}n de Kempf de un haz coherente sin torsi\'{o}n $E$
coincide con la filtraci\'{o}n de Harder-Narasimhan de $E$.
\end{thmsp}

Si reemplazamos los polinomios de Hilbert por los grados de los haces, la noci\'{o}n de
estabilidad se transforma en $\mu$-estabilidad (tambi\'{e}n conocida
como estabilidad de las pendientes) y obtenemos la
$\mu$-filtraci\'{o}n de Harder-Narasimhan. En \cite{Br,BT}, Bruasse y
Teleman dan una interpretaci\'{o}n en t\'{e}rminos de teor\'{i}a gauge de la
$\mu$-filtraci\'{o}n de Harder-Narasimhan para haces sin torsi\'{o}n y
pares holomorfos. Ellos tambi\'{e}n usan las ideas de Kempf, pero en
el marco del grupo gauge, por lo que tienen que generalizar los resultados de
Kempf a grupos infinito dimensionales.

Una correspondencia similar ha sido recientemente probada por Hoskins y Kirwan (c.f. \cite{HK}) usando un m\'etodo diferente. En la referencia se comienza con una filtraci\'on
que se encuentra en un estrato de \text{tipo de Harder-Narasimhan} fijado (lo que llamamos \textit{$m$-type}, c.f. Definici\'on \ref{mtype}). Una diferencia
con nuestro tratamiento es que ellas usan la existencia previa de la filtraci\'on de Harder-Narasimhan, mientras que nosotros no lo hacemos.

A continuaci\'{o}n, brevemente resumimos otros problemas de m\'{o}duli
para los cuales mostramos la correspondencia entre la filtraci\'{o}n
de Harder-Narasimhan y la inestabilidad maximal GIT, usando un
m\'{e}todo similar al de los haces sin torsi\'{o}n.

\subsection*{$\bullet$ Pares holomorfos}

Sea $X$ una variedad proyectiva compleja no singular. Consideramos \emph{pares holomorfos}
$$(E,\varphi:E\to \mathcal{O}_{X})$$
dados por un haz coherente sin torsi\'on de rango $r$ con determinante fijo
$\det(E)\cong \Delta$ y un morfismo al haz de estructura $\mathcal{O}_{X}$.
Obs\'ervese que se trata de un caso perticular de los tensores estudiados en la secci\'on \ref{exampletensors},
particularizando para $c=1$, $b=0$ y $s=1$.

Sea $\delta$ un polinomio de grado a lo sumo $\dim X-1$ y coeficiente director positivo. Dado un subhaz $E'\subset E$,
definimos $\epsilon(E')=1$ si $\varphi|_{E'}\neq 0$ y $\epsilon(E')=0$
en caso contrario. Un par holomorfo $(E,\varphi)$ es \emph{$\delta$-semiestable} si para cada $E'\subset E$, se tiene
$$\frac{P_{E'}-\delta\epsilon(E')}{\rk E'}\leq \frac{P_{E}-\delta\epsilon(E)}{\rk E}\; .$$

Existe una construcci\'on del espacio de m\'oduli de pares holomorfos $\delta$-semiestables fijando el polinomio
de Hilbert $P$ y el determinante $\det(E)\simeq \Delta$ en \cite{HL1} siguiendo las ideas de Gieseker,
y en \cite{HL2} (donde dichos pares son llamados \textit{framed modules}) siguiendo las ideas de Simpson.

Un par holomorfo $\delta$-inestable da un punto GIT inestable para el cual obtenemos un subgrupo uniparam\'etrico m\'aximamente
desestabilizante y una filtraci\'on de subpares. Mostramos que esta filtraci\'on no depende del entero $m$ usado en la construcci\'on del
espacio de m\'oduli en el Teorema \ref{kempfstabilizespairs}, y que coincide con la filtraci\'on de Harder-Narasimhan para pares holomorfos en el Teorema
\ref{kempfisHNpairs}.

La noci\'on de par holomorfo es dual a la de un par consistente en un haz coherente junto con una secci\'on. En ambos casos,
la condici\'on de estabilidad depende de un par\'ametro que es un polinomio. \'Esta fue la primera construcci\'on de un espacio de m\'oduli con una
condici\'on de estabilidad dependiente de par\'ametros.

\subsection*{$\bullet$ Haces de Higgs}
Sea $X$ una variedad proyectiva compleja no singular. Un \emph{haz de Higgs}
es un par $(E,\varphi)$ donde $E$ es un haz coherente sobre $X$ y
$$\varphi:E\rightarrow E\otimes\Omega^{1}_{X}\; ,$$
verificando $\varphi\wedge\varphi=0$, es un morfismo llamado el
\emph{campo de Higgs}. Si el haz $E$ es localmente libre hablaremos de \emph{fibrados de Higgs}. Decimos que un haz de Higgs es \emph{semiestable} (en el
sentido de Gieseker) si para todo subhaz propio $F\subset E$
preservado por $\varphi$ (i.e. $\varphi\big|_{F}:F\rightarrow
F\otimes \Omega^{1}_{X}$) tenemos
$$\frac{P_{F}}{\rk F}\leq \frac{P_{E}}{\rk E}\; ,$$ donde $P_{E}$ y $P_{F}$ son los respectivos polinomios de Hilbert.

Podemos pensar un haz de Higgs $(E,\varphi)$ como un haz coherente
$\mathcal{E}$ en el fibrado cotangente $T^{\ast}X$ (c.f. Lema
\ref{higgsiscoherent}) de tal forma que $\pi_{\ast}\mathcal{E}=E$, donde
$\pi:T^{\ast}X\rightarrow X$. Usamos la construcci\'on de Simpson
(c.f. \cite{Si1,Si2}) del espacio de m\'oduli de haces de Higgs con este punto de vista.
Entonces, probamos que el subgrupo uniparam\'etrico que da la m\'axima inestabilidad GIT (en el sentido de \cite{Ke})
produce una filtraci\'on de subhaces
$$0\subset \pi_{\ast}\mathcal{E}_{1} \subset \pi_{\ast}\mathcal{E}_{2} \subset \cdots \subset \pi_{\ast}\mathcal{E}_{t} \subset
\pi_{\ast}\mathcal{E}_{t+1}=E$$ que no depende de los enteros $m$ y $l$ usados en la construcci\'on de Simpson (c.f. Teorema
\ref{kempfstabilizeshiggs}). Aplicando $\pi_{\ast}$, obtenemos una filtraci\'on de subhaces de Higgs
$$0 \subset (E_{1}, \varphi|_{E_{1}}) \subset
(E_{2}, \varphi|_{E_{2}}) \subset \;\cdots\; \subset (E_{t},
\varphi|_{E_{t}}) \subset (E_{t+1},
\varphi|_{E_{t+1}})=(E,\varphi)$$ y probamos que coincide con la filtraci\'on de Harder-Narasimhan para haces de Higgs (c.f.
Corlario \ref{kempfisHNhiggs}).

Este caso tiene la particularidad de usar el embebimiento de Simpson (que depende de dos enteros $m$ y $l$) en vez de el de Gieseker, lo que hace que el m\'etodo funcione,
indistintamente, en ambos casos.

\subsection*{$\bullet$ Tensores de rango $2$}

Sea $X$ una variedad proyectiva compleja no singular de dimensi\'on $n$. Sea $E$ un haz coherente sin torsi\'on
de rango $2$ sobre $X$. Llamamos \emph{tensor de rango $2$} al par
$$
(E,\varphi:\overbrace{E\otimes\cdots \otimes E}^{\text{s
veces}}\too \mathcal{O}_{X})\; .
$$
Este es otro caso particular de tensores (c.f. secci\'on
\ref{exampletensors}), haciendo $c=1$, $b=0$, $r=2$ y $s$ arbitrario. Consideramos la construcci\'on del espacio de m\'oduli para tales
tensores de rango $2$ con determinante fijo $\det(E)\cong \Delta$.
De modo similar, probamos que la filtraci\'on de Kempf no depende del entero $m$, para $m\gg 0$ (c.f. Teorema \ref{kempfstabilizesrk2})
para construir una filtraci\'on de Harder-Narasimhan para el tensor $(E,\varphi)$, que en este caso es un subhaz de rango $1$ $L$,
$$0\subset (L,\varphi|_{L})\subset (E,\varphi)\; .$$

Cuando la variedad $X$ es una curva y el morfismo $\varphi$ es sim\'etrico, podemos interpretar esta noci\'on en t\'erminos de recubrimientos.
Mirando el lugar de anulaci\'on de $\varphi$ podemos ver el tensor $(E,\varphi)$ como un recubrimiento de grado $s$, $X'\rightarrow X$, dentro
de la superficie reglada $\mathbb{P}(E)$, para definir una noci\'on de recubrimiento estable y caracterizar geom\'etricamente el subhaz m\'aximamente desestabilizante
$L\subset E$ en t\'erminos de teor\'ia de intersecci\'on para superficies regladas. Ocurre que, la expresi\'on de la condici\'on de estabilidad para tales tensores (c.f.
(\ref{rk2stabexpression})), puede interpretarse como la estabilidad de Gieseker del haz $E$ (c.f. Definici\'on \ref{giesekerstab})
y la estabilidad de una configuraci\'on de puntos en $\mathbb{P}_{\mathbb{C}}^{1}$ (c.f. Ejemplos \ref{exG} y
\ref{exG2}), ponderadas por el par\'ametro de estabilidad.

\subsection*{$\bullet$ Tensores de rango $3$ y m\'as all\'a}

El cap\'itulo \ref{chaptersheaves} finaliza con algunas observaciones acerca del caso de los tensores de rango $3$. Poniendo $s=2$, $r=3$ en la Definici\'on
\ref{stabilityfortensors}, obtenemos el caso m\'as sencillo para el cual no podemos usar las ideas previas para demostrar que el subgrupo uniparam\'etrico m\'aximamente desestabilizante
produce una filtraci\'on de subtensores que no dependa del entero usado en la construcci\'on del espacio de m\'oduli, para valores grandes del entero.

La raz\'on de esto es que los resultados en la subsecci\'on \ref{convexcones} no se pueden aplicar. No podemos ver la funci\'on de Kempf
(c.f. Definici\'on \ref{kempffunction}) como una funci\'on en el espacio Eucl\'ideo tomando valores en los pesos de la filtraci\'on
(c.f. Proposici\'on \ref{identification}) porque los pesos depender\'an de la filtraci\'on. En este caso no podemos probar los an\'alogos a
los Lemas \ref{independenceofweightspairs} y \ref{independenceofweightsrk2}, por lo que el m\'etodo que usamos no sirve en general.

La alternativa a esto es comparar filtraciones candidatas a ser la filtraci\'on de Harder-Narasimhan mirando los valores que toma la funci\'on de Kempf en ellas, por
m\'etodos elementales. La secci\'on finaliza con la observaci\'on de que existe una clase restringida de tensores (aqu\'ellos para los cuales no se producen estos hechos y podemos
probar la independencia entre los pesos y las filtraciones), tal que los pasos de la prueba de la correspondencia pueden ser llevados a cabo
(v\'ease la Definici\'on \ref{determinedmultiindex}).

\section*{Representaciones de un carcaj}


En el cap\'itulo \ref{chapterquivers} exploramos las ideas desarrolladas en el cap\'itulo \ref{chaptersheaves} para representaciones de un carcaj. Probamos
la correspondencia an\'aloga para representaciones de un carcaj en espacios vectoriales finito dimensionales y usamos la construcci\'on functorial de un espacio de m\'oduli
para haces coherentes en \cite{ACK} para dar otra demostraci\'on del Teorema \ref{kempfisHN}.

Sea $Q$ un carcaj finito, dado por un conjunto finito de v\'ertices y
flechas entre ellos, y una representaci\'on de $Q$ en $k$-espacios vectoriales finito dimensionales, donde $k$ es un cuerpo algebraicamente cerrado de
caracter\'istica arbitraria. Existe una noci\'on de estabilidad para tales representaciones (c.f. Definici\'on \ref{Qstability}) dada por King
en \cite{Ki} y, m\'as en general, por Reineke en \cite{Re} (ambas casos particulares de la noci\'on abstracta de estabilidad para una categor\'ia abeliana que podemos encontrar en
\cite{Ru}), y una noci\'on de existencia de una \'unica filtraci\'on de Harder-Narasimhan con respecto a la condici\'on de estabilidad (c.f. Teorema \ref{HNquivers}).

En la secci\'on \ref{qvectorspaces} consideramos la construcci\'on de un espacio de m\'oduli para estos objetos dada por King (c.f. \cite{Ki}), y asociamos a una
representaci\'on inestable un punto inestable, en el sentido de la Teor\'ia Geom\'etrica de Invariantes, en un espacio de par\'ametros donde act\'ua un grupo. Entonces,
el subgrupo uniparam\'etrico dado por Kempf (c.f. Teorema \ref{kempftheoremquivers}), que es m\'aximamente desestabilizante en el sentido GIT, otorga una filtraci\'on
de subrepresentaciones. Probamos que esta filtraci\'on coincide con la filtraci\'on de Harder-Narasimhan para la representaci\'on incial, en el Teorema \ref{kempfHNquivers}.

Este caso es ligeramente diferente porque, en los anteriores, el grupo por el cual est\'abamos tomando el cociente GIT en la construcci\'on del espacio de m\'oduli,
era $\SL(V)$, pero en este caso se trata de un producto de grupos generales lineales, uno para cada v\'ertice del carcaj. Entonces, la \emph{longitud} que elegimos en el espacio de
subgrupos
uniparam\'etricos al definir la funci\'on de Kempf (v\'ease Definici\'on \ref{length}) depende de ciertos par\'ametros
(uno para cada factor simple en el grupo), y mostramos
c\'omo tenemos que colocar los par\'ametros convenientemente eligiendo una longitud particular, para ser capaces de relacionar el subgrupo uniparam\'etrico GIT m\'aximamente
desestabilizante con la filtraci\'on de Harder-Narasimhan.

Finalmente, en la secci\'on \ref{sectionqsheaves} definimos los $Q$-haces, que son representaciones de un carcaj en la categor\'ia de haces coherentes, y damos una noci\'on
de estabilidad para ellos, siguiendo \cite{AC,ACGP}. Para un carcaj de un s\'olo v\'ertice, un $Q$-haz es lo mismo que un haz coherente y usamos la construcci\'on functorial
del espacio de m\'oduli para haces dada en \cite{ACK}. En esta construcci\'on, se relacionan los haces con los m\'odulos de Kronecker para reescribir la condici\'on de estabilidad
en t\'erminos de representaciones de un carcaj en espacios vectoriales. El Teorema \ref{equivstability} relaciona todas las diferentes nociones de estabilidad que aparecen envueltas en
esta tesis, a saber, la estabilidad de un haz como $Q$-haz (que equivale a la estabilidad de Gieseker para haces), la estabilidad del m\'odulo de Kronecker asociado, la estabilidad
de la representaci\'on de otro carcaj asociado $\tilde{Q}$, y la estabilidad GIT del punto correspondiente en el espacio de par\'ametros. Usando la equivalencia de las diferentes
estabilidades, podemos aplicar el teorema de Kempf (c.f. Teorema \ref{kempftheoremquivers}) para encontrar una filtraci\'on de Harder-Narasimhan para la representaci\'on asociada de
un carcaj en espacios vectoriales y, desde aqu\'i, obtener la filtraci\'on de Harder-Narasimhan para el $Q$-haz en el Teorema \ref{HNqsheaves}.

\section*{Conclusiones}

Esta tesis contiene ideas en dos direcciones. Por una parte, exploramos la relaci\'on entre las condiciones de estabilidad y las nociones de estabilidad GIT en las construcciones
de espacios de m\'oduli. Por otra parte, relacionamos el concepto natural de m\'axima inestabilidad GIT (en el sentido de Kempf) y la filtraci\'on de Harder-Narasimhan, lo cual
produce una correspondencia en varios casos donde la filtraci\'on de Harder-Narasimhan es previamente conocida, o da una nueva noci\'on de una tal filtraci\'on en otros casos.

El esquema de la prueba es similar en todos los casos. Primero, obtener una filtraci\'on de subespacios vectoriales de secciones globales que maximice la funci\'on
de Kempf (c.f. Definici\'on \ref{kempffunction}) para el problema GIT considerado, entonces evaluar las secciones para conseguir una filtraci\'on de subobjetos,
que llamamos la \emph{$m$-filtraci\'on de Kempf}. Relacionamos la funci\'on de Kempf con una funci\'on en el espacio Eucl\'ideo (c.f. Proposici\'on \ref{identification}) para aplicar
los resultados de convexidad (v\'ease subsecci\'on \ref{convexcones}). Seguidamente, demostramos propiedades de la $m$-filtraci\'on de Kempf que la caracterizan y la har\'an
independiente del entero $m$, por lo que obtendremos una filtraci\'on que llamamos la \emph{filtraci\'on de Kempf}. Finalmente, en los casos donde la filtraci\'on de
Harder-Narasimhan es previamente conocida, probamos que las propiedades de convexidad de la filtraci\'on de Kempf son, precisamente, las de la filtraci\'on de Harder-Narasimhan,
luego por unicidad ambas filtraciones coinciden (c.f. Teorema \ref{kempfisHN}). En otros casos, como los tensores de rango $2$ en la secci\'on \ref{kempfrk2}, la filtraci\'on de
Kempf (que es \'unica) define una filtraci\'on de Harder-Narasimhan. N\'otese que, en la secci\'on \ref{qvectorspaces}, la construcci\'on del m\'oduli no depende de ning\'un entero, por
tanto no tenemos que probar un an\'alogo al Teorema \ref{kempfstabilizes}, y la correspondencia es mucho m\'as sencilla y r\'apida.

La noci\'on de longitud (c.f. Definici\'on \ref{length}) que necesitamos para definir la velocidad del subgrupo uniparam\'etrico (v\'ease secci\'on \ref{kempfsection}), juega un
papel importante. En principio, diferentes longitudes dar\'ian diferentes subgrupos uniparam\'etricos GIT m\'aximamente desestabilizantes, por tanto diferentes filtraciones de
Kempf candidatas a ser la filtraci\'on de Harder-Narasimhan. En la secci\'on \ref{qvectorspaces} observamos c\'omo, diferentes elecciones de longitud corresponden a diferentes
definiciones de estabilidad en la Definici\'on \ref{Qstability}. Por tanto, dada una noci\'on de filtraci\'on de Harder-Narasimhan (dependiente de la noci\'on de estabilidad), tenemos
que colocar los par\'ametros en la linearizaci\'on en la construcci\'on del espacio de m\'oduli, y en la definici\'on de longitud en el conjunto de subgrupos uniparam\'etricos
(c.f. Proposici\'on \ref{GITstab-stabquivers}), de forma conveniente, para lograr una correspondencia entre las filtraciones de Kempf y Harder-Narasimhan.

Potencialmente, podr\'iamos esperar usar estas ideas para definir filtraciones de Harder-Narasimhan en casos donde no han sido estudiadas, y usarlas como poderosa
herramienta para importantes aplicaciones donde ha sido usada en el pasado, tales como los teoremas de restricci\'on o el c\'alculo de los n\'umeros de Betti y los puntos racionales
en espacios de m\'oduli (c.f. \cite{HN}).

Muchos de los casos donde hemos aplicado el m\'etodo caen dentro de categor\'ias abelianas, donde la filtraci\'on de Harder-Narasimhan verifica las propiedades
de convexidad que se desprenden de su definici\'on. Ser\'ia interesante entender si esto impone una condici\'on para esperar una noci\'on de filtraci\'on de Harder-Narasimhan,
tal como la conocemos.

\vspace{2cm}

Los resultados principales de esta tesis est\'an recogidos en los
preprints \cite{GSZ,Za1,Za2}. El primero contiene la
correspondencia entre la m\'axima forma de desestabilizar en el
sentido GIT y la filtraci\'on de Harder-Narasimhan para haces
coherentes sin torsi\'on sobre variedades proyectivas y para pares
holomorfos (secciones \ref{kempfsheaves} y \ref{kempfpairs}). El
segundo contiene el resultado para representaciones de un carcaj
en espacios vectoriales de dimensi\'on finita (secci\'on
\ref{qvectorspaces}). El tercero se ocupa del caso de tensores de
rango $2$ y cubiertas estables.

\chapter*{Introduction}

\markboth{\MakeUppercase{Introduction}}{}
\addcontentsline{toc}{chapter}{Introduction}

This Ph.D. thesis is devoted to the study of the relation between
the maximal unstability in the sense of Geometric Invariant Theory
and the Harder-Narasimhan filtration in different moduli problems.
Many of the moduli problems in geometry use Geometric Invariant
Theory (abbreviated GIT) in the construction of a moduli space.
We impose, from the beginning, a notion of stability on the
objects for which we want to construct a moduli space and, by the
Geometric Invariant Theory, we associate to a stable (resp.
semistable, unstable) object, a GIT stable (resp. GIT semistable,
GIT unstable) point in certain space, establishing a
correspondence between both concepts of stability. The GIT concept
of maximal unstability has been studied by several authors, and
for our purposes we consider the work of Kempf, whose paper
\cite{Ke} explores it. On the other hand, the Harder-Narasimhan
filtration, widely used in many problems in algebraic geometry, is
the geometrical object which represents the idea of maximal
unstability for the previous notion of stability imposed on the
objects.

In this thesis we prove that both notions of maximal unstability
do coincide, and show a correspondence between them in different cases. The first chapter contains general notions about moduli problems, Geometric
Invariant Theory and the Harder-Narasimhan filtration we will use in chapters $2$ and $3$,
apart from an example of the construction of a moduli space for tensors.
In the second chapter we study different moduli
problems in relation with sheaves, or sheaves with additional
structure. We develop a technique to prove the mentioned
correspondence for torsion free coherent sheaves over arbitrary
dimensional projective varieties, holomorphic pairs, Higgs
sheaves, rank $2$ tensors, and we discuss rank $3$ tensors
as the first case of tensors for which the technique we use breaks
down. In the third chapter we study representations of quivers,
proving a similar result for representations on the category of
vector spaces and, again, coherent sheaves seen as representations
of a one vertex quiver on the category of coherent sheaves.

\section*{Maximal unstability}
Consider a moduli problem where we try to classify a class of
algebro-geometric objects modulo an equivalence relation. Usually,
we have to impose a stability condition on the objects we
classify, in order to obtain a moduli space with good properties where
each point corresponds to an equivalence class of objects. Then, by
adding additional data to the objects (what is usually called in
the literature \textit{to rigidify the data}) we include them in a
parameter space easier to work with (an affine or projective scheme). To get rid
of these new data we have to quotient by the action of a group
which is precisely encoding the changes in the additional data. And to take
this quotient we use Mumford's Geometric Invariant Theory (see
\cite{Mu} for the first edition, \cite{MF} and \cite{MFK} for
second and third editions) to obtain a projective moduli space
classifying the objects in the moduli problem. The study of the orbits of the action of this group
leads to the notion of GIT stability defining points, in the parameter space, which are GIT stable and points which are GIT unstable.

In every moduli problem using GIT, at some point one has to prove
that both notions of stability do coincide, then the stable
objects correspond to the GIT stable points, and the unstable ones
are related to the GIT unstable ones. To that purpose, Mumford
states a numerical criterion (c.f. \cite[Theorem 2.1]{Mu}) based
on ideas of Hilbert in \cite{Hi}. Theorem \ref{HMcrit}, known as
the \emph{Hilbert-Mumford criterion}, characterizes the GIT stability
through $1$-parameter subgroups, where a numerical function
(which we call the \textit{minimal relevant weight}) turns out to
be positive or negative whether the $1$-parameter subgroup
destabilizes a point or not, in the sense of GIT. Besides, when a
point is GIT unstable, we are able to talk about \textit{degrees of
unstability}, or some $1$-parameter subgroups which are more
destabilizing than others.

Based on the work of Mumford, Tits and other authors, we can measure this notion by means of a rational function on the
space of $1$-parameter subgroups, whose numerator is the numerical function of the Hilbert-Mumford criterion and
whose denominator is a \emph{length} of the $1$-parameter subgroup that we choose to avoid rescaling of the numerical
function (c.f. section \ref{kempfsection}). The \emph{center's conjecture of Tits} (c.f. \cite[p. 64]{Mu}) establishes that there exists a unique
$1$-parameter subgroup giving a maximum for this function, representing the GIT maximal unstability. Kempf
explores these ideas in a paper in $1978$ (c.f. \cite{Ke}), solving what he calls \emph{the Mumford-Tits
conjecture} (referring to Tits center's conjecture as it appears on \cite{Mu}) by proving that there exists a
unique $1$-parameter subgroup with these properties in \cite[Theorem 2.2, Theorem 3.4]{Ke} (for a correspondence between
definitions in \cite{Mu} and \cite{Ke} see \cite[Appendix 2B]{MFK}).

An unstable object gives a GIT unstable point for which there
exists a unique $1$-parameter subgroup GIT maximally
destabilizing. The different $1$-parameter subgroups produce, in a
natural way, flags (giving the \emph{flag complex} of a group $G$
studied by Tits and Mumford), hence we would like to consider the
flag associated to that unique $1$-parameter subgroup GIT
maximally destabilizing, and construct a filtration by subobjects
of the original unstable object, out of this $1$-parameter
subgroup. The main results of this Ph.D. thesis consist
on translating this $1$-parameter subgroup to a
filtration of the object, proving that this translation is well
and uniquely defined (i.e. it does not depend on several choices made during the process) and, finally, proving that it coincides with
the Harder-Narasimhan filtration in cases where it is already
known, or gives a new notion of such filtration in other cases.

\section*{Geometric Invariant Theory and moduli space of tensors}
In the first chapter I collect the necessary background about moduli spaces and Geometric Invariant Theory to present
the problem studied.

Section \ref{modulis} contains a description of what a moduli problem is, with its rigorous formulation. We
provide basic examples as the moduli space of non singular complex cubics or the formulation of the problem of
the moduli space of algebraic curves of genus $g$. Then, we recall the notions of Geometric Invariant Theory we
will need in the following, the different types of quotients, linearizations of actions of groups and the
Hilbert-Mumford criterion, essential to deal with GIT stability. We also give an example (c.f. Example \ref{hyperboles}) to
realize the concept of $S$-equivalence for the orbits of the action of a group, or the classical example
of combinations of points in $\mathbb{P}_{\mathbb{C}}^{1}$ (c.f. Examples \ref{exG} and \ref{exG2})
to apply the Hilbert-Mumford criterion.

Section \ref{exampletensors} is devoted to present a complete example of a construction of a moduli space using
Geometric Invariant Theory: the moduli space of tensors over a projective variety. This construction was first
studied for tensors over curves by Schmitt (c.f. \cite{Sch}) and then by G\'{o}mez and Sols (c.f. \cite{GS1}).
In the section, we follow \cite{GS1} but using the embedding of Gieseker (c.f. subsection
\ref{giesekerembedding}) instead of Simpson's, as it is done in \cite{GS1}. With this embedding, tensors are
collected in a parameter space on which a group is acting. Theorem \ref{GIT-delta} establishes that the
semistable tensors correspond to the GIT semistable orbits under the action of the group, hence the GIT quotient
will be a moduli space for semistable tensors. In this example, many of the general features which appear in GIT
constructions of moduli spaces take place, such as the dependence of the stability notion with a parameter, the
necessity of proving the boundedness of the set of semistable objects (a proof which usually takes a big effort,
based on results in \cite{Ma1,Ma2}), or the identification of $S$-equivalent tensors (c.f. Proposition
\ref{Sequivalence}) as semistable points.

Once we have discussed the correspondence between stability and GIT stability with the example of the moduli of
tensors, we recall the notion of the Harder-Narasimhan filtration in section \ref{HNsection}. We explain why
this filtration captures the idea of the maximal way of destabilizing an object through easy cases (as in Example \ref{HNP1}) and prove
its existence and uniqueness for the case of torsion free coherent sheaves. Then we explain the notion of the Harder-Narasimhan
filtration in the abstract context of an abelian category, which appears in \cite{Ru}.

Finally, section \ref{kempfsection} closes this chapter explaining
the ideas of \cite{Ke}. There, Kempf proves the Mumford-Tits
conjecture, asserting that a GIT unstable point has a unique
$1$-parameter subgroup which maximally destabilizes it, in the
sense of GIT, by maximizing a function, the \emph{Kempf function}
(c.f. Definition \ref{kempffunction}), in Theorem
\ref{kempftheorem0}. Hence, having this $1$-parameter subgroup
giving GIT maximal way of destabilizing and the Harder-Narasimhan
filtration, we can conjecture that they do correspond to the same
notion, and formulate the question

$ $

\textbf{Is the Harder-Narasimhan related to the $1$-parameter
subgroup given by Kempf?}

$ $

In chapters \ref{chaptersheaves} and \ref{chapterquivers} we
answer positively the previous question for different cases.

\section*{Correspondence between Kempf and
Harder-Narasimhan filtrations}

First, we summarize how to prove the correspondence for the main case, torsion free
sheaves over projective varieties. The rest of cases will be
proven in an analogous way to this main case, based on the same
ideas and techniques.

\subsection*{$\bullet$ Torsion free coherent sheaves over projective
varieties}

Let $X$ be a smooth complex projective variety and let $\mathcal{O}_X(1)$ be an ample line bundle on $X$. If $E$
is a coherent sheaf on $X$, let $P_E$ be its Hilbert polynomial with respect to $\mathcal{O}_X(1)$, i.e.,
$P_E(m)=\chi(E\otimes \mathcal{O}_X(m))$. If $P$ and $Q$ are polynomials, we write $P\leq Q$ if $P(m)\leq Q(m)$
for $m\gg 0$.

A torsion free sheaf $E$ on $X$ is called \emph{semistable} if for all proper subsheaves $F\subset E$, it is
$$
\frac{P_F}{\rk F} \leq \frac{P_E}{\rk E} \; .
$$
If it is not semistable, it is called  \emph{unstable}, and it has a canonical filtration
$$0\subset E_{1} \subset E_{2} \subset \cdots \subset E_{t} \subset
E_{t+1}=E\; ,$$ satisfying the following properties, where $E^{i}:=E_{i}/E_{i-1}$:
 \begin{enumerate}
   \item The Hilbert polynomials verify
   $$\frac{P_{E^{1}}}{\rk E^{1}}>\frac{P_{E^{2}}}{\rk E^{2}}>\ldots>\frac{P_{E^{t+1}}}{\rk E^{t+1}}$$
   \item Every $E^{i}$ is semistable
 \end{enumerate}
which is called the \emph{Harder-Narasimhan filtration} of $E$ (c.f. Theorem \ref{HNunique}).

The construction of the moduli space for these objects is originally due to Gieseker for surfaces, and
generalized to higher dimension by Maruyama (c.f. \cite{Gi1,Ma1,Ma2}). To construct the moduli space of torsion
free sheaves with fixed Hilbert polynomial $P$, we choose a suitably large integer $m$ and consider the Quot
scheme parametrizing quotients
\begin{equation}
\label{quotint} V\otimes \SO_X(-m) \too E\; ,
\end{equation}
where $V$ is a fixed vector space of dimension $P(m)$ and $E$ is a sheaf with $P_E=P$. The Quot scheme has a
canonical action by $\SL(V)$. Gieseker (c.f. \cite{Gi1}) gives a linearization of this action on a certain ample
line bundle, in order to use Geometric Invariant Theory to take the quotient by the action. The moduli space of
semistable sheaves is obtained as the GIT quotient.

Let $E$ be a torsion free sheaf which is unstable. Choosing $m$ large enough (depending on $E$), and choosing an
isomorphism $V\cong H^0(E(m))$, we obtain a quotient as in \eqref{quotint}. The corresponding point in the Quot
scheme will be GIT unstable and, by the Hilbert-Mumford criterion, there will be at least one $1$-parameter subgroup of
$\SL(V)$ which \emph{destabilizes} the point in the sense of GIT.

A \emph{$1$-parameter subgroup} of $\SL(V)$ is a non trivial homomorphism $\CC^*\to \SL(V)$. To a $1$-parameter
subgroup we associate a weighted filtration as follows. There is a basis $\{e_1,\ldots,e_p\}$ of $V$ where it
has a diagonal form
$$
t\mapsto \operatorname{diag} \big( t^{\Gamma_1},\ldots,t^{\Gamma_1}, t^{\Gamma_2},\ldots,t^{\Gamma_2}, \ldots ,
t^{\Gamma_{t+1}},\ldots,t^{\Gamma_{t+1}} \big)
$$
with $\Gamma_1<\cdots<\Gamma_{t+1}$. Among all these 1-parameter subgroups, Kempf shows that there is a
conjugacy class of maximally destabilizing $1$-parameter subgroups (i.e. maximizing the Kempf function in Definition \ref{kempffunction})
 all of them giving a unique weighted
filtration $(V_\bullet,n_\bullet)$ of $V$,
\begin{equation}
\label{filtVintro}
0\subset V_1 \subset V_2 \subset \;\cdots\; \subset V_t \subset V_{t+1}=V,
\end{equation}
and positive numbers $n_1,\, n_2,\ldots , \,n_t > 0$ (c.f. Theorem \ref{kempftheoremsheaves}).

This filtration induces a sheaf filtration of $E$ by evaluating the spaces of global sections,
$$0\subseteq E^{m}_{1} \subseteq E^{m}_{2} \subseteq \cdots \subseteq E^{m}_{t} \subseteq
E_{t+1}=E\; ,$$ which we call the \emph{$m$-Kempf filtration of
$E$}. This filtration depends on the integer $m$ which we use to
construct the moduli space and assure that the semistable sheaves
(resp. unstable) correspond to GIT semistable (resp. GIT unstable)
orbits. Then, the main point is to prove the following
\begin{thm}[c.f. Theorem \ref{kempfstabilizes}]
There exists an integer $m'\gg 0$ such that the $m$-Kempf
filtration is independent of $m$, for $m'\geq m$.
\end{thm}

Given an integer $m$, the $m$-Kempf filtration maximizes a function, called
\emph{the Kempf function},
$$
\mu(V_{\bullet},n_{\bullet})=\frac{\sum_{i=1}^{t}  n_{i} (r\dim V_{i}-r_{i}\dim V)} {\sqrt{\sum_{i=1}^{t+1}
{\dim V^{i}_{}} \Gamma_{i}^{2}}}\; ,
$$
which we identify with a geometrical function (c.f. Proposition \ref{identification})
$$\mu_{v}(\Gamma)=\frac{(\Gamma,v)}{\|\Gamma\|}\; ,$$
where $(\; , \;)$ is an inner product in $\mathbb{R}^{t+1}$ given by a diagonal matrix with elements $\dim V^{i}$
in the diagonal, and the vector $v$ having coordinates
$$v_{i}=\frac{1}{\dim V^{i}\dim V}\big[r^{i}\dim V-r\dim V^{i}\big]\; .$$

The vector $v$ is related with the flag $V_{\bullet}\subset V$ and the vector $\Gamma$ is related with the numbers $n_{\bullet}$ (by setting
$n_{i}=\frac{\Gamma_{i+1}-\Gamma_{i}}{\dim V}$) in the Kempf filtration (\ref{filtVintro}). Then, fixing a vector $v$ in an Euclidean space, consider the
function $\mu_{v}(\Gamma)$ and ask for the vector $\Gamma$
which gives maximum for $\mu_{v}$. It turns out that the answer is given by the convex envelope of the graph produced by $v$ (c.f. Theorem \ref{maxconvexenvelope}).

The $m$-Kempf filtrations, (i.e. the filtrations of $E$ we obtain by
evaluating $(V_{\bullet},n_{\bullet})$ for different integers $m$, where $V\simeq H^{0}(E(m))$) can differ for different values of $m$. However, they are, on the one hand, maximal with respect to
the value the Kempf function achieves on them and, on the other hand, verify convexity properties with
respect to $\mu_{v}(\Gamma)$. From this, we can prove different properties satisfied by the filters appearing in the filtrations,
which characterize the \emph{Kempf filtration} and show that it is independent of $m$ (c.f. Theorem \ref{kempfstabilizes}).

The filtration we obtain, which does actually not depend on the integer $m$, is
called the \emph{Kempf filtration} of $E$. Then, we observe that the two convexity
properties which were implicit in the arguments leading to prove
Theorem \ref{kempfstabilizes}, properties which are described by
Lemmas \ref{lemmaA} and \ref{lemmaB}, are also satisfied by the Kempf filtration (c.f Propositions
\ref{descendentslopes} and \ref{blocksemistability}). And we realize that these convexity properties are precisely the properties of the Harder-Narasimhan filtration
for sheaves (c.f. Theorem \ref{HNunique}), the descending slopes
and the semistability of the quotients. Therefore, by uniqueness of the Harder-Narasimhan
filtration, we prove the following
\begin{thm}[c.f. Theorem \ref{kempfisHN}]
The Kempf filtration of an unstable torsion free coherent sheaf $E$ coincides with
the Harder-Narasimhan filtration of $E$.
\end{thm}

If we replace Hilbert polynomials with degrees, the notion of stability becomes $\mu$-stability (also
known as slope stability) and we obtain the $\mu$-Harder-Narasimhan filtration. In \cite{Br,BT}, Bruasse and
Teleman give a gauge-theoretic interpretation of the $\mu$-Harder-Narasimhan filtration for torsion free sheaves
and for holomorphic pairs. They also use Kempf's ideas, but in the setting of the gauge group, so they have to
generalize Kempf's results to infinite dimensional groups.

A similar correspondence has been proved recently by Hoskins and Kirwan (c.f. \cite{HK}) by using a different method. They start
with a filtration which lays on a stratum with fixed \textit{Harder-Narasimhan type} (what we call \emph{$m$-type}, c.f. Definition \ref{mtype}).
One difference with our approach is that they use the previous existence of the Harder-Narasimhan filtration, whereas we do not use it.

Now, we briefly summarize other moduli problems for which we show
the correspond\-ence between the Harder-Narasimhan filtration and
the GIT maximal unstability, using a similar method as in the case
of torsion free sheaves.

\subsection*{$\bullet$ Holomorphic pairs}

Let $X$ be a smooth complex projective variety. Let us consider
\emph{holomorphic pairs}
$$(E,\varphi:E\to \mathcal{O}_{X})$$
given by a coherent torsion free sheaf of rank $r$ with fixed
determinant $\det(E)\cong \Delta$ and a morphism to the
structure sheaf $\mathcal{O}_{X}$. Observe that this is a particular case of the tensors studied in section
\ref{exampletensors}, by setting $c=1$, $b=0$ and $s=1$.

Let $\delta$ be a polynomial of degree at most $\dim X-1$ and positive leading coefficient. Given a subsheaf
$E'\subset E$, let $\epsilon(E')=1$ if $\varphi|_{E'}\neq 0$ and $\epsilon(E')=0$ otherwise. A holomorphic pair
$(E,\varphi)$ is \emph{$\delta$-semistable} if for every $E'\subset E$
$$\frac{P_{E'}-\delta\epsilon(E')}{\rk E'}\leq \frac{P_{E}-\delta\epsilon(E)}{\rk E}\; .$$

There is a construction of the moduli space of $\delta$-semistable holomorphic
pairs with fixed Hilbert polynomial $P$ and fixed determinant
$\det(E)\simeq \Delta$ in \cite{HL1} following Gieseker's ideas, and in \cite{HL2} (where
these pairs are called framed modules) following Simpson's ideas.

A $\delta$-unstable
holomorphic pair give a GIT unstable point for which we obtain a
$1$-parameter subgroup GIT maximally destabilizing and a
filtration of subpairs. We show that this filtration does not
depend on the integer $m$ used in the construction of the moduli
space in Theorem \ref{kempfstabilizespairs}, and that it coincides
with the Harder-Narasimhan filtration for holomorphic pairs in
Theorem \ref{kempfisHNpairs}.

This notion of holomorphic pair is dual to the pair consisting on a coherent sheaf together
with a section. In both cases, the stability condition depends on a parameter
which is a polynomial. This was the first construction of a moduli space with a stability condition
depending on parameters.

\subsection*{$\bullet$ Higgs sheaves}
Let $X$ be a smooth complex projective variety. A
\emph{Higgs sheaf} is a pair $(E,\varphi)$ where $E$ is a coherent
sheaf over $X$ and
$$\varphi:E\rightarrow E\otimes\Omega^{1}_{X}\; ,$$
verifying $\varphi\wedge\varphi=0$, a morphism called the
\emph{Higgs field}. If the sheaf $E$ is locally free we talk about
\emph{Higgs bundles}. We say that a Higgs sheaf is \emph{semistable} (in
the sense of Gieseker) if for all proper subsheaves $F\subset E$
preserved by $\varphi$ (i.e. $\varphi\big|_{F}:F\rightarrow
F\otimes \Omega^{1}_{X}$) we have
$$\frac{P_{F}}{\rk F}\leq \frac{P_{E}}{\rk E}\; ,$$ where $P_{E}$ and $P_{F}$ are the respective Hilbert polynomials.

A Higgs sheaf $(E,\varphi)$ can be thought as a coherent sheaf
$\mathcal{E}$ on the cotangent bundle $T^{\ast}X$ (c.f. Lemma
\ref{higgsiscoherent}) such that $\pi_{\ast}\mathcal{E}=E$, where
$\pi:T^{\ast}X\rightarrow X$. We use the construction of Simpson
(c.f. \cite{Si1,Si2}) of a moduli space for Higgs sheaves with
this point of view. Then, we prove that the $1$-parameter subgroup
giving the GIT maximal unstability (in the sense  of \cite{Ke})
provides a filtration of subsheaves
$$0\subset \pi_{\ast}\mathcal{E}_{1} \subset \pi_{\ast}\mathcal{E}_{2} \subset \cdots \subset \pi_{\ast}\mathcal{E}_{t} \subset
\pi_{\ast}\mathcal{E}_{t+1}=E$$ which does not depend on the
integers $m$, $l$ used in Simpson's construction (c.f. Theorem
\ref{kempfstabilizeshiggs}). By applying $\pi_{\ast}$, we get a
filtration of Higgs subsheaves
$$0 \subset (E_{1}, \varphi|_{E_{1}}) \subset
(E_{2}, \varphi|_{E_{2}}) \subset \;\cdots\; \subset (E_{t},
\varphi|_{E_{t}}) \subset (E_{t+1},
\varphi|_{E_{t+1}})=(E,\varphi)$$ and we prove that it coincides
with the Harder-Narasimhan filtration for Higgs sheaves (c.f.
Corollary \ref{kempfisHNhiggs}).

This case has the particularity of using the embedding of Simpson (which depends on two integers $m$ and $l$) instead of Gieseker's, what makes the method work, indistinctly,
both cases.

\subsection*{$\bullet$ Rank $2$ tensors}

Let $X$ be a smooth complex projective variety of dimension $n$. Let $E$ be a rank $2$ coherent torsion free
sheaf over $X$. We call a \emph{rank $2$ tensor} the pair
$$
(E,\varphi:\overbrace{E\otimes\cdots \otimes E}^{\text{s times}}\too \mathcal{O}_{X})\; .
$$
This is another particular case of tensors (c.f. section \ref{exampletensors}), by setting $c=1$, $b=0$, $r=2$ and arbitrary
$s$. Consider the given construction of the moduli space for such
rank $2$ tensors with fixed determinant $\det(E)\cong \Delta$. Similarly we
prove that the Kempf filtration does not depend on the integer
$m$, for $m\gg 0$ (c.f. Theorem \ref{kempfstabilizesrk2}) to
construct a Harder-Narasimhan filtration for the tensor
$(E,\varphi)$, which in this case is a rank $1$ subsheaf $L$,
$$0\subset (L,\varphi|_{L})\subset (E,\varphi)\; .$$

When the variety $X$ is a curve and the morphism $\varphi$ is symmetric, we can interpret this notion
in terms of coverings. Looking at the vanishing locus of $\varphi$ we can see the
tensor $(E,\varphi)$ as a degree $s$ covering $X'\rightarrow X$ lying on the
ruled surface $\mathbb{P}(E)$, to define a notion of stable
covering and characterize geometrically the maximally
destabilizing subsheaf $L\subset E$ in terms of intersection
theory for ruled surfaces. It turns out that, the expression of
the stability condition for such tensors (c.f.
(\ref{rk2stabexpression})), can be seen as the Gieseker's stability
of the sheaf $E$ (c.f. Definition \ref{giesekerstab}) and the
stability of a configuration of points in
$\mathbb{P}_{\mathbb{C}}^{1}$ (c.f. Examples \ref{exG} and
\ref{exG2}), pondered by the stability parameter.

\subsection*{$\bullet$ Rank $3$ tensors and beyond}

Chapter \ref{chaptersheaves} finishes with some observations about
the rank $3$ tensors case. Setting $s=2$, $r=3$ on Definition
\ref{stabilityfortensors} we obtain the easiest case for which we
cannot use the previous ideas to prove that the $1$-parameter
subgroup GIT maximally destabilizing produces a filtration of
subtensors which does not depend on some integer used in the
construction of the moduli space, for large values of the integer.

The reason is that results on subsection \ref{convexcones} do not
apply. We cannot see the Kempf function (c.f. Definition
\ref{kempffunction}) as a function on the Euclidean space taking
values on the weights of the filtration (c.f. Proposition
\ref{identification}) because the weights will depend on the
filtration. In this case we cannot prove analogous to Lemmas
\ref{independenceofweightspairs} and
\ref{independenceofweightsrk2}, hence the method we use does not
apply in general.

The alternative is to compare candidates to be the
Harder-Narasimhan filtration by looking at the values the Kempf function takes at them, by elementary methods. The section finishes
with the observation that there exists a restricted class of
tensors (those for which the features discussed before do not
apply and we can prove the independence between weights and filtrations), such that the steps of the proof of the correspondence do hold (c.f.
Definition \ref{determinedmultiindex}).

\section*{Representations of quivers}


In chapter \ref{chapterquivers} we explore the ideas developed in chapter \ref{chaptersheaves} for
representations of quivers. We prove the analogous correspondence for representations of quivers on finite
dimensional vector spaces and use the functorial construction of a moduli space for coherent sheaves in
\cite{ACK} to give another proof of Theorem \ref{kempfisHN}.

Let $Q$ be a finite quiver, given by a finite set of vertices and arrows between them, and a representation of
$Q$ on finite dimensional $k$-vector spaces, where $k$ is an algebraically closed field of arbitrary
characteristic. There exists a notion of stability for such representations (c.f. Definition \ref{Qstability})
given by King in \cite{Ki} and, more generally, by Reineke in \cite{Re} (both particular cases of the abstract
notion of stability for an abelian category that we can find in \cite{Ru}), and a notion of the existence of a
unique Harder-Narasimhan filtration with respect to that stability condition (c.f. Theorem \ref{HNquivers}).

In section \ref{qvectorspaces} we consider the construction of a
moduli space for these objects by King (c.f. \cite{Ki}), and
associate to an unstable representation an unstable point, in the
sense of Geometric Invariant Theory, in a parameter space where a
group acts. Then, the $1$-parameter subgroup given by Kempf (c.f.
Theorem \ref{kempftheoremquivers}), which is maximally
destabilizing in the GIT sense, gives a filtration of
subrepresentations. We prove that it coincides with the
Harder-Narasimhan filtration for the initial representation, in Theorem
\ref{kempfHNquivers}.

This case is slightly different because, in the previous ones, the
group we were taking the GIT quotient by in the construction of
the moduli space, was $\SL(N)$, but in this case it is a product of
general linear groups, one for each vertex of the quiver. Then,
the \emph{length} we choose in the space of $1$-parameter subgroups when
defining the Kempf function (c.f. Definition \ref{length}) depends
on some parameters (one for each simple factor in the group), and
we show how we have to set the parameters conveniently by choosing
a particular length, to be able to relate the $1$-parameter
subgroup GIT maximally destabilizing with the Harder-Narasimhan
filtration.

Finally, in section \ref{sectionqsheaves} we define $Q$-sheaves,
which are representations of a quiver on the category of coherent
sheaves, and give a stability notion for them, following
\cite{AC,ACGP}. For a one vertex quiver, a $Q$-sheaf is the same
that a coherent sheaf and we use the functorial construction of a moduli space for sheaves given
in \cite{ACK}. In this construction, sheaves are related to
Kronecker modules, then the stability condition turns out to be
rewritten in terms of representations of quivers on vector spaces.
Theorem \ref{equivstability} relates all different notions of
stability involved in this thesis, say, stability of the sheaf as
a $Q$-sheaf (equivalent to Gieseker's stability for sheaves),
stability of the Kronecker module associated, stability of the
representation of another quiver associated $\tilde{Q}$, and GIT
stability of the corresponding point in the parameter space. Using
the equivalence of different stabilities, we can apply the Kempf
theorem (c.f. Theorem \ref{kempftheoremquivers}) to find a
Harder-Narasimhan filtration for the associated representation of
a quiver on vector spaces and, from it, obtain the
Harder-Narasimhan filtration for the $Q$-sheaf in Theorem
\ref{HNqsheaves}.

\section*{Conclusions}

This thesis contains ideas in two directions. On the one hand, we
explore the relation between stability conditions and
GIT stability notions in constructions of moduli spaces. On the other
hand, we relate a natural concept of GIT maximal unstability (in the sense of Kempf) and the
Harder-Narasimhan filtration, which gives a correspondence in
several cases where the Harder-Narasimhan filtration is previously
known or gives a new notion of a such filtration in other cases.

The sketch of the proof is similar in all cases. First, obtaining a filtration of vector subspaces of global sections which maximizes the Kempf function
(c.f. Definition \ref{kempffunction}) for the GIT problem
considered, then evaluate the sections to get a filtration of subobjects, called \emph{the $m$-Kempf filtration}.
We relate the Kempf function with a function in the Euclidean
space (c.f. Proposition \ref{identification}) to apply results on convexity (c.f. subsection \ref{convexcones}). Next, we prove properties of the
$m$-Kempf filtration which characterize it and will make it not to depend on the integer $m$, hence we get a filtration called
\emph{the Kempf filtration}. Finally, in cases where the Harder-Narasimhan filtration is known, we prove that the convexity properties of the Kempf filtration are,
precisely, the ones of the Harder-Narasimhan filtration, hence by uniqueness both filtrations do coincide (c.f. Theorem \ref{kempfisHN}). In other cases,
as rank $2$ tensors in section \ref{kempfrk2}, Kempf filtration (which is unique) defines a Harder-Narasimhan filtration. Note that, in section
\ref{qvectorspaces}, the moduli construction does not depend on any integer, hence we do not have to prove an analogous to Theorem
\ref{kempfstabilizes}, and the correspondence is much easier and quicker.

The notion of \emph{length} (c.f. Definition \ref{length}) we need to define the \textit{speed} of the $1$-parameter subgroups (c.f. section \ref{kempfsection}) plays an important role.
In principal, different lengths would give different $1$-parameter subgroups GIT maximally destabilizing, hence different Kempf filtrations candidates to
be the Harder-Narasimhan filtration. In section \ref{qvectorspaces} we observe how, different choices of length correspond to different definitions of
stability in Definition \ref{Qstability}. Hence, given a notion of Harder-Narasimhan filtration (depending on the notion of stability), we have to set
the parameters in the linearization in the construction of the moduli space, and in the definition of the length in the set of $1$-parameter subgroups (c.f. Proposition
\ref{GITstab-stabquivers}), conveniently, in order to achieve a correspondence between the Kempf and the Harder-Narasimhan filtrations.

Potentially, we could expect to use these ideas to
define Harder-Narasimhan filtrations in cases where it has not
been studied, and use them as a powerful tool for very important
applications where it has been used in the past, as restriction
theorems or calculation of Betti numbers and rational points of
moduli spaces (c.f. \cite{HN}).

Many of the cases where we have applied the method fall into
abelian categories, where the Harder-Narasimhan filtration verifies the
convexity properties out of its definition. It would be interesting
to understand if this imposes a condition to expect a notion of Harder-Narasimhan
filtration, as we usually know.

\vspace{2cm}

The main results of this thesis are collected on the preprints
\cite{GSZ,Za1,Za2}. First one contains the correspondence between
GIT maximal way of destabilizing and the Harder-Narasimhan
filtration for torsion free coherent sheaves over projective
varieties and holomorphic pairs (sections \ref{kempfsheaves} and
\ref{kempfpairs}). Second one contains the result for
representations of quivers on finite dimensional vector spaces
(section \ref{qvectorspaces}). Third one focuses on the case of
rank $2$ tensors and stable coverings.

\chapter{Moduli spaces and maximal unstability}

\section[Constructions of moduli spaces using Geometric Invariant Theory]
{Constructions of moduli spaces using Geometric Invariant Theory}
\sectionmark{Constructions of moduli spaces using GIT}
\label{modulis}

\subsection{Moduli problems}

Since decades, the study of moduli spaces seems to be the right answer to various classification problems in
algebra and geometry. A general classification problem should consist on a collection of objects $\mathcal{A}$
and an equivalence relation $\sim$ on $\mathcal{A}$. The problem is, then, to describe the set of equivalence classes $\mathcal{A}/\sim$.
We usually refer to $\mathcal{A}/\sim$ as the \emph{quotient space}.

In principle, we can think of the solution to our problem as just
the quotient set where each equivalence class corresponds to a
point. But in the field of algebraic geometry, the objects we are
dealing with have rich algebraic and geometric structures, so we
would like this quotient set to have similar properties. Besides,
it is usual to have \textit{continuous families} of objects in
$\mathcal{A}$ and we want to reflect this fact in the quotient
space. In other words, if two objects in $\mathcal{A}$ are very
close, or are very similar (more similar than other objects in
$\mathcal{A}$), we want them to be also very close in the
quotient.

Thus, the ingredients of a moduli problem are three: the class of objects $\mathcal{A}$ we are trying to classify,
the equivalence relation $\sim$ and a notion of \emph{family} and equivalence of families. The object of the theory of moduli spaces is to provide
\textit{good spaces} (to be defined later, meaning spaces with good algebraic and geometric properties) for the quotient space
$\mathcal{A}/\sim$.

Sometimes there are discrete invariants which divide $\mathcal{A}/\sim$ into a countable number of subsets, but
this does not give a complete solution, usually. However, in many cases, in
order to consider a useful moduli space we fix these invariants and try to classify the subclass of these objects.
Examples of this are fixing rank and degree when studying the moduli space of vector bundles over a Riemann
surface, or fixing dimension and degree to consider the moduli space of hypersurfaces in a projective space.

First basic examples of moduli spaces can be the complex projective space $\mathbb{P}_{\mathbb{C}}^{n}$ as the space of lines in $\mathbb{C}^{n+1}$ which pass through the origin or,
more generally, the Grassmannian $\mathcal{GR}(k,n)$ as the moduli space of all $k$-dimensional linear subspaces of $\mathbb{C}^{n}$.

Another classic example is to construct a moduli space for the
collection $\mathcal{A}$ of all the non-singular complex cubics.
Two curves $X$ and $X'$ are equivalent, $X\sim X'$, if they are
isomorphic. By a change of coordinates we consider that all of them are of the
form $y^{2}=x(x-1)(x-\lambda)$, where $\lambda \in \mathbb{C}$.
Then we define
$$j(X)=j(\lambda)=2^{8}\dfrac{(\lambda^{2}-\lambda+1)^{3}}{\lambda^{2}(\lambda-1)^{2}}\; ,$$
called the $j$\emph{-invariant} of the curve $X$. It can be proved that two cubics are equivalent, i.e. $X\sim X'$, if and only if $j(X)=j(X')$ (c.f. \cite[Theorem 4.1]{Ha}).
We see that all the non-singular complex cubics are parametrized by the affine
complex line (an algebraic variety), the corresponding points in the line given by the $j$-invariant of each isomorphism class of curves.
Hence, to classify cubics up to isomorphism is the same that to give a $1$-dimensional variety where each point corresponds to a class of cubics. 

The fact of having many non-trivial automorphisms for some of the objects being classified makes it difficult to
have a moduli space as the set of isomorphism classes. This will be the object of study of the theory of stacks which we will not face here.
Stacks can give a different answer for the classification problem. Indeed, a stack problem is formulated as a $2$-functor problem, whose
answer falls in a more general category of spaces. To avoid that, in many cases we restrict the class of
objects $\mathcal{A}$ we are trying to classify to some subclass for which we will be able to give a moduli space.
The best example of this is the notion of stability for vector bundles or sheaves, where
we can give a solution for the moduli problem when restricting to the \emph{semistable} objects.

In the same direction, it is often possible to consider a modified moduli problem, meaning to
classify the original objects together with additional data, chosen in such a way that the identity is the only
automorphism also respecting the additional data. This choice of additional data is usually called \emph{to
rigidify the objects} or \emph{to rigidify the data}. With a suitable choice of the rigidifying data, the modified moduli problem will have a
moduli space. One of the most successful approaches to construct moduli spaces is this rigidifying the data. Consider an object
$A\in \mathcal{A}$ and suppose it to be enriched to $(A,\alpha)$, where $\alpha$ represents an additional data. In this situation,
an action by a group $G$ appears, taking $(A,\alpha)$ to $(A,\alpha')$, this is, changing the additional data for a given
object $A\in \mathcal{A}$. Hence, in our moduli problem, two objects will be equivalent if they lay on the same orbit by the
action of the group $G$. Then, in order to get rid of this choice of data, we have to quotient by the action of $G$. This is the
object of the Geometric Invariant Theory, developed by David Mumford, which provides moduli spaces as quotients of affine or
projective spaces by the action of groups.

The origin of the theory of moduli spaces started with the theory of elliptic functions, where one can show that there exists
a continuous family of these functions parametrized by the complex numbers, as in the previous example. Riemann showed in a famous
article in 1857 (c.f. \cite{Ri}) that there is a $3g-3$ dimensional family of complex structures a compact topological
surface of genus $g\geq 2$ can be endowed with. In this paper, it was coined the term \textit{moduli}, referring to the number of
parameters for the complex structure.

The modern formulation of moduli problems and definition of moduli spaces dates back to Alexander Grothendieck, (1960/1961), 
\textit{\textquotedblleft Techniques de construction en g\'eom\' etrie analytique. I. Description axiomatique de l'espace de Teichm\"{u}ller et de ses variantes\textquotedblright} 
(c.f. \cite{Gr}) in which he described the general framework, approaches and main problems using Teichm\"{u}ller spaces in complex analytic 
geometry as an example. The text describes a general method to
construct moduli spaces.

Another general approach is primarily associated with Michael Artin. Here the idea is to start with any object of the kind to be 
classified and study those objects which are closer to it, in the sense that they can be seen as deformations of the object. This is called \emph{deformation theory}.

\subsection{Formulation of moduli problems}
Given a moduli problem, i.e. a class of objects $\mathcal{A}$, an equivalence relation $\sim$ between objects and a notion of
family and equivalence of families, we want to give an algebraic structure or geometric structure to the set $\mathcal{A}/\sim$. This structure will
depend on the category we are working on and the precise context (it can be an algebraic variety, an scheme or an algebraic space, for example). In the following, we
will consider the category $\Sch_{k}$ of schemes over a field $k$, and recall that this category has fiber products. Let us denote by $\Sets$ the category of sets.

By a family of objects in $\mathcal{A}$ we understand a proper flat morphism of $k$-schemes $f:X\rightarrow S$, where fibers $X_{s}$ (i.e. $X_{s}$ is the pull back of $f$ along
the inclusion $s\in S$)
of the morphism $f$ are objects in $\mathcal{A}$. We say that $X$ is a family of objects in $\mathcal{A}$ parametrized by $S$.

To formulate a moduli problem we need that the equivalence relation $\sim$ verifies certain conditions (c.f. \cite[Conditions 1.4]{Ne})
\begin{itemize}
\item A family parameterized by a one point scheme $\{p\}$ is a single object of $\mathcal{A}$.
\item There exists a notion of equivalence between families reducing to $\sim$ for single objects in $\mathcal{A}$.
Then equivalence of objects turns out to be equivalence of families parametrized by $\{p\}$.
\item The equivalence for families is functorial, i.e. for any morphism $\varphi: S'\rightarrow S$
and a family $X$ parameterized by $S$ (i.e. $f:X\rightarrow S$), there is an induced family $\varphi^{\ast}X$ parameterized by $S'$ and
this operation satisfies functorial properties.
\end{itemize}

\begin{dfn}
\label{moduliproblem}
Let $\mathcal{A}$ be a class and let $\sim$ be an equivalence relation for families in $\mathcal{A}$. A \emph{moduli functor} is a contravariant functor
$$\mathcal{F}:\Sch_{k}\rightarrow \Sets$$
where $\mathcal{F}(S)$ denotes the set of equivalence classes of
families parameterized by $S$. The triple $(\mathcal{A},\sim,\mathcal{F})$ is called a \emph{moduli problem}.
\end{dfn}

Suppose that $M$ is a $k$-scheme with underlying set
$\mathcal{A}/\sim$. To have a family $X$ of objects in
$\mathcal{A}/\sim$ parameterized by a $k$-scheme $S$ is the same
that a map $\nu_{[X]}:S\rightarrow M$ and we would like all the
different morphisms $\nu_{[X]}:S\rightarrow M$ to be in
correspondence with the different equivalence classes of families
$[X]$ parameterized by $S$. In the language of categories and
functors this is expressed with the moduli functor in Definition
\ref{moduliproblem}. Let $\Hom(-,M)$ be the functor
of points of $M$. Recall that the functor of points of a
$k$-scheme $M$ is the contravariant functor from the category of
$k$-schemes to the category of sets, which sends a $k$-scheme $S$
to the set of morphisms from $S$ to $M$. There is a natural
transformation
$$\Phi:\mathcal{F}\rightarrow \Hom(-,M)$$
where $\Phi_{S}:\mathcal{F}(S)\rightarrow \Hom(S,M)$ is the natural map given by $\Phi_{S}([X])=\nu_{[X]}$.

To give a \emph{moduli problem} is to give a functor $\mathcal{F}$ as in Definition \ref{moduliproblem} and ask if there exists any
$k$-scheme $M$ such that $\mathcal{F}$ and the functor of points of $M$ are related, meaning that the set of
equivalence classes of families parameterized by $S$, $\mathcal{F}(S)$, is related with the set of different
morphisms from $S$ to $M$. In particular, for a one point scheme $\{p\}$, $\mathcal{F}(p)$ will be the set of equivalence classes of objects, so will
be in correspondence with the \emph{points} of $M$, $\Hom(\{p\},M)$. Hence, such $M$ will be the moduli space we are seeking.

\begin{dfn}
\label{represents}
A moduli functor $\mathcal{F}:\Sch_{k}\rightarrow \Sets$ is \emph{representable} if there exists a $k$-scheme $M$ such that $\mathcal{F}$ is
isomorphic to the functor of points of $M$
$\Hom(-,M)$. Denote such isomorphism by $\Phi$ and say that the pair $(M,\Phi)$ \emph{represents} the functor $\mathcal{F}$. A \emph{fine moduli space} for
the moduli problem considered
is a pair $(M,\Phi)$ which represents the functor $\mathcal{F}$.
\end{dfn}

Note that, by Definition \ref{represents}, if $(M,\Phi)$ represents $\mathcal{F}$ we have a natural bijection
$$\Phi(p): \mathcal{A}/\sim=\mathcal{F}(p)\rightarrow \Hom(p,M)=M$$
where $p$ is a one point $k$-scheme. Moreover, the identity morphism
$1_{M}$ determines, up to equivalence, a family $U$ parameterized
by $M$ such that every family $X$ parameterized by a $k$-scheme
$S$ is equivalent to $\nu^{\ast}_{X}U$, where
$\nu^{\ast}_{X}:S\rightarrow M$ is the morphism corresponding to
the family. The family $U$ is called a \emph{universal family} for
the moduli problem considered. Therefore, we can define a
\emph{fine moduli space} as a $k$-scheme $M$ together with a
\emph{universal family} $U$ parameterized by $M$ such that every
family is given as the pull back from $U$ by the corresponding
morphism.

\begin{dfn}
\label{corepresents}
A moduli functor $\mathcal{F}:\Sch_{k}\rightarrow \Sets$ is \emph{corepresentable} if there exists a $k$-scheme $M$ and a natural transformation $\Phi:\mathcal{F}\rightarrow \Hom(-,M)$ to the
functor of points of $M$ such that, for every $k$-scheme $N$ and a natural transformation $\Phi':\mathcal{F}\rightarrow \Hom(-,N)$, there exists a unique natural transformation
$\Psi:\Hom(-,M)\rightarrow \Hom(-N)$ such that $\Phi$ factors through $\Psi$. Such pair $(\Phi,M)$ is said to \emph{corepresent} the functor $\mathcal{F}$ and, if it exists, it is
unique up to unique isomorphism. If furthermore, $\Phi(p):\mathcal{F}(p)\rightarrow M$ is bijective, where $p$ is a one point $k$-scheme, 
we say that $(M,\Phi)$ is a \emph{coarse moduli space}
for the moduli problem considered.
\end{dfn}

There are many moduli problems for which we cannot find a fine moduli space. The existence of a coarse moduli space turns out to be a weaker solution. Note that, if $(M,\Phi)$ is
a fine moduli space, it is automatically a coarse moduli space.

One reason for the non existence of a moduli space with good
properties, easy to explain, is the \emph{jump phenomenon}. It
happens when there exists a family $X$ parametrized by a scheme
$S$ of dimension $\geq 0$ for which there is a point $s_{0}\in S$
such that
\begin{itemize}
 \item $X_{s}\sim X_{t}$ for all $s,t\in S-\{s_{0}\}$
 \item $X_{s}\nsim X_{s_{0}}$ for all $s\in S-\{s_{0}\}$
\end{itemize}
With this feature, if we include in the hypothetical moduli space the equivalence or the isomorphism class of $X_{s_{0}}$, the moduli space would be non separated. This is
the usual property shared by the \emph{unstable} objects (those which behave like $X_{s_{0}}$). The notion of stability was introduced first by Mumford, in order to construct moduli
spaces for the subclass of \emph{semistable} objects.

As an example, we can formulate the problem of finding a moduli space of algebraic curves of genus $g$. Consider the class $\mathcal{A}$ of
smooth projective curves of genus $g$ over an algebraically closed field $k$, and the equivalence relation $\sim$ being the isomorphism between curves. A family of curves parametrized by
$S$ is a proper flat morphism $f:X\rightarrow S$ between algebraic varieties where fibers are curves of genus $g$. There exists a moduli space, denoted $\mathcal{M}_{g}$,
for this moduli problem. Define a curve to be \emph{stable} if it is complete, connected,
has no singularities other than double points, and has only a finite group of automorphisms. The moduli space of stable curves of genus $g$ is usually denoted by
 $\overline{\mathcal{M}_{g}}$. The space $\overline{\mathcal{M}_{g}}$ is projective and it is a compactification of $\mathcal{M}_{g}$.

\subsection{Results on Geometric Invariant Theory}

In this section we recall the basic results of Geometric Invariant Theory we need when taking quotients by the action of groups in moduli problems.

Let $G$ be an algebraic group over an algebraically closed field
$k$. A right action on an scheme $X$ is a morphism $\sigma:X\times
G\to X$, where $\sigma(x,g)= x\cdot g,\; \forall x\in X$, such that
$x\cdot (gh)=(x\cdot g)\cdot h$ and $x\cdot e=x$, $e$ being the
identity element of $G$. A left action is defined by $(hg)\cdot x=
h\cdot (g\cdot x)$.

We denote by $x \cdot G$ the orbit of $x\in X$ by a right action of $G$ (resp. $G\cdot x$ for a left action). A morphism $f:X \to Y$ between two varieties endowed with
$G$-actions is called \emph{$G$-equivariant} if it commutes with the actions, that is $f(x)\cdot g= f(x\cdot
g)$. In the case that the action on $Y$ is trivial (i.e. $y\cdot g=y$, for all $g\in G$ and $y\in Y$), then a morphism $f$ which is $G$-equivariant is called 
\emph{$G$-invariant}.

If $X$ is an affine scheme, to construct affine
quotients is much simpler when the group $G$ is reductive. Recall
that $G$ is \emph{reductive} if its radical is isomorphic to a
direct product of copies of $k^{\ast}$. On the other hand, $G$ is
\emph{geometrically reductive} if, for every linear action of $G$
on $k^{n}$, and every $G$-invariant point $v$ of $k^{n}$, $v\neq 0$,
there exists a $G$-invariant homogeneous polynomial $f$ of degree
$\geq 1$ such that $f(v)\neq 0$. Due to results of Weil, Nagata,
Mumford and Haboush, it turns out that every reductive group is
geometrically reductive and, if a reductive group $G$ is acting on
a finitely generated $k$-algebra $R$ (as it is the ring of functions of an
affine variety $X$, $R=A(X)$), the ring of invariants $R^{G}$ is
finitely generated. Therefore, we define the quotient of an affine
variety $X$ by the action of a reductive group $G$, as the affine
variety whose ring of functions is $A(X)^{G}$.

The following example shows that the quotient of an affine scheme $X$ by the action of a reductive group $G$ can differ of 
an orbit space (c.f. Definition \ref{geometricquotient}), because the quotient $A(X)^{G}$ can possibly identify different orbits in the same point in the quotient space. 

\begin{ex}
\label{hyperboles} Consider the action
$$\xymatrix{ \sigma:\mathbb{C}^{\ast}\times \mathbb{C}^{2}\ar[r] & \mathbb{C}^{2}\\
 (\lambda,(x,y))\ar@{|->}[r] & (\lambda x,\lambda^{-1} y)}$$
whose orbits are represented in Figure \ref{hypfig}.
\begin{figure}[h]
   \begin{center}
   \includegraphics[width=6cm]{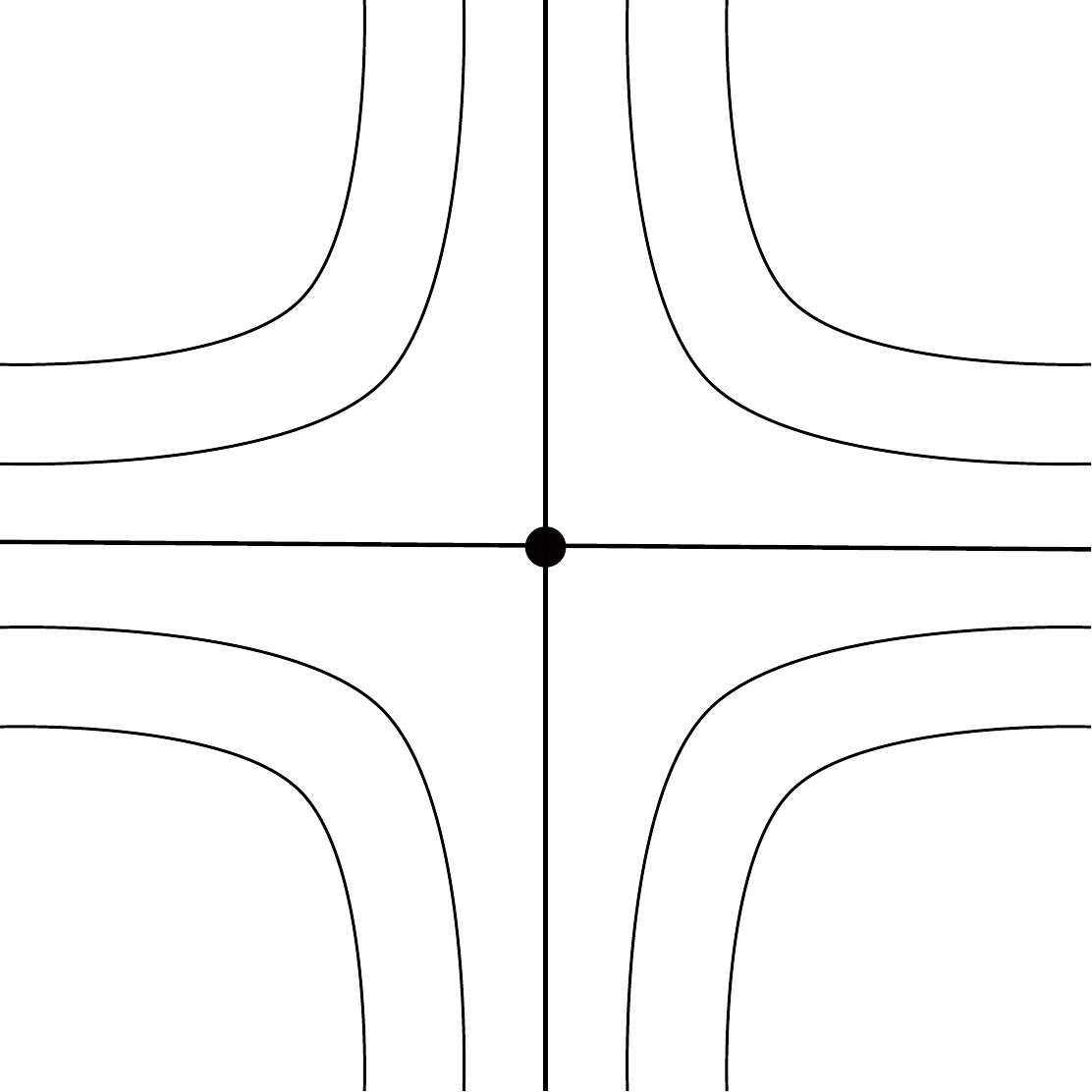}
   \end{center}
\caption{Orbits of the action in Example \ref{hyperboles}}
\label{hypfig}
\end{figure}
The orbits are the hyperboles $xy=constant$, plus three special
orbits, the $x$-axis, the $y$-axis and the origin. Observe that
the origin is in the closure of the $x$-axis and the $y$-axis.

The ring of functions of $\mathbb{C}^{2}$ is $\mathbb{C}[X,Y]$ and
the ring of invariants is
$\mathbb{C}[X,Y]^{\mathbb{C}^{\ast}}\simeq \mathbb{C}[XY]\simeq
\mathbb{C}[Z]$. So, the ring of invariants does not distinguish
between the three special orbits, and identifies them in a unique
single point in the quotient space. Hence, the orbit space (the
space where each point corresponds to an orbit) would be non
separated, but the quotient space whose ring of functions is
$\mathbb{C}[X,Y]^{\mathbb{C}^{\ast}}\simeq \mathbb{C}[Z]$ is the
affine line, which is separated.
\end{ex}

The case when $G$ acts on a projective scheme $X$ is more
complicated. We call  $ $   $\psi:G\times X\rightarrow X$, a
\emph{linearization} of the action on an ample line bundle
$\mathcal{O}_{X}(1)$. It consists of giving an action on the total
space $L$ of the line bundle $\mathcal{O}_{X}(1)$, $\sigma:G\times
L\longrightarrow L$, such that for every $g\in G$ and $x\in X$,
there exists a isomorphism which takes one fiber onto another
$L_{x}\longrightarrow L_{g\cdot x}$ (i.e. $\sigma$ is linear along
the fibers and the projection $L\rightarrow X$ is $G$-equivariant). A
linearization is the same thing as giving, for each $g\in G$, an
isomorphism of line bundles
$\widetilde{g}:\mathcal{O}_{X}(1)\longrightarrow \varphi_{g}^{*}
\mathcal{O}_{X}(1)$, ($\varphi_{g}=\psi(g,\cdot)$) which also
satisfies the previous associative property. We say also that
$\sigma=\widetilde{\psi}$ is a lifting to $L$ of the action
$\psi$:
$$\xymatrix{
G\times L \ar^(.6){\sigma=\widetilde{\psi}}[r]\ar[d] & L\ar[d]\\
G\times X \ar^(.6){\psi}[r] & X}$$

If $\mathcal{O}_{X}(1)$ is very ample, then a linearization is the same thing as a representation of $G$ on the
vector space $H^{0}(\mathcal{O}_{X}(1))$ such that the natural embedding $$X\hookrightarrow
\mathbb{P}(H^{0}(\mathcal{O}_{X}(1))^{\vee})$$ is $G$-equivariant.

Then, if we have a group $G$ acting on a projective scheme $X$ and
consider the set of orbits $X/G$, when can we define $X/G$ as a
scheme $M$, i.e., the points of $X/G$ correspond, in a natural way,
to the points of $M$?

The next example (c.f. \cite{Gi2}) illustrates some of the features which can arise when trying to
define $X/G$.

\begin{ex}\cite{Gi2}
\label{exG} Let $N$ be an integer and consider the set of all
homogeneous polynomials of degree $N$ in two variables,
$V_{N}=\{\underset{i}{\sum}a_{i}X_{0}^{i}X_{1}^{N-i}\}$. Let
$\mathbb{P}(V_{N})$ be its projectivization. The group
$G=SL(2,\mathbb{C})$ acts on $V_{N}$ as
$$P^{g}(X_{0},X_{1})=P(g^{-1}\left(
                               \begin{array}{c}
                                 X_{0} \\
                                 X_{1} \\
                               \end{array}
                             \right))$$ where $P\in V_{N}$.
The vanishing locus of each $P\in V_{N}$ consists of a finite set
of points in $\mathbb{P}^{1}$ where their multiplicities are the
orders as zeroes of $P$. Then, we can think of $\{P=0\}$ as a
divisor $D_{f}$ on $\mathbb{P}^{1}$, and $\mathbb{P}(V_{N})$ as
the space of divisors of degree $N$ on $\mathbb{P}$. Observe that
$G$ acts on divisors moving them by linear fractional
transformations.

The orbit space $\mathbb{P}(V_{N})/G$ is not a variety, because it is
not Hausdorff. To see this, let $\overline{P}\in
\mathbb{P}(V_{N})$ and let $P\in V_{N}$ be a polynomial in the corresponding line. We look for an element $Q$ in the orbit of $P$ so that
$X_{1}^{N}$ occurs
  in $Q(X_{0},X_{1})$
(i.e.,
$Q(X_{0},X_{1})=a_{0}X_{0}^{N}+a_{1}X_{0}^{N-1}X_{1}+...+a_{N-1}X_{0}X_{1}^{N-1}+X_{1}^{N}$).
Let $Q_{t}(X_{0},X_{1})=t^{N}Q(tX_{0},t^{-1}X_{1})$ and note that
$Q_{t}(X_{0},X_{1})$ lays in the orbit of $P$ and $Q$ for every
$t\neq 0$, since $Q_{t}(X_{0},X_{1})=t^{N}\cdot
Q^{g_{t}}(X_{0},X_{1})$, with $g_{t}=\left(
                               \begin{array}{cc}
                               t^{-1} & 0 \\
                               0 & t \\                                                        \end{array}
                               \right)$, and all of them give the same point in $\mathbb{P}(V_{N})$.
Therefore, $\mathbb{P}(V_{N})/G$ cannot be given a Hausdorff
topology so that $\phi:\mathbb{P}(V_{N})\longrightarrow
\mathbb{P}(V_{N})/G$ is continuous. Indeed, if $\phi$ were
continuous, it would be
$$\underset{t\rightarrow 0}{\lim Q_{t}(X_{0},X_{1})}=Q_{0}(X_{0},X_{1})=X_{1}^{N},$$
we would have $$\phi(P)=\phi(Q)=\underset{t\rightarrow 0}{\lim
\phi(Q_{t})}=\phi(X_{1}^{N})$$ and the image of $\phi$ would be
one single element. The reason of this is that the polynomial
$X_{1}^{N}$ is not in the orbit of $f$ and $g$, but it is in its
adherence. Then, when we try to define a continuous quotient map,
the adherent orbits have to go to the same point.
\end{ex}

As we have seen in Examples \ref{hyperboles} and \ref{exG}, in
order to obtain a quotient space with good properties (for example,
being Hausdorff), we have to make some considerations about
the orbits of the action of the group $G$, putting together in the
quotient space all orbits whose closures have non empty
intersection. We will call two of these orbits \emph{S-equivalent}
(c.f. Remark \ref{Seq}).

Geometric Invariant Theory, abbreviate GIT, will be a technique to
construct such quotients with good properties.

\begin{dfn}
\label{categoricalquotient} Let $X$ be a scheme endowed with a
$G$-action. A \emph{categorical quotient} is a scheme $M$ with a
$G$-invariant morphism $p:X\longrightarrow M$, such that for every 
scheme $M'$, and every $G$-invariant morphism $p':X\longrightarrow M'$, there is a
unique morphism $\varphi$ with $p'=\varphi\circ p$
$$\xymatrix{
{X} \ar_{p}[d] \ar^{p'}[rd]\\
{M} \ar@{-->}_{\exists !\, \varphi}[r]& {M'}}$$
\end{dfn}

\begin{dfn}
\label{goodquotient} Let $X$ be a scheme endowed with a
$G$-action. A \emph{good quotient} is a scheme $M$ with a
$G$-invariant morphism $p:X\longrightarrow M$ such that
\begin{enumerate}
\item $p$ is surjective and affine.
\item $p_{*}(\mathcal{O}^{G}_{X})=\mathcal{O}_{M}$, where $\mathcal{O}^{G}_{X}$ is the sheaf of $G$-invariant functions on $X$.
\item If $Z$ is a closed $G$-invariant subset of $X$, then $p(Z)$ is closed in $M$. Furthermore, if $Z_{1}$ and $Z_{2}$ are two closed $G$-invariant subsets of
$X$ with $Z_{1}\cap Z_{2}=\emptyset$, then $f(Z_{1})\cap f(Z_{2})=\emptyset$.
\end{enumerate}
\end{dfn}

\begin{dfn}
\label{geometricquotient}
A \emph{geometric quotient} is a good quotient $p:X\to M$ such that $p(x_{1})=p(x_{2})$ if and only if the orbit of $x_{1}$ is equal to the orbit of $x_{2}$.
\end{dfn}

Note that a geometric quotient is a good quotient, and a good quotient is a categorical quotient.

Let $X$ be a projective scheme, let $G$ be a reductive algebraic
group and an action $\sigma:G\times X\longrightarrow X$ of $G$ on
$X$. We call $\tilde{\sigma}$ a linearization of the action on an
ample line bundle $\mathcal{O}_{X}(1)$.

\begin{dfn}
A closed point $x\in X$ is called \emph{GIT semistable} if, for some $m>0$, there is a $G$-invariant section $s$
of $\mathcal{O}_{X}(m)$, $s\in H^{0}(X,\mathcal{O}_{X}(m))$, such that $s(x)\neq 0$. If, moreover, the orbit of
$x$ is closed in the open set of all GIT semistable points, it is called \emph{GIT polystable} and, if
furthermore, this closed orbit has the same dimension as $G$ (i.e. if $x$ has finite stabilizer), then $x$ is
called a \emph{GIT stable} point. We say that a closed point of $X$ is \emph{GIT unstable} if it is not GIT
semistable.
\end{dfn}

With this definition, the stable points are precisely the polystable points with finite stabilizer.

\begin{rem}
We consider $X$ embedded in a projective space by the ample line bundle $\mathcal{O}_{X}(1)$,
$$X\hookrightarrow \mathbb{P}(H^{0}(\mathcal{O}_{X}(1))^{\vee})=\mathbb{P}(V)\; .$$
Then, we can see a section $s\in H^{0}(\mathcal{O}_{X}(m))$ as a
homogeneous polynomial of degree $m$ in $V$. Then, the GIT
unstable points are those for which, for all $m>0$, all
$G$-invariant homogeneous polynomials vanish at the point. As all
homogeneous polynomials (in particular the $G$-invariant ones)
vanish at zero, the points which contains zero in the closure of
their orbits will be GIT unstable.

This idea of considering invariant homogeneous polynomials comes
from Hilbert who calls them \emph{nullforms} in \cite{Hi}.
\end{rem}

The central result of Mumford's Geometric Invariant Theory is the following theorem:

\begin{thm}\cite[Proposition 1.9, Theorem 1.10]{Mu}
\label{GIT} Let $X^{ss}$ (respectively, $X^{s}$) be the subset of GIT semistable points (respectively,
GIT stable). Both $X^{ss}$ and $X^{s}$ are open subsets. There is a good quotient $X^{ss}\longrightarrow
X^{ss}/\!\!/G$ (where closed points are in one-to-one correspondence to the orbits of GIT polystable points),
the image $X^{s}/\!\!/G$ of $X^{s}$ is open,
 $X/\!\!/G$ is projective, and the restriction $X^{s}\to X^{s}/\!\!/G$ is a geometric quotient.
\end{thm}

Hence, to construct good quotients, first we have to get rid of the unstable points. To find these unstable
points there exists a numerical criterion based on the use of $1$-parameter subgroups of $G$. It was first used
by Hilbert and later by Mumford, to characterize the GIT stability.

\begin{dfn}
Let G be an algebraic group over the field $k$. 
A \emph{$1$-parameter subgroup} of $G$, $\Gamma$, is a non-trivial algebraic homomorphism $\Gamma: k^{\ast}\longrightarrow G$.
\end{dfn}

Let $X$ be a projective scheme where the group $G$ acts. Suppose
that this action is linearized on a line bundle
$\mathcal{O}_{X}(1)$ and call the linearization $\sigma$. Then,
given $\Gamma$, a $1$-parameter subgroup of $G$, and given $x\in X$, we
can define $\Phi:k^{*}\longrightarrow X$ by
$\Phi(t)=\Gamma(t)\cdot x$. We say $\underset{t\rightarrow 0}{\lim
\Gamma(t)\cdot x}=\infty$ if $\Phi$ cannot be extended to a map
$\widetilde{\Phi}:k\longrightarrow X$. If $\Phi$ can be extended,
we write $\underset{t\rightarrow 0}{\lim \Gamma(t)\cdot x}=x_{0}$.

Then, the criterion is the following:

\begin{thm}
\label{HMcrit0} Let $\tilde{x}$ be a point in the affine cone over
$X$, lying over $x\in X$. With the previous notations:
    \begin{itemize}
    \item $x$ is \emph{semistable} if for all $1$-parameter subgroups $\Gamma$, $\exists\underset{t\rightarrow 0}{\lim
    \Gamma(t)\cdot \tilde{x}}\neq 0$
        or $\underset{t\rightarrow 0}{\lim \Gamma(t)\cdot \tilde{x}}=\infty$.
    \item $x$ is \emph{polystable} if it is semistable and the orbit of $\tilde{x}$ is closed.
    \item $x$ is \emph{stable} if for all $1$-parameter subgroups $\Gamma$,
    $\underset{t\rightarrow 0}{\lim \Gamma(t)\cdot \tilde{x}}=\infty$ (then the stabilizer of $x$ is finite).
    \item $x$ is \emph{unstable} if there exists a $1$-parameter subgroup $\Gamma$ such that
    $\underset{t\rightarrow 0}{\lim \Gamma(t)\cdot \tilde{x}}=0$.
    \end{itemize}
\end{thm}

The point $x_{0}$ is, clearly, a fixed point for the action of $k^{*}$ on $X$ induced by $\Gamma$. Thus, $k^{*}$
acts on the fiber of $\mathcal{O}_{X}(1)$ over $x_{0}$, say, with weight $\gamma$. One defines the numerical
function
$$\mu(\Gamma,x):=\gamma\; .$$
We will call this $\gamma$ the \textit{minimum relevant exponent} of the action of $\Gamma$ over $x$.

With the definition of $\mu(\Gamma,x)$ we can state the Hilbert-Mumford criterion of GIT stability:

\begin{thm}[\textbf{Hilbert-Mumford criterion}]\cite[Theorem
2.1]{Mu}, \cite[Theorem 4.9]{Ne}
\label{HMcrit} With the previous notations:
    \begin{itemize}
    \item $x$ is \emph{semistable} if for all $1$-parameter subgroups $\Gamma$, $\mu(\Gamma,x)\leq
    0$.
    \item $x$ is \emph{stable} if for all $1$-parameter subgroups $\Gamma$,
    $\mu(\Gamma,x)<0$.
    \item $x$ is \emph{unstable} if there exists a $1$-parameter subgroup $\Gamma$ such that
    $\mu(\Gamma,x)>0$.
    \end{itemize}
\end{thm}

\begin{ex}
\label{exG2} Returning to Example \ref{exG} (c.f. \cite{Gi2}),
we apply the Hilbert-Mumford criterion in Theorem \ref{HMcrit}. Let $G=SL(2,\mathbb{C})$
and $V_{N}=\{\underset{i}{\sum}a_{i}X_{0}^{i}X_{1}^{N-i}\}$.
Consider the following $1$-parameter subgroup of $G$
$$\Gamma(t)=\left(
               \begin{array}{cc}
                 t^{-r} & 0 \\
                 0 & t^{r} \\
               \end{array}
             \right),r>0\; .$$
Let $P(X_{0},X_{1})=\sum a_{ij}X_{0}^{i}X_{1}^{j}$ be a polynomial in $V_{N}$. We want to know when
$\underset{t\rightarrow 0}{\lim
P^{\Gamma(t)}}=\underset{t\rightarrow 0}{\lim \Gamma(t)\cdot
P}=0$. It is
$P^{\Gamma(t)}(X_{0},X_{1})=\sum
a_{ij}X_{0}^{i}X_{1}^{j}t^{r(i-j)}$, hence,
$\underset{t\rightarrow 0}{\lim P^{\Gamma(t)}}=0$ implies that
$a_{ij}=0$, if $j\geq i$. This means that, if $P$ has a factor of
$X_{0}^{k}$ with degree $k>\frac{N}{2}$, in that case, $P$ is
unstable. For a general $1$-parameter subgroup, it turns out that
$P$ is semistable if and only if $P$ has no linear factors of
degree greater that $\frac{N}{2}$.
\end{ex}

\begin{rem}
\label{Seq} A theorem of Geometric Invariant Theory (c.f.
\cite[Lemma 1.10]{Si1}) says that, if $G\cdot v$ is the orbit of a
point $v\in V$, in its closure $\overline{G\cdot v}$ there is a
unique orbit $Y\subset \overline{G\cdot v}$ such that $Y$ is
closed in $\overline{G\cdot v}$, so it is closed also in the whole
space $V$. The GIT polystable points are in correspondence with
these closed orbits. Two orbits, $G\cdot v$ and $G\cdot w$, with
the same closed orbit $Y$ in their closures $\overline{G\cdot v}$,
$\overline{G\cdot w}$, are called \emph{S-equivalent}. The points of the
moduli space are in correspondence with these distinguished closed
orbits, so the moduli space we obtain classifies polystable points, or points modulo $S$-equivalence.
\end{rem}

\begin{rem}
\label{GIT_measure} Geometric Invariant Theory states that we can reach every point in the closure of an orbit through
$1$-parameter subgroups. It can be proved (c.f. \cite[Proposition 4.3]{Ne}, \cite[Proposition 2.2]{Mu}) that a point $x$ is GIT semistable if 
$0\notin \overline{G\cdot \hat{x}}$, where $\hat{x}$ lies 
over $x$ in the affine cone. Then, GIT stability measures whether $0$ belongs to the closure of the lifted orbit or not, belonging which can be 
checked through $1$-parameter subgroups. 
\end{rem}

\section[Example: moduli of tensors]{Example of a construction of a moduli space using GIT: Moduli of tensors}
\label{exampletensors}

Here, we present a complete example of the construction of a
moduli space through Geometric Invariant Theory. We construct a
moduli space for tensors over higher dimensional projective
varieties following the Gieseker-Maruyama method. This was
constructed by Alexander Schmitt in \cite{Sch} for curves.

This section follows the paper \cite{GS1}, where the authors carry
out the same construction, but using the method of Simpson.

\subsection{Definitions and stability of tensors}
Let $X$ be a smooth projective variety of dimension $n$ over
$\mathbb{C}$. Fix an ample line bundle $\mathcal{O}_{X}(1)$ on
$X$. Fix a polynomial $P$ of degree $n$, and integers $s,c,b$. Let
$R$ be an scheme and fix a locally free sheaf $\mathcal{D}$ on
$X\times R$, i.e. a family $\{D_{u}\}_{u\in R}$ of locally free
sheaves on $X$ parametrized by $R$, where given a point $u\in R$,
we denote by $D_{u}$ the restriction of $\mathcal{D}$ to the slice
$X\times u$.

\begin{dfn}\cite[Definition 1.1]{GS1}
\label{deftensor}
A \emph{tensor} is a triple $(E,\varphi,u)$, where $E$ is a coherent
sheaf on $X$ with Hilbert polynomial $P_{E}=P$, $u$ is a point in
$R$, and $\varphi$ is a homomorphism
$$\varphi:(E^{\otimes s})^{\oplus c}\longrightarrow (\det
E)^{\otimes b}\otimes D_{u}\; ,$$ that is not identically zero.
Let $(E,\varphi,u)$ and $(F,\psi,v)$ be two tensors with
$P_{E}=P_{F}$, $\det E\simeq \det F$, and $u=v$. A homomorphism
between $(E,\varphi,u)$ and $(F,\psi,v)$ is a pair $(f,\alpha)$
where $f:E\rightarrow F$ is a homomorphism of sheaves,
$\alpha\in\mathbb{C}$, and the following diagram commutes
\begin{equation}
\label{homtensors} \xymatrix{(E^{\otimes s})^{\oplus
c}\ar[d]^{\varphi}\ar[r]^{(f^{\otimes s})^{\oplus c}} &
(F^{\otimes s})^{\oplus
c}\ar[d]^{\psi}\\
(\det E)^{\otimes b}\otimes D_{u}\ar[r]^{\hat{f}\otimes\alpha} &
(\det F)^{\otimes b}\otimes D_{v}}
\end{equation}
where $\hat{f}:\det E\longrightarrow \det F$ is the homomorphism
induced by $f$. In particular, $(E,\varphi,u)$ and $(E,
\lambda\varphi,u)$ are isomorphic for $\lambda\in
\mathbb{C}^{\ast}$.
\end{dfn}

\begin{rem}
This notion of isomorphism can be restricted by considering only
isomorphisms for which $\alpha=1$. In this case we would obtain
another category where, for example, if $E$ is simple, the set of
automorphism of $(E,\varphi,u)$ is $\mathbb{C}^{\ast}$, but if
$\alpha=1$, the set of automorphisms is
$\mathbb{Z}/(rb-s)\mathbb{Z}$ (provided $rb-s\neq 0$). If
$rb-s\neq 0$, the set of isomorphism classes will be the same
(changing $f$ into $\alpha^{1/(rb-s)}f)$), and then the moduli
spaces will be the same. If $rb-s=0$, the set of isomorphism
classes is not the same.
\end{rem}

Let $\delta$ be a polynomial with $\deg (\delta) < n=\dim X$
\begin{equation}
\label{delta} \delta(t)=\delta_{1}t^{n-1}+\delta_{2}t^{n-2}+\cdots
+\delta_{n}\in \mathbb{Q}[t],
\end{equation}
and $\delta(m)>0$ for $m\gg 0$. We denote $\tau =\delta_{j}(n-j)!$
where $\delta_{j}$ is the leading coefficient of $\delta$. We will
define a notion of stability for these tensors, which depends on
the polarization $\mathcal{O}_{X}(1)$ and $\delta$, and we will
construct, using Geometric Invariant Theory, a moduli space for
semistable tensors.

A weighted filtration $(E_\bullet,n_{\bullet})$ of a sheaf $E$ is
a filtration
\begin{equation}
\label{filtE} 0 \subset E_1 \subset E_2 \subset\cdots \subset E_t
\subset E_{t+1}=E,
\end{equation}
and rational positive numbers $n_{1},\, n_{2},\ldots , \,n_{t} >
0$. We denote $r_{i}=\rk(E_{i})$. If $t=1$ (what we will call
\emph{one-step filtration}), then we set $n_{1}=1$. The filtration
is called \emph{saturated} if all sheaves $E_{i}$ are saturated in
$E$, i.e. if $E/E_{i}$ is torsion free for all $i$.

Let $\gamma$ be a vector of $\mathbb{C}^r$ defined as
$\gamma=\sum_{i=1}^{t}n_{i} \gamma^{(\rk E_i)}$ where
$$\gamma^{(k)}:=\big( \overbrace{k-r,\ldots,k-r}^k,
 \overbrace{k,\ldots,k}^{r-k} \big)
\qquad (1\leq k < r) \, .$$ Hence, the vector is of the form
$$\gamma=(\overbrace{\gamma_{r_{1}},\ldots,\gamma_{r_{1}}}^{\rk E^1},
\overbrace{\gamma_{r_{2}},\ldots,\gamma_{r_{2}}}^{\rk E^2},
\ldots, \overbrace{\gamma_{r_{t+1}},\ldots,\gamma_{r_{t+1}}}^{\rk
E^{t+1}}) \; ,$$ where
$n_{i}=\dfrac{\gamma_{r_{i+1}}-\gamma_{r_{i}}}{r}$.

Now let $\mathcal{I}=\{1,...,t+1\}^{\times s}$ be the set of all
multi-indexes $I=(i_{i},...,i_{s})$ and define
\begin{equation}
\label{rightmu} \mu(\varphi,E_{\bullet},n_{\bullet})=\min_{I\in
\mathcal{I}} \{\gamma_{r_{i_1}}+\cdots+\gamma_{r_{i_s}}:
\,\varphi|_{(E_{i_1}\otimes\cdots\otimes E_{i_s})^{\oplus c}}\neq
0 \}.
\end{equation}

If $P_{1}$ and $P_{2}$ are two
polynomials, we write $P_{1}\prec P_{2}$ if $P_{1}(m)<P_{2}(m)$
for $m\gg 0$, and analogously for "$\leq$" and "$\preceq$".

\begin{dfn}\cite[Definition 1.3]{GS1}
\label{stabilityfortensors}
Let $\delta$ be a polynomial as in (\ref{delta}). We say that
$(E,\varphi,u)$ is \emph{$\delta$-semistable} if for all weighted
filtrations $(E_{\bullet},n_{\bullet})$ of $E$, it is
\begin{equation}
\label{stabtensors} \big(
\sum_{i=1}^{t}n_{i}(rP_{E_{i}}-r_{i}P_{E})\big) + \delta
\mu(\varphi,E_{\bullet},n_{\bullet})\preceq 0
\end{equation}
We say that $(E,\varphi,u)$ is \emph{$\delta$-stable} if we have a
strict inequality in (\ref{stabtensors}) for every weighted
filtration. If $(E,\varphi,u)$ is not $\delta$-semistable we say
that it is \emph{$\delta$-unstable}.
\end{dfn}

We assume that $\varphi$ is not identically zero, then (\ref{rightmu}) is well defined.

\begin{rem}
\label{saturatedfiltrations}
It is enough to consider saturated filtrations in Definition \ref{stabilityfortensors}. Indeed, it is $P_{E_{i}}\leq P_{\overline{E_{i}}}$ for
Hilbert polynomials, if $\overline{E_{i}}$ is the saturation of a subsheaf $E_{i}\subset E$.

Also it suffices to consider filtrations with
$\rk(E_{i})<\rk(E_{i+1})$. If not, suppose $E_{i}\subsetneq
E_{i+1}$ and $\rk E_{i}=\rk E_{i+1}$, then $E_{i}$ is not
saturated in $E_{i+1}$ and $E_{i+1}/E_{i}$ has torsion. Therefore
$E/E_{i}$ has torsion and $E_{i}$ is not saturated in $E$. Note
that the definition of (\ref{rightmu}) coincides for $E_{i}$ and
$\overline{E_{i}}$.
\end{rem}

\begin{dfn}\cite[Definition 1.7]{GS1}
\label{slopestabilitytensors}
We say that $(E,\varphi,u)$ is \emph{slope-$\tau$-semistable} if $E$ is
torsion free, and for all weighted filtrations
$(E_{\bullet},n_{\bullet})$ of $E$, it is
\begin{equation}
\label{taustab} \big( \sum_{i=1}^{t}n_{i}(r\deg E_{i}-r_{i}\deg
E)\big) + \tau \mu(\varphi,E_{\bullet},n_{\bullet})\leq 0
\end{equation}
We say that $(E,\varphi,u)$ is \emph{slope-$\tau$-stable} if we
have a strict inequality in (\ref{taustab}) for every weighted
filtration. If $(E,\varphi,u)$ is not slope-$\tau$-semistable we
say that it is slope-$\tau$-unstable.
\end{dfn}

Recall the relation
$$\tau=\delta_{j}(n-j)!$$
between the parameter $\tau$ and the leading coefficient of polynomial
$\delta$ and note that we have the following implications
$$\text{slope}-\tau-\text{stable}\Rightarrow
\delta-\text{stable}\Rightarrow
\delta-\text{semistable}\Rightarrow
\text{slope}-\tau-\text{semistable}$$ Note that, if the dimension
of the variety $X$ is $n=1$, Definitions \ref{stabilityfortensors}
and \ref{slopestabilitytensors} do coincide.

Let $\mathcal{I}=\{1,...,t+1\}^{\times s}$ be the set of all
multi-indexes $I=(i_{1},...,i_{s})$. Let us call $\nu^{i}(I)$ the
number of times that $i$ appears on the multi-index $I$ and
$\nu_{i}(I)$ the number of elements $k$ in $I$ with $k\leq i$.
Note that $\nu^{i+1}(I)=\nu_{i+1}(I)-\nu_{i}(I)$. Given a
multi-index $I\in \mathcal{I}$, we have
$$\gamma_{r_{i_1}}+\cdots+\gamma_{r_{i_s}}=\sum_{i=1}^{t}\gamma_{i}\nu^{i}(I)=
\sum_{i=1}^{t}\gamma_{r_{i}}(\nu_{i+1}(I)-\nu_{i}(I))$$
$$=s\gamma_{r_{t+1}}-\sum_{i=1}^{t}(\gamma_{r_{i+1}}-\gamma_{r_{i}})\nu_{i}(I)=
s\gamma_{r}-\sum_{i=1}^{t}n_{i}r\nu_{i}(I)$$
$$=s(\sum_{i=1}^{t}n_{i}r_{i})-\sum_{i=1}^{t}n_{i}r\nu_{i}(I)=\sum_{i=1}^{t}n_{i}(sr_{i}-\nu_{i}(I)r)\; .$$
Now let $I_{0}$ be the multi-index giving minimum in
(\ref{rightmu}). We will denote by
$\epsilon_i(\varphi,E_\bullet,n_{\bullet})$ (or just
$\epsilon_i(E_{\bullet})$ if the rest of the data is clear from
the context) the number of elements $k$ of the multi-index $I_{0}$
such that $r_{k}\leq r_{i}$. Let us call
$\epsilon^i(E_{\bullet})=\epsilon_{i+1}(E_{\bullet})-\epsilon_{i}(E_{\bullet})$.
Therefore, we can rewrite (\ref{rightmu}) as
\begin{equation}
\label{rightmu2}
\mu(\varphi,E_\bullet,n_{\bullet})=\sum_{i=1}^{t}n_{i}(sr_{i}-\epsilon_{i}(E_{\bullet})r)\;
.
\end{equation}
In the following we will consider the stability and
slope-stability conditions, (\ref{stabtensors}) and (\ref{taustab}),
with the calculation made in (\ref{rightmu2}).

\begin{rem}
\label{semistableistorsionfree} Note that, if $(E,\varphi,u)$ is
$\delta$-semistable, then it is torsion free. Indeed, consider the
filtration $0\subsetneq T(E) \subsetneq E$ where $T(E)$ is the
torsion subsheaf, and apply (\ref{stabtensors}). Then we obtain
this inequality of polynomials
$$rP_{T(E)}-\rk(T(E))P_{E}+\delta\mu(0\subsetneq T(E)
\subsetneq E)=rP_{T(E)}+\delta\mu(0\subsetneq T(E) \subsetneq
E)\preceq 0$$ which is a contradiction, because we have that the
leading coefficient of $P_{T(E)}$ is positive and $\mu(0\subsetneq
T(E) \subsetneq E)=0$ by (\ref{rightmu2}). Then $T(E)=0$ and $E$
is torsion free.
\end{rem}

\begin{lem}\cite[Lemma 1.4]{GS1}
\label{boundonmi} There is an integer $A_{1}$ (depending only on
$P$, $s$, $c$, $b$ and $\mathcal{D}$) such that it is enough to
check the stability condition (\ref{stabtensors}) for weighted
filtrations with $n_{i}\leq A_{1}$ for all $i$.
\end{lem}
\begin{pr}
Let $\mathcal{I}=\{1,\ldots,t+1\}^{\times s}$ be the set of
multi-indexes $I=(i_{1},\ldots,i_{s})$ and note that
(\ref{rightmu}) is a piece-wise linear function of $\gamma\in
\mathcal{C}$, where $\mathcal{C}\subset \mathbb{Z}^{r}$ is the
cone defined by $\gamma_{1}\leq\ldots\leq\gamma_{r}$. This is
due that it is defined as the minimum among a finite set of the linear
functions $\gamma_{r_{i_{1}}}+\cdots+\gamma_{r_{i_{s}}}$ for those $I\in
\mathcal{I}$ giving a non-zero restriction of morphism $\varphi$,
i.e. $\varphi|_{(E_{i_1}\otimes\cdots\otimes E_{i_r})^{\oplus
c}}\neq 0$. Decompose
$\mathcal{C}=\bigcup_{I\in\mathcal{I}}\mathcal{C}_{I}$ into a
finite number of subcones
$$\mathcal{C}_{I}:=\{\gamma\in\mathcal{C}:\gamma_{r_{i_{1}}}+\cdots+\gamma_{r_{i_{s}}}\leq
\gamma_{r_{i'_{1}}}+\cdots+\gamma_{r_{i'_{s}}} \text{ for all
}I'\in\mathcal{I}\}\; ,$$ such that (\ref{rightmu}) is linear on
each cone $\mathcal{C}_{I}$. Maybe some subcones $I$ are
irrelevant, meaning that $\varphi$ vanishes on them, then we set
$\mu(\varphi,E_{\bullet},n_{\bullet})|_{I}=0$. Choose one vector
$\gamma\in\mathbb{Z}^{r}$ in each edge of each cone
$\mathcal{C}_{I}$ and multiply all these vectors by $r$, so that
all their coordinates are divisible by $r$, and call this set of
vectors $S$. The vectors in $S$ come from weights $n_{i}>0$,
$i=1,\ldots,t+1$, given by the formula
$\gamma=\sum_{i=1}^{t}n_{i}\gamma^{(r_{i})}$. Hence, to obtain the
finite set $S$ of vectors it is enough to consider a finite set of
values for $n_{i}$, therefore there is a maximum value $A_{1}$.

Finally, we will show that it is enough to check
(\ref{stabtensors}) for the weights associated to the vectors in
$S$. Indeed, since the first term in
(\ref{stabtensors}) is linear on $\mathcal{C}$, then it is also
linear on each $\mathcal{C}_{I}$. Then the expression in the left
side of (\ref{stabtensors}) is linear on each subcone
$\mathcal{C}_{I}$, and hence, it is enough to check its
non-positivity on all the edges of all the cones
$\mathcal{C}_{I}$, then it is enough to check it for weights
associated to vectors in $S$.
\end{pr}

Note that the reason why we have to consider filtrations instead
of just subsheaves is that (\ref{rightmu}) is not linear as a
function of the weights $\{n_{i}\}$. But, nevertheless, we can
compare (\ref{rightmu}) for subfiltrations of a given filtration
with the following Lemma. It will be used in the proof of Theorem
2.5.

\begin{lem}\cite[Lemma 1.6]{GS1}
\label{subfiltration} Let $(E_{\bullet},n_{\bullet})$ be a
weighted filtration and let $\mathcal{T}'$ be a subset
of $\mathcal{T}=\{1,...,t\}$. Let $(E'_{\bullet},n'_{\bullet})$ be
the subfiltration obtained by considering only those terms $E_{i}$
for which $i\in \mathcal{T}'$. Then
$$\mu(\varphi,E_{\bullet},n_{\bullet})\leq
\mu(\varphi,E'_{\bullet},n'_{\bullet})+\sum_{i\in\mathcal{T}-\mathcal{T}'}n_{i}sr_{i}\;
.$$
\end{lem}
\begin{pr}
We index the filtration $(E'_{\bullet},n'_{\bullet})$ with
$\mathcal{T}'$. Let
$I'=(i'_{1},...,i'_{s})\in\{\mathcal{T}'\cup\{t+1\}\}^{\times s}$
be the multi-index giving minimum for the filtration
$(E'_{\bullet},n'_{\bullet})$. In particular, we have
$\varphi|_{(E_{i'_{1}}\otimes \cdots \otimes E_{i'_{s}})^{\oplus
c}}\neq 0$. Then
$$\mu(\varphi,E_{\bullet},n_{\bullet})=\min_{I\in
\mathcal{I}} \{\gamma_{r_{i_1}}+\cdots+\gamma_{r_{i_s}}:
\,\varphi|_{(E_{i_1}\otimes\cdots\otimes E_{i_r})^{\oplus c}}\neq
0\}\leq$$
$$\gamma_{r_{i'_1}}+\cdots+\gamma_{r_{i'_s}}=\sum_{i=1}^{t}n_{i}(sr_{i}-\nu_{i}(I')r)=\sum_{i=1}^{t}n_{i}(sr_{i}-\epsilon_{i}(E'_{\bullet})r)=$$
$$\sum_{i\in
\mathcal{T}'}n_{i}(sr_{i}-\epsilon_{i}(E_{\bullet})r)+\sum_{i\in
\mathcal{T}-\mathcal{T}'}n_{i}(sr_{i}-\epsilon_{i}(E'_{\bullet})r)
\leq \mu(\varphi,E'_{\bullet},n'_{\bullet})+\sum_{i\in
\mathcal{T}-\mathcal{T}'}n_{i}sr_{i}\; .$$\end{pr}

A \emph{family of coherent sheaves parametrized by a scheme $T$} is a coherent
sheaf $E_{T}$ on $X\times T$ which is flat over $T$, such that,
$E_{t}:=E_{T}|_{X\times \{t\}}$ is a coherent sheaf over $X$ for
every point $t\in T$. Let us define the ingredients of our moduli problem.

A \emph{family of $\delta$-semistable tensors parametrized by a
scheme $T$} is a tuple $(E_{T},\varphi_{T},u_{T},N)$, consisting
of a torsion free sheaf $E_{T}$ on $X\times T$, flat over $T$,
that restricts to a torsion free sheaf with Hilbert polynomial $P$
on every slice $X\times\{t\}$, a morphism $u_{T}:T\longrightarrow
R$, a line bundle $N$ on $T$ and a homomorphism $\varphi_{T}$,
\begin{equation}
\label{tensoronfamily} \varphi_{T}:(E_{T}^{\otimes s})^{\oplus
c}\longrightarrow (\det E_{T})^{\otimes b}\otimes
\overline{u_{T}}^{\ast}\mathcal{D}\otimes \pi_{T}^{\ast}N\; ,
\end{equation}
(where we define $\overline{u_{T}}=\id_{X}\times u_{T}$) such that
if we consider the restriction of this homomorphism on every slice
$X\times \{t\}$,
$$ \varphi_{t}:(E_{t}^{\otimes s})^{\oplus
c}\longrightarrow (\det E_{t})^{\otimes b}\otimes D_{u_{T}(t)}\;
,$$ the triple $(E_{t},\varphi_{t},u_{T}(t))$ is a
$\delta$-semistable tensor for every $t$. Particularly,
$\varphi_{t}$ is not identically zero. Two families
$(E_{T},\varphi_{T},u_{T},N)$ and
$(E'_{T},\varphi'_{T},u'_{T},N')$ parametrized by $T$ are
\emph{isomorphic} if $u_{T}=u_{T'}$ and there are isomorphisms
$f:E_{T}\longrightarrow E'_{T}$, $\alpha:N\longrightarrow N'$,
such that the induced diagram
\begin{equation}
\label{equivtensorfamily}
\xymatrix{(E_{T}^{\otimes s})^{\oplus
c}\ar[d]^{\varphi_{T}}\ar[r]^{(f^{\otimes s})^{\oplus c}} &
(E_{T}^{'\otimes s})^{\oplus
c}\ar[d]^{\varphi'_{T}}\\
(\det E_{T})^{\otimes
b}\otimes\overline{u_{T}}^{\ast}\mathcal{D}\otimes
\pi_{T}^{\ast}N\ar[r]^{\hat{f}\otimes \pi_{T}^{\ast}\alpha} &
(\det E'_{T})^{\otimes
b}\otimes\overline{u'_{T}}^{\ast}\mathcal{D}\otimes
\pi_{T}^{\ast}N'}
\end{equation} commutes, where $\pi_{T}:X\times T\rightarrow T$ is the
natural projection.

Therefore, we have a category of objects, o notion of stability, a
notion of isomorphism between objects and a notion of family and
equivalence of families. We are ready to define the \emph{functor}
for our \emph{moduli problem}.

Let $\mathcal{M}_{\delta}$ (respectively
$\mathcal{M}_{\delta}^{s}$) be the contravariant functor from the
category of schemes over $\mathbb{C}$, locally of finite type,
($\Sch/\mathbb{C}$) to the category of sets ($\Sets$) which sends
a scheme $T$ to the  set of isomorphism classes of families of
$\delta$-semistable (respectively stable) tensors parametrized by
$T$, and send a morphism $f:T'\longrightarrow T$ to the morphism
of sheaves $\tilde{f}:E'_{T'}\longrightarrow E_{T}$ given by the
pullback diagram
$$\xymatrix{
E'_{T'}\ar[r]^{\tilde{f}}\ar[d] & E_{T}\ar[d]\\
T'\times X\ar[r]^{f\times \id} & T\times X}$$ and similarly for
the morphism $\Phi:\varphi'_{T'}\longrightarrow \varphi_{T}$.

We will construct schemes $\mathfrak{M}_{\delta}$,
$\mathfrak{M}_{\delta}^{s}$ corepresenting the functors
$\mathcal{M}_{\delta}$ and $\mathcal{M}_{\delta}^{s}$ (c.f.
Definition \ref{corepresents}). In general $\mathfrak{M}_{\delta}$
will not be a coarse moduli space, because non-isomorphic tensors
can correspond to the same point in $\mathfrak{M}_{\delta}$. Then,
we will declare two such tensors $S$-equivalent, and
$\mathfrak{M}_{\delta}$ will become a coarse moduli space for the
functor of $S$-equivalence classes of tensors (c.f. Remark
\ref{Seq}). This is the main theorem (c.f. Theorem \cite[Theorem
1.8]{GS1}):

\begin{thm}
\label{maintheorem} Fix $P$, $s$, $c$, $b$ and a family
$\mathcal{D}$ of locally free sheaves on $X$ parametrized by a
scheme $R$. Let $d$ be the degree of a coherent sheaf whose
Hilbert polynomial is $P$. Let $\delta$ be a polynomial as in
(\ref{delta}).

There exists a coarse moduli space
$\mathfrak{M}_{\delta}$, projective over $\Pic^{d}(X)\times R$, of
$S$-equivalence classes of $\delta$-semistable tensors. There is
an open set $\mathfrak{M}_{\delta}^{s}$ corresponding to
$\delta$-stable tensors. Points in this open set correspond to
isomorphism classes of $\delta$-stable tensors.
\end{thm}

In Proposition \ref{S-equivalence} we will give a criterion to decide
when two tensors are $S$-equivalent. We will prove Theorem \ref{maintheorem}
in subsection \ref{proofthmtensors}.

Therefore, in the language of Section \ref{modulis}, we have our
moduli problem stated where $\mathcal{A}$ is the class of
$\delta$-semistable (resp. $\delta$-stable) tensors, the
equivalence relation $\sim$ is given by the notion of
$S$-equivalence (c.f. Remark \ref{Seq}) for which we will give a
criterion in Proposition \ref{S-equivalence} (resp. isomorphism of
tensors in Definition \ref{deftensor}), and the notion of
equivalence of families given by (\ref{equivtensorfamily}). See
\cite[Remark 1.9]{GS1} and \cite[p. 60]{Si1} for a comment on the
notion of equivalence of families giving, as a result, moduli
functors which are not sheaves.

\begin{rem}
\label{modulinotfine}
The obtention of a fine moduli space also requires the existence
of a universal family $E_{\mathfrak{M}_{\delta}}$ over
$\mathfrak{M}_{\delta}$ such that for every family $E_{T}$ of
tensors over $T$, there is a unique morphism $E_{T}\longrightarrow
\mathfrak{M}_{\delta}$ induced by pulling back the universal
family (c.f. Definition \ref{represents}). As we will see, this cannot be done for the case of
tensors as in the case of sheaves.
\end{rem}

\subsection{Results on boundedness}
\label{resultsonboundedness}
In this section we reformulate the stability for tensors using some
boundedness results to prove Theorem \ref{conditions_stability}. First we recall definitions and well known results by Simpson, Grothendieck and
Maruyama.

\begin{dfn}
A set $\mathcal{E}=\{E_{i}\}_{i\in I}$ of coherent sheaves is
\emph{bounded} if there exists a family $E_{T}\longrightarrow
X\times T$ parametrized by $T$, a scheme of finite type over
$\mathbb{C}$, such that for every $i\in I$ there exists $t\in T$
with $E_{i}\simeq E_{t}$.
\end{dfn}

Recall that a scheme $T$ is \emph{of finite type over $\mathbb{C}$} if $T$ can be covered by a finite number of open
affine subsets $\Spec A_{i}$, where each $A_{i}$ is a finitely generated $\mathbb{C}$-algebra.

\begin{dfn}
\label{regdef} A sheaf $E$ is called \emph{$m$-regular} if $h^{i}E(m-i))=0$ for $i>0$.
\end{dfn}

\begin{lem}
\label{mregularity}
If $E$ is $m$-regular then the following holds
\begin{enumerate}
\item $E$ is $m'$-regular for $m'>m$.
\item $E(m)$ is globally generated.
\item For all $m'\geq 0$ the following homomorphisms are surjective
$$
H^0(E(m))\otimes H^0(\SO_X(m'))\too H^0(E(m+m')) \; .
$$
\end{enumerate}
\end{lem}

\begin{prop}
The following properties for a family of sheaves
$\mathcal{E}=\{E_{i}\}_{i\in I}$ are equivalent:
\begin{enumerate}
    \item $\mathcal{E}$ is bounded.
    \item  The set of Hilbert polynomials
$\{P_{E}\}_{E\in\mathcal{E}}$ is finite and there exists a uniform
bound $m_{0}\in \mathbb{Z}$ such that all $E\in \mathcal{E}$ is
$m_{0}$-regular.
    \item The set of Hilbert polynomials
$\{P_{E}\}_{E\in\mathcal{E}}$ is finite and there is a coherent
sheaf $F$ such that all $E_{i}\in\mathcal{E}$ admit surjective homomorphisms
$F\longrightarrow E_{i}$.
\end{enumerate}
\end{prop}

Note that, given a bounded set $\mathcal{E}$ of coherent sheaves,
the set of Hilbert polynomials $\{P_{E}\}_{E\in\mathcal{E}}$ is
finite and hence, $\{\rk E\}_{E\in\mathcal{E}}$ and $\{\deg
E\}_{E\in\mathcal{E}}$, are bounded as sets of numbers.

We denote
\begin{equation}
\label{HpolO}
P_{\mathcal{O}_{X}}(m)=\frac{\alpha_{n}}{n!}m^{n}+\frac{\alpha_{n-1}}{(n-1)!}m^{n-1}+...
+\frac{\alpha_{1}}{1!}m+\frac{\alpha_{0}}{0!}\; ,
\end{equation}
the \emph{Hilbert polynomial} of $\mathcal{O}_{X}$, where
$\alpha_{n}=g=\deg \mathcal{O}_{X}(1)$, and
\begin{equation}
\label{HpolE}
P_{E}(m)=\frac{rg}{n!}m^{n}+\frac{d+r\alpha_{n-1}}{(n-1)!}m^{n-1}+...\;
,
\end{equation}
the \emph{Hilbert polynomial} of a sheaf $E$ with rank $r$ and degree $d$. For the following Lemma, see
\cite[Lemma 1.5 and Corollary 1.7]{Si1} and \cite[Lemme 2.4]{LeP}.
\begin{lem}\cite[Lemma 2.2]{HL2}
\label{Simpson}
Let $r>0$ be an integer. Then there exists a
constant $B$ with the following property: for every torsion free
sheaf $E$ with $0<\rk (E)\leq r$, we have
$$h^{0}(E)\leq \frac{1}{g^{n-1}n!}\big ((\rk(E)-1)([\mu_{max}(E)+B]_{+})^{n}+([\mu_{mim}(E)+B]_{+})^{n}\big)\; ,$$
where $g=\deg \mathcal{O}_{X}(1)$, $[x]_{+}=\max\{0,x\}$, and
$\mu_{max}(E)$ (resp. $\mu_{min}(E)$) is the maximum (resp.
minimum) slope of the semistable factors of the
Mumford-Harder-Narasimhan filtration of $E$ (c.f. Definition
\ref{HNdef} and Remark \ref{GHN-MHN}).
\end{lem}

\begin{lem}\cite[Lemma 2.5]{Gr}
\label{Grothendieck} Let $\mathcal{E}$ be a bounded set of sheaves
$E$ and fix a constant $C$. The set of torsion free quotients
$E\twoheadrightarrow E''$ of the sheaves $E\in \mathcal{E}$, with
$|\deg(E'')|\leq C$, is bounded.
\end{lem}

\begin{thm}\cite{Ma3}
\label{Maruyama} Fix a constant $C$. The family of sheaves $E$
with fixed Hilbert polynomial $P$ and such that $\mu_{\max}(E)\leq
C$, is bounded.
\end{thm}

Recall that we denote by $\mu(E)=\frac{\deg E}{\rk E}$ the slope of a
sheaf. As a consequence of Maruyama's result in Theorem \ref{Maruyama}, we can prove the
boundedness of the set of $\delta$-semistable tensors:

\begin{cor}
The set of $\delta$-semistable tensors $(E,\varphi,u)$ with fixed
Hilbert polynomial $P$ is bounded.
\end{cor}
\begin{pr}
Let $(E,\varphi,u)$ a $\delta$-semistable tensor. Then, we have
seen that $(E,\varphi,u)$ is $\tau$-slope-semistable, where
$\tau=\delta_{j}(n-j)!$. Then, for every weighted filtration and,
in particular, for every one-step filtration $0\subset E'\subset
E$, the expression (\ref{taustab}) holds, hence
$$(r\deg(E')-\rk(E')\deg E)\big)+\tau(s\rk(E')-\epsilon(E'\subset
E)r)\leq 0\; .$$ Then, dividing by $r\cdot \rk(E')$ we obtain this
condition for the slopes $$\mu(E')\leq
\mu(E)\tau(\frac{s}{r}+\frac{\epsilon(E'\subset E)}{\rk{E'}})\leq
\mu(E)-\tau(\frac{s}{r}+s)=C\; ,$$ where $C$ is a constant
depending only on $P$, $\tau$ and $s$, which are fixed. Hence, we
apply Theorem \ref{Maruyama}, provided $\mu(E')\leq
\mu_{\max}(E)$, for every subsheaf $E'\subset E$.
\end{pr}

This is the main theorem of this section, whose proof we will give
after some preliminary results.

\begin{thm}(c.f. \cite[Theorem 2.5]{GS1})
\label{conditions_stability} There is an integer $N_{0}$ such that
if $m\geq N_{0}$, the following properties of tensors
$(E,\varphi,u)$, with $E$ torsion free and $P_{E}=P$, are
equivalent.
\begin{enumerate}
    \item $(E,\varphi,u)$ is semistable (resp. stable).
    \item $P(m)\leq h^{0}(E(m))$ and for every weighted filtration
    $(E_{\bullet},n_{\bullet})$ as in (\ref{filtE}),
    $$(\sum_{i=1}^{t}n_{i}\big(rh^{0}(E_{i}(m))-r_{i}P(m))\big)+\delta(m)\mu(\varphi,E_{\bullet},n_{\bullet})\leq 0$$
(resp. $<$).
    \item For every weighted filtration
    $(E_{\bullet},n_{\bullet})$ as in (\ref{filtE}),
    $$\big(\sum_{i=1}^{t}n_{i}(r^{i}P(m)-rh^{0}(E^{i}(m)))\big)+\delta(m)\mu(\varphi,E_{\bullet},n_{\bullet})\leq 0$$
(resp. $<$).
\end{enumerate}
Furthermore, for any tensor $(E,\varphi,u)$ satisfying these
conditions, $E$ is $m$-regular.
\end{thm}

The set of tensors $(E,\varphi,u)$, with $E$ torsion free and $P_{E}=P$, satisfying the weak version of
conditions $1-3$ will be called $\mathcal{S}^{s}$, $\mathcal{S}'_{m}$ and $\mathcal{S}''_{m}$, respectively.

\begin{lem}
\label{bound} There is an integer $N_{1}$ and a positive constant $D$, such that if $(E,\varphi,u)$ belongs to
$\mathcal{S}=\mathcal{S}^{s}\cup \bigcup_{m\geq N_{1}} \mathcal{S}''_{m}$, then for all saturated weighted
filtrations $(E_{\bullet},n_{\bullet})$, the following holds for all $i$:
\begin{equation}
\label{leqslope} \mu(E_{i})\leq \mu(E)+D
\end{equation}
and, either $\mu(E)-D\leq\mu(E_{i})$, or
\begin{enumerate}
    \item $rh^{0}(E_{i}(m))<r_{i}(P(m)-s\delta(m))$, if
    $(E,\varphi,u)\in \mathcal{S}^{s}$ and $m\geq N_{1}$
    \item $rP_{E_{i}}-r_{i}P-r_{i}s\delta\prec 0$, if
    $(E,\varphi,u)\in \bigcup_{m\geq N_{1}}\mathcal{S}''_{m}$
\end{enumerate}
\end{lem}

\begin{pr}
Let $(E_{\bullet},n_{\bullet})$ be a weighted filtration of $E$
and let $B$ be as in Lemma \ref{Simpson}. Let $G$ be the following
polynomial,
$$G(m)=\frac{1}{g^{n-1}n!}\big ((r-1)(\mu(E)+\tau s(1-\frac{1}{r})+gm+B)^{n}+(\mu(E)-D+gm+B)^{n}\big)$$
$$=\frac{1}{g^{n-1}n!}\big
[rg^{n}m^{n}+ng^{n-1}(r\mu(E)+\tau
s\frac{(r-1)^{2}}{r}-D+rB)m^{n-1}+\cdots \big ]\; .$$ Then, the
leading coefficient of $G-(P-s\delta)$ (i.e. the coefficient of
$m^{n}$) is
$$r\frac{g}{n!}-r\frac{g}{n!}=0\; ,$$
but the coefficient of $m^{n-1}$ is
$$\big[\frac{1}{(n-1)!}(r\mu(E)+\tau
s\frac{(r-1)^{2}}{r}-D+rB)\big]-\big[\frac{1}{(n-1)!}(d+r\alpha_{n-1})+\tau
s\big]\; ,$$ so we can choose $D$ large enough so that the leading
coefficient of $G-(P-s\delta)$ is negative. We choose $D$ also to
verify $D>\tau s$.

Let $N_{1}$ be large enough so that, for $m\geq N_{1}$, the
following three expressions hold:
\begin{equation}
\label{deltapos} \delta(m)\geq 0
\end{equation}
\begin{equation}
\label{bracketpos}
\mu(E)-D+gm+B>0
\end{equation}
\begin{equation}
\label{polneg} G(m)-(P(m)-s\delta(m))<0\; .
\end{equation}

Given that the filtration is supposed to be saturated, and $E$ to be torsion free, we have $0<r_{i}<r$.

\textbf{Case 1}. Suppose that $(E,\varphi,u)\in \mathcal{S}^{s}$.
Then, $(E,\varphi,u)$ is $\tau$-slope-semistable hence, for each
$i$, consider the one-step filtration $E_{i}\subsetneq E$ and
apply (\ref{taustab}),
$$r\deg E_{i}-r_{i}\deg E +\tau(sr_{i}-\epsilon_{i}(E_{i}\subsetneq
E)r)\leq 0\; .$$ Dividing by $r_{i}\cdot r$ we get $$\mu(E_{i})\leq
\mu(E)-\tau(\frac{s}{r}+\frac{\epsilon_{i}(E_{i}\subsetneq E)}{r_{i}})\leq \mu(E)+\tau
s(1-\frac{1}{r})<\mu(E)+D$$ using $D>\tau s$, hence (\ref{leqslope}).

Let $E_{i,\max}\subset E_{i}$ be the term in the Harder-Narasimhan
filtration of $E_{i}$ with maximal slope (c.f. Definition
\ref{HNdef} and Remark \ref{GHN-MHN}). Then, the same argument as
before, applied to the filtration $E_{i,\max}\subsetneq E$, gives
\begin{equation}
\label{Eimax} \mu_{\max}(E_{i})=\mu(E_{i,\max})\leq \mu(E)+\tau s(1-\frac{1}{r}).
\end{equation}

Suppose that the first alternative does not hold, i.e.
$$\mu(E_{i})<\mu(E)-D\; .$$
Then, by Lemma \ref{Simpson},
\begin{equation}
\label{Simpsonaplic} h^{0}(E_{i}(m))\leq \frac{1}{g^{n-1}n!}\big
((r_{i}-1)([\mu_{max}(E_{i})+gm+B]_{+})^{n}+([\mu_{min}(E_{i})+gm+B]_{+})^{n}\big)\; ,
\end{equation} where note
that $\mu_{max}(E_{i}(m))=\mu_{max}(E_{i})+gm$ and
$\mu_{min}(E_{i}(m))=\mu_{min}(E_{i})+gm$. Combining the
hypothesis, the expression
$$\mu_{min}(E_{i})\leq \mu(E_{i})<\mu(E)-D\; ,$$
and the expressions (\ref{Eimax}) and (\ref{bracketpos}), the formula (\ref{Simpsonaplic}) becomes
$$h^{0}(E_{i}(m))\leq \frac{1}{g^{n-1}n!}\big ((r_{i}-1)(\mu(E)+\tau s(1-\frac{1}{r})+gm+B)^{n}+(\mu(E)-D+gm+B)^{n}\big)$$
$$\leq\frac{r_{i}}{rg^{n-1}n!}\big( (r-1)(\mu(E)+\tau
s(1-\frac{1}{r})+gm+B)^{n}+(\mu(E)-D+gm+B)^{n}\big)=\frac{r_{i}}{r}G(m)\;
.$$ Now, by using (\ref{polneg}), it is
$$rh^{0}(E_{i}(m))<r_{i}G(m)<r_{i}(P(m)-s\delta(m))\; ,$$
hence we obtain condition $1$.

\textbf{Case 2}. Suppose that $(E,\varphi,u)\in \mathcal{S}''_{m}$ for some $m\geq N_{1}$. For each $i$, consider the
quotient $E^{i}=E/E_{i}$. Let $E^{i}_{\min}$ be the last factor of the Mumford-Harder-Narasimhan filtration of
$E^{i}$ (c.f. Definition \ref{HNdef} and Remark \ref{GHN-MHN}), i.e. $\mu(E^{i}_{\min})=\mu_{\min}(E^{i})$. Let
$E'$ be the kernel
$$0\longrightarrow E'\longrightarrow E \longrightarrow
E^{i}_{\min}\longrightarrow 0\; ,$$ and consider (\ref{polneg}) in the form
$$\frac{G(m)}{r}-\frac{P(m)-s\delta(m)}{r}<0\; .$$
Using (\ref{deltapos}), we apply condition $3$ in Theorem
\ref{conditions_stability} to the one-step filtration
$E'\subsetneq E$ (note that $(E,\varphi,u)\in \mathcal{S}''_{m}$)
to get
$$\frac{G(m)}{r}<\frac{h^{0}(E^{i}_{\min}(m))}{\rk
E^{i}_{\min}}-\delta(m)\frac{(s\rk E'-\epsilon(E'\subsetneq E)r)}{r\rk E^{i}_{\min}}-\delta(m)\frac{s}{r}$$
$$\leq \frac{h^{0}(E^{i}_{\min}(m))}{\rk
E^{i}_{\min}}+\delta(m)s(-\frac{1}{r}-\frac{1}{\rk
E^{i}_{\min}}(\frac{\rk
E'}{r}-1))=\frac{h^{0}(E^{i}_{\min}(m))}{\rk E^{i}_{\min}}\; .$$

Applying Lemma \ref{Simpson},
$$\frac{G(m)}{r}<\frac{1}{g^{n-1}n!\rk E^{i}_{\min}}\big ((\rk
E^{i}_{\min}-1)[\mu_{max}(E_{\min}^{i})+gm+B]_{+}^{n}+[\mu_{\min}(E^{i}_{\min})+gm+B]_{+}^{n}\big)\;
.$$ By definition, $E^{i}_{\min}$ is semistable, then
\begin{equation}
\label{mu_max_min}
\mu(E^{i}_{\min})=\mu_{max}(E_{\min}^{i})=\mu_{min}(E_{\min}^{i}),
\end{equation}
hence
$$\frac{G(m)}{r}<\frac{1}{g^{n-1}n!}\big ([\mu(E_{\min}^{i})+gm+B]_{+}^{n}\big )\; .$$ By
definition of $G$ and (\ref{bracketpos}), we have $0<G(m)$, hence
$\mu(E_{\min}^{i})+gm+B>0$ and, then,
$$\frac{G(m)}{r}<\frac{1}{g^{n-1}n!}\big (\mu(E^{i}_{\min})+gm+B)^{n}\big)\; .$$
This inequality of polynomials holds for some $m\geq N_{1}$, therefore, it holds for larger values of $m$, and
hence, we will have this inequality between the second coefficients (note that leading coefficients are equal),
$$\frac{1}{(n-1)!}(\mu(E)+\tau s \frac{(r-1)^{2}}{r^{2}}-\frac{D}{r}+B)<\frac{1}{(n-1)!}(\mu(E^{i}_{\min})+B)\; .$$
Now, $\mu(E^{i}_{\min})\leq \mu(E^{i})=\frac{d-\deg E_{i}}{r-\rk
E_{i}}$ and, using $\frac{r-\rk E_{i}}{\rk E_{i}}<r$ and $D>\tau
s$, previous inequality gives (\ref{leqslope}).

Now, suppose that the first alternative does not hold, i.e.
$$\mu(E_{i})<\mu(E)-D\; ,$$
then
$$\frac{\deg E_{i}}{\rk E_{i}}<\frac{d}{r}-D<\frac{d}{r}-\tau s<\frac{d-\tau s}{r}\; ,$$
which is equivalent to
$$r\deg E_{i}-\rk E_{i}r+\tau s\rk E_{i}<0\; ,$$
and hence, the leading coefficient of the polynomial
$\;rP_{E_{i}}-r_{i}P-r_{i}s\delta\;$ is negative, therefore,
condition $2$ holds.
\end{pr}

\begin{lem}
\label{Sbounded}
The set $\mathcal{S}=\mathcal{S}^{s}\cup
\bigcup_{m\geq N_{1}} \mathcal{S}''_{m}$ is bounded.
\end{lem}
\begin{pr}
Let $(E,\varphi,u)\in\mathcal{S}$. Let $E'$ be a subsheaf of $E$,
and $\overline{E'}$ the saturation of $E'$ on $E$. Then we have
this exact sequence
$$0\longrightarrow E'\longrightarrow \overline{E'}\longrightarrow
T(E')\longrightarrow 0$$ and, by additivity of the Hilbert
polynomial on exact sequences, and $\rk E'=\rk \overline{E'}$, we
get $\deg (E')\leq \deg (\overline{E'})$. Then, Lemma \ref{bound}
gives (\ref{leqslope}), then
$$\mu(E')\leq \mu(\overline{E'})\leq \mu(E)+D$$
and, therefore, by Theorem \ref{Maruyama}, the set $\mathcal{S}$
is bounded.
\end{pr}

\begin{lem}
\label{S0bounded} Let $\mathcal{S}_{0}$ be the set of sheaves $E'$
such that $E'$ is a saturated subsheaf of $E$ for some
$(E,\varphi,u)\in\mathcal{S}$, and furthermore
\begin{equation}
\label{leqgeqslope} |\mu(E')-\mu(E)|\leq D.
\end{equation}
Then, $\mathcal{S}_{0}$ is bounded.
\end{lem}

\begin{pr}
Let $E'\in\mathcal{S}_{0}$. The sheaf $E''=E/E'$ is torsion free and
$$|\deg (E'')|=|\deg (E)-\deg(E')|\leq |\deg(E)|+|\deg (E')|\leq
2|\deg(E)|+rD\; ,$$ where in last inequality we use
(\ref{leqgeqslope}). Then, as $\deg E$ is fixed, by Lemma
\ref{Grothendieck}, the set of sheaves $E''$ obtained in this way
is bounded, and hence, also $\mathcal{S}_{0}$ is bounded.
\end{pr}

\begin{lem}
\label{asintotic} There exists an integer $N_{2}$ such that for
every weighted filtration $(E_{\bullet},n_{\bullet})$ as in
(\ref{filtE}), with $E_{i}\in \mathcal{S}_{0},\;\forall i$, the
inequality of polynomials (\ref{stabtensors}) in Definition
\ref{stabilityfortensors},
$$\big(\sum_{i=1}^{t}n_{i}(rP_{E_{i}}-r_{i}P)\big)+\delta\mu(\varphi,E_{\bullet},n_{\bullet})\preceq 0\; ,$$
holds if and only if it holds for a particular value of $m\geq N_{2}$. The same holds for $\prec$.
\end{lem}
\begin{pr}
Since $\mathcal{S}_{0}$ is bounded by Lemma \ref{Sbounded}, the
set of polynomials $P_{E'}$, for $E'\in\mathcal{S}_{0}$, is
finite. Lemma \ref{boundonmi} implies that we only need to
consider a finite number of values for $n_{i}$, hence the result
follows from it.
\end{pr}

\begin{pr}[\textbf{Proof of Theorem \ref{conditions_stability}}]
Given that $\mathcal{S}$ and $\mathcal{S}_{0}$ are bounded, let $N_{0}>\max\{N_{1},N_{2}\}$ (c.f. Lemmas
\ref{bound} and \ref{asintotic}) and such that all sheaves in $\mathcal{S}$ and $\mathcal{S}_{0}$ are
$N_{0}$-regular, and $E_{1}\otimes \cdots\otimes E_{s}$ is $sN_{0}$-regular for all $E_{1},...,E_{s}$ in
$\mathcal{S}_{0}$. Let $m\geq N_{0}$.

$(2.\Rightarrow 3.)$ Let $(E,\varphi,u)\in \mathcal{S}'_{m}$ and consider a weighted filtration
$(E_{\bullet},n_{\bullet})$ as in (\ref{filtE}). Note that the functor of global sections is only left exact,
then applying it to the exact sequence
$$0\longrightarrow E_{i}(m)\longrightarrow E(m)\longrightarrow
E^{i}(m)\longrightarrow 0\; ,$$ we obtain
$$0\longrightarrow H^{0}(E_{i}(m))\longrightarrow
H^{0}(E(m))\longrightarrow H^{0}(E^{i}(m))\; ,$$ and we have this
inequality for the dimensions of the vector spaces
\begin{equation}
\label{inequalityh0} h^{0}(E(m))\leq
h^{0}(E_{i}(m))+h^{0}(E^{i}(m))\; .
\end{equation}
Then, using hypothesis $P(m)\leq h^{0}(E(m))$ and (\ref{inequalityh0}), we get
$$\big(\sum_{i=1}^{t}n_{i}(r^{i}P(m)-rh^{0}(E^{i}(m))\big)+\delta(m)\mu(\varphi,E_{\bullet},n_{\bullet})=$$
$$\big(\sum_{i=1}^{t}n_{i}(r(P(m)-h^{0}(E^{i}(m)))-r_{i}P(m))\big)+\delta(m)\mu(\varphi,E_{\bullet},n_{\bullet})\leq$$
$$\big(\sum_{i=1}^{t}n_{i}(r(h^{0}(E(m))-h^{0}(E^{i}(m)))-r_{i}P(m))\big)+\delta(m)\mu(\varphi,E_{\bullet},n_{\bullet})\leq$$
$$\big(\sum_{i=1}^{t}n_{i}(rh^{0}(E_{i}(m))-r_{i}P(m))\big)+\delta(m)\mu(\varphi,E_{\bullet},n_{\bullet})\leq0\; ,$$
therefore $(E,\varphi,u)\in \mathcal{S}''_{m}$, and similarly for the strict inequality.

$(1.\Rightarrow 2.)$ Let $(E,\varphi,u)\in \mathcal{S}^{s}$ and
consider a saturated weighted filtration
$(E_{\bullet},n_{\bullet})$. Since $E$ is $N_{0}$-regular,
$P(m)=h^{0}(E(m))$. If $E_{i}\in \mathcal{S}_{0}$ (meaning that
(\ref{leqgeqslope}) holds), by choice of $N_{0}$, it is also
$P_{E_{i}}(m)=h^{0}(E_{i}(m))$. If $E_{i}\notin \mathcal{S}_{0}$
then, by definition of $\mathcal{S}_{0}$ (c.f. \ref{leqgeqslope}),
the second alternative of Lemma \ref{bound} holds, hence, as
$(E,\varphi,u)\in \mathcal{S}^{s}$, we have assertion $1.$,
$$rh^{0}(E_{i}(m))<r_{i}(P(m)-s\delta(m))\; .$$
Let $\mathcal{T}'\subset \mathcal{T}=\{1,...,t\}$ be the subset of
those $i$ for which $E_{i}\in\mathcal{S}_{0}$. Let
$(E'_{\bullet},n'_{\bullet})$ be the corresponding subfiltration.
Hence, previous argument and Lemma \ref{subfiltration} shows that
\begin{equation}
\label{1implies2}
\big(\sum_{i=1}^{t}n_{i}(rh^{0}(E_{i}(m))-r_{i}P(m))\big)+\delta(m)\mu(\varphi,E_{\bullet},n_{\bullet})\leq
\end{equation}
$$\big(\sum_{i\in\mathcal{T'}}n'_{i}(rP_{E_{i}}(m)-r_{i}P(m))\big)+\delta(m)\mu(\varphi,E'_{\bullet},n'_{\bullet})+$$
$$\big(\sum_{i\in\mathcal{T}-\mathcal{T}'}n_{i}(rh^{0}(E_{i}(m))-r_{i}P(m))+\delta(m)sr_{i}\big)\leq$$
$$\big(\sum_{i\in\mathcal{T'}}n_{i}(rP_{E_{i}}(m)-r_{i}P(m))\big)+\delta(m)\mu(\varphi,E'_{\bullet},n'_{\bullet})\leq 0\; ,$$
where last inequality follows from $(E,\varphi,u)\in \mathcal{S}^{s}$. The condition that $E_{i}$ is saturated
can be dropped, since $h^{0}(E_{i}(m))\leq h^{0}(\overline{E_{i}}(m))$ and
$\mu(\varphi,E_{\bullet},n_{\bullet})=\mu(\varphi,\overline{E_{\bullet}},n_{\bullet})$, where $\overline{E_{i}}$
is the saturated subsheaf generated by $\overline{E_{i}}$ in $E$. Therefore, $(E,\varphi,u)\in \mathcal{S}'_{m}$
and similarly for the strict inequality.

$(3.\Rightarrow 1.)$ Let $(E,\varphi,u)\in \mathcal{S}''_{m}$ and consider a saturated weighted filtration
$(E_{\bullet},n_{\bullet})$. Since $E$ is $N_{0}$-regular, $P(m)=h^{0}(E(m))$. If $E_{i}\in\mathcal{S}_{0}$,
then also $P_{E_{i}}(m)=h^{0}(E_{i}(m))$. Hence, $h^{1}(E_{i}(m))=0$ and (\ref{inequalityh0}) becomes an
equality. Then, using the additivity of the Hilbert polynomial on exact sequences, we get
$h^{0}(E^{i}(m))=P_{E^{i}}(m)$. Now, hypothesis 3. applied to the subfiltration $(E'_{\bullet},n_{\bullet})$
consisting on those terms such that $E_{i}\in \mathcal{S}_{0}$, implies
$$\big(\sum_{E_{i}\in\mathcal{S}_{0}}n_{i}(r^{i}P(m)-rP_{E^{i}}(m))\big)+\delta(m)
\mu(\varphi,E'_{\bullet},n'_{\bullet})\leq0$$ and, using Lemma \ref{asintotic}, this is equivalent to
$$\big(\sum_{E_{i}\in\mathcal{S}_{0}}n_{i}(rP_{E_{i}}-r_{i}P)\big)+\delta\mu
(\varphi,E'_{\bullet},n'_{\bullet})\preceq 0\; .$$
If $E_{i}\notin
\mathcal{S}_{0}$ then, by definition of $\mathcal{S}_{0}$ (c.f.
(\ref{leqgeqslope})) the second alternative of Lemma \ref{bound}
holds, hence, as $(E,\varphi,u)\in \mathcal{S}''_{m}$, we have
$$rP_{E_{i}}-r_{i}P+r_{i}s\delta\prec 0\; .$$ Therefore, previous arguments together with Lemma
\ref{subfiltration} give
$$\big(\sum_{i=1}^{t}n_{i}(rP_{E_{i}}-r_{i}P)\big)+\delta\mu(\varphi,E_{\bullet},n_{\bullet})\preceq$$
$$\big(\sum_{E_{i}\in \mathcal{S}_{0}}n'_{i}(rP_{E_{i}}-r_{i}P)\big)+\delta\mu(\varphi,E'_{\bullet},n'_{\bullet})+
\big(\sum_{E_{i}\notin \mathcal{S}_{0}}n_{i}(rP_{E_{i}}-r_{i}P+\delta s r_{i})\big)\preceq0\; .$$ We proceed
similarly for the strict inequality. As before, the condition that the filtration is saturated can be dropped, and
this finishes the proof of the Theorem.
\end{pr}

\begin{cor}
\label{equalto0} Let $(E,\varphi,u)$ be $\delta$-semistable,
$m\geq N_{0}$, and suppose that there exists a weighted filtration
$(E_{\bullet},n_{\bullet})$ with
$$\big(\sum_{i=1}^{t}n_{i}(rh^{0}(E_{i}(m))-r_{i}P(m))\big)+\delta(m)\mu(\varphi,E_{\bullet},n_{\bullet})=0\; .$$
Then $E_{i}\in \mathcal{S}_{0}$ and, by choice of $N_{0}$,
$h^{0}(E_{i}(m))=P_{E_{i}}(m)$ for all $i$.
\end{cor}
\begin{pr}
By the proof of the part $(1.\Rightarrow 2.)$ of Theorem
\ref{conditions_stability}, if we have this equality then all
inequalities in (\ref{1implies2}) are equalities, hence
$\mathcal{T}=\mathcal{T}'$ and $E_{i}\in \mathcal{S}_{0}$, for all
$i$.
\end{pr}

Note that in Theorem \ref{conditions_stability} we are assuming
that $E$ is torsion free. To apply the Lemma in the general case,
we will use the following Lemma.

\begin{lem}\cite[Lemma 2.11]{GS1}
\label{torsionfreeator} Fix $u\in R$. Let $(E,\varphi,u)$ be a
tensor. Assume that there exists a family
$(E_{t},\varphi_{t},u)_{t\in C}$ parametrized by a smooth curve
$C$ such that $(E_{0},\varphi_{0},u)=(E,\varphi,u)$ and $E_{t}$ is
torsion free for $t\neq 0$. Then there exists a tensor
$(F,\psi,u)$, a homomorphism
$$(E,\varphi,u)\longrightarrow (F,\psi,u)$$
such that $F$ is torsion free with $P_{E}=P_{F}$, and an exact
sequence
$$0\longrightarrow T(E)\longrightarrow
E\overset{\beta}{\longrightarrow} F\; ,$$ where $T(E)$ is the
torsion subsheaf of $E$.
\end{lem}
\begin{pr}
The family is given by a tuple $(E_{C},\varphi_{C},u_{C},N)$ as in (\ref{tensoronfamily}), where $u_{C}$ is the
constant map from $C$ to $R$ with constant value $u$. Shrinking $C$ if necessary, assume that $N$ is trivial. Let
$U=(X\times C)-\supp(T(E_{0}))$. Let $F_{C}=j_{\ast}(E_{C}|_{U})$. Since $F_{C}$ has no elements of torsion
supported on a fiber of the projection $X\times C\rightarrow C$, then $F_{C}$ is flat over $C$ (c.f.
\cite[Proposition 9.7]{Ha}). The natural map $\tilde{\beta}:E_{C}\rightarrow F_{C}$ is an isomorphism on $U$,
hence we have a homomorphism $\psi_{U}:=\varphi_{C}|_{U}$ on $U$ which extends to a homomorphism $\phi_{C}$
on $X\times C$ because $\overline{u_{C}}^{\ast}\mathcal{D}$ is locally free, where
$\overline{u_{C}}=\id_{X}\times u_{C}$. Finally define $(F,\psi)=(F_{0},\psi_{0})$, and let $\beta$ be the
homomorphism induced by $\tilde{\beta}$.
\end{pr}

\subsection{Gieseker's embedding}
\label{giesekerembedding}

In this subsection we are going to change the embedding used in
\cite{GS1} by the one given by Gieseker for the construction of
the moduli space of sheaves over algebraic surfaces.

By Serre Vanishing Theorem, choose $N\geq N_{0}$ (c.f. proof of
Theorem \ref{conditions_stability}) to be large enough so that for
all $m\geq N$, all $i>0$, all line bundles $L$ of degree $d$ and
all locally free sheaves $\{D_{u}\}_{u\in R}$, we have
$$h^{i}(L^{\otimes b}\otimes D_{u}(sm))=0\; ,$$ and
$L^{\otimes b}\otimes D_{u}(sm)$ is generated by global sections.

Fix $m\geq N$ and let $V$ be a vector space of dimension
$p:=P(m)$. The choice of $m$ implies that if $(E, \varphi, u)$ is
$\delta$-semistable then, by Theorem \ref{conditions_stability},
$E(m)$ is $m$-regular and hence, $h^{i}(E(m))=0\;,\forall i>0,$
and it is generated by global sections. Consider a tuple
$(g,E,\varphi,u)$, where $(E, \varphi, u)$ is a
$\delta$-semistable tensor and
$$g:V\longrightarrow H^{0}(E(m))$$ is an isomorphism. This induces
a quotient, as morphism of sheaves, by the evaluation map
\begin{equation}
\label{evaluation}
    \begin{array}{cccccc}
q: & V & \otimes & \mathcal{O}_{X}(-m) & \twoheadrightarrow & E\\
   & s|_{U} & \otimes & a|_{U}         & \mapsto &   s\cdot a|_{U}
   \end{array}
\end{equation}
where $s$ a section of $E$, given as the image by $g$ of an
element of $V$. Let $\mathcal{H}$ be Quot-scheme of Grothendieck
which is the scheme of such quotients with Hilbert polynomial $P$,
$$\mathcal{H}:=\Quot_{\mathcal{F},X,P}=\{q:\mathcal{F}\twoheadrightarrow
E,P_{E}=P\}\; ,$$ where $\mathcal{F}=V\otimes\mathcal{O}_{X}(-m)$.
Each quotient $q$ induces the following homomorphisms
\begin{equation}
    \begin{array}{ccc}
q(m): V\otimes\mathcal{O}_{X}(-m)\otimes\mathcal{O}_{X}(m) &
\longrightarrow & E\otimes\mathcal{O}_{X}(m)\\
 &  & \\
q(m):V\otimes\mathcal{O}_{X} & \longrightarrow & E(m)
   \end{array}
\end{equation}

\begin{equation}
    \begin{array}{ccc}
H^{0}(q(m)): V\otimes H^{0}(\mathcal{O}_{X}) &
\longrightarrow & H^{0}(E(m))\\
 & & \\
H^{0}(q(m)):V\otimes\mathbb{C} & \longrightarrow & H^{0}(E(m))\\
 & & \\
H^{0}(q(m)):V & \longrightarrow & H^{0}(E(m))
   \end{array}
\end{equation}

\begin{equation}
\wedge^{r}H^{0}(q(m)):\wedge^{r}V\longrightarrow
\wedge^{r}H^{0}(E(m))\longrightarrow H^{0}(\wedge^{r}E(m))\simeq
H^{0}(\det(E)(rm))
\end{equation}
where note that we have to choose an isomorphism between
$\wedge^{r}E$ and $\det E$. We call $Q:=\wedge^{r}H^{0}(q(m))$ and
$A:=H^{0}(\det(E)(rm))$, then we have
$$Q\in \Hom(\wedge^{r}V,A)\; .$$
Given that two of these points $Q$ differing by a scalar correspond to the
same morphism (because the isomorphism $\wedge^{r}E\simeq \det E$
is well defined up to a scalar) we get a well defined point in a
projective space
$$\overline{Q}\in \mathbb{P}(\Hom(\wedge^{r}V,A))\; .$$
Therefore we get a Grothendieck's embedding
$$\mathcal{H}\longrightarrow \mathbb{P}(\Hom(\wedge^{r}V,A))$$
and, hence, a very ample line bundle
$\mathcal{O}_{\mathcal{H}}(1)$ on $\mathcal{H}$ (depending on
$m$).

The tuple $(g,E,\varphi,u)$, where recall
$$\varphi:(E^{\otimes s})^{\oplus c}\longrightarrow (\det E)^{\otimes b}\otimes
D_{u}\; ,$$ also induces this linear map
\begin{equation}
\label{Phi} \Phi: (V^{\otimes s})^{\oplus c}\longrightarrow
H^{0}(E(m)^{\otimes s})^{\oplus c}\longrightarrow H^{0}((\det
E)^{\otimes b}\otimes D_{u}(sm))\; ,
\end{equation}
by tensoring each $E$ with $\mathcal{O}_{X}(m)$ and taking global
sections.

A Poincare bundle on $J\times X$, where $J=Pic^{d}(X)$, is a
universal family such that
\begin{equation}
\label{Poincare} \xymatrix{ \det E\simeq\mathcal{P}|_{\{det
E\}\times X}\ar[r]\ar[d] &
\mathcal{P}\ar[d]\\
\{\det E\}\times X\ar@{^{(}->}[r] & J\times X}
\end{equation}
Then, we fix an isomorphism
$$\beta:\det E\longrightarrow \mathcal{P}|_{\{\det E\}\times X}$$
and hence, $\Phi$ induces a quotient
$$(V^{\otimes s})^{\oplus c}\otimes H^{0}(\mathcal{P}|^{\otimes b}_{\{\det E\}\times
X}\otimes D_{u}(sm))^{\vee}\longrightarrow \mathbb{C}\; .$$ Note
that, if we choose a different isomorphism $\beta'$, this quotient
will only change by a scalar, so we get a well defined point
$[\Phi]$ in $\mathcal{W}$, where $\mathcal{W}$ is the projective
bundle over $J\times R$ defined as
$$\mathcal{W}=\mathbb{P}\big(((V^{\otimes s})^{\oplus
c})^{\vee}\otimes \pi_{J\times R,\ast}(\pi^{\ast}_{X\times
J}\mathcal{P}^{\otimes b}\otimes \pi^{\ast}_{X\times
R}\mathcal{D}(sm))\big )\longrightarrow J\times R\; ,$$ where
$\pi_{X\times J}$ (resp. $\pi_{J\times R}$,...) denotes the
natural projection from $X\times J\times R$ to $X\times J$ (resp.
$J\times R$,...) and we denote $\mathcal{D}(m):=\mathcal{D}\otimes
\pi_{X}^{\ast}\mathcal{O}_{X}(m)$. Note that $\pi_{J\times
R,\ast}$ $(\pi^{\ast}_{X\times J}\mathcal{P}^{\otimes b}\otimes
\pi^{\ast}_{X\times R}\mathcal{D}(sm))$ is locally free because of
the choice of $m$. Replacing $\mathcal{P}$ with another Poincare
bundle defined by tensoring with the pullback of a sufficiently
positive line bundle on $J$, we can assume that
$\mathcal{O}_{\mathcal{W}}(1)$ is very ample (this line bundle
will also depend on $m$).

A point $(\overline{Q},[\Phi])\in \mathcal{H}\times \mathcal{W}$
associated to a tuple $(g,E,\varphi,u)$ verifies that the
homomorphism $\Phi$ in (\ref{Phi}), composed with evaluation,
factors as in the diagram
\begin{equation}
\label{factorthrough} \xymatrix{
(V^{\otimes s})^{\oplus c}\otimes
\mathcal{O}_{X}\ar@{->>}[r]^{(q^{\otimes s})^{\oplus
c}(m)}\ar[d]^{\Phi}
& (E(m)^{\otimes s})^{\oplus c}\ar@(d,dl)[ddl]^{\varphi}\\
H^{0}((det E)^{\otimes b}\otimes D_{u}(sm))\otimes
\mathcal{O}_{X}\ar[d]^{ev} & \\
(det E)^{\otimes b}\otimes D_{u}(sm)}
\end{equation}

Consider the relative version of the homomorphisms in
(\ref{factorthrough}), i.e. the commutative diagram on $X\times
\mathcal{H} \times \mathcal{W}$,
\begin{equation}
\label{diagramonfamilies} \xymatrix{ 0\ar[r] &
\mathcal{K}\ar[r]\ar[rd]^{f} & (V^{\otimes s})^{\oplus c}\otimes
\mathcal{O}_{X\times\mathcal{H}\times\mathcal{W}}\ar[d]^{\Phi_{\mathcal{H}\times\mathcal{W}}}\ar[r]
& (\pi_{X\times\mathcal{H}}^{\ast}E_{\mathcal{H}}(m)^{\otimes
s})^{\oplus c}\ar[r] & 0\\
 & & \mathcal{A}:=\pi^{\ast}_{X\times J}\mathcal{P}^{\otimes
 b}\otimes \pi^{\ast}_{X\times R}\mathcal{D}\otimes
 \pi^{\ast}_{X}\mathcal{O}_{X}(sm) & &}
\end{equation}
where again $\pi_{X\times J}$ (resp. $\pi_{X}$,...) denotes the
natural projection from $X\times\mathcal{H}\times\mathcal{W}$ to
$X\times J$ (resp. $X$,...). We denote by $E_{\mathcal{H}}$ the
tautological sheaf on $X\times\mathcal{H}$, and
$\Phi_{\mathcal{H}\times\mathcal{W}}$ is the relative version of
the composition $ev\circ \Phi$ in diagram (\ref{factorthrough}).

The points $(\overline{Q},[\Phi])$ where the restriction
$\Phi_{\mathcal{H}\times\mathcal{W}}|_{X\times(\overline{Q},[\Phi])}$
factors through  $(E(m)^{\otimes s})^{\oplus c}$ as in diagram
(\ref{factorthrough}) are exactly the points where
$f_{X\times(\overline{Q},[\Phi])}$ is identically zero. Hence,
points $(\overline{Q},[\Phi])$ corresponding to tuples
$(g,E,\varphi,u)$ have to verify
$f_{X\times(\overline{Q},[\Phi])}=0$ identically, a closed
condition, then we will look for them in a closed subscheme of
$\mathcal{H}\times\mathcal{W}$. We will need the following technical Lemma.
\begin{lem}\cite[Lemma 3.1]{GS1}
\label{schemefactor} Let $Y$ be a scheme, and let
$f:\mathcal{G}\longrightarrow \mathcal{F}$ be a homomorphism of
coherent sheaves on $X\times Y$. Assume that $\mathcal{F}$ is flat
over $Y$. Then there exists a unique closed subscheme $Z'\subset
Y$ satisfying the following universal property: given a Cartesian
diagram
$$\xymatrix{
X\times Z'\ar[r]^{\overline{i}}\ar[d]^{p_{Z'}} & X\times Y\ar[d]^{p}\\
Z'\ar[r]^{i} & Y}$$ $\overline{i}^{\ast}f=0$ if and only if $h$
factors through $Z$.
\end{lem}
\begin{pr}
See \cite[Lemma 3.1]{GS1}.
\end{pr}

Let $Z'$ be the scheme given by this lemma setting
$Y=\mathcal{H}\times\mathcal{W}$ and the homomorphism
$f:\mathcal{K} \rightarrow \mathcal{A}$. Let $i:Z'\inj
\mathcal{H}\times\mathcal{W}$ and $\overline{i}=id_{X}\times i$.
Then $\overline{i}^{\ast}f=0$, and we get a commutative diagram on
$X\times Z'$,
\begin{equation}
\label{diagramrestricted}
\xymatrix{\overline{i}^{\ast}\mathcal{K}\ar[r]\ar[rd]_{\overline{i}^{\ast}f}
& (V^{\otimes s})^{\oplus c}\otimes \mathcal{O}_{X\times
Z'}\ar[d]^{\overline{i}^{\ast}\tilde{\Phi}}\ar[r] &
(\overline{i}^{\ast}\pi_{X\times\mathcal{H}}^{\ast}E_{\mathcal{H}}(m)^{\otimes
s})^{\oplus c}\ar[ld]^{\tilde{\varphi}}\ar[r] & 0\\
& \overline{i}^{\ast}\mathcal{A}& &}
\end{equation}
and, hence, there is a universal family of based tensors
parametrized by $Z'$,
\begin{equation}
\label{tautologicalfamily} \varphi_{Z'}:E_{Z'}^{\otimes
s}\longrightarrow (\det E_{Z'})^{\otimes b}\otimes
\pi_{Z'}^{\ast}\mathcal{D}\; .
\end{equation}

Thanks to the tautological family (\ref{tautologicalfamily}),
given a point $(\overline{Q},[\Phi])$ in $Z'$, we get a tuple
$(g,E,\varphi,u)$ up to isomorphism. Moreover, if
$H^{0}(q(m)):V\rightarrow H^{0}(E(m))$ is an isomorphism, then we
recover exactly the original tuple $(g=H^{0}(q(m)),E,\varphi,u)$
up to isomorphism, i.e. if $(g',E',\varphi',u')$ is another tuple
corresponding to the same point $(\overline{Q},[\Phi])$, then
there exists an isomorphism $(f,\alpha)$ between $(E,\varphi,u)$
and $(E'\varphi'u')$ as in (\ref{homtensors}), and
$H^{0}(f(m))\circ q=q'$.

Let $Z\subset Z'$ be the Zariski closure of the points associated
to $\delta$-semistable tensors. Let $\pi_{\mathcal{H}}$ and
$\pi_{\mathcal{W}}$ be the projections of $Z$ to $\mathcal{H}$ and
$\mathcal{W}$, and define a \emph{polarization} on $Z$ by
\begin{equation}
\label{polarization}
\mathcal{O}_{Z}(a_{1},a_{2}):=\pi_{\mathcal{H}}^{\ast}\mathcal{O}_{\mathcal{H}}(a_{1})\otimes
\pi_{\mathcal{W}}^{\ast}\mathcal{O}_{\mathcal{W}}(a_{2})
\end{equation}
where we choose integers $a_{1}$ and $a_{2}$ for their ratio to
verify
\begin{equation}
\label{ratio}
\frac{a_{2}}{a_{1}}=\frac{r\delta(m)}{P(m)-s\delta(m)}\; .
\end{equation}
The projective scheme $Z$ is preserved by the natural $SL(V)$
action, and this action has a natural linearization en
$\mathcal{O}_{Z}(a_{1},a_{2})$, using the linearizations
on $\mathcal{O}_{\mathcal{H}}(1)$ and
$\mathcal{O}_{\mathcal{W}}(1)$.

Recall that the points of $Z$ for which $H^{0}(q(m))$ is an
isomorphism correspond (up to isomorphism) to the tuples
$(g,E,\varphi,u)$, where $g: V\simeq H^{0}(E(m))$. To get rid of
the choice of $g$, we have to take the quotient by $GL(V)$, but if
$\lambda\in\mathbb{C}^{\ast}$, $(g,E,\varphi,u)$ and $(\lambda
g,E,\varphi,u)$ correspond to the same point, and hence it
suffices to take the quotient by the action of $SL(V)$. We
construct this quotient by using Geometric Invariant Theory.

In the following, in Proposition \ref{GITstab} we identify the GIT
semistable points in $Z$ using the Hilbert-Mumford criterion (c.f.
Theorem \ref{HMcrit}). In Theorem \ref{GIT-delta} we relate
filtrations of sheaves with filtrations of the vector space $V$,
to prove that GIT semistable points of $(\overline{Q},[\Phi])\in
Z$ coincide with the points associated to $\delta$-semistable
tensors $(E,\varphi,u)$ plus an isomorphism $g$. Therefore, we
will have $Z^{ss}=Z$.

The moduli space of $\delta$-semistable tensors,
$\mathfrak{M}_{\delta}$, will be the GIT quotient of $Z^{ss}=Z$ by
$SL(V)$,
$$\mathfrak{M}_{\delta}=Z/\!\!/SL(V)\; ,$$
which is \emph{good quotient} by Theorem \ref{GIT}.

\subsection{Application of Geometric Invariant Theory}

Recall that a $1$-parameter subgroup of $G$ is a non-trivial homomorphism $\Gamma:\mathbb{C}^{*}\longrightarrow G$. In our
case, the group is $G=SL(V)=SL(p,\mathbb{C})$. It follows from
elementary representation theory that there exists a basis
$v_{1},...,v_{p}$ of $\mathbb{C}^{p}$ and $\Gamma_{i}\in
\mathbb{Z}$ such that
$$\Gamma(t)=\begin{pmatrix}
t^{\Gamma _{1}} &  & 0 \\
& \ddots &  \\
0 &  & t^{\Gamma _{p}}
\end{pmatrix}\; ,$$ so we will refer to $\Gamma$ by giving the exponents of the diagonal, $(\Gamma_{1},...,\Gamma_{p})$.
Note that $\sum_{i=1}^{p} \Gamma_{i}=0$ because $\Gamma(t)\in
SL(V)$.

The group $SL(V)$ acts on $Z$ and this action is linearized by
means of the line bundle $\mathcal{O}_{Z}(a_{1},a_{2})$. The point
$z_{0}$ is a fixed point for the $\mathbb{C}^{*}$-action on $X$
induced by $\Gamma$. Thus, $\mathbb{C}^{*}$ acts on the fiber of
$L=\mathcal{O}_{Z}(a_{1},a_{2})$ over $z_{0}$ with weight $\gamma$
and recall the definition of the numerical function in
Theorem \ref{HMcrit}, $\mu(\Gamma,x):=\gamma$, the \textit{minimum
relevant exponent} of the action of $\Gamma$ on $x\in Z$, i.e. the
minimum exponent of the diagonal of the one parameter subgroup
which acts on a non-zero coordinate of the point $x$.

A weighted filtration $(V_\bullet,n_{\bullet})$ of $V$ is a
filtration
\begin{equation}
\label{filtV} 0 \subset V_1 \subset V_2 \subset \;\cdots\; \subset
V_t \subset V_{t+1}=V,
\end{equation}
and positive numbers $n_{1},\, n_{2},\ldots , \,n_{t} > 0$. If
$t=1$ (one-step filtration), then we will take $n_{1}=1$. To a
weighted filtration we associate a vector of $\mathbb{C}^p$
defined as $\Gamma=\sum_{i=1}^{t}n_{i} \Gamma^{(\dim V_i)}$, where
\begin{equation}
\label{semistandard} \Gamma^{(k)}:=\big(
\overbrace{k-p,\ldots,k-p}^k,
 \overbrace{k,\ldots,k}^{p-k} \big)
\qquad (1\leq k < p) \, .
\end{equation}
Hence, the vector is of the form
$$\Gamma=(\overbrace{\Gamma_1,\ldots,\Gamma_1}^{\dim V^1},
\overbrace{\Gamma_2,\ldots,\Gamma_2}^{\dim V^2}, \ldots,
\overbrace{\Gamma_{t+1},\ldots,\Gamma_{t+1}}^{\dim V^{t+1}}) \;
.$$ Giving the numbers $n_1,\ldots,n_t$ is clearly equivalent to
giving the numbers $\Gamma_1,\ldots,\Gamma_{t+1}$, because
\begin{equation}
\label{1PSFW} n_{i}=\frac{\Gamma_{i+1}-\Gamma_i}{p} \quad
\text{and} \quad \sum_{i=1}^{t+1}\Gamma_i\dim V^i = 0\; .
\end{equation}
Given a $1$-parameter subgroup $\Gamma$, we associate a weighted
filtration as follows. There is a basis $\{e_1,\ldots,e_p\}$ of
$V$ where it has a diagonal form. Let
$\Gamma_1<\cdots<\Gamma_{t+1}$ be the distinct exponents, and let
$$
0\subset V_1 \subset \cdots \subset V_{t+1}=V
$$
be the associated filtration, where each $V_{i}$ is generated by
the vectors $e_{j}$ associated to exponents $\Gamma_{j}\leq
\Gamma_{i}$. Note that two $1$-parameter subgroups define the same
filtration if and only if they are conjugate by an element of the
parabolic subgroup $P\subset \SL(V)$ defined by the filtration.

Now, let $\mathcal{I}=\{1,...,t+1\}^{\times s}$ be the set of all
multi-indexes $I=(i_{1},...,i_{s})$. Define
\begin{equation}
\label{rightmuV} \mu(\Phi,V_{\bullet},n_{\bullet})=\min_{I\in
\mathcal{I}} \{\Gamma_{\dim V_{i_1}}+\cdots+\Gamma_{\dim V_{i_s}}:
\,\Phi|_{(V_{i_1}\otimes\cdots\otimes V_{i_s})^{\oplus c}}\neq 0
\}.
\end{equation}
If $I_{0}=(i_1,\ldots,i_s)$ is the multi-index giving minimum in
(\ref{rightmuV}), we will denote by
$\epsilon_i(\overline{\Phi},V_\bullet,n_{\bullet})$ (or just
$\epsilon_i(V_{\bullet})$ if the rest of the data is clear from
the context) the number of elements $k$ of the multi-index $I_{0}$
such that $\dim V_k\leq \dim V_i$.

Given a quotient $q:V\otimes \mathcal{O}_{X}(-m)\twoheadrightarrow
E$, for each subspace $V'\subset V$, we define the subsheaf
$E_{V'}\subset E$ as the image of the restriction of $q$ to $V'$,
\begin{equation}
\label{EV} \xymatrix{ V\otimes \mathcal{O}_{X}(-m)\ar@{->>}[r]
&E\\
V'\otimes \mathcal{O}_{X}(-m)\ar@{^{(}->}[u]\ar@{->>}[r] &
E_{V'}\ar@{^{(}->}[u]}
\end{equation}
Note that, in particular, $E_{V'}(m)$ is generated by global
sections.

If the quotient $q:V\otimes \mathcal{O}_{X}(-m)\twoheadrightarrow
E$ induces an injection $V\hookrightarrow H^{0}(E(m))$, and if
$E'\subset E$ is a subsheaf, we can define
\begin{equation}
\label{VE} V_{E'}=V\cap H^{0}(E'(m))\; .
\end{equation}
We will show in Proposition \ref{GITstab2} that all quotients coming from GIT semistable points
$(\overline{Q},[\Phi])\in Z$ satisfy this injectivity property, then filtrations of subsheaves will define
filtrations of vector subspaces and viceversa. Here, there are two Lemmas relating both processes.

\begin{lem}
\label{EVE} Given a point $(\overline{Q},[\Phi])\in Z$ such that
$q$ induces an injection $V\hookrightarrow H^{0}(E(m))$, and a
weighted filtration $(E_{\bullet},n_{\bullet})$ of $E$, we have:
\begin{enumerate}
    \item \label{EVE1} $E_{V_{E_{i}}}\subset E_{i}$
    \item \label{EVE2} If $\varphi|_{(E_{i_{1}}\otimes \cdots \otimes
    E_{i_{s}})^{\oplus c}}=0$, then $\Phi|_{(V_{E_{i_{1}}}\otimes \cdots \otimes
    V_{E_{i_{s}}})^{\oplus c}}=0$
    \item \label{EVE3} $\sum_{i=1}^{t}-n_{i}\epsilon_{i}(\varphi,E_{\bullet},n_{\bullet})\leq
    \sum_{i=1}^{t}-n_{i}\epsilon_{i}(\Phi,V_{E_{\bullet}},n_{\bullet})$
\end{enumerate}
Furthermore, if $q$ induces an isomorphism $V\simeq H^{0}(E(m))$,
all $E_{i}$ are $m$-regular and all $E_{i_{1}}\otimes \cdots
\otimes E_{i_{s}}$ are $sm$-regular, then \ref{EVE1} becomes an
equality, \ref{EVE2} becomes "if and only if" and \ref{EVE3}
becomes an equality.
\end{lem}
\begin{pr}
By definition, $E_{V_{E_{i}}}=q|_{V_{E_{i}}}=q(V\cap
H^{0}(E_{i}(m))\otimes \mathcal{O}_{X}(-m))\subset E_{i}$ because
we are evaluating sections that, in particular, are sections of
$E_{i}$. And, if $E_{i}$ is $m$-regular then it will be generated
by their global sections. Hence, every point in the total space of
$E_{i}$ will be the image of any section in $H^{0}(E_{i}(m))\otimes
\mathcal{O}_{X}(-m)=V\cap H^{0}(E_{i}(m))\otimes
\mathcal{O}_{X}(-m)$, provided $V\simeq H^{0}(E(m))$. This proves
\ref{EVE1}.

To prove \ref{EVE2}, note that $\Phi$, as a morphism of sections, is given at each point by the morphism of
sheaves $\varphi$. A section of $V_{E_{i_{j}}}$ is, in particular, a section of $H^{0}(E_{i_{j}}(m))$, so if
$\varphi$ vanishes on a factor $E_{i_{j}}$, then it vanishes on $E_{i_{j}}(m)$, and therefore $\Phi$ vanishes on
$V_{E_{i_{j}}}$. If $V\simeq H^{0}(E(m))$ and $E_{i}$ is $m$-regular, by \ref{EVE1} we have
$E_{V_{E_{i}}}=E_{i}$. If $E_{i_{1}}\otimes \cdots \otimes E_{i_{s}}$ is $sm$-regular, then
$$E_{i_{1}}\otimes \cdots \otimes E_{i_{s}}\otimes
\mathcal{O}_{X}(sm)=E_{i_{1}}(m)\otimes \cdots \otimes
E_{i_{s}}(m)$$  is generated by global sections and every element
of $E_{i_{1}}(m)\otimes \cdots \otimes E_{i_{s}}(m)$ comes from a
section of
$$V_{E_{i_{1}}}\otimes \cdots \otimes V_{E_{i_{s}}}=H^{0}(E_{i_{1}}(m))\otimes\cdots\otimes H^{0}(E_{i_{s}}(m)\; .$$
Therefore, if $\Phi$ vanishes, then $\varphi$ does.

Recall that if $I_{0}$ is the multi-index giving minimum in (\ref{rightmuV}),
$\epsilon_i(\overline{\Phi},V_{E_{\bullet}},n_{\bullet})$ is the number of elements $k$ of $I_{0}$ such that
$\dim V_{E_{i_{k}}}\leq \dim V_{E_{i}}$, and similarly for $\epsilon_i(\varphi,E_\bullet,n_{\bullet})$ with $\rk
E_{i_{k}}\leq \rk E_{i}$ (c.f. (\ref{rightmu})). Then, if $\varphi$ (resp. $\Phi$) vanishes on a filter $V_{i}$
(resp. $E_{i}$), the index $i$ cannot be taken into account in the calculation of
$\epsilon_i(\overline{\Phi},V_{E_{\bullet}},n_{\bullet})$ (resp. $\epsilon_i(\varphi,E_\bullet,n_{\bullet})$).
Also note that, by \ref{VEV2}, if $\varphi$ vanishes on a filter $E_{i_{k}}$, then $\Phi$ vanishes on
$V_{E_{i_{k}}}$. Hence, all $V_{E_{i_{k}}}$ are not counted in
$\epsilon_i(\overline{\Phi},V_{E_{\bullet}},n_{\bullet})$ if they were not counted in
$\epsilon_i(\varphi,E_\bullet,n_{\bullet})$. Therefore, $\epsilon_i(\varphi,E_\bullet,n_{\bullet})\geq
\epsilon_i(\overline{\Phi},V_{E_{\bullet}},n_{\bullet})$, and this proves \ref{EVE3}.
\end{pr}

\begin{lem}
\label{VEV} Given a point $(\overline{Q},[\Phi])\in Z$ such that
$q$ induces an injection $V\hookrightarrow H^{0}(E(m))$, and a
weighted filtration $(V_{\bullet},n_{\bullet})$ of $V$, we have:
\begin{enumerate}
    \item \label{VEV1} $V_{i}\subset V_{E_{V_{i}}}$
    \item \label{VEV2} $\varphi|_{(E_{V_{i_{1}}}\otimes \cdots \otimes
    E_{V_{i_{s}}})^{\oplus c}}=0$ if and only if $\Phi|_{(V_{i_{1}}\otimes \cdots \otimes
    V_{i_{s}})^{\oplus c}}=0$
    \item \label{VEV3}
    $\sum_{i=1}^{t}-n_{i}\epsilon_{i}(\varphi,E_{V_{\bullet}},n_{\bullet})=
    \sum_{i=1}^{t}-n_{i}\epsilon_{i}(\Phi,V_{\bullet},n_{\bullet})$
\end{enumerate}
\end{lem}
\begin{pr}
Similarly to Lemma \ref{EVE}, $E_{V_{i}}=q|_{V_{i}}$, then
$V_{i}\subset V\cap H^{0}(E_{V_{i}}(m))=V_{E_{V_{i}}}$, by
definition, hence we prove \ref{VEV1}.

To prove \ref{VEV2}, as $\Phi$ is given at each point by the morphism of sheaves $\varphi$, and $E_{V_{i}}$ is
generated by the global sections of $V_{i}$, if $\Phi$ vanishes on sections $V_{i_{i}}\otimes \cdots \otimes
V_{i_{s}}$, then $\varphi$ vanishes on the respective sheaves $E_{V_{i_{1}}}\otimes \cdots \otimes
E_{V_{i_{s}}}$ and viceversa.

Statement \ref{VEV3} follows from \ref{VEV2} and the same argument that in the proof of \ref{EVE3} in Lemma
\ref{EVE}.
\end{pr}

\begin{prop}
\label{GITstab} The point $(\overline{Q},[\Phi])\in Z$ is \emph{GIT semistable} with respect to
$\mathcal{O}_{Z}(a_{1},a_{2})$ if and only if for every weighted filtration $(V_{\bullet},n_{\bullet})$ of $V$
\begin{equation}
\label{Vcrit} a_{1}\sum_{i=1}^{t}n_{i}(r\dim V_{i}-\rk E_{V_{i}}\dim V)+a_{2}\sum_{i=1}^{t}n_{i}(s\dim
V_{i}-\epsilon_{i}(V_{\bullet})\dim V)\leq 0
\end{equation}
The point $(\overline{Q},[\Phi])$ is \emph{GIT stable} if we get a strict inequality for every weighted filtration. In any case,
there exists an integer $A_{2}$ (depending only on $m$, $P$, $s$, $b$, $c$, $\mathcal{D}$) such that it is
enough to consider weighted filtrations with $n_{i}\leq A_{2}$.
\end{prop}

\begin{pr}
By the Hilbert-Mumford criterion, Theorem \ref{HMcrit}, a point is GIT semistable (resp. GIT stable) if and only
if for all $1$-parameter subgroups $\Gamma$ of $SL(V)$,
$$\mu((\overline{Q},[\Phi]),\Gamma)=a_{1}\mu(\overline{Q},\Gamma)+a_{2}\mu([\Phi],\Gamma)\leq 0$$
(resp. $<$). We have seen that, given a $1$-parameter subgroup $\Gamma$ of $\SL(V)$, we associate a weighted
filtration $(V_{\bullet},n_{\bullet})$ where each exponent $\Gamma_{i}$ corresponds to the action of $\Gamma$ on
$V^{i}=V_{i}/V_{i-1}$. Denote $\mathcal{I}'=\{1,...,t+1\}^{\times r}$. Then, the minimum weight of the action of
$\Gamma$ on $\overline{Q}\in \mathbb{P}(\Hom(\wedge^{r}V,A)$ is (c.f. \cite{Si1} and \cite{HL2}),
$$\mu(\overline{Q},\Gamma)=\min_{I\in \mathcal{I}'}\{\Gamma_{i_{1}}+\cdots+\Gamma_{i_{r}}
:Q|_{V_{i_{1}}\wedge\cdots \wedge V_{i_{r}}}\neq 0\}\; .$$ Note that $\Gamma$ acts trivially on $A=H^{0}(\det
(E)(rm))$ and observe that the evaluation $Q$ does not vanish on a wedge of $r$ sections,
$e_{i_{i}}\wedge\cdots\wedge e_{i_{r}}$, whenever $e_{i_{1}},...,e_{i_{r}}$ span the fiber of $E$ over the generic
point $x\in X$. Recall that $\Gamma_{1}<\ldots<\Gamma_{t+1}$, then to achieve the minimum we have to take the
minimum exponent $\Gamma_{1}$ as many times as possible, then take $\Gamma_{2}$ as many times as possible, and
so on, while $Q\neq 0$. Therefore, this occurs when we take $\Gamma_{1}$ a number $\rk E_{V_{1}}$ of times,
then we take $\Gamma_{2}$ a number $(\rk E_{V_{2}}-\rk E_{V_{1}})$ of times, etc, and finally we take
$\Gamma_{t+1}$ a number $(\rk E_{V_{t+1}}-\rk E_{V_{t}})$ of times, hence
$$\mu(\overline{Q},\Gamma)=\sum_{i=1}^{t+1}\Gamma_{i}(\rk
E_{V_{i}}-\rk E_{V_{i-1}})\; .$$ Making calculations
$$\mu(\overline{Q},\Gamma)=r\Gamma_{t+1}-\sum_{i=1}^{t}(\Gamma_{i+1}-\Gamma_{i})\rk
E_{V_{i}}=r\Gamma_{t+1}-\sum_{i=1}^{t}n_{i}\dim V\rk E_{V_{i}}$$
$$=r(\sum_{i=1}^{t}n_{i}\dim
V_{i})-\sum_{i=1}^{t}n_{i}\dim V\rk
E_{V_{i}}=\sum_{i=1}^{t}n_{i}(r\dim V_{i}-\rk E_{V_{i}}\dim V)\;
.$$

Now we calculate the minimum weight of the action of $\Gamma$ on
$$[\Phi]\in \mathcal{W}=\mathbb{P}\big(((V^{\otimes s})^{\oplus
c})^{\vee}\otimes \pi_{J\times R,\ast}(\pi^{\ast}_{X\times
J}\mathcal{P}^{\otimes b}\otimes \pi^{\ast}_{X\times
R}\mathcal{D}(sm))\big )\; ,$$ where note that $\Gamma$ only acts non
trivially on $V$. Then,
$$\mu([\Phi],\Gamma)=\min_{I\in \mathcal{I}}\{\Gamma_{i_{1}}+\cdots+\Gamma_{i_{s}}
:\Phi|_{(V_{i_{1}}\otimes\cdots \otimes V_{i_{s}})^{\oplus c}}\neq 0\}\; .$$ Similarly, to achieve the minimum,
we have to take $\Gamma_{1}$ as many times as $V_{1}$ can appear in the multi-index $I$ while the restriction of
$\Phi$ does not vanish, then take $\Gamma_{2}$ as many times as $V_{2}$ can appear minus the number of times
$V_{1}$ appears, and so on. This can be written in terms of the symbols $\epsilon^i(V_{\bullet})$, the number of
times that each index $i$ appears on $I$:
$$\mu([\Phi],\Gamma)=\sum_{i=1}^{t}\Gamma_{i}\epsilon^i(V_{\bullet})=
\sum_{i=1}^{t}\Gamma_{i}(\epsilon_{i+1}(V_{\bullet})-\epsilon_{i}(V_{\bullet}))$$
$$=s\Gamma_{t+1}-\sum_{i=1}^{t}(\Gamma_{i+1}-\Gamma_{i})\epsilon_{i}(V_{\bullet})=
s\Gamma_{t+1}-\sum_{i=1}^{t}n_{i}\dim V\epsilon_{i}(V_{\bullet})$$
$$=s(\sum_{i=1}^{t}n_{i}\dim
V_{i})-\sum_{i=1}^{t}n_{i}\dim
V\epsilon_{i}(V_{\bullet})=\sum_{i=1}^{t}n_{i}(s\dim
V_{i}-\epsilon_{i}(V_{\bullet})\dim V)\; .$$

The last statement follows from an argument similar to the proof
of Lemma \ref{boundonmi}, with $\mathbb{Z}^{r}$ replaced by
$\mathbb{Z}^{p}$.
\end{pr}

\begin{prop}
\label{GITstab2} The point $(\overline{Q},[\Phi])\in Z$ is GIT semistable if and only if for every weighted
filtration $(E_{\bullet},n_{\bullet})$ of $E$, it is
\begin{equation}
\label{Ecrit} \sum_{i=1}^{t}n_{i}\big((r\dim V_{E_{i}}-\rk
E_{i}\dim V)+\delta(m)(s\rk
E_{i}-\epsilon_{i}(E_{\bullet})r)\big)\leq 0\; .
\end{equation}
If $(\overline{Q},[\Phi])$ is GIT stable we get a strict inequality for every weighted filtration. Moreover, if
$(\overline{Q},[\Phi])$ is GIT semistable, then the induced map $f_{q}=H^{0}(q(m))$, $f_{q}:V\rightarrow H^{0}(E(m))$
is injective.
\end{prop}
\begin{pr}
Let us show first that if $(\overline{Q},[\Phi])$ is GIT
semistable, then the induced map $f_{q}$ is injective. Let $V'$ be
its kernel and consider the one-step filtration $V'\subset V$. We
have $E_{V'}=0$ by definition and, if we calculate
$\mu([\Phi],\Gamma)$ for the $1$-parameter subgroup $\Gamma$
associated to the one-step filtration $V'\subset V$, it is
$$\mu([\Phi],\Gamma)=s\dim V'\; ,$$
because $\epsilon_{1}=0$ ($V'$ does not appear in the multi-index giving minimum in (\ref{rightmuV}) because
$E_{V'}=0$). Therefore, applying Proposition \ref{GITstab},
$$a_{1}r\dim V'+a_{2}s\dim V'\leq 0\; ,$$
and $\dim V'=0$, hence $f_{q}$ is injective.

Using (\ref{ratio}), the inequality (\ref{Vcrit}) becomes
$$\sum_{i=1}^{t}n_{i}\big((r\dim V_{i}-\rk
E_{V_{i}}\dim V)(P(m)-s\delta(m))+r\delta(m)(s\dim
V_{i}-\epsilon_{i}(V_{\bullet})\dim V)\big)\leq 0$$ which, setting
$P(m)=\dim V$, is equivalent to
\begin{equation}
\label{Vcrit2} \sum_{i=1}^{t}n_{i}\big((r\dim V_{i}-\rk E_{V_{i}}\dim V)+\delta(m)(s\rk
E_{V_{i}}-\epsilon_{i}(V_{\bullet})r)\big)\leq 0\; .
\end{equation}

Now let $(\overline{Q},[\Phi])$ be a GIT semistable point. Take a weighted filtration $(E_{\bullet},n_{\bullet})$
of $E$. Consider the induced weighted filtration $(V_{E_{\bullet}},n_{\bullet})$ of $V$. By Proposition
\ref{GITstab} and using (\ref{ratio}) we have
$$\sum_{i=1}^{t}n_{i}\big((r\dim V_{E_{i}}-\rk
E_{V_{E_{i}}}\dim V)(P(m)-s\delta(m))+r\delta(m)(s\dim
V_{E_{i}}-\epsilon_{i}(V_{E_{\bullet}})\dim V)\big)\leq 0\; ,$$
which is equivalent to
$$\sum_{i=1}^{t}n_{i}\big((r\dim V_{E_{i}}-\rk
E_{V_{E_{i}}}\dim V)+\delta(m)(s\rk E_{V_{E_{i}}}-\epsilon_{i}(V_{E_{\bullet}})r)\big)\leq 0\; .$$

Then, by Lemma \ref{EVE}, using statement \ref{EVE1} we have
$E_{V_{E_{i}}}\subset E_{i}$, then $\rk E_{V_{E_{i}}}\leq \rk
E_{i}$, and statement \ref{EVE3} gives
$-\epsilon_{i}(E_{\bullet})\leq -\epsilon_{i}(V_{E_{\bullet}})$,
therefore
$$\sum_{i=1}^{t}n_{i}\big((r\dim V_{E_{i}}-\rk
E_{i}\dim V)+\delta(m)(s\rk
E_{i}-\epsilon_{i}(E_{\bullet})r)\big)\leq$$
$$\sum_{i=1}^{t}n_{i}\big((r\dim V_{E_{i}}-\rk
E_{V_{E_{i}}}\dim V)+\delta(m)(s\rk
E_{V_{E_{i}}}-\epsilon_{i}(V_{E_{\bullet}})r)\big)\leq 0\; ,$$
hence (\ref{Ecrit}) holds. Note that, if we start with a GIT
stable point, we can substitute the inequalities by strict
inequalities.

On the other hand, suppose that (\ref{Ecrit}) holds. Take a weighted filtration $(V_{\bullet},n_{\bullet})$ of
$V$. Then we get an induced weighted filtration $(E_{V_{\bullet}},n_{\bullet})$ of $E$. For this filtration,
(\ref{Ecrit}) becomes
$$\sum_{i=1}^{t}n_{i}\big((r\dim V_{E_{V_{i}}}-\rk E_{V_{i}}\dim
V)+\delta(m)(s\rk E_{V_{i}}-\epsilon_{i}(E_{V_{\bullet}})r)\big)\leq 0\; .$$ By Lemma \ref{VEV}, using statement
\ref{VEV1} we have $V_{i}\subset V_{E_{V_{i}}}$, then $\dim V_{i}\leq \dim V_{E_{V_{i}}}$, and statement
\ref{VEV2} gives $-\epsilon_{i}(E_{V_{\bullet}})= -\epsilon_{i}(V_{\bullet})$. Hence, applying (\ref{Vcrit2})
(which we have seen that is equivalent to (\ref{Vcrit})) for the filtration $(V_{\bullet},n_{\bullet})$ we get
$$\sum_{i=1}^{t}n_{i}\big((r\dim V_{i}-\rk
E_{V_{i}}\dim V)+\delta(m)(s\rk
E_{V_{i}}-\epsilon_{i}(V_{\bullet})r)\big)\leq$$
$$\sum_{i=1}^{t}n_{i}\big((r\dim V_{E_{V_{i}}}-\rk E_{V_{i}}\dim
V)+\delta(m)(s\rk E_{V_{i}}-\epsilon_{i}(E_{V_{\bullet}})r)\big)\leq 0\; ,$$ and therefore the point
$(\overline{Q},[\Phi])$ is GIT semistable by Proposition \ref{GITstab}. If we start with a strict inequality in
(\ref{Ecrit}), we get a GIT stable point.
\end{pr}

\begin{thm}
\label{GIT-delta} Assume $m>N$. A point $(\overline{Q},[\Phi])\in Z$ is GIT semistable (resp. GIT stable) if and
only if the corresponding tensor $(E,\varphi,u)$ is $\delta$-semistable (resp. $\delta$-stable) and the linear
map $f_{q}:V\simeq H^{0}(E(m))$ induced by $q$ is an isomorphism.
\end{thm}

\begin{pr}
$\Rightarrow)$ Suppose $(\overline{Q},[\Phi])\in Z$ is
GIT semistable and let $(E_{\bullet},n_{\bullet})$ be a weighted
filtration of $E$. We will use Theorem \ref{conditions_stability}
to prove that $(E,\varphi,u)$ is $\delta$-semistable, and similarly for stable.

Recall that, as we have seen in $(2.\Rightarrow 3.)$ in the proof of
Theorem \ref{conditions_stability}, we have this inequality for
the dimensions of the vector spaces
$$h^{0}(E(m))\leq h^{0}(E_{i}(m))+h^{0}(E^{i}(m))\; .$$ Then, by definition, $V_{E_{i}}:=V\cap
H^{0}(E_{i}(m))$ and, by dimensions formula,
$$\dim V_{E_{i}}=\dim V+h^{0}(E_{i}(m))-\dim (V\cup H^{0}(E_{i}(m)))\geq$$
$$\dim V+h^{0}(E_{i}(m))-h^{0}(E(m))\geq \dim V-h^{0}(E^{i}(m))\; .$$

Recall that $P(m)=\dim V$. Therefore we obtain the inequality of condition 3 in Theorem
\ref{conditions_stability},
$$\big(\sum_{i=1}^{t}n_{i}(\rk
E^{i}P(m)-rh^{0}(E^{i}(m)))\big)+\delta(m)\mu(\varphi,E_{\bullet},n_{\bullet})=$$
$$\sum_{i=1}^{t}n_{i}\big((\rk
E^{i}P(m)-rh^{0}(E^{i}(m)))+\delta(m)(s\rk
E_{i}-\epsilon_{i}(E_{\bullet})r)\big)=$$
$$\sum_{i=1}^{t}n_{i}\big((r(\dim V-h^{0}(E^{i}(m)))-\rk E_{i}\dim V)+\delta(m)(s\rk
E_{i}-\epsilon_{i}(E_{\bullet})r)\big)\leq$$
$$\sum_{i=1}^{t}n_{i}\big((r\dim V_{E_{i}}-\rk E_{i}\dim
V)+\delta(m)(s\rk E_{i}-\epsilon_{i}(E_{\bullet})r)\big)\leq 0\;
,$$ by Proposition \ref{GITstab2}. If we start with a GIT stable
point we get a strict inequality.

To apply Theorem \ref{conditions_stability} we need to show that
$E$ is torsion free. By Lemma \ref{torsionfreeator}, there exists
a tensor $(F,\psi, u)$ with $F$ torsion free such that
$P_{E}=P_{F}$ and a exact sequence
$$0\longrightarrow T(E)\longrightarrow E \overset{\beta}{\longrightarrow}
F$$ where define $E'':=\beta(E)$. Consider a weighted filtration
$(F_{\bullet},n_{\bullet})$ of $F$. Let $F^{i}=F/F_{i}$, and let
$E^{i}$ be the image of $E$ in $F^{i}$, $E^{i}=E''/F_{i}$. Let
$E_{i}$ be the kernel of $E\longrightarrow E^{i}$. Then
$\rk(F_{i})=\rk(E_{i})=r_{i}$, because $\rk E=\rk E^{i}+\rk E_{i}$
and $\rk E=\rk E''$. Also, $E^{i}=E''/F_{i}\subset F/F_{i}=F^{i}$,
then $h^{0}(E^{i}(m))\leq h^{0}(F^{i}(m))$. Moreover,
$\epsilon(\psi,F_{\bullet},n_{\bullet})=\epsilon(\varphi,E_{\bullet},n_{\bullet})$,
because the difference between filters of $F$ and $E$ occurs in
the $0$-rank torsion subsheaf. Using this and applying condition
$3$ in Theorem \ref{conditions_stability} to
$(F_{\bullet},n_{\bullet})$, we get
$$\big(\sum_{i=1}^{t}n_{i}\big((\rk
F^{i}P(m)-rh^{0}(F^{i}(m))\big)+\delta(m)\mu(\psi,F_{\bullet},n_{\bullet})=$$
$$\sum_{i=1}^{t}n_{i}\big((\rk
F^{i}P(m)-rh^{0}(F^{i}(m)))+\delta(m)(s\rk
F_{i}-\epsilon_{i}(\psi,F_{\bullet},n_{\bullet})r)\big)\leq$$
$$\sum_{i=1}^{t}n_{i}\big((\rk
E^{i}P(m)-rh^{0}(E^{i}(m)))+\delta(m)(s\rk
E_{i}-\epsilon_{i}(\varphi,E_{\bullet},n_{\bullet})r)\big)\leq 0\;
,$$ and hence Theorem \ref{conditions_stability} implies that
$(F,\psi,u)$ is $\delta$-semistable.

Apply condition $3$ of Theorem \ref{conditions_stability} to the
one-step filtration $T(E)\subset E$, then $E/T(E)\simeq E''$, and
$$\rk (E'')P(m)-rh^{0}(E'')(m)+\delta(m)(s\rk(T(E))-\epsilon_{1}(T(E)\subset
E)r)\leq 0$$
$$\Leftrightarrow P(m)-h^{0}(E'')(m)\leq 0\; ,$$
where note that $\rk T(E)=0$ and $\epsilon_{1}(T(E)\subset E)=0$.
Hence,
$$P(m)\leq h^{0}(E''(m))\leq h^{0}(F(m))=P(m)\; ,$$
where the second inequality follows from $E''\subset F$ and the third equality does from the conclusion of
Theorem \ref{conditions_stability} about $m$-regularity of $\delta$-semistable tensors. Hence, equality holds at
all places and $h^{0}(F(m))=h^{0}(E''(m))$. Since $F$ is globally generated, $F=E''$ and, therefore, $T(E)=0$
and $E$ is torsion free. Then, by Theorem \ref{conditions_stability}, $(E,\varphi, u)$ is $\delta$-semistable.

Finally, we have seen that if $(\overline{Q},[\Phi])\in Z$ is GIT
semistable, then the linear map $f_{q}:V\longrightarrow
H^{0}(E(m))$ is injective by Proposition \ref{GITstab}, and since
$(E,\varphi,u)$ is $\delta$-semistable, then $E$ is $m$-regular by
Proposition \ref{conditions_stability}. Given that $\dim
V=P(m)=h^{0}(E(m))$, therefore $f_{q}$ is an isomorphism.

$\Leftarrow )$ Suppose $(E,\varphi,u)$ is $\delta$-semistable, and
$q$ induces an isomorphism in the linear map
$f_{q}:V\longrightarrow H^{0}(E(m))$. Then we have
$V_{E'}=H^{0}(E'(m))$ for any subsheaf $E'\subset E$ and Theorem
\ref{conditions_stability} condition 2 says that for all weighted
filtrations $(E_{\bullet},n_{\bullet})$ of $E$,
$$\big(\sum_{i=1}^{t}n_{i}(rh^{0}(E_{i})-\rk E_{i}P(m))\big)
+\delta(m)\mu(\varphi,E_{\bullet},n_{\bullet})=$$
$$\sum_{i=1}^{t}n_{i}\big((r\dim V_{E_{i}}-\rk E_{i}P(m))
+\delta(m)(s\rk E_{i}-\epsilon_{i}(E_{\bullet})r)\big)\leq 0\; ,$$
which is exactly (\ref{Ecrit}) in Proposition \ref{GITstab2},
therefore $(\overline{Q},[\Phi])$ is GIT semistable. Similarly, if
$(E,\varphi,u)$ is $\delta$-stable we obtain a strict inequality
and, hence, $(\overline{Q},[\Phi])$ is GIT stable.
\end{pr}

\begin{cor}
Let $(E,\varphi,u)$ be a $\delta$-semistable tensor and let
$(E_{\bullet},n_{\bullet})$ be a weighted filtration of $E$. Then
the induced morphism $f_{q}:V\rightarrow H^{0}(E(m))$ is an
isomorphism and, therefore, $V=H^{0}(E(m))$ and
$V_{E_{i}}=H^{0}(E_{i}(m))$, for all $i$.
\end{cor}
\begin{pr}
It follows from the proof of Theorem \ref{GIT-delta}
\end{pr}

Now, recall that given $V$, a vector space of dimension $P(m)$, and a $1$-parameter subgroup $\Gamma$ of $SL(V)$
given in its diagonal form
$$\Gamma=(\overbrace{\Gamma_1,\ldots,\Gamma_1}^{\dim V^1},
\overbrace{\Gamma_2,\ldots,\Gamma_2}^{\dim V^2}, \ldots, \overbrace{\Gamma_{t+1},\ldots,\Gamma_{t+1}}^{\dim
V^{t+1}}) \; ,$$ we get a weighted filtration $(V_{\bullet},n_{\bullet})$ of $V$ and a splitting $V=\oplus_{i}
V^{i}$ of this filtration. Defining $E_{V_{i}}=q(V_{i}\otimes \mathcal{O}_{X}(-m))$ we obtain a weighted
filtration $(E_{\bullet},n_{\bullet})$ of $E$.

Conversely, let $(E_{\bullet},n_{\bullet})$ be a weighted filtration of $E$ and $V=\oplus_{i} V^{i}$ a splitting
of the filtration $V_{i}=H^{0}(E_{i}(m))$. This gives a $1$-parameter subgroup $\Gamma$ of $SL(V)$ defined as
$v^{i}\mapsto t^{\lambda_{i}}v^{i}$, for $v^{i}\in V^{i}$, with relations (\ref{1PSFW}).

We will use the following proposition to prove the criterion for
S-equivalence.

\begin{prop}
\label{S-equivalence} Suppose that $m>N$. Let $(E,\varphi,u)$ be a $\delta$-semistable tensor, $f:V\simeq
H^{0}(E(m))$ an isomorphism, and let $(\overline{Q},[\Phi])\in Z$ be the corresponding GIT semistable point (c.f. Theorem \ref{GIT-delta}). The
above construction gives a bijection between $1$-parameter subgroups of $SL(V)$ with
$\mu((\overline{Q},[\Phi]),\Gamma)=0$ on the one hand, and weighted filtrations $(E_{\bullet},n_{\bullet})$ of
$E$ with
\begin{equation}
\label{Sequivcrit} \sum_{i=1}^{t}n_{i}\big((rP_{E_{i}}-\rk
E_{i}P)\big)+\delta\mu(\varphi,E_{\bullet},n_{\bullet})=0
\end{equation}
together with a splitting of the filtration $H^{0}(E_{\bullet}(m))$ of $V\simeq H^{0}(E(m))$ on the other hand.
\end{prop}
\begin{pr}
Let $\Gamma$ be a $1$-parameter subgroup of $SL(V)$ with
$\mu((\overline{Q},[\Phi]),\Gamma)=0$. Then we get a weighted
filtration $(V_{\bullet},n_{\bullet})$ of $V$ and, by evaluating,
a weighted filtration $(E_{V_{\bullet}},n_{\bullet})$ of $E$. By
hypothesis, the proof of Proposition \ref{GITstab} gives equality
in (\ref{Vcrit}) applied to $(V_{\bullet},n_{\bullet})$, and in
the proof of Proposition \ref{GITstab2} we have seen that (\ref{Vcrit})
is equivalent to (\ref{Vcrit2}). Therefore, using Lemma \ref{VEV},
statement \ref{VEV3}, we get
$$\sum_{i=1}^{t}n_{i}\big((r\dim V_{i}-\rk E_{V_{i}}P(m))+\delta(m)(s\rk
E_{V_{i}}-\epsilon_{i}(E_{V_{\bullet}})r)\big)=$$
$$\sum_{i=1}^{t}n_{i}\big((r\dim V_{i}-\rk E_{V_{i}}P(m))\big)+\delta(m)\mu(\varphi,E_{V_{\bullet}},n_{\bullet})=0\; ,$$
where we also use that $\dim V=P(m)$ and (\ref{rightmu2}).
Statement \ref{VEV1} of Lemma \ref{VEV} gives $V_{i}\subset
V_{E_{V_{i}}}=H^{0}(E_{V_{i}}(m))$, hence
$$\sum_{i=1}^{t}n_{i}\big((rh^{0}(E_{V_{i}}(m))-\rk E_{V_{i}}P(m))\big)+\delta(m)\mu(\varphi,E_{V_{\bullet}},n_{\bullet})\geq 0$$
but, by Theorem \ref{conditions_stability}, condition 2, this must be non-positive, hence
$V_{i}=H^{0}(E_{V_{i}}(m))=V_{E_{V_{i}}}$, and last inequality is an equality. By Corollary \ref{equalto0},
$E_{V_{i}}\in \mathcal{S}_{0}$, and hence $h^{0}(E_{V_{i}}(m))=P_{E_{V_{i}}}(m)$ for all $i$. Therefore, as
equality holds for $m$, Lemma \ref{asintotic} gives the equality of polynomials
$$\sum_{i=1}^{t}n_{i}\big((rP_{E_{V_{i}}}-\rk E_{V_{i}}P)\big)+\delta\mu(\varphi,E_{V_{\bullet}},n_{\bullet})=0\; .$$

Conversely, let $(E_{\bullet},n_{\bullet})$ be a weighted filtration of $E$ such that (\ref{Sequivcrit}) holds,
together with a splitting of the filtration $H^{0}(E_{\bullet}(m))$ of $V\simeq H^{0}(E(m))$, and let $\Gamma$
be the associated $1$-parameter subgroup of $SL(V)$. Evaluating the expression (\ref{Sequivcrit}) in $m$ gives,
in particular,
$$\sum_{i=1}^{t}n_{i}\big((rP_{E_{i}}(m)-\rk E_{i}P(m))\big)+\delta(m)\mu(\varphi,E_{\bullet},n_{\bullet})=$$
$$\sum_{i=1}^{t}n_{i}\big((rP_{E_{i}}(m)-\rk E_{i}P(m))+\delta(m)(s\rk{E_{i}}-\epsilon_{i}(E_{\bullet})r)\big)=0\; ,$$
using (\ref{rightmu2}). By the proof of implication $(3.\Rightarrow 1.)$ in Theorem \ref{conditions_stability},
since we get an equality, it is $E_{i}\in \mathcal{S}_{0}$ for all $i$, hence $\dim
V_{E_{i}}=h^{0}(E_{i}(m))=P_{E_{i}}(m)$ for all $i$, and the previous equality becomes
$$\sum_{i=1}^{t}n_{i}\big((r\dim V_{E_{i}}-\rk E_{i}P(m))+\delta(m)(s\rk{E_{i}}-\epsilon_{i}(E_{\bullet})r)\big)=0\; .$$
Using the strong version of Lemma \ref{EVE}, $E_{V_{E_{i}}}=E_{i}$ and
$\epsilon_{i}(E_{\bullet})=\epsilon_{i}(V_{E_{\bullet}})$, then
$$\sum_{i=1}^{t}n_{i}\big((r\dim V_{E_{i}}-\rk E_{V_{E_{i}}}\dim
V)+\delta(m)(s\rk E_{V_{E_{i}}}-\epsilon_{i}(V_{E_{\bullet}})r)\big)=0$$ which is (\ref{Vcrit2}) applied to the
weighted filtration $(V_{E_{\bullet}},n_{\bullet})$ of $V$ and, by the proof of Proposition \ref{GITstab2},
equivalent to equality in (\ref{Vcrit}), therefore $\mu((\overline{Q},[\Phi]),\Gamma)=0$.

We have that $V_{i}$ generates $E_{V_{i}}(m)$, then take
$H^{0}(E_{V_{i}}(m))=V\cap H^{0}(E_{V_{i}}(m))=V_{E_{V_{i}}}$,
which we have seen that is equal to $V_{i}$. Conversely, take
$E_{i}$, then we have that $H^{0}(E_{i}(m))=V\cap
H^{0}(E_{i}(m))=V_{E_{i}}$ and, evaluating, we obtain
$E_{V_{E_{i}}}$, which we have seen that is equal to $E_{i}$.
Therefore, this gives the bijection.
\end{pr}

\subsection{Proof of Theorem \ref{maintheorem}}
\label{proofthmtensors}

\begin{pr}[\textbf{Proof of the Theorem \ref{maintheorem}}]
We follow \cite[Proof of Theorem 1.8]{GS1} which also follows closely \cite[Proof of Main Theorem
0.1]{HL2}.

Recall notation from section \ref{giesekerembedding} and Theorem
\ref{GIT}. We will use Theorem \ref{GIT-delta}, where we show that
GIT semistable points correspond to $\delta$-semistable tensors.
Let $\mathfrak{M}_{\delta}$ (respectively
$\mathfrak{M}^{s}_{\delta}$) be the GIT quotient of $Z$
(respectively $Z^{s}$) by $SL(V)$. Since $Z$ is projective,
$\mathfrak{M}_{\delta}$ is also projective and, by Theorem
\ref{GIT}, $\mathfrak{M}^{s}_{\delta}$ is an open subset of the
projective scheme $\mathfrak{M}_{\delta}$. The restriction
$Z^{s}\longrightarrow \mathfrak{M}^{s}_{\delta}$ to the stable
points is a geometric quotient where the fibers are
$SL(V)$-orbits, and hence the points of
$\mathfrak{M}^{s}_{\delta}$ correspond to isomorphism classes of
$\delta$-stable tensors. We have to show that
$\mathfrak{M}_{\delta}$ corepresents the functor
$\mathcal{M}_{\delta}$ (c.f. Definition \ref{corepresents}).

Let $(E_{T},\varphi_{T},u_{T},N)$ be a family of $\delta$-semistable tensors parametrized by a scheme $T$, as in
(\ref{tensoronfamily}). Then, $\mathcal{V}:=\pi_{T,\ast}(E_{T}\otimes \pi_{X}^{\ast}\mathcal{O}_{X}(m))$ is
locally free on $T$. The family $E_{T}$ induces a map $\Delta: T\longrightarrow Pic^{d}(X)$, sending $t\in T$ to
$\det E_{t}$. We can cover $T$ with small open sets $T_{i}$ such that for each $i$ we can find an isomorphism
$$\beta_{T_{i}}:\det E_{T_{i}}\longrightarrow
\overline{\Delta_{i}}^{\ast}\mathcal{P}\; ,$$ where $\mathcal{P}$
is the Poincare bundle defined in (\ref{Poincare}), and a
trivialization
$$g_{T_{i}}:V\otimes \mathcal{O}_{T_{i}}\longrightarrow
\mathcal{V}|_{T_{i}}\; ,$$ where note that $\mathcal{V}|_{t\in
T_{i}}\simeq E_{t\in T_{i}}(m)$ and $H^{0}(\mathcal{V}|_{t\in
T_{i}})\simeq H^{0}(E|_{t\in T_{i}}(m))\simeq V$. Using this
trivialization we obtain a family of quotients parametrized by
$T_{i}$,
$$q_{T_{i}}:V\otimes
\pi_{X}^{\ast}\mathcal{O}_{X}(-m)\twoheadrightarrow E_{T_{i}}\;
,$$ giving a map $T_{i}\longrightarrow \mathcal{H}$. And, using
the quotient $q_{T_{i}}$ and the isomorphism $\beta_{T_{i}}$, we
have another family of quotients parametrized by $T_{i}$,
$$(V^{\otimes s})^{\oplus c}\otimes \big(
\pi_{T_{i},\ast}(\overline{\Delta_{i}}^{\ast}\mathcal{P}^{\otimes b}\otimes
\overline{u_{T_{i}}}^{\ast}\mathcal{D}\otimes \pi_{X}^{\ast}\mathcal{O}_{X}(sm))\big)^{\vee}\twoheadrightarrow
N_{T_{i}}\; ,$$ giving an element of $\mathcal{W}$ for each $t\in T_{i}$. Then, using the representability
properties of $\mathcal{H}$ and $\mathcal{W}$, we obtain a morphism $T_{i}\longrightarrow \mathcal{H}\times
\mathcal{W}$. By Lemma \ref{schemefactor}, this morphism factors through $Z'$ and its image is in
$Z^{ss}$, because a $\delta$-semistable
tensor gives a GIT-semistable point (c.f. Theorem \ref{GIT-delta}). Compose with the geometric quotient to $\mathfrak{M}_{\delta}$ to obtain maps
$$\hat{f}_{i}:T_{i}\overset{f_{i}}{\longrightarrow}
Z^{ss}\longrightarrow \mathfrak{M}_{\delta}\; .$$ Note that the
morphism $f_{i}$ is independent of the choice of isomorphism
$\beta_{T_{i}}$, because of the universal property of the Poincare
bundle $\mathcal{P}$. A different choice of isomorphism
$g_{T_{i}}$ will change $f_{i}$ to $h_{i}\cdot f_{i}$, where
$h_{i}:T_{i}\longrightarrow GL(V)$, then $\hat{f}_{i}$ is
independent of the choice of $g_{T_{i}}$. Glue the morphisms
$\hat{f}_{i}$ to give a morphism
$$\hat{f}:T\longrightarrow \mathfrak{M}_{\delta}\; ,$$
and hence we have a natural transformation from the moduli functor
to the functor of points of $\mathfrak{M}_{\delta}$,
$$\mathcal{M}_{\delta}\longrightarrow \underline{\mathfrak{M}_{\delta}}\; .$$

Recall that there is a tautological family
(\ref{tautologicalfamily}) of tensors parametrized by $Z'$. By
restriction to $Z^{ss}$, we obtain a tautological family of
$\delta$-semistable tensors parametrized by $Z^{ss}$. If $Y$ is
another scheme with a natural transformation
$\mathcal{M}_{\delta}\longrightarrow \underline{Y}$, then the
tautological family defines a $SL(V)$-invariant morphism
$Z^{ss}\longrightarrow Y$, hence this factors through the quotient
$\mathfrak{M}_{\delta}$. Then, the natural transformation
$\mathcal{M}_{\delta}\longrightarrow \underline{Y}$ factors
through $\underline{\mathfrak{M}_{\delta}}$ and this proves that
$\mathfrak{M}_{\delta}$ corepresents the functor
$\mathcal{M}_{\delta}$.
\end{pr}

\begin{rem}
Note that this is not a fine moduli space because the analog of the uniqueness result of \cite[Lemma 1.6]{HL2}
does not hold in general for tensors.
\end{rem}

Now let us give a criterion for S-equivalence. If $(E,\varphi,u)$
and $(F,\psi,u)$ are two $\delta$-stable tensors then we have seen
that they correspond to the same point in the moduli space if and
only if they are isomorphic. But if they are strictly
$\delta$-semistable (i.e. $\delta$-semistable but not
$\delta$-stable), they can be $S$-equivalent (i.e. they correspond
to the same point in the moduli space), but not isomorphic. Hence,
given a tensor $(E,\varphi,u)$, we will show that there exists a
canonical representative of its S-equivalence class
$(E^{S},\varphi^{S},u)$, such that two tensors $(E,\varphi,u)$ and
$(F,\psi,u)$ will be S-equivalent if and only if
$(E^{S},\varphi^{S},u)$ $(F^{S},\psi^{S},u)$ are isomorphic.

Let $(E,\varphi,u)$ be a strictly $\delta$-semistable. Then, by
Proposition \ref{GIT-delta}, the corresponding point
$(\overline{Q},[\Phi])$ is strictly GIT semistable, by Theorem
\ref{GITstab} there exists at least one $1$-parameter subgroup
$\Gamma$ of $SL(V)$ with $\mu((\overline{Q},[\Phi]),\Gamma)=0$
and, by Proposition \ref{S-equivalence}, $\Gamma$ corresponds to a
weighted filtration $(E_{\bullet},n_{\bullet})$ with
$$\sum_{i=1}^{t}n_{i}\big((rP_{E_{i}}-\rk
E_{i}P)\big)+\delta\mu(\varphi,E_{\bullet},n_{\bullet})=0\; ,$$
which we will call an \emph{admissible} weighted filtration for a
strictly $\delta$-semistable tensor.

Let $\mathcal{I}_{0}$ be the set of pairs $(k,I)$ where $k$ is an
integer with $1\leq k \leq c$, and $I=(i_{1},...,i_{s})$ is a
multi-index with $1\leq i_{j}\leq t+1$, such that the restriction
of $\varphi$
$$\varphi_{k,I}:\overbrace{0\oplus \cdots\oplus 0}^{k-1}\oplus
(E_{i_{1}}\otimes \cdots \otimes E_{i_{s}})\oplus
\overbrace{0\oplus \cdots\oplus 0}^{c-k}\longrightarrow (\det
E)^{\otimes b}\otimes D_{u}$$ is nonzero and
$$\gamma_{r_{i_{1}}}+\cdots +\gamma_{r_{i_{s}}}=\mu(\varphi,E_{\bullet},n_{\bullet})\; .$$

Note that, if $(k,I)\in \mathcal{I}_{0}$ and
$I'=(i'_{1},...,i'_{s})$ is a multi-index with $I'\neq I$ and
$i'_{j}\leq i_{j}$ for all $j$, then $\varphi_{k,I'}=0$, by
definition of $\mu(\varphi,E_{\bullet},n_{\bullet})$. Hence, if
$(k,I)\in\mathcal{I}_{0}$, the restriction $\varphi_{k,I}$ defines
a homomorphism in the quotient
$$\varphi'_{k,I}:\overbrace{0\oplus \cdots\oplus 0}^{k-1}\oplus
(E'_{i_{1}}\otimes \cdots \otimes E'_{i_{s}})\oplus \overbrace{0\oplus \cdots\oplus 0}^{c-k}\longrightarrow
(\det E)^{\otimes b}\otimes D_{u}\; ,$$ where $E'_{i}=E_{i}/E_{i+1}$. If $(k,I)$ is not in $\mathcal{I}_{0}$,
then define $\varphi'_{k,I}:=0$. Therefore, we can define
$$\big(E'=E'_{1}\oplus \ldots \oplus E'_{t+1}\text{   ,
}\varphi'=\bigoplus_{(k,I)}\varphi'_{k,I}\big)$$ in which we are
using that $\det E\simeq \det E'$, hence
$(E',\varphi',u)$ is well-defined up to isomorphism and we call it
the \emph{admissible deformation} associated to the admissible
filtration $(E_{\bullet},n_{\bullet})$ of $E$. Observe that this notion depends on the weighted filtration chosen.

\begin{prop}\cite[Proposition 4.1]{GS1}
\label{Sequivalence} The tensor $(E',\varphi',u)$ is strictly $\delta$-semistable and it is $S$-equivalent to
$(E,\varphi,u)$. If we repeat this process, after a finite number of iterations, the process will stop, i.e. we
will obtain tensors isomorphic to each other. We call this tensor $(E^{S},\varphi^{S},u)$ and it verifies
\begin{enumerate}
    \item The isomorphism class of $(E^{S},\varphi^{S},u)$ is
    independent of the choices made, i.e. the weighted filtrations
    chosen.
    \item Two tensors $(E,\varphi,u)$ and $(F,\psi,u)$ are
    S-equivalent if and only if $(E^{S},\varphi^{S},u)$ is
    isomorphic to $(F^{S},\psi^{S},u)$.
\end{enumerate}
\end{prop}
\begin{pr}
First, we recall some observations about GIT quotients. Let $Z$ be a projective variety with an action of a
group $G$ linearized on an ample line bundle $\mathcal{O}_{Z}(1)$. Two points in the open subset $Z^{ss}$ of
semistable points are GIT equivalent, or they give the same point in the moduli space, if the closures (in $Z^{ss}$) of their orbits do intersect
(c.f. Remark \ref{Seq}). Let $z\in Z^{ss}$ and let $B(z)$ be the unique closed orbit in the closure $\overline{G\cdot z}$ in $Z^{ss}$ of its orbit $G\cdot z$.
If $z$ is not in $B(z)$, there exists a $1$-parameter subgroup $\Gamma$ such that the limit
$z_{0}=\lim_{t\rightarrow 0}\Gamma(t)\cdot z$ is in $\overline{G\cdot z}\backslash G\cdot z$ (for instance, we
can take the $1$-parameter subgroup given by \cite[Lemma 1.25]{Si1}). Note that it is $\mu(z,\Gamma)=0$
(by semistability of $z$, $\mu(z,\Gamma)\leq 0$, and if it were negative, then we would have $G\cdot z=B(z)$).
Conversely, if $\Gamma$ is a $1$-parameter subgroup with $\mu(z,\Gamma)=0$, then the limit, $z_{0}$, is GIT
semistable (\cite[Proposition 2.14]{GS0}). Observe that $G\cdot z_{0}\subset \overline{G\cdot z}\backslash
G\cdot z$, therefore, $\dim G\cdot z_{0}<\dim G\cdot z$. Repeating the process with $z_{0}$ instead of $z$, we
get a sequence of points which stops after a finite number of steps, and gives $\tilde{z}\in B(z)$. Two points $z_{1}$ and $z_{2}$ in
$Z^{ss}$ are S-equivalent if and only if $B(z_{1})=B(z_{2})$.

Let $(E,\varphi,u)$ be a $\delta$-semistable tensor with an isomorphism $f:V\simeq H^{0}(E(m))$, and let
$z=(\overline{Q},[\Phi])\in Z$ be the corresponding GIT semistable point. Recall from Proposition
\ref{S-equivalence} the bijection between $1$-parameter subgroups $\Gamma$ of $SL(V)$ with
$\mu(z,\Gamma)=0$ on the one hand, and weighted filtrations $(E_{\bullet},n_{\bullet})$ of $E$ with
$$\big(\sum_{i=1}^{t}n_{i}(rP_{E_{i}}-\rk
E_{i}P)\big)+\delta\mu(\varphi,E_{\bullet},n_{\bullet})\big)=0$$ together with a splitting of the filtration
$H^{0}(E_{\bullet}(m))$ of $V=H^{0}(E(m))$ on the other hand. A $1$-parameter subgroup $\Gamma$ acting on $z$ defines a
morphism $\mathbb{C}^{\ast}\rightarrow Z$ which extends to
$$h:\mathbb{C}\longrightarrow Z\; ,$$
with $h(t)=\Gamma(t)\cdot z$ for $t\neq 0$ and whose limit is $h(0)=\lim_{t\rightarrow 0}\Gamma(t)\cdot
z=z_{0}$.

If we pull back by $h$ the universal family parametrized by $Z$, we obtain another family
$(q_{T},E_{T},\varphi_{T},u)$, where
$$E_{T}=\bigoplus_{i=1}^{t+1}E^{i}\otimes t^{\gamma_{r_{i}}}\subset
E\otimes_{\mathbb{C}}t^{-N}\mathbb{C}[t]\subset E\otimes_{\mathbb{C}}\mathbb{C}[t,t^{-1}]\; ,$$ where recall
that $t^{\Gamma_{r_{i}}}$ acts on each $E^{i}$. We get the morphisms,
$$\xymatrix{
q_{T}:V\otimes \mathcal{O}_{X}(-m)\otimes
\mathbb{C}[t]\ar[r]^{\xi} & \oplus_{i}V^{i}\otimes
\mathcal{O}_{X}(-m)\otimes t^{\Gamma_{r_{i}}}\ar[r] & E_{T}\\
v^{i}\otimes 1\ar@{|->}[r] & v^{i}\otimes
t^{\Gamma_{r_{i}}}\ar@{|->}[r] & q(v^{i})\otimes
t^{\Gamma_{r_{i}}}}$$ and
$$\xymatrix{\varphi_{T}:(E_{T}^{\otimes s})^{\oplus
c}\ar[r] & (\det E_{T})^{\otimes b}\otimes
\overline{u_{T}}^{\ast}\mathcal{D}\otimes \pi_{T}^{\ast}N\\
(\underbrace{0,\ldots,0}_{k-1},w_{i_{1}}t^{\Gamma_{r_{i_{1}}}}\cdot
\ldots \cdot
w_{i_{s}}t^{\Gamma_{r_{i_{s}}}},\underbrace{0,\ldots,0}_{c-k})\ar@{|->}[r]
& \varphi(\underbrace{0,\ldots,0}_{k-1},w_{i_{1}}\cdot \ldots
\cdot w_{i_{s}},\underbrace{0,\ldots,0}_{c-k})\otimes
t^{\Gamma_{r_{i_{1}}}+\cdots +\Gamma_{r_{i_{s}}}}}$$

Then, $(q_{t},E_{t},\varphi_{t},u)$ corresponds to $h(t)$ (in
particular, if $t\neq 0$, then $(E_{t},\varphi_{t},u)$ is
canonically isomorphic to $(E,\varphi,u)$), and
$(E_{0},\varphi_{0},u)$ is the admissible deformation associated
to $(E_{\bullet},n_{\bullet})$. Note that $1$ follows from the
universality of the construction and $2$ follows from the previous
discussion.
\end{pr}

\section{The Harder-Narasimhan filtration}
\label{HNsection}

In section \ref{exampletensors} we have constructed a moduli space for tensors, by restricting the class of
objects that we classify, the $\delta$-semistable tensors. This is the usual situation when constructing a
moduli space, to restrict the original moduli problem by introducing a stability condition.

In some sense, the construction of a moduli space answers the classification problem for the class of the
semistable objects. For the rest, the unstable objects, there is a main tool in algebraic geometry, called the
Harder-Narasimhan filtration, to study them.

We will recall the original Harder-Narasimhan filtration for vector bundles and torsion free sheaves in this
section. At the end, we will discuss the abstract generalization of this notion for an abelian category.

\subsection{Harder-Narasimhan filtration for sheaves}

We consider, first, the case of vector bundles over curves. Let $E$ be a holomorphic vector bundle over a smooth
projective complex curve $X$. Let
$$\mu(E):=\frac{\deg E}{\rk E}$$ be its \emph{slope}.

\begin{dfn}
\label{slopestability} $E$ is \emph{semistable} if for every proper holomorphic subbundle $F\subset E$, it is
$\mu(F)\leq\mu(E)$. If the inequality is strict for every proper subbundle we say that $E$ is \emph{stable}. A
holomorphic vector bundle is \emph{unstable} if it exists a proper subbundle verifying $\mu(F)>\mu(E)$.
\end{dfn}

The holomorphic vector bundles of fixed rank $r$ and degree $d$ over an algebraic curve $X$ of genus $g$ were
studied by Grothendieck for $g=0$ and Atiyah for $g=1$. With the previous definition of stability, Narasimhan
and Seshadri constructed a moduli space for holomorphic bundles over algebraic curves of genus $g$. They did it,
first for bundles of degree $0$ (c.f. \cite{NS1}) and later in general (c.f. \cite{NS2}), through
representations of the fundamental group, and putting in correspondence the semistable bundles with the
semistable points, in the sense of GIT, defined by Mumford. Then, Gieseker (c.f. \cite{Gi1}) gave an algebraic
construction for a moduli space of torsion free sheaves over an algebraic surface and Maruyama (c.f. \cite{Ma1})
extended the construction to higher dimensional varieties.

As we announced in Section \ref{modulis}, we impose a condition on the objects we are trying to classify, the
notion of stability, and restrict our classification problem to the semistable objects. What we have to do then
is, in all the moduli problems which arise as GIT quotients of a space by the action of a group, to show that
the semistable objects, with respect to the definition of stability we give from the beginning, correspond to
the semistable orbits in the sense of Geometric Invariant Theory. Therefore, we obtain a moduli space for the
class of semistable objects which is a good quotient (c.f. Definition \ref{goodquotient}) where each point
corresponds to an $S$-equivalence class of semistable objects (c.f. Remark \ref{Seq}). If we restrict the moduli
problem to the stable objects, we get a geometrical quotient (c.f. Definition \ref{geometricquotient}) which is
an orbit space where, indeed, each point corresponds to a stable orbit and represents an isomorphism class of
stable objects.

Harder and Narasimhan prove the existence of a canonical filtration for a holomorphic vector bundle over a
smooth algebraic curve (c.f. \cite{HN}). The construction of the filtration is based on the existence of a
unique subbundle which maximally contradicts the stability in Definition \ref{slopestability} (in
\cite[Proposition 1.3.4]{HN} this subbundle is called SCSS, a subbundle which \textit{``strongly contradicts the
semistability''}), taking the quotient by this subbundle and repeating the process by recursion (c.f.
\cite[Lemma 1.3.7]{HN}). In that article, Harder and Narasimhan use the filtration to decompose an unstable
vector bundle in blocks and calculate some numbers in relation with the cohomology groups of the moduli space.
Within the years, the so-called \emph{Harder-Narasimhan filtration} has been proved to be extremely useful in
the study of properties of moduli spaces in algebraic geometry.

Let us show how to construct the Harder-Narasimhan filtration in
an easy case, where $E$ is a holomorphic vector bundle of rank $r$
and degree $d$ over a smooth projective complex curve $X$ of genus
$g$.

Suppose that $E$ is unstable and let $\mu(E)=\frac{d}{r}$ be its slope. By definition of stability there are
subbundles $E^{\prime }$ of rank $r^{\prime }<r$ and degree $d$, $0\subsetneqq
E^{\prime}\subsetneqq E$ such that $\mu(E^{\prime})=\frac{d^{\prime }}{r^{\prime }}>\mu(E) =\frac{d}{r}$. We
choose $E_{1}$ with $\mu (E_{1})>\mu(E) $ to be maximal and of maximal rank among those of maximal slope (i.e.
if $\exists E_{1}^{\prime }$ with $\mu(E_{1}^{\prime })=\mu(E _{1})$, then $E_{1}^{\prime }\subseteq E_{1}$). We
will call $E_{1}$ \emph{the maximal destabilizing subbundle of $E$}. Now we consider the subbundle $F=E/E_{1}$.
If it is semistable we are done, and the Harder-Narasimhan filtration is $0\subsetneqq E_{1}\subsetneqq E$.
If not, in analogy with the previous case, there exists $0\subsetneqq F_{1}\subsetneqq F$, of maximal slope
and of maximal rank among those of maximal slope, hence we have
$$\begin{array}{ccccccc}
        0 & \subsetneqq & E_{1} & \subsetneqq & E_{2} & \subsetneqq & E \\
        & & \downarrow & & \downarrow & & \downarrow\\
        & & 0 & \subsetneqq & E_{2}/E_{1}=F_{1} & \subsetneqq & E/E_{1}=F
        \end{array}$$
Call $r_{1},r_{2},d_{1},d_{2}$ the ranks and degrees of $E_{1}$ and $E_{2}$ respectively. These two properties
hold:
\begin{itemize}
    \item The quotient $E_{2}/E_{1}$ is semistable. Indeed, if $E_{2}/E_{1}=F_{1}$ were not semistable, there would
    exists $0 \subsetneqq F_{2}\subsetneqq F_{1}$ with
$\mu(F_{2})>\mu(F_{1})$, contradicting the choice of $F_{1}$.
    \item It is $\mu(E_{1}/0)=\mu(E_{1})>\mu(E_{2}/E_{1})$, because if we had
$\mu(E_{1})\leq \mu(E_{2}/E_{1})\Longleftrightarrow \frac{d_{1}}{r_{1}}\leq
\frac{d_{2}-d_{1}}{r_{2}-r_{1}}\Longleftrightarrow d_{1}r_{2}-d_{1}r_{1} \leq
d_{2}r_{1}-d_{1}r_{1}\Longleftrightarrow \frac{d_{1}}{r_{1}}\leq \frac{d_{2}}{r_{2}}\Longleftrightarrow
\mu(E_{1})\leq \mu(E_{2})$, and we have chosen $E_{1}$ of maximal slope among the subbundles of $E$ and
$E_{1}\subsetneqq E_{2}$.
\end{itemize}

Repeating the process, if the quotient $G=E/E_{2}$ is not semistable, we can choose $0\subsetneqq
G_{1}\subsetneqq G$ with maximal slope and rank, and we obtain
$$\begin{array}{ccccccccc}
      0 & \subsetneqq & E_{1} & \subsetneqq & E_{2} & \subsetneqq & E_{3} & \subsetneqq & E\\
       & & \downarrow & & \downarrow & & \downarrow & & \downarrow \\
       & & 0 & \subsetneqq & E_{2}/E_{1}=F_{1} & \subsetneqq & E_{3}/E_{1}=F_{2} & \subsetneqq
       & E/E_{1}=F\\
      & & & & \downarrow & & \downarrow & & \downarrow \\
     & & & & 0=F_{1}/E_{2} & \subsetneqq & E_{3}/E_{2}=G_{1} & \subsetneqq & E/E_{2}=G
    \end{array}\; .$$ By analogy, we get that $F_{2}/F_{1}=\frac{E_{3}/E_{1}}{E_{2}/E_{1}}\simeq E_{3}/E_{2}$ is
    semistable and
$$\mu(F_{1})>\mu(G_{1})\Longleftrightarrow \mu(E_{2}/E_{1})>\mu(E_{3}/E_{2})\; .$$

By iterating until we get a semistable quotient $E/E_{t}$, we obtain the Harder-Narasimhan filtration:
$$0\subset E_{1}\subset E_{2}\subset \dots \subset E_{t}\subset E_{t+1}=E$$
which verifies
\begin{itemize}
    \item $\mu(E^{1})>\mu(E^{2})>\mu(E^{3})>...>\mu(E^{t})>\mu(E^{t+1})=\mu $, where $\mu(E^{i})=
    \frac{\deg E^{i}}{\rk E^{i}}$
    \item $E^{i}:=E_{i}/E_{i-1}$ is semistable, $\forall i\in \left\{1,...,t+1\right\} $ where $E_{0}=0 $
\end{itemize}
And the process has to stop by finiteness of the rank of $E$.

Therefore, note that we can exhibit unstable vector bundles as extensions of semistable ones in this way. Given an
unstable vector bundle we have its Harder-Narasimhan filtration
$$0\subset E_{1}\subset E_{2}\subset \dots \subset E_{t}\subset E_{t+1}=E\; .$$ This breaks into short exact sequences
$$\begin{array}{ccccccccc}
0 & \rightarrow & \underset{semistable}{E_{1}} & \rightarrow & E_{2} & \rightarrow & \underset{semistable}{E_{2}/E_{1}} & \rightarrow & 0 \\
0 & \rightarrow & E_{2} & \rightarrow & E_{3} & \rightarrow & \underset{semistable}{E_{3}/E_{2}} & \rightarrow & 0 \\
 & &  & & \dots & & & & \\
0 & \rightarrow & E_{t} & \rightarrow & E & \rightarrow &
\underset{semistable}{E/E_{t}} & \rightarrow & 0
\end{array}\; ,$$
where vector bundles on the right are semistable. Using the
Harder-Narasimhan filtration we can think of semistable bundles as
building blocks for holomorphic vector bundles.

Now we give the definition and the proof of the existence and uniqueness of the Harder-Narasimhan filtration for
torsion free sheaves over smooth projective varieties.

Let $X$ be a smooth projective variety and fix an ample line
bundle $\mathcal{O}_{X}(1)$. For every coherent sheaf over $X$,
$E$, let $P_{E}$ its Hilbert polynomial with respect to
$\mathcal{O}_{X}(1)$, i.e $P_{E}(m)=\chi(E\otimes
\mathcal{O}_{X}(m))$. If $P$ and $Q$ are polynomials, we write
$P\leq Q$ if $P(m)\leq Q(m)$ for $m\gg 0$.
\begin{dfn}\cite[Definition 0.1]{Gi1}
\label{polstability} Let $E$ be a torsion free sheaf over $X$. We
say that $E$ is \emph{semistable} if for all proper subsheaves
$F\subset E$, it is
$$\frac{P_{F}}{\rk F}\leq\frac{P_{E}}{\rk E}\; .$$ If the inequality is strict for every proper subsheaf we say that $E$ is
\emph{stable}.
\end{dfn}

Note that, if $E$ is a holomorphic vector bundle of rank $r$ and degree $d$ over an algebraic curve $X$ of
genus $g$, the Hilbert polynomial of $E$ is $P_{E}(m)=rm+d+r(1-g)$ and Definition \ref{polstability} is
equivalent to Definition \ref{slopestability}. We often refer to Definition \ref{polstability} as \emph{Gieseker or Maruyama stability},
whereas Definition \ref{slopestability} is usually
called \emph{Mumford, Takemoto or slope stability}, both definitions coinciding for curves.

\begin{dfn}
\label{HNdef} Let $E$ be a torsion free sheaf of rank $r$ over a
smooth projective algebraic variety $X$. A \emph{Harder-Narasimhan
filtration} for $E$ is a sequence
$$0\subset E_{1}\subset E_{2}\subset \ldots \subset E_{t}\subset E_{t+1}=E$$
verifying
\begin{itemize}
   \item The Hilbert polynomials verify
   $$\frac{P_{E^{1}}}{\rk E^{1}}>\frac{P_{E^{2}}}{\rk E^{2}}>\ldots>\frac{P_{E^{t+1}}}{\rk E^{t+1}}$$
   \item Every $E^{i}$ is semistable
 \end{itemize}
where
$E^{i}:=E_{i}/E_{i-1}$.
\end{dfn}

\begin{rem}
\label{GHN-MHN} Note that the Harder-Narasimhan filtration with
Gieseker stability is a refinement of the one with Mumford
stability, with the inequalities holding between the Hilbert
polynomials in one case, or their leading coefficients in the other.
\end{rem}

A sheaf $E$ is \emph{pure of dimension $n$} if its support has dimension $n$ and it has no subsheaves supported
on a locus of lower dimension.

\begin{thm}\cite[Proposition 1.3.9]{HN}, \cite[Theorem 1.3.4]{HL3}
\label{HNunique} Every pure sheaf $E$ of dimension $n$ over a smooth projective variety $X$ has a unique
Harder-Narasimhan filtration.
\end{thm}
\begin{lem}
\label{maxdes} Let $E$ be a torsion free sheaf. Then, there exists a subsheaf $F\subset E$ such that for all
subsheaves $G\subset E$, one has $\frac{P_{F}}{\rk F}\geq \frac{P_{G}}{\rk G}$ and, in case of equality
$G\subset F$. Moreover, $F$ is uniquely determined and $F$ is semistable, called the \emph{maximal destabilizing
subsheaf of $E$}.
\end{lem}
\begin{pr}
Note that $F$ to be semistable and uniquely determined follows from the first property.

Define an order relation on the set of subsheaves of $E$
by $F_{1}\leq F_{2}$ if and
only if $F_{1}\subset F_{2}$ and $\frac{P_{F_{1}}}{\rk F_{1}}\leq
\frac{P_{F_{2}}}{\rk F_{2}}$. Every ascending chain is bounded by
$E$, then by Zorn's Lemma, for every subsheaf
$F$ there exists $F\subset F'\subset E$ such that
$F'$ is maximal with respect to $\leq$. Let
$F$ be $\leq$-maximal with $F$ of minimal rank
among all maximal subsheaves of $E$. Let us show that $F$ has the required properties.

Suppose there exists $G\subset E$ with $\frac{P_{G}}{\rk G}\geq
\frac{P_{F}}{\rk F}$. First, note that we can assume $G\subset
F$ by replacing $G$ by $F\cap G$. Suppose that $G\nsubseteq F$, then
$F$ is a proper subsheaf of $F+G$ and hence $\frac{P_{F}}{\rk
F}>\frac{P_{F+G}}{\rk F+G}$, by definition of $F$. From the sequence
$$0\rightarrow F\cap G\rightarrow F\oplus G\rightarrow F+G\rightarrow 0$$ we get
$$P_{F}+P_{G}=P_{F\oplus G}=P_{F\cap G}+P_{F+G}$$
and
$$\rk F+\rk G=\rk (F\oplus G)=\rk (F\cap G)+\rk (F+G)\; .$$
Calculating we have
$$\rk (F\cap G)(\frac{P_{G}}{\rk G}-\frac{P_{F\cap G}}{\rk(F\cap G)})=$$
$$\rk(F+G)(\frac{P_{F+G}}{\rk (F+G)}-
\frac{P_{F}}{\rk F})+(\rk G-\rk(F\cap G))(\frac{P_{F}}{\rk F}-\frac{P_{G}}{\rk G})\; .$$
Then, together with the two inequalities $\frac{P_{F}}{\rk F}\leq
\frac{P_{G}}{\rk G}$ and $\frac{P_{F}}{\rk F}> \frac{P_{F+G}}{\rk (F+G)}$ we
obtain
$$\frac{P_{G}}{\rk G}-\frac{P_{F\cap G}}{\rk (F\cap G)}<0$$
and hence
$$\frac{P_{F}}{\rk F}<\frac{P_{F\cap G}}{\rk (F\cap G)}\; ,$$
which proves the assert that we can suppose $G\subset F$.

Now, let $G\subset F$ with $\frac{P_{G}}{\rk G}>\frac{P_{F}}{\rk
F}$ such that $G$ is $\leq$-maximal in $F$. Then let $G'\geq G$,
$\leq$-maximal in $E$. We obtain the inequalities $\frac{P_ {F}}{\rk
F}<\frac{P_{G}}{\rk G}\leq \frac{P_{G'}}{\rk G'}$. Because of the maximality
of $G'$ and $F$ it is $G'\nsubseteq F$, because otherwise
$\rk(G')<\rk(F)$ but $\rk(F)$ is minimal by hypothesis. Therefore,
$F$ is a proper subsheaf of $F+G'$ and $\frac{P_{F}}{\rk
F}>\frac{P_{F+G'}}{\rk (F+G')}$. The previous inequalities
$\frac{P_{F}}{\rk F}<\frac{P_{G'}}{\rk G'}$ and $\frac{P_{F}}{\rk
F}>\frac{P_{F+G'}}{\rk (F+G')}$ give
$$\frac{P_{F\cap G'}}{\rk (F\cap G')}>\frac{P_{G'}}{\rk G'}\geq\frac{P_{G}}{\rk G}\; .$$
Given that $G\subset F\cap G'\subset F$, we get a contradiction with the hypothesis on $G$.
\end{pr}

\begin{pr}[Proof of the Theorem]
Lemma \ref{maxdes} shows the existence of a
Harder-Narasimhan filtration for $E$. Let $E_{1}$ the maximal
destabilizing subsheaf and suppose that the corresponding quotient
$E/E_{1}$ has a Harder-Narasimhan filtration,
$$0\subset G_{0}\subset G_{1}\subset\ldots\subset G_{t-1}=E/E_{1}\; ,$$
by induction hypothesis. We define $E_{i+1}$ as the pre-image of $G_{i}$ and it is $\frac{P_ {E_{1}}}{\rk
E_{1}}>\frac{P_{E_{2}/E_{1}}}{\rk E_{2}/E_{1}}$ because, if not, we get $\frac{P_{E_{1}}}{\rk
E_{1}}\leq\frac{P_{E_{2}}}{\rk E_{2}}$, which contradicts the maximality of $E_{1}$.

Next we prove the uniqueness. Let $E_{\bullet}$ and $E'_{\bullet}$ be two Harder-Narasimhan filtrations of the
same sheaf $E$. We consider, without loss of generality, $\frac{P_{E'_{1}}}{\rk E'_{1}}\geq \frac{P_{E_{1}}}{\rk
E_{1}}$. Let $j$ be the minimal index verifying $E'_{1}\subset E_{j}$. The composition
$$E'_{1}\rightarrow E_{j}\rightarrow E_{j}/E_{j-1}$$ is a
non-trivial homomorphism of semistable sheaves which implies
$$\frac{P_{E_{j}/E_{j-1}}}{\rk E_{j}/E_{j-1}}\geq \frac{P_{E'_{1}}}{\rk E'_{1}}\geq
\frac{P_{E_{1}}}{\rk E_{1}}\geq \frac{P_{E_{j}/E_{j-1}}}{\rk E_{j}/E_{j-1}}$$ where first inequality comes from
the fact that, if there exists a non-trivial homomorphism between semistable sheaves, then the Hilbert
polynomial of the target is greater or equal than the one of the first sheaf. Therefore, equality holds everywhere,
and this implies that the index $j$ is equal to $1$, so that $E'_{1}\subset E_{1}$. Then, by semistability of $E_{1}$, it is $\frac{P_{E'_{1}}}{\rk
E'_{1}}\leq \frac{P_{E_{1}}}{\rk E_{1}}$, and we can repeat the argument interchanging the roles of $E_{1}$ and
$E'_{1}$ to show that $E_{1}=E'_{1}$. By induction we can assume that uniqueness holds for the Harder-Narasimhan
filtration of $E/E_{1}$ to show that $E'_{i}/E_{1}=E_{i}/E_{1}$, which completes the proof.
\end{pr}

\begin{rem}
If a torsion free sheaf $E$ is already semistable, we can still talk about its Harder-Narasimhan filtration
which is the trivial filtration $0\subset E$.
\end{rem}

Next we show how the Harder-Narasimhan filtration looks like in the easiest case, for an unstable vector bundle
over $X=\mathbb{P}_{\mathbb{C}}^{1}$.
\begin{ex}
\label{HNP1} Let $X=\mathbb{P}_{\mathbb{C}}^{1}$. We know, by a theorem of Grothendieck, that a vector bundle
$E$ over $\mathbb{P}_{\mathbb{C}}^{1}$ splits on line bundles
$$E=\mathcal{O}_{\mathbb{P}_{\mathbb{C}}^{1}}(a_{1})\oplus \mathcal{O}_{\mathbb{P}_{\mathbb{C}}^{1}}(a_{1})\oplus
\cdots \oplus\mathcal{O}_{\mathbb{P}_{ \mathbb{C}}^{1}}(a_{1}) \oplus \mathcal{O}_{\mathbb{P}_{\mathbb{C}
}^{1}}(a_{2})\oplus \cdots \oplus \mathcal{O}_{\mathbb{P}_{\mathbb{C} }^{1}}(a_{2})\oplus
\mathcal{O}_{\mathbb{P}_{\mathbb{C}}^{1}}(a_{3}) \oplus \cdots\oplus
\mathcal{O}_{\mathbb{P}_{\mathbb{C}}^{1}}(a_{s})$$ with $a_{1}>a_{2}>...>a_{s}$, and we call $b_{i}$ the number
of times each line bundle $\mathcal{O}_{\mathbb{P}_{\mathbb{C}}^{1}}(a_{i})$ appears (c.f. \cite[Theorem
1.3.1]{HL3}). Thus, the slope of $E$ is the average of the degrees $a_{i}$ of the line bundles appearing in the
decomposition of $E$,
$$\mu(E)=\frac{\deg E}{rk\text{ }E}=\frac{a_{1}b_{1}+\cdots +a_{s}b_{s}}{b_{1}+\cdots +b_{s}}\; .$$

With the notation of Theorem \ref{HNunique}, it is clear that
$$E^{1}=\mathcal{O}_{\mathbb{P}_{\mathbb{C}}^{1}}(a_{1})\oplus
\mathcal{O}_{\mathbb{P}_{\mathbb{C}}^{1}}(a_{1})\oplus \dots\oplus
\mathcal{O}_{\mathbb{P}_{\mathbb{C}}^{1}}(a_{1})$$ with
$$\mu(E^{1})=\frac{\overbrace{a_{1}+\cdots+a_{1}}^{b_{1}\text{ times }}}{b_{1}}=a_{1}>\mu(E)$$
and
$$F=E/E_{1}=\mathcal{O}_{\mathbb{P}_{\mathbb{C}}^{1}}(a_{2})\oplus\cdots\oplus
\mathcal{O}_{\mathbb{P}_{\mathbb{C}}^{1}}(a_{2})\oplus \mathcal{O}_{\mathbb{P}_{\mathbb{C}}^{1}}(a_{3}) \oplus
\cdots\oplus \mathcal{O}_{\mathbb{P}_{\mathbb{C}}^{1}}(a_{s})\; .$$ Then it is also
$$E^{2}=\mathcal{O}_{\mathbb{P}_{\mathbb{C}}^{1}}(a_{2})\oplus
\cdots\oplus \mathcal{O}_{\mathbb{P} _{\mathbb{C}}^{1}}(a_{2})$$ which lifts to
$$E_{2}=\mathcal{O}_{\mathbb{P}_{\mathbb{C} }^{1}}(a_{1})\oplus \cdots\oplus
\mathcal{O}_{\mathbb{P}_{\mathbb{C}}^{1}}(a_{1})\oplus
\mathcal{O}_{\mathbb{P}_{\mathbb{C}}^{1}}(a_{2})\oplus \cdots\oplus \mathcal{O}_{\mathbb{P}
_{\mathbb{C}}^{1}}(a_{2})\; .$$ Repeating the process we obtain a filtration
$$0\subset E_{1}\subset E_{2}\subset \cdots\subset E_{s-1}\subset E_{s}=E$$
where
$$E_{i}=\mathcal{O}_{\mathbb{P}_{\mathbb{C}}^{1}}(a_{1})\oplus \ldots \oplus
\mathcal{O}_{\mathbb{P}_{ \mathbb{C}}^{1}}(a_{1}) \oplus \ldots \oplus \mathcal{O}_{\mathbb{P}_{\mathbb{C}
}^{1}}(a_{i})\oplus \cdots\oplus \mathcal{O}_{\mathbb{P}_{\mathbb{C} }^{1}}(a_{i})$$ which is the
Harder-Narasimhan filtration. It is clear that each $E^{i}=E_{i}/E_{i-1}$ is semistable and
$\mu(E^{1})>\mu(E^{2})>\mu(E^{3})>...>\mu(E^{s-1})>\mu(E^{s})$.
\end{ex}

\subsection{Harder-Narasimhan filtration in an abelian category}
Finally, we would like to close this section with some comments about stability notions and the concept of the
Harder-Narasimhan filtration in a more general context. Rudakov defines in \cite{Ru} a notion of stability for
objects in an abelian category.

Let $\mathcal{C}$ be an abelian category. We define a \emph{preorder} on the objects making possible to compare
two nonzero objects, i.e. if $A\neq 0$, $B\neq 0$ are objects of $\mathcal{C}$ it is one of the following
$A\prec B$, $A\succ B$ or $A\asymp B$, being possible to have $A\asymp B$ although $A\neq B$.
\begin{dfn}\cite[Definition 1.1]{Ru}
 A \emph{stability structure} on $\mathcal{C}$ is a preorder on $\mathcal{C}$ such that for every
 short exact sequence $0\rightarrow A\rightarrow B\rightarrow C\rightarrow 0$ it happens one
of the following
\begin{itemize}
 \item $A\prec B\Leftrightarrow A\prec C\Leftrightarrow B\prec C$
 \item $A\succ B\Leftrightarrow A\succ C\Leftrightarrow B\succ C$
 \item $A\asymp B\Leftrightarrow A\asymp C\Leftrightarrow B\asymp C$
\end{itemize}
\end{dfn}

Note that this property is satisfied by the category of holomorphic vector bundles over curves, if we associate
to each object $E$ the numerical function given by its slope
$$\mu(E)=\frac{\deg E}{\rk E}$$
and define the preorder as $$E\prec F\Leftrightarrow \mu(E)<\mu(F)\;\;\;\; E\asymp F\Leftrightarrow
\mu(E)=\mu(F)\; .$$ Also, for torsion free sheaves over projective varieties, if we associate to each sheaf $E$
the polynomial function given by $\frac{P_{E}}{\rk E}$, where $P_{E}$ is the Hilbert polynomial of $E$, we have
a stability structure in the corresponding category by defining the preorder with the obvious relations between
the polynomial functions.

\begin{dfn}
An object $A\in \mathcal{C}$ is \emph{semistable} if it is nonzero
and for every nontrivial subobject $B\subset A$, we have
$B\preccurlyeq A$. We say that $A$ is \emph{stable} if we have a
strict inequality for every nontrivial subobject.
\end{dfn}

Let us mention three properties for an abelian category
$\mathcal{C}$ to have, in order to assure that a Harder-Narasimhan
filtration exists for an unstable object in $\mathcal{C}$ (these
properties appear in \cite{Ru}).

\begin{dfn}
An object $A\in \mathcal{C}$ is \emph{quasi-noetherian} if a chain verifying
$$A_{1}\subset A_{2}\subset \ldots \subset A$$
and
$$A_{1}\preccurlyeq A_{2}\preccurlyeq \ldots \preccurlyeq A$$
stabilizes. We say that $\mathcal{C}$ is quasi-noetherian if every $A\in \mathcal{C}$ is.
\end{dfn}

\begin{dfn}
An object $A\in \mathcal{C}$ is \emph{weakly-noetherian} if it is quasi-noetherian and a chain verifying
$$A_{1}\supset A_{2}\supset \ldots \supset A$$
and
$$A_{1}\succcurlyeq A_{2}\succcurlyeq \ldots \succcurlyeq A$$
stabilizes. We say that $\mathcal{C}$ is weakly-noetherian if every $A\in \mathcal{C}$ is.
\end{dfn}

\begin{dfn}
An object $A\in \mathcal{C}$ is \emph{weakly-artinian} if a chain verifying
$$A_{1}\supset A_{2}\supset \ldots \supset A$$
and
$$A_{1}\preccurlyeq A_{2}\preccurlyeq \ldots \preccurlyeq A$$
stabilizes. We say that $\mathcal{C}$ is weakly-artinian if every $A\in \mathcal{C}$ is.
\end{dfn}

\begin{thm}\cite[Theorem 2]{Ru}
\label{HNabelian} Let $\mathcal{C}$ be an abelian category with a
given stability structure, which is weakly-noetherian and
weakly-artinian. For every object $A\in \mathcal{C}$ there exists
a filtration
$$0\subset A_{1}\subset A_{2}\subset \ldots \subset A_{t}\subset A_{t+1}=A$$
such that
\begin{itemize}
   \item $A^{1}\succ A^{2}\succ \ldots \succ A^{t}\succ A^{t+1}$
   \item Every $A^{i}$ is semistable
 \end{itemize}
where $A^{i}:=A_{i}/A_{i-1}$.
\end{thm}

For an object to be quasi-noetherian it is needed to prove the
existence and uniqueness of a maximally destabilizing subobject
(c.f. Lemma \ref{maxdes}) and the stronger weakly-noetherian is
used when lifting the filtration of the quotient $A/A_{1}$, which
exists by hypothesis in the recursion (c.f. Proof of Theorem
\ref{HNunique}). The condition of being
 weakly-artinian assures that the recursive process when constructing the Harder-Narasimhan filtration
 of an object finishes in a finite number of steps, i.e. the Harder-Narasimhan
 filtration is finite.

An object is called \emph{noetherian} if every ascending chain on it stabilizes. Clearly, being noetherian implies being weakly-noetherian.

Note that the category of coherent sheaves over a projective variety is abelian and noetherian, as well as the
category of finite dimensional representations of quivers (which we will see in Chapter 3), hence the existence
of a Harder-Narasimhan filtration in these cases can be seen as a particular case of Theorem \ref{HNabelian}.

\section{Kempf theorem}
\label{kempfsection}

In the previous sections we have studied moduli problems for which we have to impose a stability condition in
order to have a moduli space with good properties. By \textit{rigidifying the data}, we add extra data to the
objects we are classifying, and this leads us to an action of a group in a space which changes the extra data,
for a given object in the moduli problem. Then, we use Geometric Invariant Theory to take the quotient by the
group and obtain a moduli space with the desired properties.

In the example of the construction of a moduli space for tensors (c.f. section \ref{exampletensors}), the extra
data we add to a tensor is the isomorphism between a vector space and a space of global sections of the
(twisted) sheaf of the tensor. Different isomorphisms differ by an element of a general linear group, and this
is the group we take the quotient by, using GIT (c.f. subsection \ref{giesekerembedding}).

In this kind of constructions, one of the main points appears to be the correspondence between semistable
objects and semistable orbits or GIT semistable points (c.f. Theorem \ref{GIT-delta}). Recall that GIT stability
can be checked by $1$-parameter subgroups (c.f. Hilbert-Mumford criterion, Theorem \ref{HMcrit}): a point $x$ is
unstable if there exists any $1$-parameter subgroup $\Gamma$ which makes some numerical function, the so-called
\textit{minimal relevant weight} $\mu(x,\Gamma)$, positive.

The GIT stability criterion exposed in Theorem \ref{HMcrit0} asserts that a point $x$ is GIT unstable if there exists a $1$-parameter subgroup $\Gamma$ such that
    $$\underset{t\rightarrow 0}{\lim \Gamma(t)\cdot \tilde{x}}=0\; ,$$ where $\tilde{x}$ is a point in the affine cone,
lying over $x$. This is, the Hilbert-Mumford criterion says that the fact of $0$ appearing as the limit point in
the orbit of the linearized action of the group $G$ can be checked through $1$-parameter subgroups. Then,
Theorem \ref{HMcrit0}, the numerical criterion, expresses that fact with the positivity of the numerical
function $\mu(x,\Gamma)$.

The function $\mu(x,\Gamma)$ can be thought as a measure of how rapidly we can reach $0$ from a point
$\tilde{x}$ in the affine cone, lying over $x$, through different $1$-parameter subgroups. Let us see this with
an easy example.

\begin{ex}
\label{speeduns} Consider the group $G=SL(3,\mathbb{C})$ and let $\Gamma:\mathbb{C}^{\ast}\rightarrow
SL(3,\mathbb{C})$ be a $1$-parameter subgroup. There exists a basis of $\mathbb{C}^{3}$ where $\Gamma$ takes the
diagonal form
$$\begin{pmatrix}
t^{\Gamma _{1}} & 0 & 0 \\
0 & t^{\Gamma _{1}} & 0 \\
0 & 0 & t^{\Gamma _{3}}
\end{pmatrix}\; ,$$
where we order the exponents as $\Gamma_{1}<\Gamma_{2}<\Gamma_{3}$ and they verify
$\Gamma_{1}+\Gamma_{2}+\Gamma_{3}=0$. Now consider that it acts on $\mathbb{P}_{\mathbb{C}}^{3}$ and let
$x=[0:x_{2}:x_{3}]$ be a point in homogeneous coordinates, $x_{2}\neq 0, x_{3}\neq 0$. Let
$\tilde{x}=(0,x_{2},x_{3})$ be a point in the affine cone lying over $x$. Then, $\lim_{t\rightarrow 0}
\Gamma(t)\cdot \tilde{x}=\lim_{t\rightarrow 0}(0,t^{\Gamma_{2}}\cdot x_{2},t^{\Gamma_{3}}\cdot x_{3})$ and we
say that $\Gamma$ acts on the limit with weight $\Gamma_{2}$, the minimal relevant exponent of the action of
$\Gamma$ over $x$, i.e. $\mu(x,\Gamma)=\Gamma_{2}$. A point is GIT unstable if $0$ can be reached in the closure
of the orbit of the linearized action through $1$-parameter subgroups. And $\lim_{t\rightarrow 0} \Gamma(t)\cdot
\tilde{x}=0$ if and only if $\Gamma_{2}>0$. Hence, $x$ is GIT unstable if there exists any $\Gamma$ with
$\Gamma_{2}>0$.

Observe that the specific value $\Gamma_{2}$ can be thought as a measure of how rapidly we can move from
$\tilde{x}=(0,x_{2},0)$ to $0$. The greater $\Gamma_{2}$ is, the faster $\lim_{t\rightarrow 0} \Gamma(t)\cdot
\tilde{x}$ takes $\tilde{x}$ to $0$. Hence, $\mu(x,\Gamma)$ encodes, in this sense, the \textit{speed of
unstability}.
\end{ex}


A first question which arises is, could we possibly find a $1$-parameter subgroup giving the greatest
\textit{speed of unstability} as in Example \ref{speeduns}?

The answer would be: not yet. Note that if we multiply the exponents appearing in the diagonal of $\Gamma$ by
the same integer, we still obtain a $1$-parameter subgroup of $SL(3,\mathbb{C})$ giving a positive value for
$\mu(x,\Gamma)$, hence it also destabilizes the point $x$. But the value $\mu(x,\Gamma)$ is multiplied by this
integer, hence we cannot yet well define a unique $1$-parameter subgroup $\Gamma$ giving maximum for
$\mu(x,\Gamma)$. We have to introduce a notion of \textit{length} in the set of $1$-parameter subgroups, to be
able to calibrate this kind of features.

Let $G$ be a connected reductive algebraic group over $k$ and let $T$ be a maximal torus. Let $N$ be the
normalizer of $T$ and let $N/T$ be the Weyl group. Let $\Gamma(G)$ be the set of $1$-parameter subgroups of $G$.
For a $k$-point $g\in G$ and $\Gamma\in \Gamma(G)$, define $g\ast \Gamma$ as the $1$-parameter subgroup $g\ast
\Gamma(t)=g\cdot \Gamma(t)\cdot g^{-1}$. We define a notion of \emph{length} for $\Gamma\in \Gamma(G)$ (c.f.
\cite[p. 305]{Ke}).
\begin{dfn}
\label{length}
A \emph{length} is a non-negative real function on $\Gamma(G)$ verifying
\begin{itemize}
 \item If $g\in G$ is $k$-rational, $\|g\ast \Gamma\|=\|\Gamma\|$ for any $\Gamma\in \Gamma(G)$.
 \item For any maximal torus $T$ of $G$, there is a positive definite integral-valued bilinear
 form $(\;,\;)$ on $\Gamma(T)$ such that $(\Gamma,\Gamma)=\|\Gamma\|^{2}$, for any $\Gamma\in
\Gamma(T)$.
\end{itemize}
\end{dfn}
 As it is pointed out in \cite{Ke}, the first property is the invariance of the length by the
 action of the Weyl group of $G$ with respect to $T$. And, given a positive definite
integral-valued bilinear form $(\;,\;)$ on $\Gamma(T)$ invariant by the Weyl group, where $T$ is a maximal
torus, it corresponds to a unique length $\|\cdot\|$ on $\Gamma (G)$, verifying the property
$(\Gamma,\Gamma)=\|\Gamma\|^{2}$, for any $\Gamma\in \Gamma(G)$.

\begin{rem}
\label{killinglengths} If $G$ is simple in characteristic zero all choices of length will be multiples of
the Killing form in the Lie algebra $\mathfrak{g}$ (note that in this case $\Gamma(G)\subset \mathfrak{g}$) and,
in general, for an almost simple group in arbitrary characteristic, all lengths differ also by a scalar (c.f. \cite[p. 305]{Ke}).

However, if $G$ has different simple factors, there are more choices of lengths. We can obtain different lengths by choosing
a linear combination of the Killing forms in each simple factor with positive coefficients. 
\end{rem}

Given a choice of length in $G$, we can define the function appearing on \cite[Theorem 2.2]{Ke}.

\begin{dfn}
 \label{kempffunction}
Let $G$ be a reductive algebraic group over an algebraically closed field of arbitrary characteristic. Let
$G\times X\rightarrow X$ be an action of $G$ on a $k$-scheme $X$. Consider a length in $\Gamma(G)$, as in
Definition \ref{length}. For a point $x\in X$ and a $1$-parameter subgroup $\Gamma\in\Gamma(G)$, let
$\mu(x,\Gamma)$ be the numerical function of the Hilbert-Mumford criterion as in Theorem \ref{HMcrit}. We define
the following function
$$K(x,\Gamma)=\frac{\mu(x,\Gamma)}{\|\Gamma\|}\; .$$
We call this function the \emph{Kempf function}.
\end{dfn}

\begin{rem}
The numerator of the Kempf function is precisely the \textit{speed of unstability} we discussed about in Example
\ref{speeduns}, and the denominator serves for normalizing that quantity with respect to scalar multiples of $\Gamma$.
Then, we will refer to the $1$-parameter subgroup which maximally contradicts the stability condition in the
sense of GIT by talking of that $\Gamma$ which gives maximum for the Kempf function $K(x,\Gamma)$.
\end{rem}

Geometric Invariant Theory contains a study of the dependence of $\mu(x,\Gamma)$ with the $1$-parameter subgroup
$\Gamma$ (c.f. \cite[Section 2.2]{Mu}). It is based on the previous study of a metric space called the
\emph{flag complex}, by Tits, which is the space of $1$-parameter subgroups modulo certain equivalence
relation for which, the values of $\mu(x,\Gamma)$ that we obtain are multiples. Then, a new function can be
defined on this flag complex, whose positivity or negativity coincides with the one of $\mu(x,\Gamma)$, encoding
GIT stability but forgetting about rescaling the minimal relevant weight $\mu(x,\Gamma)$ with multiples of the
$1$-parameter subgroups.

The \textit{conjecture of Mumford-Tits} (as it is stated in the introduction of \cite{Ke}, or \textit{Tits'
center conjecture} in \cite[Appendix 2B]{MFK}, see \cite[p. 64]{Mu}) says that, if a $k$-rational point $x$ is
unstable with respect to the action of $G$, we can find a special $1$-parameter subgroup giving maximum for the
Kempf function $K(x,\Gamma)$. Kempf explores this idea in \cite{Ke} and solves positively the Mumford-Tits
conjecture, finding that there exists an special class of $1$-parameter subgroups which moves most rapidly
toward the origin. Kempf shows it in more generality, using a closed $G$-invariant set $S$, instead of just the
one point set $\{0\}$, to define a point to be $S$-unstable if the closure of its orbit intersects $S$. Kempf
uses this to prove that the Hilbert-Mumford criterion (i.e. the checking of the GIT stability of a point by
$1$-parameter subgroups) holds for actions of algebraic groups over algebraically closed fields of arbitrary characteristic first,
and then, the analogous result for perfect fields. For a correspondence between Kempf and Mumford's GIT
language, see \cite[Appendix 2B]{MFK}.

The precise statement of the Kempf's result is the following:

\begin{thm}\cite[Theorem 2.2]{Ke}
\label{kempftheorem0} Let $G$ be a reductive algebraic group over an algebraically closed field of arbitrary
characteristic. Let $G\times X\rightarrow X$ be an action of $G$ on a $k$-scheme $X$. Let $x\in X$ be a
$k$-point and suppose that $x$ is GIT unstable, i.e. there exists a $1$-parameter subgroup $\Gamma$ such that
$\mu(x,\Gamma)>0$. Define a length in $\Gamma(G)$ as in Definition \ref{length} and consider the Kempf function
$K(x,\Gamma)=\frac{\mu(x,\Gamma)}{\|\Gamma\|}$. Then, the function $K(x,\Gamma)$ achieves a maximum $B$, taken
over all $\Gamma\in \Gamma(G)$  and there exists a parabolic subgroup $P\subset G$ such that in each maximal
torus $T$ conjugated by $P$, there exists a unique $1$-parameter subgroup $\Gamma\in \Gamma(T)$ achieving the
maximal value $K(x,\Gamma)=B$.
\end{thm}

We can say that the Harder-Narasimhan filtration (c.f. section \ref{HNsection}) is the best filtration which destabilizes
an unstable object, with respect to the given definition of stability, among all possible filtrations by subobjects. Its construction (c.f. Theorem \ref{HNunique})
 is based on the existence of a maximally destabilizing subobject (c.f. Lemma \ref{maxdes} in the case of sheaves), which tells us that there is no better choice for
the first element of the filtration. Then, the recursive process of the construction implies that, at every
step, we do the best possible, finding maximally destabilizing subobjects for the quotients which successively
appear. Given an unstable object, we can obtain a maximally destabilizing subobject, and follow the construction
to complete it until the Harder-Narasimhan filtration.

On the other hand, Theorem \ref{kempftheorem0} implies that, whenever we have a GIT unstable point, we can find
a special $1$-parameter subgroup giving maximal unstability in the sense of Geometric Invariant Theory, with respect to maximizing the Kempf function in Definition \ref{kempffunction}.

Then, consider a notion of stability for a category such that there exists a construction of a moduli space
of semistable (or stable) objects. Consider that the construction of the moduli space is given through Geometric
Invariant Theory, by means of rigidifying the data and taking the quotient of a space by a group, to get rid of
the extra data. In that case, we have a correspondence between unstable objects and GIT unstable objects (as in 
Theorem \ref{GIT-delta} in the construction of the moduli of tensors). In some cases, to give a notion of
maximal unstability for an unstable object we have the Harder-Narasimhan filtration. And, to give a notion of
GIT maximal unstability we have the $1$-parameter subgroup given by Kempf in Theorem \ref{kempftheorem0}.
$$ $$
Therefore, the natural question which arises is,

$$ $$

\textbf{Is the Harder-Narasimhan filtration related to the $1$-parameter subgroup given by Kempf?}

$$ $$

The main purpose of this thesis will be to explore this idea by establishing a correspondence between both
notions which answers positively the question in different cases.

\chapter[Correspondence between Kempf and Harder-Narsimhan filtrations]{Correspondence between Kempf and Harder-Narsimhan filtrations}
\chaptermark{Corresp. Kempf and Harder-Narasimhan}
\label{chaptersheaves}

\section{Torsion free sheaves over projective varieties}
\label{kempfsheaves}

In this first section of the chapter, we describe the main case of the correspondence between the $1$-parameter
subgroup giving the GIT maximal unstability in the sense of Kempf (c.f. Theorem \ref{kempftheorem0}) and the
Harder-Narasimhan filtration (c.f. Theorem \ref{HNunique}). The machinery and the ideas described here will
serve, in the remaining sections of the chapter, to prove the analogous result for other other moduli problems.

\vspace{2cm}

Let $X$ be a smooth complex projective variety, and let $\SO_X(1)$ be an ample
line bundle on $X$. If $E$ is a coherent sheaf on $X$, let $P_E$ be its Hilbert
polynomial with respect to $\SO_X(1)$, i.e.,
$P_E(m)=\chi(E\otimes \SO_X(m))$.

We will briefly describe the construction of the moduli space for these objects. This is originally due to
Gieseker for surfaces (c.f. \cite{Gi1}), and it was generalized to higher dimension by Maruyama (c.f.
\cite{Ma1,Ma2}). First, we give Giesekers's definition of stability for torsion free sheaves. Recall that, if
$P$ and $Q$ are polynomials, we write $P\leq Q$ if $P(m)\leq Q(m)$ for $m\gg 0$.


\begin{dfn}
\label{giesekerstab} \cite[Definition 0.1]{Gi1} A torsion free sheaf $E$ on $X$ is called \emph{semistable} if
for all proper subsheaves $F\subset E$, the following inequality between polynomials hold,
$$
\frac{P_F}{\rk F} \leq \frac{P_E}{\rk E} \; .$$ If strict inequality holds for every proper subsheaf, we say
that $E$ is \emph{stable}.
\end{dfn}

To construct the moduli space of torsion free sheaves with fixed
Hilbert polynomial $P$, we choose a suitably large integer $m$ and
consider the Quot scheme parametrizing quotients
\begin{equation}
\label{quot}
V\otimes \SO_X(-m) \too E
\end{equation}
where $V$ is a fixed vector space of dimension $P(m)$ and $E$ is a
sheaf with $P_E=P$.
The Quot scheme has a canonical action by
$\SL(V)$. Gieseker (c.f. \cite{Gi1}) gives a linearization of this action
on a certain ample line bundle, in order to use Geometric
Invariant Theory to take the quotient by the action.
The moduli space of semistable sheaves is obtained
as the GIT quotient.

As we said, at the beginning of section \ref{HNsection}, the construction of a moduli space for semistable
torsion free sheaves solves the classification problem partially. If a sheaf $E$ is not semistable, it is called
\emph{unstable}, and it has a canonical filtration:

\begin{thm}\cite[Proposition 1.3.9]{HN}
\label{HNdefthm}
Given a torsion free sheaf $E$, there exists a unique filtration
$$0\subset E_{1} \subset E_{2} \subset \cdots \subset E_{t} \subset
E_{t+1}=E\; ,$$ which satisfies the following properties, where
$E^{i}:=E_{i}/E_{i-1}$:
 \begin{enumerate}
   \item The Hilbert polynomials verify
   $$\frac{P_{E^{1}}}{\rk E^{1}}>\frac{P_{E^{2}}}{\rk E^{2}}>\ldots>\frac{P_{E^{t+1}}}{\rk E^{t+1}}$$
   \item Every $E^{i}$ is semistable
 \end{enumerate}
This filtration is called the \emph{Harder-Narasimhan filtration} of
$E$
\end{thm}
\begin{pr}
C.f. Theorem \ref{HNunique}.
\end{pr}

In this section we will develop a series of arguments to establish a correspondence between the $1$-parameter
subgroup given by Kempf in Theorem \ref{kempftheorem0} and the Harder-Narasimhan filtration in Theorem
\ref{HNdefthm} to show that both notions do coincide.

\subsection{Moduli space and Kempf theorem}
\label{modulisheaves}

We will recall Gieseker's construction (c.f. \cite{Gi1}) of the moduli space of semistable torsion free sheaves
with fixed Hilbert polynomial $P$ and fixed determinant $\det(E)\cong \Delta$.

Recall that a coherent sheaf is called $m$-regular if $h^i(E(m-i))=0$ for all $i>0$ (c.f. Definition
\ref{regdef} and Lemma \ref{mregularity}). Let $m$ be a suitable large integer, so that $E$ is $m$-regular for
all semistable $E$ (c.f. \cite[Corollary 3.3.1 and Proposition 3.6]{Ma1}). Let $V$ be a vector space of
dimension $p:=P(m)$. Given an isomorphism $V\cong H^0(E(m))$ we obtain a quotient
$$
q:V\otimes \SO_X(-m) \surj E\; ,
$$
hence a homomorphism
$$
Q:\wedge^r V \cong \wedge^r H^0(E(m)) \too H^0(\wedge^r(E(m)))
\cong H^0(\Delta(rm))=:A
$$
and points
$$
Q\in \Hom(\wedge^r V , A) \qquad\overline{Q} \in
\PP(\Hom(\wedge^r V , A) )\; ,
$$
where $Q$ is well defined up to a scalar because the isomorphism $\det(E)\cong \Delta$ is well defined up to a
scalar, and hence $\overline{Q}$ is a well defined point the the projective space. Two different isomorphisms
between $V$ and $H^0(E(m))$ differ by the action of an element of $\GL(V)$, but, since an homothecy does not
change the point $\overline{Q}$, to get rid of the choice of isomorphism it is enough to take the quotient by
the action of $\SL(V)$.

We recall from section \ref{exampletensors} the correspondence between weighted filtrations and $1$-parameter subgroups.
A \emph{weighted filtration} $(V_\bullet,n_\bullet)$
of $V$ is a filtration
\begin{equation}
\label{filtV} 0\subset V_1 \subset V_2 \subset \;\cdots\; \subset V_t \subset V_{t+1}=V,
\end{equation}
and rational numbers $n_1,\, n_2,\ldots , \,n_t > 0$.
To a weighted filtration we associate a vector of $\mathbb{C}^p$
defined as $\Gamma=\sum_{i=1}^{t}n_i \Gamma^{(\dim V_i)}$
where
\begin{equation}
\label{semistandard}
\Gamma^{(k)}:=\big( \overbrace{k-p,\ldots,k-p}^k,
 \overbrace{k,\ldots,k}^{p-k} \big)
\qquad (1\leq k < p)
\, .
\end{equation}
Hence, the vector is of the form
$$\Gamma=(\overbrace{\Gamma_1,\ldots,\Gamma_1}^{\dim V^1},
\overbrace{\Gamma_2,\ldots,\Gamma_2}^{\dim V^2}, \ldots, \overbrace{\Gamma_{t+1},\ldots,\Gamma_{t+1}}^{\dim
V^{t+1}}) \; ,$$ where $V^i=V_i/V_{i-1}$. Giving the numbers $n_1,\ldots,n_t$ is equivalent to giving the
numbers $\Gamma_1,\ldots,\Gamma_{t+1}$ by setting
$$n_i=\frac{\Gamma_{i+1}-\Gamma_i}{p}
\qquad \text{and} \quad \sum_{i=1}^{t+1}\Gamma_i\dim V^i = 0\; .$$ A $1$-parameter subgroup of $\SL(V)$ is a
non-trivial homomorphism
$$\Gamma:\mathbb{C}^{\ast}\to \SL(V)\; .$$
To a $1$-parameter subgroup we associate a weighted
filtration as follows. There is a basis $\{e_1,\ldots,e_p\}$ of
$V$ where it has a diagonal form
$$
t\mapsto \operatorname{diag} \big(
t^{\Gamma_1},\ldots,t^{\Gamma_1},
t^{\Gamma_2},\ldots,t^{\Gamma_2}, \ldots ,
t^{\Gamma_{t+1}},\ldots,t^{\Gamma_{t+1}} \big)
$$
with
$\Gamma_1<\cdots<\Gamma_{t+1}$. Let
$$
0\subset V_1 \subset \cdots \subset V_{t+1}=V
$$
be the associated filtration. Finally recall that two $1$-parameter subgroups give the same
filtration if and only if they are
conjugate by an element of the parabolic subgroup of
$\SL(V)$ defined by the filtration.

The basis $\{e_1,\ldots,e_p\}$, together with a basis $\{w_j\}$ of $A$, induces a basis of $\Hom(\wedge^r V,A)$
indexed in a natural way by tuples $(i_1,\ldots,i_r,j)$ with $i_1<\cdots <i_{r}$, and the coordinate
corresponding to such an index is acted by the $1$-parameter subgroup as
$$
Q_{i_1,\cdots,i_r,j} \mapsto t^{\Gamma_{i_1}+\cdots+\Gamma_{i_r}} Q_{i_1,\cdots,i_r,j}\; .
$$
The coordinate $(i_1,\ldots,i_r,j)$ of the point corresponding to $E$ is non-zero if and only if the evaluations
of the sections $e_1,\ldots,e_r$ are linearly independent for generic $x\in X$. Therefore, the numerical
function (i.e. the \textit{minimal relevant weight}) which has to be calculated to apply Hilbert-Mumford
criterion for GIT stability (c.f. Theorem \ref{HMcrit}) is
\begin{eqnarray}
  \label{eq:muleft}
\mu(\overline{Q},V_\bullet,n_\bullet)&=&\min \{\Gamma_{i_1}+\cdots+\Gamma_{i_r}: \,
Q_{i_1,\ldots ,i_r,j}\neq 0 \} \notag\\
&=&\min \{\Gamma_{i_1}+\cdots+\Gamma_{i_r}: \,
Q(e_{i_1}\wedge\cdots \wedge e_{i_r})\neq 0 \} \notag\\
&=&\min \{\Gamma_{i_1}+\cdots+\Gamma_{i_r}: \,
e_{i_1}(x),\ldots,e_{i_r}(x)  \\
& & \qquad\text{linearly independent for generic
  $x\in X$}\}   \notag
\end{eqnarray}

After a short calculation (originally due to Gieseker) we obtain
\begin{equation}
\label{mrw}
\mu(\overline{Q},V_\bullet,n_\bullet)= \sum_{i=1}^{t} n_i ( r \dim V_i - r_i \dim V)
= \sum_{i=1}^{t+1} \frac{\Gamma_i}{\dim V} ( r^i \dim V - r\dim
V^i)
\end{equation}
 (recall $n_i=\frac{\Gamma_{i+1}-\Gamma_i}{p}$), where
$r_{i}=\rk E_{i}$, $E_{i}$ is the sheaf generated by evaluation of the sections of $V_{i}$ and $r^{i}=\rk
E^{i}$, being $E^{i}=E_{i}/E_{i-1}$.

By the Hilbert-Mumford criterion in Theorem \ref{HMcrit}, a point
$$
\overline{Q}\in \PP(\Hom(\wedge^r V , A) )
$$
is \emph{GIT semistable} if and only if for all weighted filtrations, it is
$$
\mu(\overline{Q},V_\bullet,n_\bullet)\leq 0\; .
$$
A point $\overline{Q}$ is \emph{GIT stable} if we get a strict inequality for all weighted filtrations. Using
the previous calculation, this can be stated as follows:

\begin{lem}
\label{HMcritsheaves}
A point $\overline{Q}$ is GIT semistable (resp. GIT stable) if
for all weighted filtrations $(V_\bullet,n_\bullet)$
$$\sum_{i=1}^{t}  n_i ( r \dim V_i - r_i \dim V)\leq 0$$
(resp. $<$).
\end{lem}

A \emph{weighted filtration} $(E_\bullet,n_\bullet)$ of a sheaf $E$ of rank $r$ is a filtration
\begin{equation}
\label{filtE} 0\subset E_1 \subset E_2 \subset \;\cdots\;\subset E_t \subset E_{t+1}=E,
\end{equation}
and rational numbers $n_1,\, n_2,\ldots , \,n_t > 0$.
To a weighted filtration we associate a vector of $\CC^r$
defined as $\gamma=\sum_{i=1}^{t}n_i \gamma^{(\rk E_i)}$ where
$$
\gamma^{(k)}:=\big( \overbrace{k-r,\ldots,k-r}^k,
\overbrace{k,\ldots,k}^{r-k} \big)\qquad (1\leq k < r)
\, .
$$
Hence, the vector is of the form
$$
\gamma=(\overbrace{\gamma_1,\ldots,\gamma_1}^{\rk E^1},
\overbrace{\gamma_2,\ldots,\gamma_2}^{\rk E^2},
\ldots,\overbrace{\gamma_{t+1},\ldots,\gamma_{t+1}}^{\rk E^{t+1}}) \;
,
$$
where $n_{i}=\dfrac{\gamma_{i+1}-\gamma_{i}}{r}$, and
$E^{i}=E_{i}/E_{i-1}$.

The following theorem follows from \cite{Gi1,Ma1,Ma2}.

\begin{thm}
\label{GITequiv} Let $E$ be a sheaf. There exists an integer
$m_{0}(E)$ such that, for $m>m_{0}(E)$, the associated point
$\overline{Q}$ is GIT semistable if and only if the sheaf is
semistable.
\end{thm}

From this property, Gieseker, in the case of algebraic surfaces and, later, Maruyama, for higher dimensional
varieties, constructs a moduli space of semistable torsion free sheaves as the GIT quotient, which is a
projective scheme (c.f. \cite[Theorem 0.3]{Gi1}, \cite[Theorem 4.11]{Ma2}).

Let $E$ be an unstable torsion free sheaf over $X$ of Hilbert polynomial $P$. We choose an integer $m_{0}$
larger than $m_0(E)$ (c.f. Theorem \ref{GITequiv}), also larger than the integer used in Gieseker's construction
of the moduli space, and such that $E$ is $m$-regular. Let $V$ be a vector space of dimension $P(m)=h^{0}(E(m))$
and fix an isomorphism $V\simeq H^{0}(E(m))$.

Recall that, through Geometric Invariant Theory, stability of a point in the parameter space can be checked by
$1$-parameter subgroups (c.f. Hilbert-Mumford criterion, Theorem \ref{HMcrit}). In other words, a point is
unstable if there exists any $1$-parameter subgroup which makes the quantity (\ref{mrw}) positive. It is a
natural question to ask if there exists a best way of destabilizing a GIT unstable point in this sense, i.e. a
$1$-parameter subgroup which gives maximum for (\ref{mrw}).

As we showed in Section \ref{kempfsection}, Kempf explores this idea in \cite{Ke} and answers yes to the
question, finding that there exists an special class of $1$-parameter subgroups which moves most rapidly toward
the origin.

We have seen that giving a weighted filtration, i.e. a filtration of vector subspaces $V_{1}\subset \cdots
\subset V_{t}\subset V$ and rational numbers $n_{1},\cdots,n_{t}>0$, is equivalent to giving a parabolic
subgroup with weights, which determines uniquely the vector $\Gamma$ of a $1$-parameter subgroup and two of
these $1$-parameter subgroup are conjugated by the parabolic and come from the same weighted filtration. Hence,
the data of $\Gamma$ is equivalent to the data of $(V_{\bullet},n_{\bullet})$.

Define the function in Definition \ref{kempffunction},
$$K(x,\Gamma)=\frac{\mu(x,\Gamma)}{\|\Gamma\|}\; .$$
as the following function
\begin{equation}
 \label{kempffunctionsheaves}
K(x,\Gamma)=\frac{\sum_{i=1}^{t}  n_{i} (r\dim V_{i}-r_{i}\dim V)}
{\sqrt{\sum_{i=1}^{t+1} {\dim V^{i}_{}} \Gamma_{i}^{2}}}=\mu(V_{\bullet},n_{\bullet})\; ,
\end{equation}
which we call \emph{Kempf function}. The numerator of the function coincides
with the calculation of the minimal relevant weight by
Hilbert-Mumford criterion for GIT stability (c.f. \eqref{mrw}), and the denominator is a function $||\cdot ||$ in the set $\Gamma(SL(V))$ of $1$-parameter subgroups of $SL(V)$,
which is precisely the norm of the vector
$$
\Gamma=(\overbrace{\Gamma_1,\ldots,\Gamma_1}^{\dim V^1},\overbrace{\Gamma_2,\ldots,\Gamma_2}^{\dim V^2},\ldots,\overbrace{\Gamma_{t+1},
\ldots,\Gamma_{t+1}}^{\dim V^{t+1}})
$$
associated to each $1$-parameter subgroup $\Gamma$.

To define the Kempf function we need to choose a \textit{length} in $\Gamma(SL(V))$ (c.f. Definition
\ref{length}). Recall that for a simple group $G$ (as it is the case of $G=\SL(V)$) every bilinear symmetric
invariant form is a multiple of the Killing form (c.f. Remark \ref{killinglengths}), and this norm $||\Gamma ||$
we choose verifies these properties. Hence, the function we defined in (\ref{kempffunctionsheaves}) is a Kempf
function as in Definition \ref{kempffunction}.

We take the GIT quotient by the group $G=\SL(V)$, for which, Theorem \ref{kempftheorem0} (c.f. \cite[Theorem
2.2]{Ke}) states that whenever there exists any $\Gamma$ giving a positive value for the numerator of the
function (i.e. whenever there exists a 1-parameter subgroup whose numerical function (\ref{mrw}) is positive,
which is equivalent to the sheaf $E$ to be unstable), there exists a unique parabolic subgroup containing a
unique $1$-parameter subgroup in each maximal torus, giving maximum for the Kempf function i.e., there exists a
unique weighted filtration for which the Kempf function achieves a maximum.

Note that $\mu(V_{\bullet},n_{\bullet})=\mu(V_{\bullet},\alpha n_{\bullet})$, for every $\alpha>0$, hence by
multiplying each $n_{i}$ by the same scalar $\alpha$, which we call \textit{rescaling the weights}, we get
another 1-parameter subgroup but the same value for the Kempf function. Hence, we divide by the norm in the
Kempf function to get a well defined maximal weighted filtration, i.e. defined up to rescaling.

Therefore, Theorem \ref{kempftheorem0} rewritten in our case asserts
the following:

\begin{thm}
\label{kempftheoremsheaves} There exists a unique weighted filtration
$$0\subset V_{1}\subset \cdots \subset V_{t+1}= V$$
and rational numbers $n_{1},\cdots,n_{t}>0$, up to multiplication by a scalar, called the \emph{Kempf filtration
of V}, such that the Kempf function $\mu(V_{\bullet},n_{\bullet})$ achieves the maximum among all filtrations
and positive weights $n_{i}>0$.
\end{thm}

We construct a filtration by subsheaves of $E$ out of the Kempf filtration of $V$ in Theorem
\ref{kempftheoremsheaves}. Recall that $E$ is an unstable torsion free sheaf over $X$ of Hilbert polynomial $P$.
Let $m$ be an integer, $m\geq m_{0}$ and let $V$ be a vector space of dimension $P(m)=h^{0}(E(m))$ (recall that
$m_0$ was defined before). We fix an isomorphism $V\simeq H^{0}(E(m))$ and let $V_{1}\subset \cdots \subset
V_{t+1}= V$ be the filtration of vector spaces given by Theorem \ref{kempftheoremsheaves}, called the
\emph{Kempf filtration of V}. For each index $i$, let $E^{m}_{i}\subset E$ be the subsheaf generated by $V_{i}$
under the evaluation map. We call this filtration
$$
0\subseteq E^{m}_1 \subseteq E^{m}_2 \subseteq \;\cdots\; \subseteq E^{m}_t \subseteq E^{m}_{t+1}=E\; ,
$$
the \emph{$m$-Kempf filtration of E}. Note that it depends on the integer $m$ we choose in the process.

The question we finished Section \ref{kempfsection} with was

\vspace{0.7cm} \textbf{Is the Harder-Narasimhan filtration related to the $1$-parameter subgroup given by
Kempf?} \vspace{0.7cm}

The maximal unstability with respect to Definition \ref{giesekerstab} is given by the Harder-Narasimhan
filtration (c.f. Theorem \ref{HNdefthm}) and the GIT maximal unstability is encoded in the Kempf filtration of
$V$, by Theorem \ref{kempftheoremsheaves}. This filtration of vector subspaces can be evaluated to get a
filtration of subsheaves, the $m$-Kempf filtration of $E$, depending on an integer $m$. Therefore, the previous
question turns out to be more concrete:

\vspace{0.7cm} \textbf{Does the $m$-Kempf filtration coincide with the Harder-Narasimhan filtration?}
\vspace{0.7cm}

The answer will be yes. In the following pages we will develop a technique to prove the following two theorems:

\begin{thm}
\label{kempfstabilizes} There exists an integer $m'\gg 0$ such that the $m$-Kempf filtration of $E$ is
independent of $m$, for $m'\geq m$.
\end{thm}

This filtration we obtain, independent of the integer $m$, will be called the \emph{Kempf filtration} of $E$.

\begin{thm}
\label{kempfisHN}
The Kempf filtration of an unstable torsion free coherent sheaf $E$ coincides with
the Harder-Narasimhan filtration of $E$.
\end{thm}

The method we use to prove Theorem \ref{kempfstabilizes} and Theorem \ref{kempfisHN} will be translated to other
moduli problems to prove an analogous result in the subsequent sections of this chapter.

\subsection{Results on convexity}
\label{convexcones}

In this subsection we define the machinery which will serve us in
the following. We study a function on a convex set, and how to
maximize it. It will turn out to be that this function will be in
correspondence with the Kempf function and we will use the results
of this subsection to figure out properties about the Kempf
filtration.

Endow $\mathbb{R}^{t+1}$ with an inner product
$(\cdot,\cdot)$ defined by a diagonal matrix
 $$
 \left(
 \begin{array}{ccc}
 b^1 & & 0 \\
  & \ddots & \\
 0 & & b^{t+1}\\
 \end{array}
 \right)
 $$
where $b^i$ are positive integers. Let
$$
\mathcal{C}= \big\{ x\in \mathbb{R}^{t+1} : x_1<x_2<\cdots <x_{t+1}
\big\}\; ,
$$
$$
\overline{\mathcal{C}}= \big\{ x\in \mathbb{R}^{t+1} :
x_1\leq x_2\leq \cdots \leq x_{t+1} \big\}\; ,
$$
and let $v= (v_1,\cdots,v_{t+1})\in \mathbb{R}^{t+1}-\{0\}$ verifying
$\sum_{i=0}^{t+1} v_{i}b^{i}=0$. Define the function
\begin{eqnarray*}
  \label{eq:mu}
\mu_{v}:\overline{\mathcal{C}}-\{0\} & \to & \mathbb{R}\\
\Gamma & \mapsto & \mu_{v}(\Gamma)=\frac{(\Gamma,v)}{||\Gamma||}
\end{eqnarray*}
and note that  $\mu_{v}(\Gamma)=||v||\cdot \cos\beta(\Gamma,v)$,
where $\beta(\Gamma, v)$ is the angle between $\Gamma$ and $v$.
Then, the function $\mu_{v}(\Gamma)$ does not depend on the norm
of $\Gamma$ and takes the same value on every point of the ray
spanned by each $\Gamma$.

Assume that there exists $\Gamma\in\overline{\mathcal{C}}$ with
$\mu_v(\Gamma)>0$. In that case, we want to find a vector
$\Gamma\in \overline{\mathcal{C}}$ which maximizes the function
defined before.

Let $w^i=-b^iv_i$, $w_0=0$, $w_i=w^1+\cdots+w^i$, $b_0=0$, and $b_i=b^1+\cdots+b^i$. Note that $w_{t+1}=0$, by
construction. We draw a graph joining the points with coordinates $(b_i,w_i)$. Note that this graph has $t+1$
segments, each segment has slope $-v_i$ and width $b^i$. This is the graph drawn with a thin line in the figure.
Now draw the convex envelope of this graph (thick line in Figure \ref{convexfigure}), whose coordinates we
denote by $(b_i,\widetilde{w_i})$, and let us define
$\Gamma_{i}=-\frac{\widetilde{w_{i}}-\widetilde{w_{i-1}}}{b^i}$. In other words, the quantities $-\Gamma_i$ are
the slopes of the convex envelope graph. We call the vector defined in this way $\Gamma_{v}$. Note that the
vector $\Gamma_{v}=(\Gamma_1,\cdots,\Gamma_{t+1})$ belongs to $\overline{\mathcal{C}}$ by construction and
$\Gamma_{v}\neq 0$.

$$ $$

\setlength{\unitlength}{1cm}

\begin{figure}[h]
   \begin{center}
\begin{picture}(13,7)(-1,-1)
\thicklines
\put(0,0){\line(1,0){11.5}}
\put(0,0){\line(0,1){5.5}}
\put(6,-0.3){\makebox(0,0)[c]{$b_i$}} \put(-0.4,3){\makebox(0,0)[c]{\rotatebox{90}{$w_i,\widetilde{w_i}$}}}
\put(0,0){\makebox(0,0){$\circ$}} \put(0,0){\line(1,1){4}} \put(2,2){\makebox(0,0){$\circ$}}
\put(3,3){\makebox(0,0){$\circ$}} \put(4,4){\makebox(0,0){$\circ$}} \put(4,4){\line(4,1){4}}
\put(5,4.25){\makebox(0,0){$\circ$}} \put(8,5){\makebox(0,0){$\circ$}} \put(8,5){\line(1,-1){2}}
\put(9.5,3.5){\makebox(0,0){$\circ$}} \put(10,3){\makebox(0,0){$\circ$}} \put(10,3){\line(1,-3){1}}
\put(11,0){\makebox(0,0){$\circ$}} \thinlines \put(0,0){\makebox(0,0){$\circ$}} \put(0,0){\line(2,1){2}}
\put(2,1){\makebox(0,0){$\circ$}} \put(2,1){\line(1,1){1}} \put(3,2){\makebox(0,0){$\circ$}}
\put(3,2){\line(1,2){1}} \put(4,4){\line(1,-1){1}} \put(5,3){\makebox(0,0){$\circ$}} \put(5,3){\line(3,2){3}}
\put(9,4){\makebox(0,0){$\circ$}} \put(9,4){\line(1,-4){0.5}} \put(9.5,2){\makebox(0,0){$\circ$}}
\put(9.5,2){\line(1,2){0.5}} \put(2,0.6){\makebox(0,0)[l]{$(b_1,w_1)$}}
\put(1.8,2.2){\makebox(0,0)[r]{$(b_1,\widetilde{w_1})$}} \put(3,1.6){\makebox(0,0)[l]{$(b_2,w_2)$}}
\put(2.8,3.2){\makebox(0,0)[r]{$(b_2,\widetilde{w_2})$}}
\put(3.8,4.2){\makebox(0,0)[r]{$(b_3,\widetilde{w_3}=w_3)$}}
\put(5,4.7){\makebox(0,0)[c]{$(b_4,\widetilde{w_4})$}} \put(5,2.6){\makebox(0,0)[c]{$(b_4,w_4)$}}
\put(8,5.4){\makebox(0,0)[c]{$(b_5,\widetilde{w_5}=w_5)$}}
\end{picture}
\end{center}
\caption{Convex envelope $\Gamma_{v}$ of $v$} \label{convexfigure}
\end{figure}
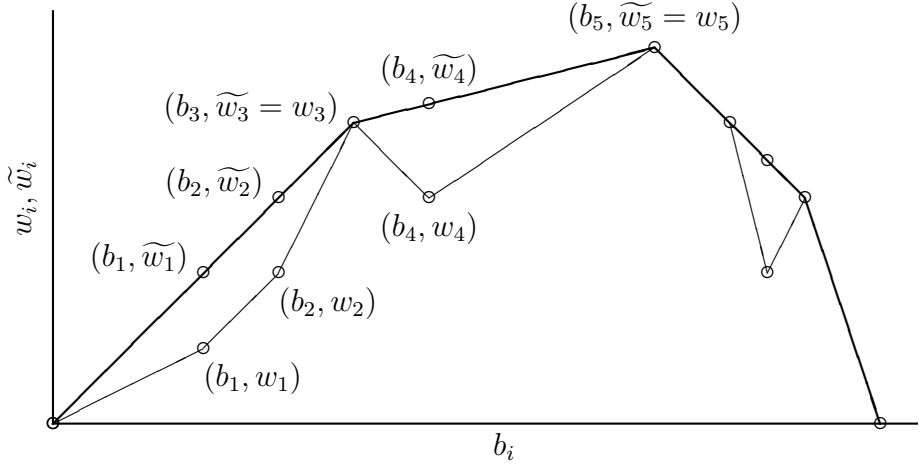

\begin{rem}
\label{rmk} Observe that $\widetilde{w_{i}}> w_i$, then $\Gamma_{i}=\Gamma_{i+1}$. Indeed, if
$\widetilde{w_{i}}> w_i$, there will be a segment in the convex envelope joining two vertices such that
$w_{j}=\widetilde{w_{j}}$ and $w_{k}=\widetilde{w_{k}}$, with $j<i$ and $i<k$. Then, it is clear that all
segments joining the intermediate heights $\widetilde{w_{l}}$, $j<l<k$, will have the same slope, in particular
$\Gamma_{i}=\Gamma_{i+1}$.
\end{rem}

\begin{thm}
\label{maxconvexenvelope} The vector $\Gamma_{v}$ defined in this way (c.f. Figure \ref{convexfigure}) gives a
maximum for the function $\mu_v$ on its domain.
\end{thm}

Before proving the theorem we need some lemmas.

\begin{lem}
\label{close}
Let $v= (v_1,\cdots,v_{t+1})\in \mathbb{R}^{t+1}-\{0\}$ verifying
$\sum_{i=0}^{t+1} v_{i}b^{i}=0$.
Let $\Gamma$ be the point in $\overline{\mathcal{C}}$ which is closest
to $v$. Then $\Gamma$ achieves the maximum of $\mu_v$.
\end{lem}

\begin{pr}
For any $\alpha\in \mathbb{R}^{>0}$, the vector $\alpha \Gamma$ is also in $\overline{\mathcal{C}}$, so in
particular $\Gamma$ is the closest point to $v$ in the line $\alpha \Gamma$. This point is the orthogonal
projection of $v$ into the line $\alpha \Gamma$, and the distance is
\begin{equation}
  ||v|| \sin \beta(v,\Gamma)\label{eq:sin}\; ,
\end{equation}
where $\beta(\Gamma,v)$ is the angle between $\Gamma$ and $v$. But, a vector $\Gamma\in \overline{\mathcal{C}}$
minimizes \eqref{eq:sin} if and only if it maximizes
$$
||v|| \cos \beta(\Gamma,v) = \frac{(\Gamma,v)}{||\Gamma||}\; ,
$$
so the lemma is proved.
\end{pr}

We say that an affine hyperplane in $\mathbb{R}^{t+1}$ separates a
point $v$ from $\mathcal{C}$ if $v$ is on one side of the
hyperplane and all the points of $\mathcal{C}$ are on
the other side of the hyperplane.

\begin{lem}
\label{hyperplane} Let $v\notin \overline{\mathcal{C}}$. A point
$\Gamma\in \overline{\mathcal{C}}-\{0\}$ gives minimum distance to $v$
if and only if the hyperplane $\Gamma+(v-\Gamma)^\perp$ separates
$v$ from $\mathcal{C}$.
\end{lem}

\begin{pr}
$\Rightarrow$) Let $\Gamma\in \overline{\mathcal{C}}$ and assume that there is a point $w\in \mathcal{C}$ on the
same side of the hyperplane as $v$. The segment going from $\Gamma$ to $w$ is in $\overline{\mathcal{C}}$ (by
convexity of $\overline{\mathcal{C}}$), but there are points in this segment (near $\Gamma$), which are closer
to $v$ than $\Gamma$.

$\Leftarrow)$ Let $\Gamma$ be a point in $\overline{\mathcal{C}}$
such that $\Gamma+(v-\Gamma)^\perp$ separates $v$ from
$\mathcal{C}$. Let $w \in\overline{\mathcal{C}}$ be another point.
Let $w'$ be the intersection of the hyperplane and the segment
which goes from $w$ to $v$. Since the hyperplane separates
$\mathcal{C}$ from $v$, either $w'=w$ or $w'$ is in the interior
of the segment. Therefore
$$
d(w,v) \geq d(w',v) \geq d(\Gamma,v)\; ,
$$
where the last inequality follows from the fact that $\Gamma$
is the orthogonal projection of $v$ to the hyperplane.
\end{pr}

\begin{pr}[Proof of the Theorem \ref{maxconvexenvelope}]
Let $\Gamma_{v}=(\Gamma_{1},...,\Gamma_{t+1})$ be the vector in
the hypothesis of the theorem. If $v\in \overline{\mathcal{C}}$,
then $\Gamma_{v}=v$, and use Lemma \ref{close} to conclude. If
$v\notin \overline{\mathcal{C}}$, by Lemmas \ref{close} and
\ref{hyperplane}, it is enough to check that the hyperplane
$\Gamma_{v} + (v-\Gamma_{v})^\perp$ separates $v$ from
$\mathcal{C}$.

Let $\Gamma_{v}+\epsilon\in \mathcal{C}$, $\epsilon\in\mathbb{R}^{t+1}$. The condition that
$\Gamma_{v}+\epsilon$ belongs to $\mathcal{C}$ means that
\begin{equation}
  \label{eq:leq}
\epsilon_i-\epsilon_{i+1} < \Gamma_{i+1}-\Gamma_{i}
\end{equation}
The hyperplane separates $v$ from $\mathcal{C}$ if and only if
$(v-\Gamma_{v},\epsilon)<0$ for all such $\epsilon$. Therefore we
calculate (using the convention $\widetilde{w_0}=0$, $w_0=0$, and
$\widetilde{w_{t+1}}=w_{t+1}=0$)
$$
(v-\Gamma_{v},\epsilon) = \sum_{i=1}^{t+1}
b^i(v_i-\Gamma_i)\epsilon_i = \sum_{i=1}^{t+1}
(-w^i+(\widetilde{w_i}-\widetilde{w_{i-1}}))\epsilon_i =
$$
$$
=\sum_{i=1}^{t+1} \big( (\widetilde{w_i}-\widetilde{w_{i-1}})-(w_i-w_{i-1})\big)\epsilon_i
= \sum_{i=1}^{t+1} (\widetilde{w_{i}}- w_i)(\epsilon_i-\epsilon_{i+1})\; .
$$
If $\widetilde{w_i}=w_i$, then the corresponding summand is zero. On the other hand, if $\widetilde{w_i}> w_i$,
then $\Gamma_{i+1}=\Gamma_i$ (c.f. Remark \ref{rmk}), and \eqref{eq:leq} implies $\epsilon_i-\epsilon_{i+1}<0$.
In any case, the summands are always non-positive, and there is at least one which is negative (because $v\notin
\overline{\mathcal{C}}$ and then $v\neq \Gamma_{v}$ and $\widetilde{w_i}> w_i$ for at least one $i$). Hence
$$(v-\Gamma_{v},\epsilon)<0\; .$$
\end{pr}

Therefore, the function $\mu_{v}(\Gamma)$ achieves its maximum for
the value $\Gamma_{v}\in \overline{\mathcal{C}}-\{0\}$ (or any other
point on the ray $\alpha\Gamma_{v}$) defined as the convex envelope
of the graph associated to $v$.

\subsection{Graph and identification}

In the last section we studied a geometrical function,
$\mu_{v}(\Gamma)$, very similar to the Kempf function. This new
function depends on two arguments, one is a vector $\Gamma\in\overline{\mathcal{C}}-\{0\}$, where
$$
\overline{\mathcal{C}}= \big\{ x\in \mathbb{R}^{t+1} :
x_1\leq x_2\leq \cdots \leq x_{t+1} \big\}\; ,
$$
and the other is $v=(v_1,\cdots,v_{t+1})\in \mathbb{R}^{t+1}-\{0\}$ verifying $\sum_{i=0}^{t+1} v_{i}b^{i}=0$,
for certain coefficients $b^{i}$ of an inner product in an Euclidean space. We will relate both functions where
the first argument $\Gamma$ will be associated to a $1$-parameter subgroup (or to a weighted filtration
$(V_{\bullet},n_{\bullet})$ which we recall that is equivalent), and the second one will be associated to the
numerical invariants of the Kempf filtration of $V$,
$$0\subset V_{1}\subset \cdots \subset V_{t+1}= V$$
(c.f. Theorem \ref{kempftheoremsheaves}) and the $m$-Kempf filtration of $E$
$$
0\subseteq E^{m}_1 \subseteq E^{m}_2 \subseteq \;\cdots\; \subseteq E^{m}_t \subseteq E^{m}_{t+1}=E
$$ obtained by evaluating. With this, we will be able to prove
properties of the filters appearing on the different $m$-Kempf filtrations for each $m$, out from convexity
properties of the function $\mu_{v}$ (c.f. Theorem \ref{maxconvexenvelope}). Both functions have to be maximized
by the convex envelope of the graph defined by $v$, or the Kempf filtration of $V$, therefore both notions have
to correspond to the same filtrations. And to make precise that relation, we have to encode the $m$-Kempf
filtration as a graph.

\begin{dfn}
\label{graph} Let $m\geq m_{0}$. Given $0\subset V_{1}\subset \cdots \subset V_{t+1}= V$, a filtration of vector
spaces of $V$, define
$$v_{m,i}=m^{n+1}\cdot \frac{1}{\dim V^{i}\dim V}\big[r^{i}\dim V-r\dim V^{i}\big]\; ,$$
$$b_{m}^{i}=\dfrac{1}{m^{n}}\dim V^{i}>0\; ,$$
$$w_{m}^{i}=-b_{m}^{i}\cdot v_{m,i}=m\cdot \frac{1}{\dim V}\big[r\dim V^{i}-r^{i}\dim V\big]\; .$$
Also let
$$b_{m,i}=b_{m}^{1}+\ldots +b_{m}^{i}=\dfrac{1}{m^{n}}\dim V_{i}\; ,$$
$$w_{m,i}=w_{m}^{1}+\ldots +w_{m}^{i}=m\cdot \frac{1}{\dim V}\big[r\dim V_{i}-r_{i}\dim V\big]\; .$$
We call the graph defined by points $(b_{m,i},w_{m,i})$ the \emph{graph associated to the filtration} $V_{\bullet}\subset V$.
\end{dfn}

Now we can identify the Kempf function (\ref{kempffunctionsheaves}) in Theorem \ref{kempftheoremsheaves}
$$\mu(V_{\bullet},n_{\bullet})=\frac{\sum_{i=1}^{t} n_{i}(r\dim V_{i}-r_{i}\dim V)}{\sqrt{\sum_{i=1}^{t+1}
{\dim V^{i}} \Gamma_{i}^{2}}}\; ,$$ with the function in Theorem
\ref{maxconvexenvelope} up to a factor which is a power of $m$, by
defining $v_{m,i}$, the coordinates of vector $v_{m}$, and
$b_{m}^{i}$, the eigenvalues of the inner product, as in
Definition \ref{graph}. Note that $-v_{m,i}$ are the slopes of the
graph associated to the filtration $V_{\bullet}\subset V$. To give the weights $n_{i}$ is the same that to give the coordinates $\Gamma_{i}$
(recall the discussion about the correspondence between $1$-parameter subgroups of $SL(V)$ and weighted filtrations). Also note that
$\sum_{i=1}^{t+1}v_{m,i}b^{i}_{m}=0$. Then, an easy calculation
shows that

\begin{prop}
\label{identification}
For every integer $m$, the following equality holds
$$
\mu(V_{\bullet},n_{\bullet})=m^{(-\frac{n}{2}-1)}\cdot
\mu_{v_{m}}(\Gamma)
$$
between the Kempf function (\ref{kempffunctionsheaves}) in Theorem \ref{kempftheoremsheaves} and the function
in Theorem \ref{maxconvexenvelope}.
\end{prop}

In the following, we will omit the subindex $m$ for the numbers
$v_{m,i}$, $b_{m,i}$, $w_{m,i}$ in the definition of the graph
associated to a filtration of vector spaces, where it is clear
from the context. Hence, given $V\simeq H^{0}(E(m))$ we will refer
to a filtration $V_{\bullet}\subset V$ and a vector
$v=(v_{1},\ldots,v_{t+1})$ as the vector of the graph associated
to the filtration.

\begin{rem}
\label{growth}
We introduce the factor $m^{n+1}$ in Definition \ref{graph} for
convenience, so that $v_{m,i}$ and $b_{m}^{i}$ have order zero on
$m$, because $\dim V=P(m)$ appears in their expressions. Then, the
size of the graph does not change when $m$ grows.
\end{rem}

Now, let us prove two lemmas encoding the convexity properties of the graph associated to the Kempf filtration.
They will be strongly used in the following, to show properties shared by the possible filters $E_{i}^{m}$
appearing in the different $m$-Kempf filtrations and, finally, to prove Theorems \ref{kempfstabilizes} and
Theorem \ref{kempfisHN}.

\begin{lem}
\label{lemmaA}
Let $0\subset V_{1}\subset \cdots \subset V_{t+1}= V$ be the Kempf filtration of $V$ (cf. Theorem
\ref{kempftheoremsheaves}). Let $v=(v_{1},...,v_{t+1})$ be the vector of the graph associated to this filtration
by Definition \ref{graph}. Then
$$
v_{1}<v_{2}<\ldots<v_{t}<v_{t+1}\; ,
$$
i.e., \emph{the graph is convex}.
\end{lem}

\begin{pr}
By Theorem \ref{kempftheoremsheaves} the maximum of $\mu(V_{\bullet},n_{\bullet})$ among all filtrations
$V_{\bullet}\subset V$ and weights $n_{i}>0,\forall i$ is achieved by a unique weighted filtration
$(V_{\bullet},n_{\bullet})$, $n_{i}>0,\forall i$, up to rescaling. Let $V_{\bullet}\subset V$ be this
filtration, and allow $n_{i}$ to vary. By Proposition \ref{identification} $\mu(V_{\bullet},n_{\bullet})$ is
equal to $\mu_v$ up to a constant factor. By Theorem \ref{maxconvexenvelope}, $\mu_v$ achieves the maximum on
$\Gamma_{v}$. The vector $\Gamma_{v}$ corresponds to the weights $n_{i}$ given by Theorem
\ref{kempftheoremsheaves}. Summing up, if $V_\bullet\subset V$ is Kempf filtration of $V$, then the vector
$\Gamma_{v}=(\Gamma_{1},\ldots,\Gamma_{t+1})$ verifies $\Gamma_{i}<\Gamma_{i+1},\forall i$.

Assume that, for the Kempf filtration of $V$, there exists some $i$ such
that $v_{i}\geq v_{i+1}$. Then $v\notin \mathcal{C}$ and, by Lemma
\ref{close}, $\Gamma_{v}\in \overline{\mathcal{C}}\backslash
\mathcal{C}$, which means that there exists some $j$ with
$\Gamma_{j}=\Gamma_{j+1}$, but we have just seen that this is
impossible.
\end{pr}

\begin{lem}
\label{lemmaB} Let $0\subset V_{1}\subset \cdots \subset V_{t+1}= V$ be the Kempf filtration of $V$ (cf. Theorem
\ref{kempftheoremsheaves}). Let $W$ be a vector space with $V_{i}\subset W\subset V_{i+1}$ and consider the new
filtration $V'_{\bullet}\subset V$
\begin{equation}
    \begin{array}{ccccccccccccccccc}
    0 & \subset & V'_{1} & \subset & \cdots & \subset & V'_{i} & \subset & V'_{i+1} & \subset &
    V'_{i+2} & \subset & \cdots & \subset & V'_{t+2} & = & V\\
    || & & || & & & & || & & || & & || & & & & & & ||  \\
    0 & \subset & V_{1}& & & \subset & V_{i} & \subset & W & \subset & V_{i+1} & \subset & \cdots & \subset & V_{t+1} & = & V
    \end{array}
\end{equation}
Then, $v'_{i+1}\geq v_{i+1}$. We say that \emph{the Kempf filtration is the convex envelope of every
refinement}.
\end{lem}

\begin{pr}
The graph associated to $V'_{\bullet}\subset V$ has one more point
than the graph associated to $V_{\bullet}\subset V$, hence it is a
refinement of the graph associated to Kempf filtration of $V$.
Therefore the convex envelope of the graph associated to $v'$ has
to be equal to the graph associated to $v$, and this happens only
when the extra point associated to $W$ is not above the graph
associated to $v$, which means that the slope $-v'_{i+1}$ has to
be less or equal than $-v_{i+1}$.
\end{pr}

\begin{rem}
Note that Lemmas \ref{lemmaA} and \ref{lemmaB} assert two properties
similar to the ones of the Harder-Narasimhan filtration (c.f.
Theorem \ref{HNdefthm}). This will be the key point in the proof
of Theorem \ref{kempfisHN}.
\end{rem}

Hence, we will prove that, for $m$ large enough, the $m$-Kempf filtration stabilizes in the sense
$E^{m}_{i}=E^{m+l}_{i},\forall i,\forall l>0$, in Theorem \ref{kempfstabilizes}. The $m$-Kempf filtration for
$m\gg 0$ will be called the Kempf filtration of $E$, and the goal is to show that it coincides with the
Harder-Narasimhan filtration of $E$ in Theorem \ref{kempfisHN}.

\subsection{Properties of the $m$-Kempf filtration}

We will show that the filters
appearing in the different $m$-Kempf filtrations form a bounded
family.

First recall Lemma \ref{Simpson} in subsection \ref{resultsonboundedness}. Also recall the definition of the
Hilbert polynomials of $\mathcal{O}_{X}$ in (\ref{HpolO}) and $E$ in (\ref{HpolE}). Then, let us define
\begin{equation}
\label{C_constant}
C=\max\{r|\mu _{\max }(E)|+\frac{d}{r}+r|B|+|A|+1\;,\;1\},
\end{equation}
a positive constant.

\begin{prop}
\label{boundedness} Given an integer $m$ and a vector space
$V\simeq H^{0}(E(m))$, we have the Kempf filtration
$V_{\bullet}\subset V\simeq H^{0}(E(m))$ (c.f. Theorem \ref{kempftheoremsheaves}) and, by evaluation, the
$m$-Kempf filtration $E_{\bullet}^{m}\subset E$. There exists an
integer $m_{2}$ such that for $m\geq m_{2}$, each filter in the
$m$-Kempf filtration of $E$ has slope $\mu(E^{m}_{i})\geq
\dfrac{d}{r}-C$.
\end{prop}

\begin{pr}
Choose an $m_{1}\geq m_{0}$ such that for $m\geq m_{1}$
$$[\mu_{max}(E)+gm+B]_{+}=\mu_{max}(E)+gm+B$$
and
$$[\frac{d}{r}-C+gm+B]_{+}=\frac{d}{r}-C+gm+B\; .$$
Now, let $m\geq m_{1}$ and let
$$0 \subseteq E^{m}_1 \subseteq E^{m}_2 \subseteq \;\cdots\; \subseteq E^{m}_t \subseteq E^{m}_{t+1}=E$$
be the $m$-Kempf filtration of $E$.

Suppose we have a filter $E^{m}_{i}\subseteq E$, of rank $r_{i}$
and degree $d_{i}$, such that $\mu (E^{m}_{i})<\frac{d}{r}-C$. The subsheaf $E_{i}^{m}(m)\subset E(m)$ satisfies the estimate in Lemma \ref{Simpson},
$$h^{0}(E^{m}_{i}(m))\leq \frac{1}{g^{n-1}n!}\big ((r_{i}-1)([\mu_{max}(E^{m}_{i})+gm+B]_{+})^{n}+([\mu_{min}(E^{m}_{i})+gm+B]_{+})^{n}\big)\; ,$$
where $\mu_{max}(E^{m}_{i}(m))=\mu_{max}(E^{m}_{i})+gm$ and similarly for $\mu_{min}$.

Note that $\mu_{max}(E^{m}_{i})\leq \mu_{max}(E)$ and
$\mu_{min}(E^{m}_{i})\leq \mu(E^{m}_{i})<\frac{d}{r}-C$, so
$$h^{0}(E^{m}_{i}(m))\leq \frac{1}{g^{n-1}n!}\big ((r_{i}-1)([\mu_{max}(E)+gm+B]_{+})^{n}+([\frac{d}{r}-C+gm+B]_{+})^{n}\big)\; ,$$
and, by choice of $m$,
$$h^{0}(E^{m}_{i}(m))\leq \frac{1}{g^{n-1}n!}\big ((r_{i}-1)(\mu_{max}(E)+gm+B)^{n}+(\frac{d}{r}-C+gm+B)^{n}\big)=G(m)\; ,$$
where
$$G(m)=\frac{1}{g^{n-1}n!}\big
[r_{i}g^{n}m^{n}+ng^{n-1}\big((r_{i}-1)\mu_{max}(E)+\frac{d}{r}-C+r_{i}B\big)m^{n-1}+\cdots
\big ]\; .$$

Recall that, by Definition \ref{graph}, to such filtration we
associate a graph with heights, for each $j$,
$$w_{j}=w^{1}+\ldots +w^{j}=m\cdot \frac{1}{\dim V}\big[r\dim V_{j}-r_{j}\dim V\big]\; .$$
To reach a contradiction, it is enough to show that $w_{i}<0$. In that case, the graph has to be convex by Lemma
\ref{lemmaA}. If $w_{i}<0$ there is a $j<i$ such that $-v_{j}<0$, because the graph starts at the origin. Hence,
the rest of the slopes of the graph are negative, $-v_{k}<0$, $k\geq i$, because the slopes have to be
decreasing. Then $w_{i}>w_{i+1}>\ldots w_{t+1}$, and $w_{t+1}<0$. But it is
$$w_{t+1}=m\cdot \frac{1}{\dim V}\big[r\dim V_{t+1}-r_{t+1}\dim V\big]=0\; ,$$
because $r_{t+1}=r$ and $V_{t+1}=V$, then the contradiction.

Let us show that $w_{i}<0$. Since $E^{m}_{i}(m)$ is generated by
$V_{i}$ under the evaluation map, it is $\dim V_{i}\leq
h^{0}(E^{m}_{i}(m))$, hence
$$
w_{i}=\frac{m}{\dim V}\big[r\dim V_{i}-r_{i}\dim V\big]\leq
$$
$$
\leq \frac{m}{P(m)}\big[rh^{0}(E^{m}_{i}(m))-r_{i}P(m)\big]\leq
\frac{m}{P(m)}\big[rG(m)-r_{i}P(m)\big]\; .
$$

Hence, $w_{i}<0$ is equivalent to
$$\Psi(m)=rG(m)-r_{i}P_{E}(m)<0\; ,$$
where $\Psi(m)=\xi_{n}m^{n}+\xi_{n-1}m^{n-1}+\cdots +\xi_{1}m+\xi_{0}$ is an $n^{th}$-order polynomial. Let us
calculate the $n^{th}$-coefficient:
$$\xi_{n}=(rG(m)-r_{i}P(m))_{n}=r\frac{r_{i}g}{n!}-r_{i}\frac{rg}{n!}=0\; .$$
Then, $\Psi(m)$ has no coefficient in order $n^{th}$. Let us calculate the $(n-1)^{th}$-coefficient:
$$\xi_{n-1}=(rG(m)-r_{i}P(m))_{n-1}=(rG_{n-1}-r_{i}\frac{A}{(n-1)!})\; ,$$
where $G_{n-1}$ is the $(n-1)^{th}$-coefficient of the polynomial $G(m)$,
$$G_{n-1}=\frac{1}{g^{n-1}n!}ng^{n-1}((r_{i}-1)\mu_{max}(E)+\frac{d}{r}-C+r_{i}B)=$$
$$\frac{1}{(n-1)!}((r_{i}-1)\mu_{max}(E)+\frac{d}{r}-C+r_{i}B)\leq$$
$$\frac{1}{(n-1)!}((r_{i}-1)|\mu_{max}(E)|+\frac{d}{r}-C+r_{i}|B|)\leq$$
$$\frac{1}{(n-1)!}(r|\mu_{max}(E)|+\frac{d}{r}-C+r|B|)<\frac{-|A|}{(n-1)!}\; ,$$
last inequality coming from the definition of $C$ in \eqref{C_constant}.
Then
$$\xi_{n-1}<r\big(\frac{-|A|}{(n-1)!}\big)-r_{i}\frac{A}{(n-1)!}=
\frac{-r|A|-r_{i}A}{(n-1)!}<0$$ because $-r|A|-r_{i}A<0$.

Therefore $\Psi(m)=\xi_{n-1}m^{n-1}+\cdots +\xi_{1}m+\xi_{0}$
with $\xi_{n-1}<0$, so there exists $m_{2}\geq m_{1}$ such that for $m\geq
m_{2}$ we will have $\Psi(m)<0$ and
$w_{i}<0$, then the contradiction.
\end{pr}

\begin{prop}
\label{regular} There exists an integer $m_{3}$ such that for
$m\geq m_{3}$ the sheaves $E^{m}_{i}$ and
$E^{m,i}=E^m_i/E^m_{i-1}$ are $m_{3}$-regular. In particular their
higher cohomology groups, after twisting with
$\mathcal{O}_{X}(m_{3})$, vanish and they are generated by global
sections.
\end{prop}
\begin{pr}
Note that $\mu(E^{m}_{i})\leq \mu_{\max}(E)$. Then, although $E^{m}_{i}$ depends on $m$, its slope is bounded above and below by numbers which do not
depend on $m$, (cf. Proposition \ref{boundedness}) and furthermore it is a subsheaf of $E$. Hence, the set of possible isomorphism classes for $E^{m}_{i}$ is bounded.
Apply Serre Vanishing Theorem choosing $m_{3}\geq m_{2}$.
\end{pr}

\begin{prop}
Let $m\geq m_{3}$. For each filter $E_{i}^{m}$ in the $m$-Kempf filtration, we have $\dim V_{i}=h^{0}(E_{i}^{m}(m))$, therefore $V_{i}\cong H^{0}(E_{i}^{m}(m))$.
\label{task}
\end{prop}

\begin{pr}
Let $V_{\bullet}\subseteq V$ be the Kempf filtration of $V$ (cf. Theorem \ref{kempftheoremsheaves}) and let $E_{\bullet}^{m}\subseteq E$ be the $m$-Kempf filtration of $E$. We know that each $V_{i}$ generates the subsheaf $E_{i}^{m}$, by definition, then we have the following diagram:

$$\begin{array}{ccccccccccc}
    0 & \subset & V_{1} & \subset & V_{2} & \subset & \cdots & \subset & V_{t+1} & = & V \\
     & & \cap & & \cap & & & & & & ||  \\
      &     & H^{0}(E_{1}^{m}(m)) & \subset & H^{0}(E_{2}^{m}(m)) & \subset & \cdots & \subset & H^{0}(E_{t+1}^{m}(m)) & = & H^{0}(E(m))
      \end{array}$$

Suppose that there exists an index $i$ such that $V_{i}\neq H^{0}(E_{i}^{m}(m))$. Let $i$ be the index such that $V_{i}\neq H^{0}(E_{i}^{m}(m))$ and $\forall j>i$ it is $V_{j}=H^{0}(E_{j}^{m}(m))$. Then we have the diagram:

\begin{equation}
 \label{filtrationV}
    \begin{array}{ccccccccccccc}
    V_{i} & \subset & V_{i+1}\\
    \cap & & ||\\
    H^{0}(E_{i}^{m}(m)) & \subseteq & H^{0}(E_{i+1}^{m}(m))
    \end{array}
\end{equation}

Therefore $V_{i}\subsetneq H^{0}(E_{i}^{m}(m))\subseteq V_{i+1}$ and we can consider a new filtration by adding the filter $H^{0}(E_{i}^{m}(m))$:

\begin{equation}
\label{filtrationV'}
    \begin{array}{ccccccccccccc}
    V_{i} & \subset & H^{0}(E_{i}^{m}(m)) & \subset & V_{i+1}\\
    || & & || & & ||\\
    V'_{i} & \subset &  V'_{i+1} & \subset & V'_{i+2}
    \end{array}
\end{equation}

Note that we are in situation of Lemma \ref{lemmaB}, where
$W=H^{0}(E_{i}^{m}(m))$, filtration $V_{\bullet}$ is
\eqref{filtrationV} and filtration $V'_{\bullet}$ is
\eqref{filtrationV'}.

The graph associated to filtration $V_{\bullet}$, by Definition \ref{graph}, is given by the points
$$(b_{i},w_{i})=\big(\dfrac{\dim V_{i}}{m^{n}},\frac{m}{\dim V}(r\dim V_{i}-r_{i}\dim V)\big)\; ,$$
where the slopes of the graph are given by

$$-v_{i}=\frac{w^{i}}{b^{i}}=\frac{w_{i}-w_{i-1}}{b_{i}-b_{i-1}}=$$
$$\frac{m^{n+1}}{\dim V}\big(r-r^{i}\frac{\dim V}{\dim V^{i}}\big)\leq \frac{m^{n+1}}{\dim V}\cdot r:=R$$
and equality holds if and only if $r^{i}=0$.

Now, the new point which appears in the graph of the filtration $V'_{\bullet}$ is
$$Q=\big(\dfrac{h^{0}(E_{i}^{m}(m))}{m^{n}},\frac{m}{\dim V}(rh^{0}(E_{i}^{m}(m))-r_{i}\dim V)\big)\; .$$

Point $Q$ joins two new segments appearing in this new graph. The slope of the segment between $(b_{i},w_{i})$ and $Q$ is, by a
similar calculation,
$$-v'_{i+1}=\dfrac{m^{n+1}}{\dim V}\cdot r=R\; .$$

By Lemma \ref{lemmaA}, the graph is convex, so $v_{1}<v_{2}<\ldots<v_{t+1}$. As $E^{m}_{1}$ is a non-zero
torsion free sheaf, it has positive rank $r_{1}=r^{1}$ and so it follows $v_{1}>-R$. On the other hand, by Lemma
\ref{lemmaB}, $v'_{i+1}\geq v_{i+1}$. Hence
$$-R<v_{1}<v_{2}<\ldots<v_{i+1}\leq v'_{i+1}=-R\; ,$$
which is a contradiction.

Therefore, $\dim V_{i}=h^{0}(E_{i}^{m}(m))$, for every filter in
the $m$-Kempf filtration.
\end{pr}

\begin{cor}
\label{rank} For every filter $E_{i}^{m}$ in the $m$-Kempf filtration, it is $r^{i}=\rk
E_{i}^{m}/E_{i-1}^{m}>0$.
\end{cor}
\begin{pr}
In the proof of Proposition \ref{task} we have seen that $r^{i}=0$ is equivalent to $-v_{i}=R$. Then, the result
follows from that because it is $r^{1}=r_{1}>0$ and $-R<v_{1}<v_{2}<\ldots<v_{t+1}$.
\end{pr}

\subsection{Proof of Theorem \ref{kempfstabilizes}: the $m$-Kempf filtration stabilizes}
\label{sectionkempfstabilizes}

In Proposition \ref{regular} we have seen that, for any $m\geq m_{3}$,
all the filters $E^{m}_{i}$ of the $m$-Kempf filtration of $E$ are
$m_{3}$-regular. Hence, $E^{m}_{i}(m_{3})$ is generated by the
subspace $H^{0}(E^{m}_{i}(m_{3}))$ of $H^{0}(E(m_{3}))$, and the
filtration of sheaves
$$
0\subset E^{m}_{1} \subset E^{m}_{2} \subset \cdots \subset
E^{m}_{t_{m}} \subset E^{m}_{t_{m}+1}=E
$$
is the filtration associated to the filtration of vector spaces
$$
0\subset H^{0}(E^{m}_{1}(m_{3})) \subset H^{0}(E^{m}_{2}(m_{3}))
\subset \cdots \subset H^{0}(E^{m}_{t_{m}}(m_{3})) \subset
H^{0}(E^{m}_{t_{m}+1}(m_{3}))=H^{0}(E(m_{3}))
$$
by the evaluation map (c.f. Lemma \ref{mregularity}). Note that the dimension of the vector space
$H^{0}(E(m_{3}))$ does not depend on $m$ and, by Corollary
\ref{rank}, the length $t_{m}+1$ of the $m$-Kempf filtration of
$E$ is, at most, equal to $r$, the rank of $E$, a bound which does
not also depend on $m$. Note that, also because of Corollary \ref{rank}, each subsheaf in the $m$-Kempf filtration of $E$  is strictly contained
in the following one, for $m\geq m_{3}$. 

\begin{dfn}
\label{mtype}
We call \emph{$m$-type} to the tuple of different Hilbert
polynomials appearing in the $m$-Kempf filtration of $E$
$$
(P_{1}^{m},\ldots,P_{t_{m}+1}^{m})\; ,
$$
where $P_{i}^{m}:=P_{E_{i}^{m}}$.
\end{dfn}

Note that
$P^{i,m}:=P_{E^{m}_{i}/E^{m}_{i-1}}=P_{E^{m}_{i}}-P_{E^{m}_{i-1}}$,
so they are defined in terms of elements of each $m$-type.

\begin{prop}
\label{Pisfinite}
For all integers $m\geq m_{3}$, the set of
possible $m$-types
$$
\mathcal{P}=\big\{(P_{1}^{m},\ldots, P_{t_{m}+1}^{m})\big\}
$$
is finite.
\end{prop}
\begin{pr}
Once we fix $V\cong H^{0}(E(m_{3}))$ of dimension
$h^{0}(E(m_{3}))$ (which does not depend on $m$), all possible
filtrations by vector subspaces of $V$ are parametrized by a finite-type
scheme. Therefore the set of all possible $m$-Kempf filtrations of
$E$, for $m\geq m_{3}$, is bounded and $\mathcal{P}$ is finite.
\end{pr}

Recall that the vector $v$ can be recovered from the filtration $V_{\bullet}\subset V$ and the vector $\Gamma$
from the weights $n_{i}$. Then, given $m$, the $m$-Kempf filtration achieves the maximum for the Kempf function
$\mu(V_{\bullet},n_{\bullet})$ (c.f. (\ref{kempffunctionsheaves})), which is the same, by Proposition
\ref{identification}, that achieving the maximum for the function
$$
\mu_v(\Gamma)=\frac{(\Gamma,v)}{||\Gamma||}\; ,
$$
among all vectors $v$ coming from filtrations $V_{\bullet}\subset V$ and vectors $\Gamma\in \mathcal{C}-\{0\}$,
where
$$
\mathcal{C}= \big\{ x\in \mathbb{R}^{t+1} : x_1<x_2<\cdots
<x_{t+1} \big\}\; .
$$

By Definition \ref{graph} we associate a graph to the $m$-Kempf
filtration, given by $v_{m}$. Recall that, by Lemma \ref{lemmaA}
the graph is convex, meaning $v_{m}\in \mathcal{C}$, which implies
$\Gamma_{v_{m}}=v_{m}$ by Lemma \ref{close}. Then, given $v_{m}$
associated to the $m$-Kempf filtration
\begin{equation}
\label{maxvalue}
\max_{\Gamma\in \overline{\mathcal{C}}}\mu_{v_{m}}(\Gamma)=\mu_{v_{m}}(\Gamma_{v_{m}})=\frac{(\Gamma_{v_{m}},v_{m})}{||\Gamma_{v_{m}}||}=
\frac{(v_{m},v_{m})}{||v_{m}||}=||v_{m}||\; ,
\end{equation}
where recall that we defined in Definition \ref{graph}
$$
v_{m,i}=m^{n+1}\cdot \frac{1}{\dim V^{i}\dim V}\big[r^{i}\dim V-r\dim
V^{i}\big]\; ,
$$
$$
b_{m}^{i}=\frac{1}{m^{n}}\cdot \dim V^{i}\; ,
$$
and, thanks to Propositions \ref{regular} and \ref{task}, we can rewrite
$$
v_{m,i}=m^{n+1}\cdot
\frac{1}{P^{i,m}(m)P(m)}\big[r^{i}P(m)-rP^{i,m}(m)\big]\; ,
$$
$$
b_{m}^{i}=\frac{1}{m^{n}}\cdot P^{i,m}(m)\; .
$$

Let
$$
v_{m,i}(l)=m^{n+1}\cdot
\frac{1}{P^{i,m}(l)P(l)}\big[r^{i}P(l)-rP^{i,m}(l)\big]
$$
be the coordinates of the graph associated to the $m$-Kempf filtration but where the polynomials are evaluated at
another variable $l$. Let us define
$$
\Theta_{m}(l)=(\mu_{v_{m}(l)}(\Gamma_{v_{m}(l)}))^{2}=||v_{m}(l)||^{2}\;
,
$$
where the second equality follows by an argument similar to
\eqref{maxvalue}. Note that $\Theta_m(l)$ is a rational
function on $l$. Let
$$
\mathcal{A}=\{\Theta_{m}:m\geq m_{3}\}
$$
which is a finite set by Proposition \ref{Pisfinite}. We say that
$f_{1}\prec f_{2}$ for two rational functions, if the inequality
$f_{1}(l)<f_{2}(l)$ holds for $l\gg 0$, and let $K$ be the maximal
function in the finite set $\mathcal{A}$, with respect to the
defined ordering.

Note that the value $\Theta_{m}(m)$ is the square of the maximum
of the Kempf function $\mu_{v_{m}}(\Gamma)$, by \eqref{maxvalue},
achieved for the maximal filtration $V_{\bullet}\subset V\simeq
H^{0}(E(m))$ of vector spaces which gives the vector $v_{m}$. This
weighted filtration is the only one which gives the value
$\sqrt{\Theta_{m}(m)}$ for the Kempf function.

\begin{lem}
\label{uniquefunction} There exists an integer $m_{4}\geq m_{3}$
such that $\forall m\geq m_{4}$, $\Theta_{m}=K$.
\end{lem}
\begin{pr}
Choose $m_{4}$ such that $K(l) \geq \Theta_{m}(l)$, $\forall l\geq
m_{4}$ and every $\Theta_{m}\in \mathcal{A}$ with equality only
when $\Theta_{m}= K$. Let $m\geq m_{4}$. Given that the Kempf
function achieves the maximum over all possible filtrations and
weights (c.f. Theorem \ref{kempftheoremsheaves}), we have $\Theta_{m}(m)\geq
K(m)$, because $K$ is another rational function built with other
$m'$-type, i.e., other values for the polynomials appearing on the
rational function (c.f. Definition \ref{mtype}). Combining both inequalities we obtain
$\Theta_m(m)=K(m)$ for all $m\geq m_{4}$.
\end{pr}

\begin{prop}
\label{eventually} Let $l_{1}$ and $l_{2}$ be integers with
$l_{1}\geq l_{2} \geq m_{4}$. Then the $l_{1}$-Kempf filtration of
$E$ is equal to the $l_{2}$-Kempf filtration of $E$.
\end{prop}

\begin{pr}
By construction, the filtration
\begin{equation}
  \label{filt1}
  H^{0}(E^{l_{1}}_{1}(l_{1})) \subset H^{0}(E^{l_{1}}_{2}(l_{1})) \subset \cdots \subset H^{0}(E^{l_{1}}_{t_{1}}(l_{1}))
  \subset H^{0}(E^{l_{1}}_{t_{1}+1}(l_{1}))=H^{0}(E(l_{1}))
\end{equation}
is the $l_{1}$-Kempf filtration of $V\simeq H^{0}(E(l_{1}))$. Now
consider the filtration $V'_{\bullet}\subset V\simeq
H^{0}(E(l_{1}))$ defined as follows
\begin{equation}
\label{filt2}
H^{0}(E^{l_{2}}_{1}(l_{1})) \subset H^{0}(E^{l_{2}}_{2}(l_{1})) \subset \cdots \subset H^{0}(E^{l_{2}}_{t_{2}}(l_{1})) \subset
H^{0}(E^{l_{2}}_{t_{2}+1}(l_{1}))=H^{0}(E(l_{1})) \; .
\end{equation}
We have to prove that \eqref{filt2} is in fact the $l_{1}$-Kempf
filtration of $V\simeq H^{0}(E(l_{1}))$.

Since $l_{1},l_{2}\geq m_{4}$, by Lemma \ref{uniquefunction} we
have $\Theta_{l_{1}}=\Theta_{l_{2}}=K$. Then,
$\Theta_{l_{1}}(l_{1})=\Theta_{l_{2}}(l_{1})$ and, by uniqueness
of the Kempf filtration (c.f. Theorem \ref{kempftheoremsheaves}), the
filtrations
 \eqref{filt1} and \eqref{filt2} coincide. Since, in particular $l_{1},l_{2}\geq m_{3}$, $E_{i}^{l_{1}}$ and $E_{i}^{l_{2}}$
 are $l_{1}$-regular by Proposition \ref{regular}. Hence, $E_{i}^{l_{1}}(l_{1})$ and $E_{i}^{l_{2}}(l_{1})$ are
 generated by their global sections (c.f. Lemma \ref{mregularity})
 $H^{0}(E_{i}^{l_{1}}(l_{1}))$ and $H^{0}(E_{i}^{l_{2}}(l_{1}))$, respectively. By the
 previous argument, $H^{0}(E_{i}^{l_{1}}(l_{1}))=H^{0}(E_{i}^{l_{2}}(l_{1}))$, therefore $E_{i}^{l_{1}}(l_{1})=E_{i}^{l_{2}}(l_{1})$. By
tensoring with $\mathcal{O}_{X}(-l_{1})$, this implies that
 the filtrations $E^{l_{1}}_{\bullet} \subset E$
  and $E^{l_{2}}_{\bullet} \subset E$ coincide.
\end{pr}

Therefore, Theorem \ref{kempfstabilizes} follows from Proposition \ref{eventually} and it is proved that, eventually, the Kempf filtration does not depend on the integer $m$.

\begin{dfn}
\label{kempffiltration}
If $m\geq m_{4}$, the $m$-Kempf filtration of $E$ is called
\emph{the Kempf filtration of $E$},
$$0\subset E_{1} \subset E_{2} \subset \cdots \subset E_{t} \subset E_{t+1}=E\; .$$
\end{dfn}

\subsection{Proof of Theorem \ref{kempfisHN}: Kempf filtration is Harder-Narasimhan filtration}

Recall that the Kempf theorem (c.f. Theorem \ref{kempftheoremsheaves}) asserts
that given an integer $m$ and $V\simeq H^{0}(E(m))$, there exists a
unique weighted filtration of vector spaces $V_{\bullet}\subseteq V$
which gives maximum for the Kempf function
$$
\mu(V_{\bullet},n_{\bullet})=
\frac {\sum_{i=1}^{t+1} \frac{\Gamma_i}{\dim V} ( r^i \dim V - r\dim V^i)}
{\sqrt{\sum_{i=1}^{t+1} {\dim V^{i}} \Gamma_{i}^{2}}}\; .
$$
This filtration induces a filtration of sheaves, called the Kempf
filtration of $E$,
$$
0\subset E_{1} \subset E_{2} \subset \cdots \subset E_{t} \subset
E_{t+1}=E
$$
which is independent of $m$, for $m\geq m_{4}$, by Proposition
\ref{eventually}, hence it only depends on $E$. From now on, we assume
$m\geq m_4$.

Based on the fact we can rewrite the Kempf function as a certain scalar product divided by a norm (c.f.
Proposition \ref{identification}), we shave seen that the Kempf filtration is encoded by a graph with two
convexity properties (c.f. Lemmas \ref{lemmaA} and \ref{lemmaB}). We can express the data related to the
filtration of vector spaces with the data of the filtration of sheaves. Since $m\geq m_3$, the sheaves $E_{i}$
and $E^i$ are $m$-regular $\forall i$, and
\begin{equation}
\label{substitutions}
    \begin{array}{c}
    \dim V_{i}=h^{0}(E_{i}(m))=P_{E_{i}}(m) =: P_i(m)\\
\dim V^{i}=h^{0}(E^{i}(m))=P_{E^{i}}(m) =: P^i(m)
    \end{array}
\end{equation}
(c.f. Proposition \ref{regular} and Proposition \ref{task}). Recall
that the Kempf function is a function on $m$, with order
$m^{-\frac{n}{2}-1}$ at zero (c.f. Proposition \ref{identification})
then we consider the function $K$, where
$$
K(m)=m^{\frac{n}{2}+1}\cdot
\mu(V_{\bullet},m_{\bullet})=\mu_{v_{m}}(\Gamma)\; .
$$
Making the substitutions \eqref{substitutions},
and using the relation $\gamma_{i}=\frac{r}{P}\Gamma_{i}$ (c.f. (\ref{filtV}) and (\ref{filtE})),
$$
K(m)=m^{\frac{n}{2}+1}\cdot \frac {\sum_{i=1}^{t+1} \frac{\gamma_{i}}{r}[(r^i P - rP^i)]}
 {\sqrt{\sum_{i=1}^{t+1} P^{i}\frac{P^{2}}{r^{2}} \gamma_{i}^{2}}}\; ,
$$
which is a function on $m$ whose square is a rational function (since $P$ and $P^i$ are
polynomials on $m$).
Therefore we get
$$
K(m)=m^{\frac{n}{2}+1}\cdot \frac{1}{P}\frac {\sum_{i=1}^{t+1}
\gamma_{i}[r^{i}P-rP^{i}]}
 {\sqrt{\sum_{i=1}^{t+1} P^{i}\gamma_{i}^{2}}}\; .
$$

\begin{prop}
\label{finalfunction} Given a sheaf $E$, there exists a unique filtration
$$0\subset E_{1} \subset E_{2} \subset \cdots \subset E_{t} \subset
E_{t+1}=E$$ with positive weights $n_{1},\ldots,n_{t}$,
$n_{i}=\frac{\gamma_{i+1}-\gamma_{i}}{r}$, which gives maximum for
the function
$$K(m)=m^{\frac{n}{2}+1}\cdot \frac {\sum_{i=1}^{t+1}
P^{i}\gamma_{i}[\frac{r^{i}}{P^{i}}-\frac{r}{P}]}
 {\sqrt{\sum_{i=1}^{t+1} P^{i}\gamma_{i}^{2}}}\; .$$
\end{prop}

Similarly, we have defined the coordinates $v_{i}$ (slopes of segments of the graph), as
$$v_{i}=m^{n+1}\cdot \big[\frac{r^{i}}{P^{i}}-\frac{r}{P}\big]$$
(c.f. Definition \ref{graph}). Therefore we can express the function $K$ as
$$K(m)=m^{-\frac{n}{2}}\cdot \frac {\sum_{i=1}^{t+1}
P^{i}\gamma_{i}v_{i}} {\sqrt{\sum_{i=1}^{t+1}
P^{i}\gamma_{i}^{2}}}=m^{-\frac{n}{2}}\cdot
\frac{(\gamma,v)}{||\gamma||}\; ,$$
where the scalar product is
given by the diagonal matrix
$$\left(
    \begin{array}{cccc}
      P^{1} &  &  & 0 \\
       & P^{2} &  &  \\
       &  & \ddots &  \\
      0 &  &  & P^{t+1} \\
    \end{array}
  \right)$$

Finally, we use Lemmas \ref{lemmaA} and \ref{lemmaB} to show that the Kempf filtration verifies the two
properties of the Harder-Narasimhan filtration for sheaves (c.f. Theorem \ref{HNunique}) hence, by uniqueness,
both filtrations have to coincide.

\begin{prop}
\label{descendentslopes} Given the Kempf filtration of a sheaf
 $E$,
$$0\subset E_{1} \subset E_{2} \subset \cdots \subset E_{t} \subset E_{t+1}=E$$
it verifies
$$\frac{P^{1}}{r^{1}}>\frac{P^{2}}{r^{2}}>\ldots>\frac{P^{t+1}}{r^{t+1}}$$
\end{prop}
\begin{pr}
The coordinates of the vector $v$ associated to the filtration are, for $m$ large enough, $v_{i}=m^{n+1}\cdot
(\frac{r^{i}}{P^{i}}-\frac{r}{P})$. Now apply Lemma \ref{lemmaA} which says that $v$ is convex, i.e.
$v_{1}<\ldots<v_{t+1}$.
\end{pr}

\begin{prop}
\label{blocksemistability} Given the Kempf filtration of a sheaf $E$,
$$0\subset E_{1} \subset E_{2} \subset \cdots \subset E_{t} \subset E_{t+1}=E\; ,$$
each one of the blocks $E^{i}=E_{i}/E_{i-1}$ is semistable.
\end{prop}
\begin{pr}
Consider the graph associated to the Kempf filtration of $E$.
Suppose that any of the blocks has a destabilizing subsheaf. Then,
it corresponds to a point above of the graph of the filtration.
The graph obtained by adding this new point is a refinement of the
graph of the Kempf filtration, whose convex envelope is not the
original graph, which contradicts Lemma \ref{lemmaB}.
\end{pr}

\begin{cor}[c.f. Theorem \ref{kempfisHN}]
The Kempf filtration of a sheaf $E$ coincides with the Harder-Narasimhan filtration.
\end{cor}
\begin{pr}
By Propositions \ref{descendentslopes} and \ref{blocksemistability} the Kempf filtration verifies the two
properties of the Harder-Narasimhan filtration. By uniqueness of the Harder-Narasimhan filtration (c.f. Theorem
\ref{HNunique}) both filtrations do coincide.
\end{pr}

\section{Holomorphic pairs}
\label{kempfpairs}

In this section we prove the correspondence between the Kempf filtration and the Harder-Narasimhan filtration
for holomorphic pairs. It follows the scheme of the proof given for torsion free coherent sheaves in section
\ref{kempfsheaves}. First, we give some definitions and the notion of stability for the construction of the
moduli space of holomorphic pairs. It can be deduced from the construction of the moduli space of tensors in section \ref{exampletensors}. 

\vspace{1cm}

Let $X$ be a smooth complex projective variety. Let us consider
\emph{holomorphic pairs}
$$(E,\varphi:E\to \mathcal{O}_{X})$$
given by a coherent torsion free sheaf of rank $r$ with fixed
determinant $\det(E)\cong \Delta$ and a morphism to a the
structure sheaf $\mathcal{O}_{X}$. Note that the definition of
holomorphic pair coincides with the definition of tensor in
Definition \ref{deftensor}, with $s=1$, $c=1$, $b=0$, $R=\Spec
\mathbb{C}$ and $\mathcal{D}=\mathcal{O}_{X}$ is the structure
sheaf over $X\times R\simeq X$.

A \emph{weighted filtration} $(E_\bullet,n_\bullet)$ of a sheaf $E$ of rank $r$ is a filtration
\begin{equation}
\label{filtEpairs} 0 \subset E_{1} \subset E_{2} \subset
\;\cdots\; \subset E_{t} \subset E_{t+1}=E,
\end{equation}
and rational numbers $n_1,\, n_2,\ldots , \,n_t > 0$.

Let $\delta$ be a polynomial of degree at most $\dim X-1$ and
positive leading coefficient. We rephrase Definition
\ref{stabtensors} for the case of holomorphic pairs. See also the
calculation made in (\ref{rightmu2}).

\begin{dfn}
\label{stabilitypairs} A holomorphic pair $(E,\varphi)$ is \emph{$\delta$-semistable} if for all weighted
filtrations $(E_{\bullet},n_{\bullet})$ (c.f. (\ref{filtEpairs})),
$$\sum_{i=1}^{t}  n_i ( r  P_{E_i} - r_i P_E) +\delta \sum_{i=1}^t
n_i \big(r_i -\epsilon(E_{i}) r \big)\leq 0\; ,$$ where
$\epsilon(E_{i})=1$ if $\varphi|_{E_{i}}\neq 0$ and
$\epsilon(E_{i})=0$ otherwise. If the strict inequality holds for every weighted filtration, we say that
$(E,\varphi)$ is \emph{$\delta$-stable}. If $(E,\varphi)$ is not $\delta$-semistable, we say that it is \emph{$\delta$-unstable}.
\end{dfn}

\begin{dfn}
\label{subpairs} Given a holomorphic pair $(E,\varphi:E\rightarrow \mathcal{O}_{X})$, let $(E',\varphi|_{E'})$ be
a \emph{subpair} where $E'\subset E$ is a subsheaf and $\varphi|_{E'}$ is the restriction of the morphism
$\varphi$. Let $E''=E/E'$ and define the holomorphic pair $(E'',\varphi|_{E''})$ where, if
$\varphi|_{E'}\neq 0$, define $\varphi|_{E''}:=0$, and if $\varphi|_{E'}=0$, $\varphi|_{E''}$ is the induced
morphism in the quotient sheaf. We call $(E'',\varphi|_{E''})$ a \emph{quotient pair} of $(E,\varphi)$. For every pair $(E,\varphi:E\rightarrow \mathcal{O}_{X})$, define
$\epsilon(E)=1$ if $\varphi|_{E}\neq 0$ and $\epsilon(E)=0$ otherwise. Recall that we define a morphism of pairs
$(E,\varphi)\rightarrow (F,\psi)$ as a morphism of sheaves $\alpha:E\rightarrow F$ such that $\psi\circ
\alpha=\varphi$ (c.f. Definition \ref{deftensor}).
\end{dfn}

\begin{dfn}
\label{corrected} Let $(E,\varphi:E\rightarrow \mathcal{O}_{X})$ be a holomorphic pair. We define the \emph{corrected Hilbert polynomial}
of $(E,\varphi)$ as
$$\overline{P}_{E}:=P_{E}-\delta\epsilon(E)$$
\end{dfn}
Note that the exact sequence of sheaves
$$0\rightarrow E'\rightarrow E\rightarrow E''\rightarrow 0$$ verify
$$\overline{P}_{E}=\overline{P}_{E'}+\overline{P}_{E''}$$ for the corrected polynomials.

\begin{rem}
Note that the definition of quotient pair in Definition \ref{subpairs} does not imply that 
$$0\rightarrow (E',\varphi|_{E'})\rightarrow (E,\varphi)\rightarrow (E'',\varphi|_{E''})\rightarrow 0$$
is an exact sequence in the category of tensors, where $E'\subset E$ and $E'':=E/E'$. Nonetheless, we keep that definition 
for the additivity of the corrected Hilbert polynomials to hold on exact sequences of sheaves. 
\end{rem}

From Definition \ref{stabilitypairs} it can be directly deduced the
following equivalent definition, which appears on \cite[Definition 1.1]{HL2}.
\begin{prop}
\label{subobjects} A pair $(E,\varphi)$ is $\delta$-unstable if
and only if there exists a subpair $(F,\varphi|_{F})$ with
$\frac{\overline{P}_{F}}{\rk F}>\frac{\overline{P}_{E}}{\rk E}$.
\end{prop}
\begin{pr}
If $(E,\varphi)$ is $\delta$-unstable, there exists a filtration
$$0\subset E_{1} \subset E_{2} \subset \cdots \subset E_{t} \subset
E_{t+1}=E$$ and weights $n_{i}>0$ such that
$$\sum_{i=1}^{t}  n_i ( r  P_{E_i} - r_i P_E) +\delta \sum_{i=1}^t
n_i \big(r_i -\epsilon(E_{i}) r \big)=$$
$$\sum_{i=1}^{t}n_{i}\big(r(P_{E_{i}}-\delta\epsilon(E_{i}))-r_{i}(P_{E}-\delta)\big)=
\sum_{i=1}^{t}n_{i}(r\overline{P}_{E_{i}}-r_{i}\overline{P}_{E})>0\;
.$$ As the weights $n_{i}$ are positive, there exists any $i$ such
that
$$r\overline{P}_{E_{i}}-r_{i}\overline{P}_{E}>0\Leftrightarrow$$
$$\frac{\overline{P}_{E_{i}}}{\rk E_{i}}>\frac{\overline{P}_{E}}{\rk E}\; .$$

On the other hand, if there exists $(F,\varphi|_{F})$ with
$\frac{\overline{P}_{F}}{\rk F}>\frac{\overline{P}_{E}}{\rk E}$,
the one-step filtration
$$0\subset F\subset E$$ gives a positive quantity in the expression
of Definition \ref{stabilitypairs}. Therefore $(E,\varphi)$ is
$\delta$-unstable.
\end{pr}

\subsection{Moduli space of holomorphic pairs}

We recall the construction of the moduli space of
$\delta$-semistable pairs with fixed polynomial $P$ and fixed
determinant $\det(E)\simeq \Delta$. This was done in \cite{HL1}
following Gieseker's ideas, and in \cite{HL2} following Simpson's
ideas. Here, we use Gieseker's method (although \cite{HL1} assumes
that $X$ is a curve or a surface, thanks to Simpson's bound
\cite[Corollary 1.7]{Si1}, we can follow Gieseker's method for any
dimension). As we said at the beginning of the section, the
construction can be derived from the construction of a moduli
space for tensors in section \ref{exampletensors}, where $s=1$,
$c=1$, $b=0$, $R=\Spec \mathbb{C}$ and
$\mathcal{D}=\mathcal{O}_{X}$ is the structure sheaf over $X\times
R\simeq X$, (c.f. Definition \ref{deftensor}).

Let $m$ be an integer, so that $E$ is $m$-regular for all
semistable $E$ (c.f. \cite[Corollary 3.3.1 and Proposition
3.6]{Ma1}). Let $V$ be a vector space of dimension $p:=P(m)$.
Given an isomorphism $V\cong H^0(E(m))$, we obtain a quotient
$$
q:V\otimes \SO_X(-m) \surj E\; ,
$$
hence a homomorphism
$$
Q:\wedge^r V \cong \wedge^r H^0(E(m)) \too H^0(\wedge^r(E(m)))
\cong H^0(\Delta(rm))=:A
$$
and points
$$
Q\in \Hom(\wedge^r V , A) \qquad\overline{Q} \in \PP(\Hom(\wedge^r
V , A) )\; .$$ The morphism $\varphi:E\too \mathcal{O}_{X}$
induces a homomorphism
$$
\Phi:V=H^0(E(m)) \too H^0(\mathcal{O}_{X}(m))=:B
$$
and hence points
$$
\Phi\in \Hom(V,B) \qquad \overline{\Phi}\in \PP(\Hom(V,B))\; .
$$
If we change the isomorphism $V\cong H^0(E(m))$ by a homothecy, we
obtain another point in the line defined by $Q$, but the point $\overline{Q}$ does not change, and similarly for $\overline{\Phi}$.

Two different isomorphisms $V\cong
H^0(E(m))$ differ by an element of $SL(V)$, hence this group acts on the two projective spaces we have defined. We choose a polarization $\mathcal{O}(a_{1},a_{2})$
(c.f. (\ref{polarization})) to give a linearization of the action of $SL(V)$. By the Hilbert-Mumford criterion (c.f. Theorem \ref{HMcrit}), a point
$$(\overline{Q},\overline{\Phi}) \in
\PP(\Hom(\wedge^r V , A) )\times \PP(\Hom(V,B))
$$
is \emph{GIT semistable} with respect to the natural linearization on
$\mathcal{O}(a_1,a_2)$ if and only if for all weighted filtrations
it is
$$
\mu(\overline{Q},V_\bullet,n_\bullet) + \frac{a_2}{a_1}
\mu(\overline{\Phi},V_\bullet, n_\bullet) \leq 0\;  ,
$$
where each numerical function $\mu$ is the calculation of the
minimal relevant weight of the action of a $1$-parameter subgroup
$\Gamma$ on each projective space. Recall from section
\ref{exampletensors} the correspondence between $1$-parameter
subgroups and weighted filtrations $(V_{\bullet},n_{\bullet})$.

\begin{prop}
\label{HMcritpairs} A point $(\overline{Q},\overline{\Phi})$ is
\emph{GIT $a_2/a_1$-semistable} if for all weighted filtrations
$(V_\bullet,n_\bullet)$ we have
$$
\sum_{i=1}^{t}  n_i ( r \dim V_i - r_i \dim V) +\frac{a_2}{a_1}
\sum_{i=1}^t n_i \big(\dim V_i-\epsilon_i(\overline\Phi)\dim V
\big) \leq 0\; .
$$
\end{prop}
\begin{pr}
C.f. Proposition \ref{GITstab}.
\end{pr}

Recall that $r_{i}$ is the rank of the subsheaf $E_{i}\subset E$ generated by $V_{i}$ by the evaluation map. Also recall that, if $j$
is the index giving minimum in \eqref{rightmuV}, we will define
$\epsilon_i(\overline{\Phi},V_{\bullet})=1$ if $i\geq j$ and
$\epsilon_i(\overline{\Phi},V_{\bullet})=0$ otherwise. We will
denote $\epsilon_i(\overline{\Phi},V_{\bullet})$ by just
$\epsilon_i(\overline\Phi)$ if the filtration $V_{\bullet}$ is
clear from the context. Let us call
$\epsilon^i(\overline\Phi)=
\epsilon_i(\overline\Phi)-\epsilon_{i-1}(\overline\Phi)$
and note that
$\epsilon_{i}(\overline{\Phi},V_{\bullet})=\epsilon(E_{V_{i}})$,
in Definition \ref{subpairs}.

\begin{rem}
\label{independenceofweightspairs} Note that the definition of $\epsilon_i(\overline\Phi)$ is independent of the
weights $n_{\bullet}$ or the vector $\Gamma$ associated to them. Indeed,
$\epsilon_{i}(\overline{\Phi},V_{\bullet})=\epsilon(E_{V_{i}})$ only depends on the vanishing of the morphism
$\varphi$ on the subsheaves $E_{V_{i}}$ (c.f. Definition \ref{subpairs}).
\end{rem}

\begin{thm}
\label{GIT-deltapairs}
Let $(E,\varphi)$ be a holomorphic pair.
There exists an integer $m_0$ such that, for $m\geq m_0$, the
associated point $(\overline{Q},\overline{\Phi})$ is GIT
$a_2/a_1$-semistable if and only if the pair is
$\delta$-semistable, where
$$
\frac{a_2}{a_1} = \frac{r\delta(m)}{P_E(m)-\delta(m)}
$$
\end{thm}
\begin{pr}
C.f. Theorem \ref{GIT-delta}.
\end{pr}

Let $(E,\varphi)$ be a $\delta$-unstable holomorphic pair. Let
$m_{0}$ be the integer in Theorem \ref{GIT-deltapairs} (i.e. such
that the $\delta$-stability of the tensor coincides with the GIT
stability). If necessary, choose another $m_{0}$ such that the
sheaf $E$ is $m_{0}$-regular.

Let $m\geq m_{0}$ be an integer and let $V$ be a vector space of
dimension $P(m)=h^{0}(E(m))$. Fix an isomorphism $V\simeq
H^{0}(E(m))$. Given a filtration of vector subspaces $0\subset V_{1}\subset
\cdots \subset V_{t+1}= V$ and positive numbers
 $n_{1},\cdots,n_{t}>0$, i.e., given a weighted filtration, we define now the function
$$\mu(V_{\bullet},n_{\bullet})=\frac{\sum_{i=1}^{t}  n_{i} (r\dim V_{i}-r_{i}\dim V)+
\frac{a_{2}}{a_{1}}\sum_{i=1}^{t} n_{i} \big(\dim
V_{i}-\epsilon_{i}(\overline\Phi)\dim V \big)(\leq) 0}
{\sqrt{\sum_{i=1}^{t+1} {\dim V^{i}_{}} \Gamma_{i}^{2}}}\; ,$$
which is a \emph{Kempf function} for this problem, as in the case of sheaves (c.f. Definition \ref{kempffunction}).

We can apply Theorem \ref{kempftheoremsheaves} to obtain
\begin{equation}
\label{vkempf-filtpairs} 0\subset V_{1}\subset \cdots \subset
V_{t+1}= V\; ,
\end{equation}
the \emph{Kempf filtration} of $V$. Let
\begin{equation}
\label{ekempf-filtpairs}
0\subseteq (E^{m}_{1},\varphi|_{E_{1}^{m}})\subseteq
(E^{m}_{2},\varphi|_{E_{2}^{m}})\subseteq\cdots
(E^{m}_{t},\varphi|_{E_{t}^{m}})\subseteq
(E^{m}_{t+1},\varphi|_{E_{t+1}^{m}})\subseteq (E,\varphi)
\end{equation}
be the \emph{$m$-Kempf filtration of the pair $(E,\varphi)$}, where $E_{i}^{m}\subset E$ is the subsheaf generated by $V_{i}$
under the evaluation map.

We will apply the same techniques as in section \ref{kempfsheaves}
to prove the following theorem:

\begin{thm}
\label{kempfstabilizespairs}
There exists an integer $m'\gg 0$ such that the $m$-Kempf filtration of the holomorphic pair $(E,\varphi)$ is independent of $m$, for $m\geq m'$.
\end{thm}

\subsection{The $m$-Kempf filtration stabilizes with $m$}

In this subsection we give a proof of Theorem
\ref{kempfstabilizespairs}, based on the same arguments as in the
case of sheaves. As we did in section \ref{kempfsheaves}, we
associate a graph to the $m$-Kempf filtration of a
$\delta$-unstable pair $(E,\varphi)$, to relate the Kempf function
with the function $\mu_{v}(\Gamma)$ in Theorem
\ref{maxconvexenvelope}.

\begin{dfn}
\label{graphpairs} Let $m\geq m_{0}$. Given $0\subset V_{1}\subset
\cdots \subset V_{t+1}= V$ a filtration of vector spaces of $V$, let
$$v_{m,i}=m^{n+1}\cdot \frac{1}{\dim V^{i}\dim V}\big[r^{i}\dim V-r\dim V^{i}+\dfrac{a_{2}}{a_{1}}(\epsilon^{i}(\overline\Phi)\dim V-
\dim V^{i})\big]\; ,$$
$$b_{m}^{i}=\dfrac{1}{m^{n}}\dim V^{i}>0$$
$$w_{m}^{i}=-b_{m}^{i}\cdot v_{m,i}=m\cdot \frac{1}{\dim V}\big[r\dim V^{i}-r^{i}\dim V+\dfrac{a_{2}}{a_{1}}(\dim V^{i}-\epsilon^{i}
(\overline\Phi)\dim V) \big]\; .$$
Also let
$$b_{m,i}=b_{m}^{1}+\ldots +b_{m}^{i}=\dfrac{1}{m^{n}}\dim V_{i}\; ,$$
$$w_{m,i}=w_{m}^{1}+\ldots +w_{m}^{i}=m\cdot \frac{1}{\dim V}\big[r\dim V_{i}-r_{i}\dim V+\dfrac{a_{2}}{a_{1}}(\dim V_{i}-\epsilon_{i}(\overline\Phi)\dim V)\big]\; .$$
We call the graph defined by points $(b_{m,i},w_{m,i})$ the
\emph{graph associated to the filtration} $V_{\bullet}\subset V$.
\end{dfn}

Now, applying Proposition \ref{identification}, we can identify as
well the new Kempf function in Theorem \ref{kempftheoremsheaves},
$$\mu(V_{\bullet},n_{\bullet})=\frac{\sum_{i=1}^{t}  n_{i} (r\dim V_{i}-r_{i}\dim V)+\frac{a_{2}}{a_{1}}
\sum_{i=1}^{t} n_{i} \big(\dim
V_{i}-\epsilon_{i}(\overline\Phi)\dim V
\big)}{\sqrt{\sum_{i=1}^{t+1} {\dim V^{i}} \Gamma_{i}^{2}}}\; ,$$
with the function in Theorem \ref{maxconvexenvelope}, where the coordinates of the
graph are now given as in Definition \ref{graphpairs}.

We will use Lemmas \ref{lemmaA} and \ref{lemmaB}, to give the
analogous to Propositions \ref{boundedness} and \ref{task} for the
case of holomorphic pairs.

Let
\begin{equation} \label{C_constant_pairs} C=\max\{r|\mu
_{\max }(E)|+\frac{d}{r}+r|B|+|A|+\delta_{n-1}(n-1)!+1\;,\;1\}
\end{equation}
be positive constant, where $\delta_{n-1}$ is the
$(n-1)^{th}$-degree coefficient of the polynomial $\delta(m)$ (if
$\deg (\delta)<n-1$, then set $\delta_{n-1}=0$).
\begin{prop}
\label{boundednesspairs} Given a sufficiently large $m$, each
filter $E_{i}^{m}$ in the $m$-Kempf filtration of $(E,\varphi)$ (cf. \eqref{ekempf-filtpairs})  has slope
$\mu(E^{m}_{i})\geq \dfrac{d}{r}-C$.
\end{prop}

\begin{pr}
Choose an integer $m_{1}$ such that for $m\geq m_{1}$
$$[\mu_{max}(E)+gm+B]_{+}=\mu_{max}(E)+gm+B$$
and
$$[\frac{d}{r}-C+gm+B]_{+}=\frac{d}{r}-C+gm+B\; .$$
Let $m_{2}$ be such that $P_{E}(m)-\delta(m)>0$ for $m\geq m_{2}$.
Now consider $m\geq \max\{m_{0},m_{1},m_{2}\}$ and let
$$0\subseteq (E^{m}_{1},\varphi|_{E_{1}^{m}})\subseteq
(E^{m}_{2},\varphi|_{E_{2}^{m}})\subseteq\cdots
(E^{m}_{t},\varphi|_{E_{t}^{m}})\subseteq
(E^{m}_{t+1},\varphi|_{E_{t+1}^{m}})\subseteq (E,\varphi)$$
be the $m$-Kempf filtration of $(E,\varphi)$.

Suppose we have a filter $E^{m}_{i}\subseteq E$, of rank $r_{i}$
and degree $d_{i}$, such that $\mu (E^{m}_{i})<\frac{d}{r}-C$. The
subsheaf $E_{i}^{m}(m)\subset E(m)$ satisfies the estimate in
Lemma \ref{Simpson},
$$h^{0}(E^{m}_{i}(m))\leq \frac{1}{g^{n-1}n!}\big ((r_{i}-1)([\mu_{max}(E^{m}_{i})+gm+B]_{+})^{n}+([\mu_{min}(E^{m}_{i})+gm+B]_{+})^{n}\big)\; ,$$
where $\mu_{max}(E^{m}_{i}(m))=\mu_{max}(E^{m}_{i})+gm$ and
similarly for $\mu_{min}$.

Note that $\mu_{max}(E^{m}_{i})\leq \mu_{max}(E)$ and
$\mu_{min}(E^{m}_{i})\leq \mu(E^{m}_{i})<\frac{d}{r}-C$, so
$$h^{0}(E^{m}_{i}(m))\leq \frac{1}{g^{n-1}n!}\big ((r_{i}-1)([\mu_{max}(E)+gm+B]_{+})^{n}+([\frac{d}{r}-C+gm+B]_{+})^{n}\big)\; ,$$
and, by choice of $m$,
$$h^{0}(E^{m}_{i}(m))\leq \frac{1}{g^{n-1}n!}\big ((r_{i}-1)(\mu_{max}(E)+gm+B)^{n}+(\frac{d}{r}-C+gm+B)^{n}\big)=G(m)\; ,$$
where
$$G(m)=\frac{1}{g^{n-1}n!}\big
[r_{i}g^{n}m^{n}+ng^{n-1}\big((r_{i}-1)\mu_{max}(E)+\frac{d}{r}-C+r_{i}B\big)m^{n-1}+\cdots
\big ]\; .$$

Recall that, by Definition \ref{graph}, to the filtration
\eqref{vkempf-filtpairs}
we
associate a graph with heights, for each $j$
$$w_{j}=w^{1}+\ldots +w^{j}=m\cdot \frac{1}{\dim V}\big[r\dim V_{j}-r_{j}\dim V+\dfrac{a_{2}}{a_{1}}(\dim V_{j}-\epsilon_{j}(\overline\Phi)\dim V )\big]\; .$$

We will show that $w_{i}<0$ and will get a contradiction as in
Proposition \ref{boundedness}. Since $E^{m}_{i}(m)$ is generated
by $V_{i}$ under the evaluation map, it is $\dim V_{i}\leq
h^{0}(E^{m}_{i}(m))$, hence
$$w_{i}=\frac{m}{\dim V}\big[r\dim V_{i}-r_{i}\dim V+\dfrac{a_{2}}{a_{1}}(\dim V_{i}-\epsilon_{i}(\overline\Phi)\dim V )\big]\leq$$
$$\frac{m}{P_{E}(m)}\big[rh^{0}(E^{m}_{i}(m))-r_{i}P_{E}(m)+\dfrac{r\delta(m)}{P_{E}(m)-\delta(m)}(h^{0}(E^{m}_{i}(m))-
\epsilon_{i}(\overline\Phi)P_{E}(m))\big]\leq$$
$$\frac{m}{P_{E}(m)}\big[rG(m)-r_{i}P_{E}(m)+\dfrac{r\delta(m)}{P_{E}(m)-\delta(m)}
(G(m)-\epsilon_{i}(\overline\Phi)P_{E}(m))\big]=$$
$$m\cdot \frac{\big[(P_{E}(m)-\delta(m))(rG(m)-r_{i}P_{E}(m))+
(r\delta(m))(G(m)-\epsilon_{i}(\overline\Phi)P_{E}(m))\big]}{P_{E}(m)(P_{E}(m)-\delta(m))}\;
.$$

Hence, $w_{i}<0$ is equivalent to
$$\Psi(m)=(P_{E}(m)-\delta(m))(rG(m)-r_{i}P_{E}(m))+(r\delta(m))(G(m)-\epsilon_{i}(\overline\Phi)P_{E}(m))<0$$
and $\Psi(m)=\xi_{2n}m^{2n}+\xi_{2n-1}m^{2n-1}+\cdots
+\xi_{1}m+\xi_{0}$ is a $(2n)^{th}$-order polynomial. Let us
calculate the higher order coefficient:
$$\xi_{2n}=(P_{E}(m)-\delta(m))_{n}(rG(m)-r_{i}P_{E}(m))_{n}+(r\delta(m))_{n}(G(m)-
\epsilon_{i}(\overline\Phi)P_{E}(m))_{n}=$$
$$(P_{E}(m)-\delta(m))_{n}(r\frac{r_{i}g}{n!}-r_{i}\frac{rg}{n!})+0=0\; .$$
Then, $\Psi(m)$ has no coefficient in order $(2n)^{th}$. Let us
calculate the $(2n-1)^{th}$-coefficient:
$$\xi_{2n-1}=(P_{E}(m)-\delta(m))_{n}(rG(m)-r_{i}P_{E}(m))_{n-1}+(r\delta(m))_{n-1}(G(m)-
\epsilon_{i}(\overline\Phi)P_{E}(m))_{n}=$$
$$\frac{rg}{n!}(rG_{n-1}-r_{i}\frac{A}{(n-1)!})+r\delta_{n-1}(\frac{r_{i}g}{n!}-
\epsilon_{i}(\overline\Phi)\frac{rg}{n!})$$ where $G_{n-1}$ is the
$(n-1)^{th}$-coefficient of the polynomial $G(m)$,
$$G_{n-1}=\frac{1}{g^{n-1}n!}ng^{n-1}((r_{i}-1)\mu_{max}(E)+\frac{d}{r}-C+r_{i}B)=$$
$$\frac{1}{(n-1)!}((r_{i}-1)\mu_{max}(E)+\frac{d}{r}-C+r_{i}B)\leq$$
$$\frac{1}{(n-1)!}((r_{i}-1)|\mu_{max}(E)|+\frac{d}{r}-C+r_{i}|B|)\leq$$
$$\frac{1}{(n-1)!}(r|\mu_{max}(E)|+\frac{d}{r}-C+r|B|)<\frac{-|A|}{(n-1)!}-\delta_{n-1}\; ,$$
last inequality coming from the definition of $C$ in
\eqref{C_constant_pairs}. Then
$$\xi_{2n-1}<\frac{rg}{n!}\big(r(\frac{-|A|}{(n-1)!}-\delta_{n-1})-r_{i}\frac{A}{(n-1)!}\big )+r\delta_{n-1}\big(\frac{r_{i}g}{n!}-\epsilon_{i}(\overline\Phi)\frac{rg}{n!}\big)=$$
$$\frac{rg}{n!}\big [ \big (\frac{-r|A|-r_{i}A}{(n-1)!}\big )
-r\delta_{n-1}+\delta_{n-1}(r_{i}-\epsilon_{i}(\overline\Phi)r)\big
]=$$
$$\frac{rg}{n!}\big [ \big (\frac{-r|A|-r_{i}A}{(n-1)!}\big )
+\delta_{n-1}(r_{i}-(1+\epsilon_{i}(\overline\Phi))r)\big ]<0$$
because $-r|A|-r_{i}A<0$,
$r_{i}-(1+\epsilon_{i}(\overline\Phi))r<0$ and $\delta_{n-1}\geq
0$. Note that if $r_{i}=r$, then
$\epsilon_{i}(\overline\Phi)=\epsilon_{t+1}(\overline\Phi)=1$.

Therefore $\Psi(m)=\xi_{2n-1}m^{2n-1}+\cdots +\xi_{1}m+\xi_{0}$
with $\xi_{2n-1}<0$, so there exists $m_{3}$ such that for $m\geq
m_{3}$ we will have $\Psi(m)<0$ and $w_{i}<0$, then the
contradiction.
\end{pr}

Now we can prove the following proposition in a similar way as we
proved Proposition \ref{regular}.

\begin{prop}
\label{regularpairs} There exists an integer $m_{4}$ such that for
$m\geq m_{4}$ the sheaves $E^{m}_{i}$ and
$E^{m,i}=E^m_i/E^m_{i-1}$ are $m_{4}$-regular. In particular their
higher cohomology groups, after twisting with
$\mathcal{O}_{X}(m_{4})$, vanish and they are generated by global
sections.
\end{prop}

\begin{prop}
\label{taskpairs}
Let $m\geq m_{4}$. For each filter $E_{i}^{m}$
in the $m$-Kempf filtration of $(E,\varphi)$ (c.f. \eqref{ekempf-filtpairs}) we have $\dim
V_{i}=h^{0}(E_{i}^{m}(m))$, therefore $V_{i}\cong H^{0}(E_{i}^{m}(m))$.
\end{prop}

\begin{pr}
Let $V_{\bullet}\subseteq V$ be the Kempf filtration of $V$ (cf.
Theorem \ref{kempftheoremsheaves}) and let $(E_{\bullet}^{m},\varphi|_{E_{\bullet}^{m}})\subseteq (E,\varphi)$ be the
$m$-Kempf filtration of $(E,\varphi)$ (cf.
\eqref{vkempf-filtpairs} and \eqref{ekempf-filtpairs}). We know
that each $V_{i}$ generates the subsheaf $E_{i}^{m}$, by
definition, then we have the following diagram:

$$\begin{array}{ccccccccccc}
    0 & \subset & V_{1} & \subset & V_{2} & \subset & \cdots & \subset & V_{t+1} & = & V \\
     & & \cap & & \cap & & & & & & ||  \\
      &     & H^{0}(E_{1}^{m}(m)) & \subset & H^{0}(E_{2}^{m}(m)) & \subset & \cdots & \subset & H^{0}(E_{t+1}^{m}(m)) & = & H^{0}(E(m))
      \end{array}$$

Suppose there exists an index $i$ such that $V_{i}\neq
H^{0}(E_{i}^{m}(m))$. Let $i$ be the index such that $V_{i}\neq
H^{0}(E_{i}^{m}(m))$ and $\forall j>i$ it is
$V_{j}=H^{0}(E_{j}^{m}(m))$. Then we have the diagram:

\begin{equation}
 \label{filtrationV_pair}
    \begin{array}{ccc}
    V_{i} & \subset & V_{i+1}\\
    \cap & & ||\\
    H^{0}(E_{i}^{m}(m)) & \subset & H^{0}(E_{i+1}^{m}(m))
    \end{array}
\end{equation}

Therefore $V_{i}\subsetneq H^{0}(E_{i}^{m}(m))\subsetneq V_{i+1}$
and we can consider a new filtration by adding the filter
$H^{0}(E_{i}^{m}(m))$:

\begin{equation}
\label{filtrationV'_pair}
    \begin{array}{ccccc}
    V_{i} & \subset & H^{0}(E_{i}^{m}(m)) & \subset & V_{i+1}\\
    || & & || & & ||\\
    V'_{i} & &  V'_{i+1} & & V'_{i+2}
    \end{array}
\end{equation}

Note that $V_{i}$ and $H^{0}(E_{i}^{m})$ generate the same sheaf
$E_{i}^{m}$, hence we are in situation of Lemma \ref{lemmaB},
where $W=H^{0}(E_{i}^{m})$, filtration $V_{\bullet}$ is
\eqref{filtrationV_pair} and filtration $V'_{\bullet}$ is
\eqref{filtrationV'_pair}.

The graph associated to filtration $V_{\bullet}$, by Definition
\ref{graph}, is given by the points
$$(b_{i},w_{i})=(\dfrac{\dim V_{i}}{m^{n}},\frac{m}{\dim V}\big(r\dim V_{i}-r_{i}\dim V+\dfrac{a_{2}}{a_{1}}
(\dim V_{i}-\epsilon_{i}(\overline\Phi,V_{\bullet})\dim V))\big)\;
,$$ where the slopes of the graph are given by

$$-v_{i}=\frac{w^{i}}{b^{i}}=\frac{w_{i}-w_{i-1}}{b_{i}-b_{i-1}}=$$
$$\frac{m^{n+1}}{\dim V}\big(r-r^{i}\frac{\dim V}{\dim V^{i}}+\frac{a_{2}}{a_{1}}(1-\epsilon^{i}(\overline\Phi,V_{\bullet})
\frac{\dim V}{\dim V^{i}})\big)\leq$$
$$\frac{m^{n+1}}{\dim V}\big(r+\frac{a_{2}}{a_{1}}\big):=R$$
and equality holds if and only if $r^{i}=0$. Here note that
$r^{i}=0$ implies $\epsilon^{i}(\overline\Phi,V_{\bullet})=0$.

Now, the new point which appears in the graph of the filtration
$V'_{\bullet}$ is
$$Q=\big(\dfrac{h^{0}(E_{i}^{m}(m))}{m^{n}},\frac{m}{\dim V}(rh^{0}(E_{i}^{m}(m))-r_{i}\dim V+
\dfrac{a_{2}}{a_{1}}(h^{0}(E_{i}^{m}(m))- \epsilon_{i}(\overline\Phi,V_{\bullet})\dim V))\big)\; ,$$ where we
write $\epsilon_{i}(\overline\Phi,V_{\bullet})$ instead of $\epsilon_{i}(\overline\Phi,V'_{\bullet})$, because
they are equal given that $V_{i}=V'_{i}$.

Point $Q$ joins two new segments appearing in this new graph. The
slope of the segment between $(b_{i},w_{i})$ and $Q$ is, by a
similar calculation,
$$-v'_{i+1}=\dfrac{m^{n+1}}{\dim V}(r+\dfrac{a_{2}}{a_{1}})=R\; .$$

By Lemma \ref{lemmaA}, the graph is convex, so
$v_{1}<v_{2}<\ldots<v_{t+1}$. As $E_{1}^{m}$ is a non-zero torsion
free sheaf, it has positive rank $r_{1}=r^{1}$ and hence it
follows $v_{1}>-R$.

Recall that, by definition,
$\epsilon_{i}(\overline\Phi,V_{\bullet})$ is equal to $1$ if
$\overline{\Phi}|_{V_{i}}\neq 0$ and $0$ otherwise. Then, it is
clear that
\begin{equation}
\label{epsilon}
    \begin{array}{c}
\epsilon_{j}(\overline{\Phi},V'_{\bullet})=\epsilon_{j}(\overline{\Phi},V_{\bullet})\;\;\;,j\leq i\\
\epsilon_{j}(\overline{\Phi},V'_{\bullet})=\epsilon_{j-1}(\overline\Phi,V_{\bullet})\;\;\;,j>i,
    \end{array}
\end{equation}
and note that
$\epsilon_{i}(\overline{\Phi},V_{\bullet})=\epsilon_{i+1}(\overline{\Phi},V'_{\bullet})$.
Then, the graph associated to $V'_{\bullet}\subset V$ is a
refinement of the graph associated to Kempf filtration
$V_{\bullet}\subset V$, therefore by Lemma \ref{lemmaB},
$v'_{i+1}\geq v_{i+1}$. Hence,
$$-R<v_{1}<v_{2}<\ldots<v_{i+1}\leq v'_{i+1}=-R\; ,$$
which is a contradiction.

Therefore, $\dim V_{i}=h^{0}(E_{i}^{m}(m))$, for every filter in the
$m$-Kempf filtration.
\end{pr}

\begin{cor}
\label{rankpairs} For every filter $E_{i}^{m}$ in the $m$-Kempf
filtration of $(E,\varphi)$ (c.f. \eqref{ekempf-filtpairs}), it is $r^{i}>0$.
\end{cor}
\begin{pr}
C.f. Corollary \ref{rank}.
\end{pr}

Now let us recall the results on subsection
\ref{sectionkempfstabilizes}. By Proposition \ref{regularpairs},
for any $m\geq m_{4}$, all the filters $E^{m}_{i}$ of the
$m$-Kempf filtration of the pair $(E,\varphi)$ are $m_{4}$-regular
and hence, the sheaves of the $m$-Kempf filtration 
$$0\subset (E^{m}_{1},\varphi|_{E_{1}^{m}})\subset
(E^{m}_{2},\varphi|_{E_{2}^{m}})\subset\cdots
(E^{m}_{t},\varphi|_{E_{t}^{m}})\subset
(E^{m}_{t+1},\varphi|_{E_{t+1}^{m}})\subset (E,\varphi)$$
are obtained by evaluating the filtration of vector subspaces
$$0\subset H^{0}(E^{m}_{1}(m_{4})) \subset H^{0}(E^{m}_{2}(m_{4})) \subset \cdots \subset H^{0}(E^{m}_{t_{m}}(m_{4})) \subset H^{0}(E^{m}_{t_{m}+1}(m_{4}))=H^{0}(E(m_{4}))$$
(c.f. Lemma \ref{mregularity}), of a unique vector space $H^{0}(E(m_{4}))$, whose
dimension is independent of $m$. Note that, because of Corollary \ref{rankpairs}, each subpair in the $m$-Kempf filtration of $(E,\varphi)$  is strictly contained
in the following one, for $m\geq m_{3}$. Let
$$(P_{1}^{m},\ldots,P_{t_{m}+1}^{m})$$ be the \emph{$m$-type}
of the $m$-Kempf filtration of $(E,\varphi)$ (c.f. Definition \ref{mtype}) and let
$$\mathcal{P}=\big\{(P_{1}^{m},\ldots, P_{t_{m}+1}^{m})\big\}$$
be the set of possible $m$-types, which is a finite set by the same argument as in Proposition \ref{Pisfinite}.

By Definition \ref{graphpairs} we associate a graph to the
$m$-Kempf filtration, given by $v_{m}$, which, thanks to
Propositions \ref{regularpairs} and \ref{taskpairs}, can be
rewritten as
$$v_{m,i}=m^{n+1}\cdot \frac{1}{P_{m}^{i}(m)P(m)}\big[r^{i}P(m)-rP_{m}^{i}(m)+\dfrac{r\delta(m)}{P(m)-\delta(m)}(\epsilon^{i}
(\overline\Phi)P(m)-P_{m}^{i}(m))\big]\; ,$$
$$b_{m}^{i}=\frac{1}{m^{n}}\cdot P_{m}^{i}(m)\; .$$

Define
$$v_{m,i}(l)=l^{n+1}\cdot \frac{1}{P_{m}^{i}(l)P(l)}\big[r^{i}P(l)-rP_{m}^{i}(l)+\dfrac{r\delta(l)}{P(l)-\delta(l)}(\epsilon^{i}
(\overline\Phi)P(l)-P_{m}^{i}(l))\big]\; ,$$ the coordinates of
the graph where the polynomials are evaluated on $l$ and let
$$
\Theta_{m}(l)=(\mu_{v_{m}(l)}(\Gamma_{v_{m}(l)}))^{2}=||v_{m}(l)||^{2}\;
,
$$
as in \eqref{maxvalue}. Let $\mathcal{A}$ be the finite set (c.f. Proposition \ref{Pisfinite})
$$
\mathcal{A}=\{\Theta_{m}:m\geq m_{4}\}\; .
$$
Let $K$ be a rational function which is maximal in $\mathcal{A}$
and, by a similar argument as in Lemma \ref{uniquefunction}, there
exists an integer $m_{5}$ with $\Theta_{m}=K$, $\forall m\geq
m_{5}$. Finally, we can prove the following

\begin{prop}
\label{eventuallypairs} Let $l_{1}$ and $l_{2}$ be integers with
$l_{1}\geq l_{2} \geq m_{5}$. Then the $l_{1}$-Kempf filtration of
$E$ is equal to the $l_{2}$-Kempf filtration of the holomorphic
pair $(E,\varphi)$.
\end{prop}
\begin{pr}
C.f. Proposition \ref{eventually}.
\end{pr}

Therefore, Theorem \ref{kempfstabilizespairs} follows from Proposition \ref{eventuallypairs}.

\begin{dfn}
If $m\geq m_{5}$, the $m$-Kempf filtration of $(E,\varphi)$ is
called \emph{the Kempf filtration of $(E,\varphi)$},
$$0\subset (E_{1},\varphi|_{E_{1}})\subset
(E_{2},\varphi|_{E_{2}})\subset\cdots
(E_{t},\varphi|_{E_{t}})\subset
(E_{t+1},\varphi|_{E_{t+1}})\subset (E,\varphi)\; .$$
\end{dfn}

\subsection{Harder-Narasimhan filtration for holomorphic pairs}
\label{hn-pairs}

Let $m\geq m_5$. Kempf's theorem
(c.f. Theorem \ref{kempftheoremsheaves}) asserts that given
$V\simeq H^{0}(E(m))$, there exists a unique
weighted filtration of vector spaces $(V_{\bullet},n_{\bullet})$
which gives maximum for the Kempf function
$$\mu(V_{\bullet},n_{\bullet})=\frac {\sum_{i=1}^{t+1} \frac{\Gamma_i}{\dim V} ( r^i \dim V - r\dim V^i)+\frac{a_2}{a_1}
 \sum_{i=1}^{t+1} \frac{\Gamma_{i}}{\dim V} \big(\epsilon^i(\overline\Phi)\dim V -\dim V^i\big)}
 {\sqrt{\sum_{i=1}^{t+1} {\dim V^{i}} \Gamma_{i}^{2}}}\; .$$ This filtration induces a filtration of holomorphic subpairs, called
the Kempf filtration of $(E,\varphi)$,
$$0\subset (E_{1},\varphi|_{E_{1}})\subset
(E_{2},\varphi|_{E_{2}})\subset\cdots
(E_{t},\varphi|_{E_{t}})\subset
(E_{t+1},\varphi|_{E_{t+1}})\subset (E,\varphi)\; ,$$
which is independent of $m$, for $m\geq m_{5}$, by Theorem
\ref{kempfstabilizespairs}, hence it is unique.

We proceed in a similar way to Section \ref{kempfsheaves} (c.f. Proof of Theorem \ref{kempfisHN}), to rewrite the Kempf
function for holomorphic pairs in terms of Hilbert polynomials of sheaves. Let
$\epsilon^{i}:=\epsilon^i(\overline\Phi)=\epsilon^{i}(\varphi)$
and note that $\epsilon^{i}=1$ for the unique index $i$ in the Kempf filtration such that
$\varphi|_{E_{i}}\neq 0$ and $\varphi|_{E_{i-1}}=0$, and
$\epsilon^{i}=0$ otherwise. Let us call this index $j$ in the
following.

\begin{prop}
\label{finalfunctionpairs} Given a holomorphic pair $(E,\varphi:E\rightarrow
\mathcal{O}_{X})$, there exists a unique filtration
$$0\subset (E_{1},\varphi|_{E_{1}})\subset
(E_{2},\varphi|_{E_{2}})\subset\cdots
(E_{t},\varphi|_{E_{t}})\subset
(E_{t+1},\varphi|_{E_{t+1}})\subset (E,\varphi)$$ with positive weights $n_{1},\ldots,n_{t}$, which
gives maximum for the function
$$K(m)=\frac{m^{\frac{n}{2}+1}}{P}\cdot \frac {\sum_{i=1}^{t+1} \gamma_{i}[(r^i P - rP^i)+\frac{r\delta}{P-\delta}
 (\epsilon^i P -P^i)]}
 {\sqrt{\sum_{i=1}^{t+1} P^{i} \gamma_{i}^{2}}}\; .$$
\end{prop}

Similarly, we can express the function $K$ in Proposition \ref{finalfunctionpairs} as
$$K(m)=m^{-\frac{n}{2}}\cdot \frac
{\sum_{i=1}^{t+1} P^{i}\gamma_{i}v_{i}} {\sqrt{\sum_{i=1}^{t+1}
P^{i}\gamma_{i}^{2}}}=m^{-\frac{n}{2}}\cdot
\frac{(\gamma,v)}{||\gamma||}\; ,$$

where the coordinates $v_{i,m}$ (slopes of segments of the graph),
now are
$$v_{i}=m^{n+1}\cdot \frac{1}{P^{i}P}\big[r^{i}P-rP^{i}+\dfrac{r\delta}{P-\delta}
(\epsilon^{i}P-P^{i})\big]$$ and the scalar product is, again,
$$\left(
    \begin{array}{cccc}
      P^{1} &  &  &  \\
       & P^{2} &  &  \\
       &  & \ddots &  \\
       &  &  & P^{t+1} \\
    \end{array}
  \right)$$

With Definition \ref{corrected}, the coordinates of the graph are
$$v_{i}=m^{n+1}\cdot
\frac{r^{i}r}{P^{i}(P-\delta)}\big(\frac{\overline{P}_{E}}{r}-\frac{\overline{P}_{E^{i}}}{r^{i}}\big)\;
,$$ where $\overline{P}_{E^{i}}=P^{i}-\delta\epsilon^{i}$ is the corrected Hilbert polynomial of
the quotient pair $(E^{i},\varphi|_{E^{i}})$ (c.f. Definitions
\ref{subpairs} and \ref{corrected}).

Now, we define a Harder-Narasimhan filtration for a holomorphic
pair, analogously to the notion for torsion free sheaves,
substituting the Hilbert polynomials by the corrected Hilbert
polynomials.

\begin{dfn}
Given a holomorphic pair $(E,\varphi:E\rightarrow
\mathcal{O}_{X})$, a filtration
$$0\subset (E_{1},\varphi|_{E_{1}})\subset
(E_{2},\varphi|_{E_{2}})\subset\cdots
(E_{t},\varphi|_{E_{t}})\subset
(E_{t+1},\varphi|_{E_{t+1}})\subset (E,\varphi)$$ is called a \textbf{Harder-Narasimhan filtration} of
$(E,\varphi)$ if it satisfies these two properties, where $E^{i}:=(E_{i}/E_{i-1},\varphi|_{E_{i}/E_{i-1}})$,
 \begin{enumerate}
   \item The corrected Hilbert polynomials verify
   $$\frac{\overline{P}_{E^{1}}}{\rk E^{i}}>\frac{\overline{P}_{E^{2}}}{\rk E^{2}}>\ldots>\frac{\overline{P}_{E^{t+1}}}{\rk E^{t+1}}$$
   \item Every quotient pair $(E^{i},\varphi|_{E^{i}})$ is $\delta$-semistable as a holomorphic pair.
 \end{enumerate}
\end{dfn}

Next, we prove the existence and uniqueness for the Harder-Narasimhan filtration of a holomorphic pair. The proof follows similarly to Theorem \ref{HNunique}. 

\begin{thm}
\label{HNuniquepairs} Every pair $(E,\varphi)$ has a unique
Harder-Narasimhan filtration.
\end{thm}
\begin{lem}
Let $(E,\varphi)$ be a pair. Then, there exists a subsheaf
$F\subseteq E$ such that for all subsheaves $G\subset E$, one has
$\frac{\overline{P}_{F}}{\rk F}\geq \frac{\overline{P}_{G}}{\rk G}$, and in case of equality
$G\subseteq F$. Moreover, $F$ is uniquely determined and
$(F,\varphi|_{F})$ is $\delta$-semistable, called the \emph{maximal
destabilizing subpair of $(E,\varphi)$}.
\end{lem}
\begin{pr}
The last two assertions follow from the first, where note that
being $\delta$-semistable can be checked by subpairs, by Lemma
\ref{subobjects}.

Define an order relation on the set of subpairs of $(E,\varphi)$
by $(F_{1},\varphi|_{F_{1}})\leq(F_{2},\varphi|_{F_{2}})$ if and
only if $F_{1}\subset F_{2}$ and $\frac{\overline{P}_{F_{1}}}{\rk F_{1}}\leq
\frac{\overline{P}_{F_{2}}}{\rk F_{2}}$. Every ascending chain is bounded by
$(E,\varphi)$, then by Zorn's Lemma, for every subpair
$(F,\varphi|_{F})$ there exists a $F\subset F'\subset E$ such that
$(F',\varphi|_{F'})$ is maximal with respect to $\leq$. Let
$(F,\varphi|_{F})$ be $\leq$-maximal with $F$ of minimal rank
among all maximal subpairs and let us show that $(F,\varphi|_{F})$ is the maximal destabilizing subpair.

Suppose that $\exists\;\; G\subset E$ with $\frac{\overline{P}_{G}}{\rk G}\geq \frac{\overline{P}_{F}}{\rk F}$. First, we show that we can assume $G\subset F$
by replacing $G$ by $G\cap F$. Indeed, if $G\nsubseteq F$, then
$F$ is a proper subsheaf of $F+G$ and hence
$\frac{\overline{P}_{F}}{\rk F}>\frac{\overline{P}_{F+G}}{\rk F+G}$, by definition of $F$. Let the exact sequence
$$0\rightarrow F\cap G\rightarrow F\oplus G\rightarrow
F+G\rightarrow 0$$ out of which we get
$$P_{F}+P_{G}=P_{F\oplus G}=P_{F\cap
G}+P_{F+G}$$ and $$\rk(F)+\rk(G)=\rk(F\oplus G)=\rk(F\cap
G)+\rk(F+G)\; .$$ Calculating we have
$$\rk(F\cap G)(\frac{P_{G}}{\rk G}-\frac{P_{F\cap G}}{\rk (F\cap G)}=$$
$$\rk(F+G)(\frac{P_{F+G}}{\rk (F+G)}-
\frac{P_{F}}{\rk F})+(\rk(G)-\rk(F\cap G))(\frac{P_{F}}{\rk F}-\frac{P_{G}}{\rk G})\; .$$
Using
$$\epsilon(F\cap G)+\epsilon(F+G)\leq\epsilon(F)+\epsilon(G)$$
we get
$$(P_{F}-\delta\epsilon(F))+(P_{G})-\delta\epsilon(G))=(P_{F\cap G}-\delta\epsilon(F\cap G))+(P_{F+G}-\delta\epsilon(F+G))$$ and,
similarly,
$$\rk(F\cap G)\big(\frac{\overline{P}_{G}}{\rk G}-\frac{\overline{P}_{F\cap G}}{\rk (F\cap G)}\big)\leq$$
$$\rk(F+G)\big(\frac{\overline{P}_{F+G}}{\rk (F+G)}-\frac{\overline{P}_{F}}{\rk F}\big)+
(\rk(G)-\rk(F\cap G))\big(\frac{\overline{P}_{F}}{\rk F}-\frac{\overline{P}_{G}}{\rk G}\big)\; .$$ Then,
using the inequalities $\frac{\overline{P}_{F}}{\rk F}\leq
\frac{\overline{P}_{G}}{\rk G}$ and $\frac{\overline{P}_{F}}{\rk F}> \frac{\overline{P}_{F+G}}{\rk (F+G)}$, we obtain
$$\frac{\overline{P}_{G}}{\rk G}-\frac{\overline{P}_{F\cap G}}{\rk (F\cap G)}<0$$
and hence
$$\frac{\overline{P}_{F}}{\rk F}<\frac{\overline{P}_{F\cap G}}{\rk (F\cap G)}\; ,$$
hence we can suppose that $G\subset F$. 

Let $G\subset F$ with $\frac{\overline{P}_{G}}{\rk G}>\frac{\overline{P}_{F}}{\rk F}$ such
that $(G,\varphi|_{G})$ is $\leq$-maximal in $(F,\varphi|_{F})$.
Then let $(G',\varphi|_{G'})\geq (G,\varphi|_{G})$ to be $\leq$-maximal
in $(E,\varphi)$. We obtain, 
$\frac{\overline{P}_{F}}{\rk F}<\frac{\overline{P}_{G}}{\rk G}\leq \frac{\overline{P}_{G'}}{\rk G'}$ and, by
maximality of $(G',\varphi|_{G'})$ and $(F,\varphi|_{F})$ it is
$G'\nsubseteq F$, since otherwise it would be $\rk G<\rk F$ which contradicts 
the minimality of $\rk F$, therefore $F$ is a proper subsheaf of
$F+G'$. Then we obtain $\frac{\overline{P}_{F}}{\rk F}>\frac{\overline{P}_{F+G'}}{\rk (F+G')}$ and the inequalities 
$\frac{\overline{P}_{F}}{\rk F}<\frac{\overline{P}_{G'}}{\rk G'}$ and
$\frac{\overline{P}_{F}}{\rk F}>\frac{\overline{P}_{F+G'}}{\rk (F+G')}$ give
$$\frac{\overline{P}_{F\cap G'}}{\rk (F\cap G')}>\frac{\overline{P}_{G'}}{\rk G'}\geq\frac{\overline{P}_{G}}{\rk G}\; .$$
Therefore, as $G\subset F\cap G'\subset F$, we get a contradiction. 
\end{pr}

\begin{pr}[Proof of the Theorem \ref{HNuniquepairs}]
With the previous Lemma we are able to show the existence of a Harder-Narasimhan
filtration for $(E,\varphi)$. Let $(E_{1},\varphi|_{E_{1}})$ the
maximal destabilizing subpair and suppose that the corresponding
quotient $(E/E_{1},\varphi|_{E/E_{1}})$ has a Harder-Narasimhan
filtration,
$$0\subset G_{0}\subset G_{1}\subset\ldots\subset G_{t-1}=E/E_{1}\; ,$$
by induction. We define $E_{i+1}$ to be the pre-image of $G_{1}$ and it
is $\frac{\overline{P}_{E_{1}}}{\rk E_{1}}>\frac{\overline{P}_{E_{2}/E_{1}}}{\rk E_{2}/E_{1}}$ because, if
not, we would have $\frac{\overline{P}_{E_{1}}}{\rk E_{1}}\leq\frac{\overline{P}_{E_{2}}}{\rk E_{2}}$,
contradicting the maximality of $(E_{1},\varphi|_{E_{1}})$.

To show the uniqueness, suppose that $E_{\bullet}$ and $E'_{\bullet}$
are two Harder-Narasimhan filtrations of $(E,\varphi)$ and consider, without loss
of generality, that $\frac{\overline{P}_{E'_{1}}}{\rk E'_{1}}\geq \frac{\overline{P}_{E_{1}}}{\rk E_{1}}$. Call $j$ an index which is
 minimal such that $E'_{1}\subset E_{j}$. The composition
$$E'_{1}\rightarrow E_{j}\rightarrow E_{j}/E_{j-1}$$ is a
non-trivial homomorphism of semistable sheaves which implies
$$\frac{\overline{P}_{E_{j}/E_{j-1}}}{\rk E_{j}/E_{j-1}}\geq \frac{\overline{P}_{E'_{1}}}{\rk E'_{1}}\geq
\frac{\overline{P}_{E_{1}}}{\rk E_{1}}\geq \frac{\overline{P}_{E_{j}/E_{j-1}}}{\rk E_{j}/E_{j-1}}\; ,$$ where first
inequality comes from the fact that if there exists a non-trivial
homomorphism between semistable pairs, then the corrected Hilbert
polynomial of the target is greater or equal than the one of the
first pair. Hence, equality holds everywhere, implying $j=1$ so
that $E'_{1}\subset E_{1}$. Then, by semistability of the pair $(E_{1},\varphi|_{E_{1}})$, it
is $\frac{\overline{P}_{E'_{1}}}{\rk E'_{1}}\leq \frac{\overline{P}_{E_{1}}}{\rk E_{1}}$, and we can
repeat the argument interchanging the roles of $E_{1}$ and
$E'_{1}$ to show that $E_{1}=E'_{1}$. By induction we can assume that
uniqueness holds for the Harder-Narasimhan filtration of
$(E/E_{1},\varphi|_{E/E_{1}})$. This implies that $E'_{i}/E_{1}=E_{i}/E_{1}$ and completes the
proof.
\end{pr}

Now we will give the analogous to Propositions
\ref{descendentslopes} and \ref{blocksemistability}.

\begin{prop}
\label{descendentpolpairs} Given the Kempf filtration of a holomorphic pair
$(E,\varphi)$,
$$0\subset (E_{1},\varphi|_{E_{1}})\subset
(E_{2},\varphi|_{E_{2}})\subset\cdots
(E_{t},\varphi|_{E_{t}})\subset
(E_{t+1},\varphi|_{E_{t+1}})\subset (E,\varphi)$$
it verifies
$$\frac{\overline{P}_{E^{1}}}{\rk E^{i}}>\frac{\overline{P}_{E^{2}}}{\rk E^{2}}>\ldots>\frac{\overline{P}_{E^{t+1}}}{\rk E^{t+1}}\; .$$
\end{prop}
\begin{pr}
Let $j$ be the unique index such that $\epsilon^{j}=1$. By Lemma
\ref{lemmaA} it is
$$v_{1}<v_{2}<\ldots v_{j-1}<v_{j}<v_{j+1}<\ldots
<v_{t+1}\; .$$ Note that for $i\neq j$ it is $\overline{P}^{i}=P^{i}-\delta\epsilon^{i}=P^{i}$, hence
$v_{i-1}<v_{i}$ implies $\frac{\overline{P}_{E^{i-1}}}{\rk E^{i-1}}>\frac{\overline{P}_{E^{i}}}{\rk E^{i}}$ for
all $i\neq j,j+1$.

Now the inequality $v_{j-1}<v_{j}$ is
$$\frac{r^{j-1}r}{P^{j-1}(P-\delta)}(\frac{P-\delta}{r}-\frac{P^{j-1}}{r^{j-1}})<
\frac{r^{j}r}{P^{j}(P-\delta)}(\frac{P-\delta}{r}-\frac{P^{j}-\delta}{r^{j}})$$
or, equivalently,
$$-\delta \frac{rP^{j-1}}{P-\delta}<P^{j-1}r^{j}-P^{j}r^{j-1}\;
.$$ The function $\frac{rP^{j-1}}{P-\delta}$ is a homogeneous
rational function whose limit at infinity is $r^{j-1}$, hence for
large values of the variable we obtain this inequality between the polynomials
$$-\delta r^{j-1}<P^{j-1}r^{j}-P^{j}r^{j-1}\; ,$$
which is equivalent to $\frac{\overline{P}_{E^{j-1}}}{\rk E^{j-1}}>\frac{\overline{P}_{E^{j}}}{\rk E^{j}}$. A similar argument
proves that  $\frac{\overline{P}_{E^{j}}}{\rk E^{j}}>\frac{\overline{P}_{E^{j+1}}}{\rk E^{j+1}}$.
\end{pr}

\begin{prop}
\label{blocksemistabilitypairs} Given the Kempf filtration of a
holomorphic pair $(E,\varphi)$,
$$0\subset (E_{1},\varphi|_{E_{1}})\subset
(E_{2},\varphi|_{E_{2}})\subset\cdots
(E_{t},\varphi|_{E_{t}})\subset
(E_{t+1},\varphi|_{E_{t+1}})\subset (E,\varphi)$$
each one of the quotient pairs $(E^{i},\varphi|_{E^{i}})$ is semistable as a holomorphic pair.
\end{prop}
\begin{pr}
Suppose that any of the blocks has a destabilizing subpair and apply a similar argument to the one in Proposition \ref{blocksemistability}.
\end{pr}

Hence, having seen the convexity properties of the Kempf filtration in
Propositions \ref{descendentpolpairs} and
\ref{blocksemistabilitypairs}, we get that the Kempf filtration
of a holomorphic pair $(E,\varphi)$ is a Harder-Narasimhan filtration. Given that every holomorphic pair has a unique Harder-Narasimhan filtration by Theorem \ref{HNuniquepairs},
therefore it will be the same that the Kempf filtration.

\begin{cor}
\label{kempfisHNpairs}
Let $(E,\varphi)$ be a $\delta$-unstable
holomorphic pair. The Kempf filtration is the same that the
Harder-Narasimhan filtration.
\end{cor}
\begin{pr}
By Propositions \ref{descendentpolpairs} and
\ref{blocksemistabilitypairs} the Kempf filtration is a
Harder-Narasimhan filtration, which is unique by Theorem
\ref{HNuniquepairs}, hence both filtrations are the same.
\end{pr}

\section{Higgs sheaves}
\label{kempfhiggs}

Here we consider the moduli space of Higgs sheaves constructed by
Simpson in \cite{Si1,Si2} and use the techniques of the previous
sections to prove the analogous result in this case, the
correspondence between Kempf and Harder-Narasimhan filtrations.

\vspace{1cm}

 Let $X$ be a smooth complex projective variety of
dimension $n$. A \emph{Higgs sheaf} is a pair $(E,\varphi)$ where
$E$ is a coherent sheaf over $X$ and $\varphi:E\rightarrow
E\otimes\Omega^{1}_{X}$ verifying $\varphi\wedge\varphi=0$, a
morphism called the \emph{Higgs field}. We call $(E,\varphi)$ a \emph{Higgs bundle} if $E$ is a locally free sheaf. Recall that
$\Omega^{1}_{X}=T^{\ast}X$, the cotangent bundle.

We say that a Higgs sheaf $(E,\varphi)$ is \emph{semistable} (in the
sense of Gieseker) if for all proper subsheaves $F\subset E$,
preserved by $\varphi$ (i.e. $\varphi\big|_{F}:F\rightarrow
F\otimes \Omega^{1}_{X}$) we have
$$\frac{P_{F}}{\rk F}\leq\frac{P_{E}}{\rk E}\; ,$$ where $P_{E}$ and $P_{F}$ are the respective Hilbert
polynomials of $E$ and $F$. We say that $(E,\varphi)$ is
\emph{stable} if we have a strict inequality for every subsheaf
preserved by $\varphi$.

A Higgs field can be thought as a coherent sheaf $\mathcal{E}$ on the cotangent bundle $T^{\ast}X$, supported on
the spectral curve (the eigenvalues of the Higgs field). Note that, to define a sheaf of
$\mathcal{O}_{T^{\ast}X}$-modules on the total space of $T^{\ast}X$ we have to determine how to multiply a section by a
function on the vertical variables, which is given by the Higgs field, by definition.

Let $Z$ be a projective completion of $T^{\ast}X$, $Z=\mathbb{P}(T^{\ast}X\oplus\mathcal{O})$, and
$D=Z-T^{\ast}X$, the divisor at infinity.

\begin{lem}\cite[Lemma 6.8]{Si2}
\label{higgsiscoherent}
A Higgs sheaf $(E,\varphi)$ on $X$ is the same thing that a coherent sheaf
$\mathcal{E}$ on $Z$ such that $\supp(E)\cap D=\emptyset$, where $E=\pi_{\ast}\mathcal{E}$ and
$\pi:T^{\ast}X\rightarrow X$. This identification is compatible with morphisms, giving an equivalence of
categories. The condition that $E$ is torsion free is the same that $\mathcal{E}$ is of pure dimension $n=\dim
(X)$.
\end{lem}

Choose $k$ so that
$$\mathcal{O}_{Z}(1)\overset{def}{=}\pi^{\ast}\mathcal{O}_{X}(k)\otimes_{\mathcal{O}_{Z}}\mathcal{O}_{Z}(D)$$
is ample on $Z$ (c.f. \cite[Appendix A, Theorem 5.1]{Ha}). In
particular,
$\mathcal{O}_{T^{\ast}X}(1)=\pi^{\ast}\mathcal{O}_{X}(k)$. Thus, for
any coherent sheaf $\mathcal{E}$ on $Z$ with support not meeting
$D$, the Hilbert polynomials of $\mathcal{E}$ and
$\pi_{\ast}\mathcal{E}$ differ by rescaling on the variable $m$
$$P_{\mathcal{E}}(m)=P_{\pi_{\ast}\mathcal{E}}(km)\; .$$

Recall that, the condition for $E$ to be torsion free is equivalent to $\mathcal{E}$ being pure of dimension $n$. To relate the stability of a Higgs
sheaf $(E,\varphi)$ with the stability of the associated sheaf $\mathcal{E}$, we have to modify Definition \ref{giesekerstab} as in \cite{Si1,Si2}, 
which was stated only for torsion free sheaves.
Recall the expression
of the Hilbert polynomial of a sheaf $\mathcal{E}$, in (\ref{HpolE}),
$$P_{\mathcal{E}}(m)=\frac{rg}{n!}m^{n}+\frac{d+r\alpha_{n-1}}{(n-1)!}m^{n-1}+...$$
We define $r=\rk \mathcal{E}$, the rank of $\mathcal{E}$, such that the coefficient of the leading term of the Hilbert polynomial is $\frac{rg}{n!}$. We also define $d=\deg \mathcal{E}$,
the degree of $\mathcal{E}$, as the corresponding coefficient appearing in the expression. A coherent sheaf $\mathcal{E}$ is of pure dimension $n$ if it has no subsheaves supported
on a lower dimensional locus.

\begin{dfn}
\label{simpsonstab}
A coherent sheaf $\mathcal{E}$ on $X$ is called \emph{semistable} if it is pure of dimension $n$, and for all proper subsheaves $\mathcal{F}\subset \mathcal{E}$, it is
$$
\frac{P_{\mathcal{F}}}{\rk \mathcal{F}} \leq \frac{P_{\mathcal{E}}}{\rk \mathcal{E}} \; .$$ If strict inequality holds for every proper subsheaf, we say
that $\mathcal{E}$ is \emph{stable}. If $(E,\varphi)$ is not semistable, we say that it is \emph{unstable}.
\end{dfn}

Given that the Higgs subsheaves of $(E,\varphi)$ correspond to the coherent
subsheaves of $\mathcal{E}$, and since a subsheaf of $\mathcal{E}$
is the same that a subsheaf of $\pi_{\ast}\mathcal{E}$ preserved
by the action of the symmetric algebra $\Sym^{\ast}(T^{\ast}X)$,
we have the following

\begin{cor}
\label{samestability} $(E,\varphi)$ is a semistable Higgs sheaf
if and only if the corresponding sheaf $\mathcal{E}$, by Lemma
\ref{higgsiscoherent}, is semistable as a coherent sheaf (c.f. Definition \ref{simpsonstab}).
\end{cor}

\subsection{Moduli space of Higgs sheaves}

Given a polynomial $P$, we denote by $k^{\ast}P$ the polynomial
such that $k^{\ast}P(m)=P(km)$. Denote by
$\mathfrak{M}(\mathcal{O}_{T^{\ast}X},k^{\ast}P)$ be the moduli
space of coherent sheaves $\mathcal{E}$ over $\mathcal{O}_{T^{\ast}X}$ with
Hilbert polynomial $k^{\ast}P$. By \cite[$\S 1$]{Si1}, Lemma
\ref{higgsiscoherent} and Corollary \ref{samestability}, the
scheme $\mathfrak{M}(\mathcal{O}_{T^{\ast}X},k^{\ast}P)$
corepresents the functor $\mathcal{M}_{Higgs}(X,P)$
which associates a scheme $S$ to the set of isomorphism classes of
semistable Higgs sheaves $(E,\varphi)$ on $X$, over $S$, with Hilbert
polynomial $P$. Therefore, we put
$$\mathfrak{M}_{Higgs}(X,P)=\mathfrak{M}(\mathcal{O}_{T^{\ast}X},k^{\ast}P)$$
and let us construct the scheme $\mathfrak{M}(\mathcal{O}_{T^{\ast}X},k^{\ast}P)$ following Simpson's method
(c.f. \cite[$\S 1$]{Si1}).

Let $P$ be a polynomial of degree $n=\dim X$. There exists an integer $N$, sufficiently large, such that for
$m\geq N$, $\mathcal{E}(m)$ is generated by global sections and $h^{i}(\mathcal{E}(m))=0$ for $i>0$. Then,
choose $m\geq N$ and fix an isomorphism
$$\alpha: V\simeq \mathbb{C}^{k^{\ast}P(m)}=\mathbb{C}^{P(km)}$$ to obtain a quotient
$$q:V\otimes \mathcal{W}\twoheadrightarrow \mathcal{E}\; ,$$ where $\mathcal{W}=\mathcal{O}_{Z}(-m)$. Let $\mathcal{H}$
be the Hilbert scheme of quotients $$\mathcal{H}=\Hilb(V\otimes
\mathcal{W},k^{\ast}P)=\{V\otimes\mathcal{W}\rightarrow \mathcal{E}\rightarrow 0\}\; ,$$ with
$P_{\mathcal{E}}(m)=P(km)=k^{\ast}P(m)$. We define an embedding of this Hilbert scheme to a projective space. Let $l\gg m$ be an integer such
that $H^{1}(\kernel (V\otimes \mathcal{W}\twoheadrightarrow \mathcal{E})(l))=0$. Then, $q$ induces the following
homomorphisms
$$q:V\otimes \mathcal{W}(l)\twoheadrightarrow \mathcal{E}(l)$$
$$q':V\otimes H^{0}(\mathcal{W}(l))\twoheadrightarrow H^{0}(\mathcal{E}(l))$$
$$q'':\bigwedge^{P(kl)}(V\otimes H^{0}(\mathcal{W}(l))\twoheadrightarrow \bigwedge^{P(kl)}H^{0}
(\mathcal{E}(l))\simeq \mathbb{C}$$ Hence, it defines a
Grothendieck's embedding on the Grassmannian manifold
$$\mathcal{H}=\Hilb(V\otimes \mathcal{W},k^{\ast}P(m))\overset{\mathcal{L}_{l,m}}{\hookrightarrow}\mathbb{P}
(\bigwedge^{P(kl)}(V^{\vee}\otimes
H^{0}(\mathcal{W}(l))^{\vee}))\; ,$$ where $\mathcal{L}_{l,m}$ is
the very ample line bundle (depending on $l$ and $m$) given by the
pullback of the canonical line bundle on the Grassmannian by the
embedding. Note that, given a point $q\in \mathcal{H}$, if
$H^{0}(q(m)):V\rightarrow H^{0}(\mathcal{E}(m))$ is an
isomorphism, we can recover the sheaf $\mathcal{E}$ together with
the isomorphism $\alpha: V\simeq \mathbb{C}^{P(km)}$.

The group $GL(V)$ of changes of isomorphism $V\simeq
\mathbb{C}^{P(km)}$, acts on $\mathcal{H}$ and the line bundle
$\mathcal{L}_{l,m}$. Note that, if we change the isomorphism by a
homothecy we obtain a different point in the line bundle defined
by $q$, hence the point $\overline{q}$ in the projective space is
the same. Hence, as in the case of the Gieseker embedding (c.f. subsection
\ref{modulisheaves}), we can get rid of the choice of isomorphism by
dividing by the action of $SL(V)$. Let $Q\subset \mathcal{H}$ the
$SL(V)$-invariant open subset of quotients where $\mathcal{E}$ is
a semistable sheaf of pure dimension $n=\dim X$ and the induced
morphism $\alpha: V\simeq \mathbb{C}^{P(km)}$ is an isomorphism.
There exists a good quotient (c.f. Definition \ref{goodquotient})
$$\mathfrak{M}(\mathcal{O}_{Z},k^{\ast}P)=Q/SL(V)\; ,$$ (c.f. \cite[Theorem 1.19]{Si1} and \cite{Mu}). Let $Q'\subset Q$ be the
subset of those quotient sheaves $\mathcal{E}$ whose support does
not meet $D$ (which is also $SL(V)$-invariant and is,
set-theoretically, the inverse image of a subset of
$\mathfrak{M}(\mathcal{O}_{Z},k^{\ast}P)$). Therefore, a good
quotient $Q'/SL(V)$ exists and it is equal to an open subset which
we denote $\mathfrak{M}(\mathcal{O}_{T^{\ast}X},k^{\ast}P)\subset
\mathfrak{M}(\mathcal{O}_{Z},k^{\ast}P)$.

As we said before, $\mathfrak{M}(\mathcal{O}_{T^{\ast}X},k^{\ast}P)$ corepresents the functor
$\mathcal{M}_{Higgs}(X,P)$, hence,
$$\mathfrak{M}_{Higgs}(X,P)=\mathfrak{M}(\mathcal{O}_{T^{\ast}X},k^{\ast}P)\; .$$

Therefore, to construct a moduli space for Higgs sheaves
$(E,\varphi)$, where $P$ is the fixed Hilbert polynomial of $E$,
we construct a particular moduli space of sheaves. From now on,
let us consider semistable coherent sheaves $\mathcal{E}$ over $Z=\mathbb{P}(T^{\ast}X\oplus\mathcal{O})$
of pure dimension $n=\dim X$ and fixed Hilbert polynomial $P$. We
consider the construction of the moduli space of sheaves following
Simpson's method. Giving a sheaf $\mathcal{E}$ and an isomorphism $V\cong
H^{0}(\mathcal{E}(m))$, we obtain a point $\overline{q}$ in the GIT space of
parameters. We have to prove that a point $\overline{q}$ is GIT
semistable if and only if it corresponds to a semistable sheaf
$\mathcal{E}$, to conclude that a moduli space for semistable sheaves can be
obtained as the good quotient of the space of GIT semistable
points, by the group $SL(V)$.

\begin{rem}
Note that the sheaves $\mathcal{E}$ are torsion sheaves supported on a
subscheme of $Z=\mathbb{P}(T^{\ast}X\oplus\mathcal{O})$. In this
case we cannot use Gieseker's embedding in the construction of the
moduli space, because if we take $\wedge^{r}\mathcal{E}$ we get
the zero sheaf (c.f. subsection \ref{modulisheaves}). Simpson
develops his method based on Grothendieck's ideas which gives a
solution to the problem in this case (c.f. \cite[p. 53]{Si1}).
\end{rem}

Let us calculate the numerical function on the set of $1$-parameter subgroups, to apply the Hilbert-Mumford criterion (c.f. Theorem \ref{HMcrit}).

Let $l,m$ be integers as before, and let $V$ be a complex vector space of dimension $P(m)$. Recall that a weighted filtration $(V_\bullet,n_\bullet)$ of $V$ is
a filtration
\begin{equation}
\label{filtVhiggs} 0 \subset V_1 \subset V_2 \subset \;\cdots\; \subset V_t \subset V_{t+1}=V,
\end{equation}
and positive rational numbers $n_1,\, n_2,\ldots , \,n_t > 0$. Let $\Gamma$ be the $1$-parameter subgroup associated to
the weighted filtration (c.f. subsection \ref{modulisheaves}) given by
$$\Gamma=(\overbrace{\Gamma_1,\ldots,\Gamma_1}^{\dim V^1},
\overbrace{\Gamma_2,\ldots,\Gamma_2}^{\dim V^2}, \ldots, \overbrace{\Gamma_{t+1},\ldots,\Gamma_{t+1}}^{\dim
V^{t+1}}) \; ,$$ where $V^i=V_i/V_{i-1}$.

Let $W=H^{0}(\mathcal{O}(l-m))$ be a vector space where $SL(V)$ acts
trivially. The basis $\{e_1,\ldots,e_p\}$, together with a basis
$\{w_j\}$ of $W$, induces a basis of
$\bigwedge^{P(l)}(V^{\vee}\otimes W^{\vee})$ indexed in a natural
way by tuples $(i_1,\ldots,i_{P(l)},j)$ (the indexes $i_j$ being
skewsymmetric), and the coordinate corresponding to such an index
is acted by $\Gamma$ with exponent
$$\Gamma_{i_1}+\cdots+\Gamma_{i_{P(l)}}\; .$$
The coordinate $(i_1,\ldots,i_{P(l)},j)$ of the point corresponding to the sheaf $\mathcal{E}$ is non-zero if and only if the
evaluations of the sections $e_1,\ldots,e_{P(l)}$ are linearly independent for generic $x\in X$. Therefore, the
numerical function in Theorem \ref{HMcrit} is
\begin{eqnarray}
  \label{eq:mulefthiggs}
\mu(\overline{q},V_\bullet,n_\bullet)&=&\min \{\Gamma_{i_1}+\cdots+\Gamma_{i_{P(l)}}: \,
q''(e_{i_1}\wedge\cdots \wedge e_{i_{P(l)}})\neq 0 \} \notag\\
&=&\min \{\Gamma_{i_1}+\cdots+\Gamma_{i_{P(l)}}: \,
e_{i_1}(x),\ldots,e_{i_{P(l)}}(x)  \\
& & \qquad\text{linearly independent for generic
  $x\in X$}\}   \notag
\end{eqnarray}

Let $\mathcal{E}_{V_{i}}$ be the subsheaf generated by $V_{i}$ and let $\mathcal{E}_{V^{i}}$ be the subsheaf generated by
$V^{i}=V_{i}/V_{i-1}$. Let $P_{\mathcal{E}_{V_{i}}}$ and $P_{\mathcal{E}_{V^{i}}}$ be the corresponding Hilbert polynomials.
Note that $P_{\mathcal{E}_{V_{i}}}-P_{\mathcal{E}_{V_{i-1}}}=P_{\mathcal{E}_{V^{i}}}$. Given a $1$-parameter subgroup $\Gamma$, we get
$$\mu(\overline{q},V_\bullet,n_\bullet)= \sum_{i=1}^{t} n_i ( P(l) \dim V_i - P_{\mathcal{E}_{V_{i}}}(l)\dim V)=$$
$$\sum_{i=1}^{t+1} \frac{\Gamma_i}{\dim V} ( P_{\mathcal{E}_{V^{i}}}(l) \dim V - P(l)\dim V^i)\; ,$$ where recall that
$n_i=\frac{\Gamma_{i+1}-\Gamma_i}{p}$.

By the Hilbert-Mumford criterion, Theorem \ref{HMcrit}, a point
$$\overline{q}\in \PP(\bigwedge^{P(l)}(V^{\vee}\otimes W^{\vee}))$$
is \emph{GIT semistable} if and only if for all weighted filtrations it is
$$\mu(\overline{q},V_\bullet,n_\bullet)\leq 0\; .$$
Using the previous calculations, this can be stated as follows:
\begin{prop}
\label{HMcrithiggs} A point $\overline{q}$ is \emph{GIT semistable} if
for all weighted filtrations $(V_\bullet,n_\bullet)$
$$\sum_{i=1}^{t}  n_{i} ( P(l)\dim V_{i} - P_{\mathcal{E}_{V_{i}}}(l)\dim V)\leq 0\; .$$
A point $\overline{q}$ is \emph{GIT stable} if we get a strict inequality for every weighted filtration.
\end{prop}

Then, the next result completes the sketch of the construction of a moduli space for Higgs sheaves.

\begin{thm}\cite[Theorem 1.19]{Si1}
\label{GITstabisstabhiggs} Fix a polynomial $P$ of degree $n=\dim
X$. There exist integers $m_{0}$ and $l_{0}$ ($l_{0}$ depending on
$m_{0}$) such that for $m\geq m_{0}$ and $l\geq l_{0}$, a point
$q$ in $\Hilb(V\otimes \mathcal{O}_{Z}(-m),P)$ is GIT semistable
(for the action of $SL(V)$ with respect to the embedding into a
projective space), if and only if the quotient $\mathcal{E}$ is a semistable
coherent sheaf of pure dimension $n$ and the map $V\rightarrow
H^{0}(E(m))$ is an isomorphism.
\end{thm}

Let $(E,\varphi)$ be an unstable Higgs sheaf and $\mathcal{E}$ its associated coherent sheaf by
Lemma \ref{higgsiscoherent}. Let $m_{0}$, $l_{0}$ be integers as in Theorem
\ref{GITstabisstabhiggs}. Choose $m_{1}\geq m_{0}$, $l_{1}\geq
l_{0}$ such that $\mathcal{E}$ is also $m_{1}$-regular. Now choose $m\geq
m_{1}$ and fix an isomorphism $V\simeq H^{0}(\mathcal{E}(m))$. Given a
weighted filtration $(V_{\bullet},n_{\bullet})$, define for each
$l\geq l_{1}$ the function

$$\mu(V_{\bullet},n_{\bullet})=\frac{\sum_{i=1}^{t}  n_{i} (P(l)\dim V_{i}-P_{\mathcal{E}_{V_{i}}}(l)\dim V)}
{\sqrt{\sum_{i=1}^{t+1} {\dim V^{i}_{}} \Gamma_{i}^{2}}}\; ,$$
which is a \emph{Kempf function} for this problem (c.f. Definition \ref{kempffunction}). Note that, for each integer $m$, this is a
polynomial function on $l$, whose coefficients depend on $m$.

Let
\begin{equation}
\label{vkempf-filthiggs} 0\subset V_{1}\subset \cdots \subset
V_{t+1}= V
\end{equation}
be the \emph{Kempf filtration} of vector spaces
given by By Theorem \ref{kempftheoremsheaves}, and let
\begin{equation}
\label{ekempf-filthiggs} 0\subseteq \mathcal{E}^{m}_1 \subseteq
\mathcal{E}^{m}_2 \subseteq \;\cdots\; \subseteq \mathcal{E}^{m}_t \subseteq
\mathcal{E}^{m}_{t+1}=\mathcal{E}\; ,
\end{equation}
be the \emph{$m$-Kempf filtration of $\mathcal{E}$}, where $\mathcal{E}_{i}^{m}\subset \mathcal{E}$
is the subsheaf generated by $V_{i}$ under the evaluation map. Making the correspondence of Lemma \ref{higgsiscoherent}, we call the \emph{$m$-Kempf
filtration of the Higgs sheaf $(E,\varphi)$} to
\begin{equation}
\label{ekempf-filthiggs2} 0 \subseteq (E_{1}^{m}, \varphi|_{E_{1}^{m}}) \subseteq
(E_{2}^{m}, \varphi|_{E_{2}^{m}}) \subseteq \;\cdots\; \subseteq (E_{t}^{m}, \varphi|_{E_{t}^{m}}) \subseteq
(E_{t+1}^{m}, \varphi|_{E_{t+1}^{m}})=(E,\varphi)\; ,
\end{equation}
where $E_{i}^{m}=\pi_{\ast}\mathcal{E}^{m}_{i}$.

\begin{rem}
 \label{independentofl}
Recall that for two rational (in particular polynomial) functions,
we define an ordering by saying that $f_{1}\prec f_{2}$ if
$f_{1}(l)<f_{2}(l)$, for $l\gg 0$. Then, although in the
construction of the moduli space and in the definition of the
Kempf function we use the integer $l$, we view the Kempf function
as a polynomial function on $l$, having fixed previously $m$. We
define the $m$-Kempf filtration of $\mathcal{E}$ as the one which
maximizes the Kempf function having fixed $m$, seen as a
polynomial function on $l$, by Theorem \ref{kempftheoremsheaves}.
Note that we also talk about the $m$-Kempf filtration of
$\mathcal{E}$, without mentioning $l$.
\end{rem}

We will proceed as in the previous sections of the chapter to prove the
following theorem:

\begin{thm}
\label{kempfstabilizeshiggs}
There exists an integer $m'\gg 0$ such that the $m$-Kempf filtration of $\mathcal{E}$ is independent of $m$, for $m\geq m'$.
\end{thm}

\subsection{The $m$-Kempf filtration stabilizes with $m$}

We will prove Theorem \ref{kempfstabilizeshiggs} in an analogous
way to the cases of torsion free sheaves and holomorphic pairs in
sections \ref{kempfsheaves} and \ref{kempfpairs}. First, we
associate a graph to the $m$-Kempf filtration of $\mathcal{E}$.

\begin{dfn}
\label{graphhiggs} Let $m\geq m_{1}$ and $l\geq l_{1}$. Given $0\subset V_{1}\subset \cdots \subset V_{t+1}= V$, a filtration
of vector spaces of $V$, define
$$v_{m,i}(l)=m^{n+1}\cdot \frac{1}{\dim V^{i}\dim V}\big[P_{\mathcal{E}_{V^{i}}}(l)\dim V-P(l)\dim V^{i}\big]\; ,$$
$$b_{m}^{i}=\dfrac{1}{m^{n}}\dim V^{i}>0\; ,$$
$$w_{m}^{i}(l)=-b_{m}^{i}\cdot v_{m,i}(l)=m\cdot \frac{1}{\dim V}\big[P(l)\dim V^{i}-P_{\mathcal{E}_{V^{i}}}(l)\dim V\big]\; .$$
Also let
$$b_{m,i}=b_{m}^{1}+\ldots +b_{m}^{i}=\dfrac{1}{m^{n}}\dim V_{i}\; ,$$
$$w_{m,i}(l)=w_{m}^{1}(l)+\ldots +w_{m}^{i}(l)=m\cdot \frac{1}{\dim V}\big[P(l)\dim V_{i}-P_{\mathcal{E}_{V_{i}}}(l)\dim V\big]\; .$$
We call the graph defined by points $(b_{m,i},w_{m,i}(l))$ the \emph{graph associated to the filtration}
$V_{\bullet}\subset V$. Note that, having fixed $m$, the coordinates of the graph are polynomials on $l$.
\end{dfn}

We use a similar argument that in Proposition \ref{identification} to identify
the Kempf function in Theorem \ref{kempftheoremsheaves},
$$\mu(V_{\bullet},n_{\bullet})=\frac{1}{l^{n}}\frac{\sum_{i=1}^{t} n_{i}(P(l)\dim V_{i}-P_{\mathcal{E}_{V_{i}}}(l)\dim V)}
{\sqrt{\sum_{i=1}^{t+1} {\dim V^{i}} \Gamma_{i}^{2}}}\; ,$$ where $n_{i}=\frac{\Gamma_{i}-\Gamma_{i-1}}{\dim
V}$, with the function in Theorem \ref{maxconvexenvelope}, where the coordinates of the graph are as in Definition
\ref{graphhiggs}.

\begin{prop}
\label{identificationhiggs}
For every integer $m$, the following equality holds
$$\mu(V_{\bullet},n_{\bullet})=\frac{1}{m^{(\frac{n}{2}+1)}}\cdot \mu_{v_{m}(l)}(\Gamma)=
\frac{1}{m^{(\frac{n}{2}+1)}}\cdot\frac{(\Gamma,v_{m}(l))}{\|\Gamma\|}\; ,$$
between the Kempf function on Theorem \ref{kempftheoremsheaves} and the function in Theorem \ref{maxconvexenvelope}.
\end{prop}

\begin{rem}
\label{growthhiggs}
Again, we introduce factors $m^{n+1}$ in Definition \ref{graphhiggs} for $v_{m,i}(l)$ and $b_{m}^{i}$ to have
order zero on $m$ (c.f. Remark \ref{growth}). As a result, the graph keeps its dimensions when $m$ grows, the coordinates being
polynomials on the variable $l$.
\end{rem}

In the following, we will omit the integers $m$ and $l$ for the quantities $v_{m,i}(l)$, $b_{m,i}$, $w_{m,i}(l)$
in the definition of the graph associated to a filtration of vector spaces, where it is clear from the context.

We give the analogous to Propositions \ref{boundedness} and
\ref{task} for the case of Higgs sheaves, using Lemmas
\ref{lemmaA} and \ref{lemmaB}.

Define
\begin{equation}
\label{C_constant_higgs}
C=\max\{r|\mu _{\max }(E)|+\frac{d}{r}+r|B|+|A|+1\;,\;1\},
\end{equation}
a positive constant, where recall that $r=\rk \mathcal{E}=\rk E$ and $d=\deg\mathcal{E}=\deg E$, $E=\pi_{\ast}\mathcal{E}$,
are the rank and the degree of $\mathcal{E}$ as the corresponding coefficients
of the Hilbert polynomial for a pure sheaf.

\begin{prop}
\label{boundednesshiggs} Given sufficiently large integers $m$ and $l$,
each filter in the $m$-Kempf filtration of $\mathcal{E}$ has slope $\mu(\mathcal{E}^{m}_{i})\geq \dfrac{d}{r}-C$.
\end{prop}

\begin{pr}
Choose an integer $m_{2}\geq m_{1}$ such that for $m\geq m_{2}$
$$[\mu_{max}(\mathcal{E})+gm+B]_{+}=\mu_{max}(\mathcal{E})+gm+B\; ,$$
and
$$[\frac{d}{r}-C+gm+B]_{+}=\frac{d}{r}-C+gm+B\; .$$
Now let $m\geq m_{2}$. Let
$$0 \subseteq \mathcal{E}^{m}_1 \subseteq \mathcal{E}^{m}_2 \subseteq \;\cdots\; \subseteq \mathcal{E}^{m}_t \subseteq \mathcal{E}^{m}_{t+1}=\mathcal{E}$$
be the $m$-Kempf filtration of $\mathcal{E}$. 

Let $\mathcal{E}^{m}_{i}\subseteq \mathcal{E}$ a subsheaf of rank $r_{i}$ and degree $d_{i}$, such that $\mu
(\mathcal{E}^{m}_{i})<\frac{d}{r}-C$, and suppose that $\mathcal{E}_{i}^{m}(m)\subset \mathcal{E}(m)$ satisfies the estimate in Lemma
\ref{Simpson}. Analogously to Proposition \ref{boundedness},
$$h^{0}(\mathcal{E}^{m}_{i}(m))\leq \frac{1}{g^{n-1}n!}\big ((r_{i}-1)(\mu_{max}(\mathcal{E})+gm+B)^{n}+(\frac{d}{r}-C+gm+B)^{n}\big)=G(m)\; ,$$
where
$$G(m)=\frac{1}{g^{n-1}n!}\big
[r_{i}g^{n}m^{n}+ng^{n-1}\big((r_{i}-1)\mu_{max}(\mathcal{E})+\frac{d}{r}-C+r_{i}B\big)m^{n-1}+\cdots \big ]\; .$$

By Definition \ref{graphhiggs}, to such filtration we associate a graph with heights, for each $j$
$$w_{m,j}(l)=w_{m}(l)^{1}+\ldots +w_{m}(l)^{j}=m\cdot \frac{1}{\dim V}\big[P(l)\dim V_{j}-P_{\mathcal{E}^{m}_{i}}(l)\dim V\big]\; ,$$
and we will show that $w_{m,i}(l)<0$ for $m$ and $l$ large enough,
to get a contradiction as in Proposition \ref{boundedness}. Note that the coordinates $w_{m,i}(l)$ are polynomials on $l$, so
$w_{m,i}(l)<0$ for $l\gg 0$, contradicts Lemma \ref{lemmaA} for the Kempf function constructed with the integer $m$.

Given that $\mathcal{E}^{m}_{i}(m)$ is generated by $V_{i}$ under the
evaluation map, it is $\dim V_{i}\leq h^{0}(\mathcal{E}^{m}_{i}(m))$, hence
$$w_{m,i}(l)=\frac{m}{\dim V}\big[P(l)\dim V_{i}-P_{\mathcal{E}^{m}_{i}}(l)\dim V\big]\leq$$
$$\frac{m}{P_{\mathcal{E}}(m)}\big[P(l)h^{0}(\mathcal{E}^{m}_{i}(m))-P_{\mathcal{E}^{m}_{i}}(l)P_{\mathcal{E}}(m)\big]\leq
\frac{m}{P_{\mathcal{E}}(m)}\big[P(l)G(m)-P_{\mathcal{E}^{m}_{i}}(l)P_{\mathcal{E}}(m)\big]\; .$$

Hence, $w_{m,i}(l)<0$ is equivalent to
$$\Phi_{m}(l)=P(l)G(m)-P_{\mathcal{E}^{m}_{i}}(l)P(m)<0\: ,$$
where $\Phi_{m}(l)$ can be seen as an $n^{th}$-order polynomial on $l$,
$$\Phi_{m}(l)=\alpha_{n}(m)l^{n}+\alpha_{n-1}(m)l^{n-1}+\cdots
+\alpha_{1}(m)l+\alpha_{0}(m)\: .$$ Hence, it is sufficient to
show that $\alpha_{n}(m)<0$ for an integer $m$ sufficiently large.

Note that
$$\alpha_{n}(m)=rG(m)-r_{i}P(m)<0\; ,$$
is the same polynomial as $\Psi(m)$ in the proof of Proposition
\ref{boundedness}. Hence, by the same argument
$$\alpha_{n}(m)=\xi_{n-1}m^{n-1}+\cdots +\xi_{1}m+\xi_{0}$$ with
$\xi_{n-1}<0$, so there exists $m_{3}\geq m_{2}$ such that for
$m\geq m_{3}$ we will have $\alpha_{n}(m)<0$. Hence, there exists
an integer $l\gg 0$, depending on $m_{3}$, such that for $m\geq m_{3}$ we have
 $\Phi_{m}(l)\prec 0$ as a polynomial (c.f. Remark \ref{independentofl}), hence $w_{m,i}(l)<0$,
which is a contradiction.
\end{pr}

Now we can assure the $m$-regularity of the family of the
subsheaves appearing in the different $m$-Kempf filtrations of
$\mathcal{E}$, similarly to Proposition \ref{regular}.

\begin{prop}
\label{regularhiggs} There exists an integer $m_{4}$ such that for
$m\geq m_{4}$ the sheaves $\mathcal{E}^{m}_{i}$ and
$\mathcal{E}^{m,i}=\mathcal{E}^m_i/\mathcal{E}^m_{i-1}$ are
$m_{4}$-regular. In particular their higher cohomology groups,
after twisting with $\mathcal{O}_{X}(m_{4})$, vanish and they are
generated by global sections.
\end{prop}

\begin{prop}
\label{taskhiggs}
Let $m\geq m_{4}$. For each filter $\mathcal{E}_{i}^{m}$
in the $m$-Kempf filtration, we have $\dim
V_{i}=h^{0}(\mathcal{E}_{i}^{m}(m))$, therefore $V_{i}\cong
H^{0}(\mathcal{E}_{i}^{m}(m))$.
\end{prop}

\begin{pr}
The proof follows analogously to the proof of Proposition \ref{task}.

Let $m\geq m_{4}$. Let $V_{\bullet}\subseteq V$ be the Kempf
filtration of $V$ (cf. (\ref{vkempf-filthiggs})) and let
$\mathcal{E}_{\bullet}^{m}\subseteq \mathcal{E}$ be the $m$-Kempf filtration of $\mathcal{E}$
(c.f. (\ref{ekempf-filthiggs})). We know that each $V_{i}$ generates
the subsheaf $\mathcal{E}_{i}^{m}$, by definition, then we have the diagram:

$$\begin{array}{ccccccccccc}
    0 & \subset & V_{1} & \subset & V_{2} & \subset & \cdots & \subset & V_{t+1} & = & V \\
     & & \cap & & \cap & & & & & & ||  \\
      &     & H^{0}(\mathcal{E}_{1}^{m}(m)) & \subset & H^{0}(\mathcal{E}_{2}^{m}(m)) & \subset & \cdots & \subset & H^{0}(\mathcal{E}_{t+1}^{m}(m))&
      = & H^{0}(\mathcal{E}(m))
      \end{array}$$

Let $i$ be the first index such that $V_{i}\neq H^{0}(\mathcal{E}_{i}^{m}(m))$, then we have the diagram:

\begin{equation}
 \label{filtrationVhiggs}
    \begin{array}{ccc}
    V_{i} & \subset & V_{i+1}\\
    \cap & & || \\
    H^{0}(\mathcal{E}_{i}^{m}(m)) & \subset & H^{0}(\mathcal{E}_{i+1}^{m}(m))
    \end{array}
\end{equation}

Therefore we consider a new filtration by adding
the filter $H^{0}(\mathcal{E}_{i}^{m}(m))$

\begin{equation}
\label{filtrationV'higgs}
    \begin{array}{ccccccccccccc}
    V_{i} & \subset & H^{0}(\mathcal{E}_{i}^{m}(m)) & \subset & V_{i+1}\\
    || & & || & & ||\\
    V'_{i} & &  V'_{i+1} & & V'_{i+2}
    \end{array}
\end{equation}

Then, $V_{i}$ and $H^{0}(\mathcal{E}_{i}^{m}(m))$ generate the same
sheaf $\mathcal{E}_{i}^{m}$, hence we are in situation of Lemma
\ref{lemmaB}, where $W=H^{0}(\mathcal{E}_{i}^{m}(m))$, filtration
$V_{\bullet}$ is \eqref{filtrationVhiggs} and filtration
$V'_{\bullet}$ is \eqref{filtrationV'higgs}.

Now the graph associated to filtration $V_{\bullet}$ is, by Definition
\ref{graphhiggs}, given by the points
$$(b_{m,i},w_{m,i}(l))=\big(\dfrac{\dim V_{i}}{m^{n}},\frac{m}{\dim V}(P(l)\dim V_{i}-P_{\mathcal{E}^{m}_{i}}(l)\dim V)\big)\; ,$$
and the slopes $-v_{m,i}(l)$ of the graph are given by

$$-v_{m,i}(l)=\frac{w_{m}^{i}(l)}{b_{m}^{i}}=\frac{w_{m,i}(l)-w_{m,i-1}(l)}{b_{m,i}-b_{m,i-1}}=$$
$$\frac{m^{n+1}}{\dim V}\big(P(l)-P_{\mathcal{E}^{i,m}}(l)\frac{\dim V}{\dim V^{i}}\big)\; ,$$
which is an $n^{th}$-order polynomial on $l$ whose leading coefficient
is
$$\alpha^{i}(m)=\frac{m^{n+1}}{\dim V}\big(r-r^{i}\frac{\dim V}{\dim V^{i}}\big)\leq \frac{m^{n+1}}{\dim V}\cdot r:=R\; .$$
Equality holds if and only if $r^{i}=0$.

The new point which appears in graph of the filtration
$V'_{\bullet}$ is
$$Q=\big(\dfrac{h^{0}(\mathcal{E}_{i}^{m}(m))}{m^{n}},\frac{m}{\dim V}(P(l)h^{0}(\mathcal{E}_{i}^{m}(m))-P_{\mathcal{E}^{m}_{i}}(l)\dim V)\big)\; ,$$
joining two new segments appearing in this new graph. The slope of
the segment between $(b_{m,i},w_{m,i}(l))$ and $Q$ is, similarly,
$$-v'_{m,i}(l)=\dfrac{m^{n+1}}{\dim V}\cdot P(l)\; ,$$
again an $n^{th}$-order polynomial on $l$ whose leading coefficient is
$$\alpha'^{i}(m)=\dfrac{m^{n+1}}{\dim V}\cdot r=R\: .$$
By Lemma \ref{lemmaA}, the graph is convex, so
$v_{m,1}(l)<v_{m,2}(l)<\ldots<v_{m,t+1}(l)$. On the other hand, by
Lemma \ref{lemmaB}, $v'_{m,i}(l)\geq v_{m,i}(l)$. Therefore for a
sufficiently large $l$ we have the following inequalities between
the leading coefficients of the $-v'_{m,i}(l)$,
$$\alpha^{1}(m)\geq\alpha^{2}(m)\geq\ldots \geq\alpha^{t+1}(m)\; ,$$
and
$$\alpha'^{i}(m)\leq \alpha^{i}(m)\; .$$
Besides, $r^{1}=r_{1}>0$, then $R>\alpha^{1}(m)$. Indeed, $\mathcal{E}$ is
pure, then it has no torsion elements on its support, hence also the subsheaf $E_{1}^{m}$, and a rank $0$
pure sheaf should be the zero sheaf. Hence
$$R>\alpha^{1}(m)\geq\alpha^{2}(m)\geq\ldots \geq\alpha^{i}(m)\geq \alpha'^{i}(m)=R\; ,$$
which is a contradiction.

Therefore, for $m\geq m_{4}$, every filter in the $m$-Kempf filtration of $\mathcal{E}$ verifies $\dim V_{i}=h^{0}(\mathcal{E}_{i}^{m}(m))$.
\end{pr}

\begin{cor}
\label{rankhiggs} Given $m\geq m_{4}$, for every filter $\mathcal{E}_{i}^{m}$ in the $m$-Kempf
filtration, it is $r^{i}>0$.
\end{cor}
\begin{pr}
It follows from Proposition \ref{taskhiggs} the same way as in
Corollary \ref{rank}. Indeed, given that $r^{i}=0$ is equivalent
to $\alpha^{i}(m)=R$, note that it is $r^{1}=r_{1}>0$ and
$R>\alpha^{1}(m)\geq\alpha^{2}(m)\geq\ldots\geq\alpha^{t+1}(m)$.
\end{pr}

Next, we again recall the results on subsection
\ref{sectionkempfstabilizes}.

By Proposition \ref{regularhiggs}, for any $m\geq m_{4}$, all the
filters $\mathcal{E}^{m}_{i}$ of the $m$-Kempf filtration of $\mathcal{E}$ are $m_{4}$-regular and hence, the filtration of
sheaves
$$0\subset \mathcal{E}^{m}_{1} \subset \mathcal{E}^{m}_{2} \subset \cdots \subset \mathcal{E}^{m}_{t_{m}} \subset \mathcal{E}^{m}_{t_{m}+1}=\mathcal{E}$$
is the filtration associated to the filtration of vector subspaces
$$0\subset H^{0}(\mathcal{E}^{m}_{1}(m_{4})) \subset H^{0}(\mathcal{E}^{m}_{2}(m_{4})) \subset \cdots \subset H^{0}(\mathcal{E}^{m}_{t_{m}}(m_{4}))
\subset H^{0}(\mathcal{E}^{m}_{t_{m}+1}(m_{4}))=H^{0}(\mathcal{E}(m_{4}))$$ by the evaluation map (c.f. Lemma \ref{mregularity}), of
a unique vector space $H^{0}(\mathcal{E}(m_{4}))$, whose dimension does not depend on $m$. Let $P^{m}_{i}:=P_{\mathcal{E}^{m}_{i}}$ and $P^{i,m}:=P_{\mathcal{E}^{i,m}}$. Let
$$(P_{1}^{m},\ldots,P_{t_{m}+1}^{m})$$ be the \emph{$m$-type}
of the $m$-Kempf filtration of $\mathcal{E}$ (c.f. Definition
\ref{mtype}) and let
$$\mathcal{P}=\big\{(P_{1}^{m},\ldots, P_{t_{m}+1}^{m})\big\}$$
be the set of possible $m$-types, which is a finite set (c.f. Proposition
\ref{Pisfinite}).

By Definition \ref{graphhiggs} we associate a graph to the
$m$-Kempf filtration of $\mathcal{E}$, given by $v_{m}(l)$. By Propositions
\ref{regularhiggs} and \ref{taskhiggs} it can be rewritten as
$$v_{m,i}(l)=m^{n+1}\cdot \frac{1}{P^{i,m}(m)P(m)}\big[P^{i,m}(l)P(m)-P(l)P^{i,m}(m)\big]\; ,$$
$$b_{m}^{i}=\frac{1}{m^{n}}\cdot P^{i,m}(m)\; .$$

Note that, given the $m$-Kempf filtration, its $m$-type is fixed.
Hence, the coordinates of the graph, $v_{m,i}(l)$ are polynomials on $l$, whose coefficients are fixed
(c.f. Definition \ref{graphhiggs}). Then, fixing the $m$-type, $v_{m}(l)$ defines a different graph for each $l$.

We define a functional on $\mathcal{P}$ which assigns to each $m$-type (to each $m$-Kempf filtration) the
function
$$\Theta_{m}(l)=(\mu_{v_{m}(l)}(\Gamma_{v_{m}(l)}))^{2}=||v_{m}(l)||^{2}\; ,$$ which is, given $m$, a polynomial function on $l$ (c.f. \eqref{maxvalue}).
By finiteness of $\mathcal{P}$ there is a finite list of such
possible functions
$$
\mathcal{A}=\{\Theta_{m}:m\geq m_{4}\}\; ,
$$
so, analogously to Lemma \ref{uniquefunction}, we can choose $K$ to be a polynomial function such that there exists an integer $m_{5}$ with $\Theta_{m}=K$,
for all $m\geq m_{5}$, meaning that
two polynomial functions do coincide if they do for large values of the variable.

\begin{prop}
\label{eventuallyhiggs} Let $a_{1}$ and $a_{2}$ be integers with $a_{1}\geq a_{2} \geq m_{5}$.
The $a_{1}$-Kempf filtration of $\mathcal{E}$ is equal to the $a_{2}$-Kempf filtration
of $\mathcal{E}$.
\end{prop}
\begin{pr}
The proof follows analogously to the proof of Proposition \ref{eventually}.
\end{pr}

\begin{dfn}
\label{kempffiltrationhiggs}
If $m\geq m_{5}$, and for $l\geq l_{5}$, the $m$-Kempf filtration of $\mathcal{E}$ is called \emph{the Kempf filtration
of $\mathcal{E}$},
$$0\subset \mathcal{E}_{1} \subset \mathcal{E}_{2} \subset \cdots \subset \mathcal{E}_{t} \subset \mathcal{E}_{t+1}=\mathcal{E}\; .$$
\end{dfn}

By applying $\pi_{\ast}$ to the Kempf filtration we obtain the following definition.

\begin{dfn}
 \label{kempffiltrationhiggs2}
The following filtration
$$0 \subset (E_{1}, \varphi|_{E_{1}}) \subset
(E_{2}, \varphi|_{E_{2}}) \subset \;\cdots\; \subset (E_{t}, \varphi|_{E_{t}}) \subset
(E_{t+1}, \varphi|_{E_{t+1}})=(E,\varphi)\; ,$$
where $E_{i}=\pi_{\ast}\mathcal{E}_{i}$ is called \emph{the Kempf filtration of the Higgs sheaf $(E,\varphi)$}.
\end{dfn}

\subsection{Harder-Narasimhan filtration for Higgs sheaves}

Recall that, by the Kempf theorem (c.f. Theorem \ref{kempftheoremsheaves}), given an integer $m$ and $V\simeq
H^{0}(\mathcal{E}(m))$, there exists a unique weighted filtration of vector spaces $V_{\bullet}\subseteq V$ which gives
maximum for the Kempf function, which in this case is
$$\mu(V_{\bullet},n_{\bullet})=\frac {\sum_{i=1}^{t+1} \frac{\Gamma_i}{\dim V} ( P_{\mathcal{E}^i}(l) \dim V - P(l)\dim V^i)}
 {\sqrt{\sum_{i=1}^{t+1} {\dim V^{i}} \Gamma_{i}^{2}}}\; .$$ This filtration induces the Kempf filtration of $\mathcal{E}$,
$$0\subset \mathcal{E}_{1} \subset \mathcal{E}_{2} \subset \cdots \subset \mathcal{E}_{t} \subset \mathcal{E}_{t+1}=\mathcal{E}$$
which is independent of $m$, for $m\geq m_{5}$, by Proposition \ref{eventuallyhiggs}, hence it only depends on $\mathcal{E}$.

We proceed again as in Section \ref{kempfsheaves} (c.f. Proof of Theorem \ref{kempfisHN}), to rewrite the Kempf
function in terms of Hilbert polynomials of sheaves. We set $P=P_{\mathcal{E}}$, $P^{i}=P_{\mathcal{E}^{i}}$, and recall the
relation
$$\gamma_{i}=\frac{r}{P(m)}\Gamma_{i}\; .$$

\begin{prop}
\label{finalfunctionhiggs} Given $\mathcal{E}$, a pure sheaf of dimension $n$, there exists a unique filtration
$$0\subset \mathcal{E}_{1} \subset \mathcal{E}_{2} \subset \cdots \subset \mathcal{E}_{t} \subset
\mathcal{E}_{t+1}=\mathcal{E}$$ with positive weights $n_{1},\ldots,n_{t}$, which gives maximum for the function
$$K_{m}(l):=m^{\frac{n}{2}+1}\cdot \frac{1}{P(m)}\frac {\sum_{i=1}^{t+1}
\gamma_{i}[P^i(l)P(m)-P^i(l)P(m)]}
 {\sqrt{\sum_{i=1}^{t+1} P^{i}(m)\gamma_{i}^{2}}}\; .$$
\end{prop}

The coordinates of the graph $v_{m,i}(l)$ are given by
$$v_{m,i}(l)=m^{n+1}\cdot \frac{1}{P^{i}(m)P(m)}\big[P^i(l)P(m)-P(l)P^{i}(m)\big]\; .$$
hence, the function $K$ is
$$K_{m}(l)=m^{-\frac{n}{2}}\cdot l^{n}\cdot \frac {\sum_{i=1}^{t+1}
P^{i}(m)\gamma_{i}v_{i}} {\sqrt{\sum_{i=1}^{t+1} P^{i}(m)\gamma_{i}^{2}}}=m^{-\frac{n}{2}}\cdot l^{n}\cdot
\frac{(\gamma,v)}{||\gamma||}\; ,$$ where the scalar product is given by
$$\left(
    \begin{array}{cccc}
      P^{1}(m) &  &  &  \\
       & P^{2}(m) &  &  \\
       &  & \ddots &  \\
       &  &  & P^{t+1}(m) \\
    \end{array}
  \right)$$

\begin{prop}
\label{descendentslopeshiggs} Given the Kempf filtration of a sheaf $\mathcal{E}$,
$$0\subset \mathcal{E}_{1} \subset \mathcal{E}_{2} \subset \cdots \subset \mathcal{E}_{t} \subset \mathcal{E}_{t+1}=\mathcal{E}$$
it verifies
$$\frac{P_{\mathcal{E}^{1}}}{\rk \mathcal{E}^{1}}>\frac{P_{\mathcal{E}^{2}}}{\rk \mathcal{E}^{2}}>\ldots>\frac{P_{\mathcal{E}^{t+1}}}{\rk \mathcal{E}^{t+1}}\; .$$
\end{prop}
\begin{pr}
By Lemma \ref{lemmaA}, the vector $v_{m}(l)$ is convex for $m\geq m_{4}$ and $l\gg 0$. Therefore, seeing $v_{m,i}(l)$ as polynomials on $l$,
$$v_{m,i}(l)<v_{m,i+1}(l)\Leftrightarrow \frac{P_{\mathcal{E}^{i}}(l)}{P_{\mathcal{E}^{i}}(m)}-\frac{P(l)}{P(m)}<\frac{P_{\mathcal{E}^{i+1}}(l)}{P_{\mathcal{E}^{i+1}}(m)}-\frac{P(l)}{P(m)}
\Leftrightarrow$$
$$\frac{\rk(\mathcal{E}^{i})}{P_{\mathcal{E}^{i}}(m)}<\frac{\rk(\mathcal{E}^{i+1})}{P_{\mathcal{E}^{i+1}}(m)}\Leftrightarrow 
\frac{P_{\mathcal{E}^{i}}(m)}{\rk(\mathcal{E}^{i})}>\frac{P_{\mathcal{E}^{i+1}}(m)}{\rk(\mathcal{E}^{i+1})}\; ,$$
where the second equivalence holds for $l\gg 0$.
\end{pr}

\begin{prop}
\label{blocksemistabilityhiggs} Given the Kempf filtration of $\mathcal{E}$,
$$0\subset \mathcal{E}_{1} \subset \mathcal{E}_{2} \subset \cdots \subset \mathcal{E}_{t} \subset \mathcal{E}_{t+1}=\mathcal{E}$$
each one of the blocks $\mathcal{E}^{i}=\mathcal{E}_{i}/\mathcal{E}_{i-1}$ is semistable.
\end{prop}
\begin{pr}
C.f. Proposition \ref{blocksemistability}.
\end{pr}

Theorem \ref{HNunique}, which provides the construction of the
Harder-Narasimhan filtration, holds for pure sheaves (c.f.
\cite[Theorem 1.3.4]{HL3}), as it is the present case of the sheaf
$\mathcal{E}$ supported on the cotangent bundle associated to a
Higgs sheaf $(E,\varphi)$.

\begin{prop}\cite[Theorem 1.3.4]{HL3}
\label{HNpurehiggs} Given a pure sheaf $\mathcal{E}$, there exists a unique filtration
$$0\subset \mathcal{E}_{1} \subset \mathcal{E}_{2} \subset \cdots \subset \mathcal{E}_{t} \subset
\mathcal{E}_{t+1}=\mathcal{E}\; ,$$ which satisfies these two properties
 \begin{enumerate}
   \item The Hilbert polynomials verify
   $$\frac{P_{\mathcal{E}^{1}}}{\rk \mathcal{E}^{1}}>\frac{P_{\mathcal{E}^{2}}}{\rk \mathcal{E}^{2}}>\ldots>\frac{P_{\mathcal{E}^{t+1}}}{\rk \mathcal{E}^{t+1}}$$
   \item Every block $\mathcal{E}^{i}=\mathcal{E}^{i}/\mathcal{E}_{i}$ is semistable.
 \end{enumerate}
This filtration is called the \textbf{Harder-Narasimhan filtration} of $\mathcal{E}$.
\end{prop}

\begin{cor}
\label{kempfisHNhiggs}
The Kempf filtration of a sheaf
$\mathcal{E}$ is the Harder-Narasimhan filtration.
\end{cor}
\begin{pr}
See Propositions \ref{descendentslopeshiggs} and \ref{blocksemistabilityhiggs} and use uniqueness of the
Harder-Narasimhan filtration of $\mathcal{E}$ in Theorem \ref{HNpurehiggs}.
\end{pr}

Making the correspondence between Higgs sheaves $(E,\varphi)$ over $X$ and sheaves $\mathcal{E}$ of pure
dimension $n=\dim X$ over $T^{\ast}X$ we can construct a filtration of Higgs subsheaves
$$0\subset \pi_{\ast}\mathcal{E}_{1} \subset \pi_{\ast}\mathcal{E}_{2} \subset \cdots \subset \pi_{\ast}\mathcal{E}_{t} \subset
\pi_{\ast}\mathcal{E}_{t+1}=E$$
i.e.
$$0 \subset (E_{1}, \varphi|_{E_{1}}) \subset
(E_{2}, \varphi|_{E_{2}}) \subset \;\cdots\; \subset (E_{t},
\varphi|_{E_{t}}) \subset (E_{t+1},
\varphi|_{E_{t+1}})=(E,\varphi)\; ,$$ where $E_{i}=\pi_{\ast}\mathcal{E}_{i}$, which coincides with the
Harder-Narasimhan filtration of a Higgs sheaf $(E,\varphi)$, which
appears in the literature (c.f. \cite{AB}).

\section{Rank $2$ tensors}
\label{kempfrk2}

In this section we study the case of tensors where the sheaf has rank $2$. The
moduli space of tensors has been studied in section
\ref{exampletensors}. Here we prove the analogous correspondence between Kempf and Harder-Narasimhan filtrations
for rank $2$ tensors, similarly to the cases of torsion free sheaves, holomorphic pairs and Higgs sheaves.

\vspace{1cm}

Let $X$ be a smooth complex projective variety of dimension $n$.
Let $E$ be a coherent torsion free sheaf over $X$, of rank $2$. We
call a \emph{rank $2$ tensor} the pair consisting of
$$
(E,\varphi:\overbrace{E\otimes\cdots \otimes E}^{\text{s
times}}\too \mathcal{O}_{X})\; .
$$
These objects are particular cases of the ones studied in section \ref{exampletensors} for arbitrary $s$, $c=1$,
$b=0$, $R=\Spec \mathbb{C}$ and $\mathcal{D}=\mathcal{O}_{X}$, meaning the structure sheaf over $X\times R\simeq X$,
in Definition \ref{deftensor}.

Let $\delta$ be a polynomial of degree at most $\dim X-1=n-1$ and
positive leading coefficient. Recall the definition of
$\delta$-stability for tensors. Recall Definition
\ref{stabilityfortensors} and calculation made in
(\ref{rightmu2}). Recall Remark \ref{saturatedfiltrations} which
says that, in Definition \ref{stabilityfortensors} it suffices to
check the condition on filtrations with $\rk E_{i}<\rk E_{i+1}$.
Hence, as the rank of $E$ is $2$, the only filtrations we have to check are one-step filtrations, i.e.
subsheaves of rank $1$, and we can rewrite the stability condition as follows:

\begin{dfn}
\label{stabilitytensorsrk2}
A rank $2$ tensor $(E,\varphi)$ is \emph{$\delta$-semistable} if for
every rank $1$ subsheaf $L\subset E$
\begin{equation}
\label{stabilityrk2}
 (2P_{L}-P_{E})+\delta(s-2\epsilon(L))\leq 0,
\end{equation}
where $\epsilon(L)$ is the number of times that $L$ appears in the
multi-index $(i_1,\ldots,i_s)$ giving the minimum in
(\ref{rightmu}) and $P_{E}$, $P_{L}$ are the Hilbert polynomials
of $E$ and $L$ respectively. If the inequality is strict for every $L$, we say that $(E,\varphi)$ is \emph{$\delta$-stable}.
If $(E,\varphi)$ is not $\delta$-semistable, we say that it is \emph{$\delta$-unstable}.
\end{dfn}

\subsection{Moduli space of $\rk 2$ tensors}

We recall the main points of the construction of the moduli space for tensors with fixed determinant
$\det(E)\cong \Delta$ of degree $d$ and $\rk(E)=2$. The general construction was explained in section
\ref{exampletensors}, following Gieseker's method. The present case can be obtained by setting $c=1$, $b=0$,
arbitrary $s$, $R=\Spec \mathbb{C}$ and $\mathcal{D}=\mathcal{O}_{X}$, the structure sheaf over $X\times
R\simeq X$, in Definition \ref{deftensor}.

Let $V$ be a vector space of dimension $p:=h^0(E(m))$, where $m$
is a suitable large integer (in particular, $E(m)$ is generated by
global sections and $h^i(E(m))=0$ for $i>0$). Given an isomorphism
$V\cong H^0(E(m))$ we obtain a point
$$(\overline{Q},\overline{\Phi}) \in
\PP(\Hom(\wedge^r V , A) )\times \PP(\Hom(V^{\otimes s},B))\; .
$$
If we change the isomorphism $\det(E)\cong \Delta$, we obtain a different point in the line defined by $Q$.
Likewise, if we change the isomorphism $V\cong H^0(E(m))$ by a homothecy, we obtain a different point in the
line defined by $Q$. In both cases, the point $\overline{Q}$ in the projective space is the same. The same
applies for $\overline{\Phi}$. If we fix once and for all a basis of $V$, then giving an isomorphism between $V$
and $H^0(E(m))$ is equivalent to giving a basis of $H^0(E(m))$. A change of basis is given by an element of
$\GL(V)$, but, since an homothecy does not change the point $(\overline{Q},\overline{\Phi})$, when we want to
get rid of this choice it is enough to divide by the action of $\SL(V)$.

A weighted filtration $(V_\bullet,n_\bullet)$ of $V$ is a filtration
\begin{equation}
\label{filtVrk2} 0 \subset V_1 \subset V_2 \subset \;\cdots\; \subset V_t \subset V_{t+1}=V,
\end{equation}
and rational numbers $n_1,\, n_2,\ldots , \,n_t > 0$, and recall
that this is equivalent to giving a $1$-parameter subgroup
$\Gamma: \CC^*\to \SL(V)$ (c.f. subsection \ref{modulisheaves})
represented by the vector
$$
\Gamma=(\overbrace{\Gamma_1,\ldots,\Gamma_1}^{\dim V^1},
\overbrace{\Gamma_2,\ldots,\Gamma_2}^{\dim V^2},
\ldots,
\overbrace{\Gamma_{t+1},\ldots,\Gamma_{t+1}}^{\dim V^{t+1}}) \; .
$$

By the Hilbert-Mumford criterion (c.f. Theorem \ref{HMcrit}), a point
$$(\overline{Q},\overline{\Phi}) \in
\PP(\Hom(\wedge^r V , A) )\times \PP(\Hom(V^{\otimes s},B))
$$
is \emph{GIT semistable} with respect to the natural linearization on $\SO(a_1,a_2)$ if and only if for all weighted
filtrations
$$
\mu(\overline{Q},V_\bullet,n_\bullet) + \frac{a_2}{a_1} \mu(\overline{\Phi},V_\bullet,
n_\bullet) \leq 0\; ,$$
and recall the numerical function which has to be calculated to apply Mumford criterion for GIT stability (c.f.
Proposition \ref{GITstab}).
\begin{prop}
\label{GITstabrk2}
A point $(\overline{Q},\overline{\Phi})$ is \emph{GIT $a_2/a_1$-semistable} if for all
weighted filtrations $(V_\bullet,n_\bullet)$,
$$
\sum_{i=1}^{t}  n_i ( r \dim V_i - r_i \dim V) +\frac{a_2}{a_1} \sum_{i=1}^t n_i \big( s\dim
V_i-\epsilon_i(\overline\Phi)\dim V \big) \leq 0\; .
$$
\end{prop}

Here, $E_{V_i}$ is the subsheaf of $E$ generated by $V_i$ and $r_{i}=\rk E_{V_{i}}$. If $I=(i_1,\ldots,i_s)$ is the multi-index giving
minimum in (\ref{rightmuV}) (c.f. Section \ref{exampletensors}), we will denote by
$\epsilon_i(\overline{\Phi},V_\bullet,n_\bullet)$ (or just $\epsilon_i(\overline{\Phi})$ if the rest of the data
is clear from the context) the number of elements $k$ of the multi-index $I$ such that $\dim V_k\leq \dim V_i$.
Let $\epsilon^i(\overline\Phi)=\epsilon_i(\overline\Phi)-\epsilon_{i-1}(\overline\Phi)$.

Then, recall Theorem \ref{GIT-delta}:
\begin{thm}
\label{GIT-deltark2} Let $(E,\varphi)$ be a tensor. There exists
an $m_{0}$ such that, for $m\geq m_{0}$ the associated point
$(\overline{Q},\overline{\Phi})$ is GIT $a_2/a_1$-semistable if
and only if the tensor is $\delta$-semistable, where
$$
\frac{a_2}{a_1} = \frac{r\delta(m)}{P_E(m)-s\delta(m)}\; .
$$
\end{thm}

Let $X$ be a smooth complex projective variety of dimension $n$. Let us consider rank $2$ tensors
$$
(E,\varphi:\overbrace{E\otimes\cdots \otimes E}^{\text{s
times}}\too \mathcal{O}_{X})
$$
given by a coherent torsion free sheaf $E$ of rank $2$ over $X$
with fixed determinant $\det(E)\cong \Delta$ and a morphism
$\varphi$ from a tensor product of $s$ copies of $E$ to the
trivial line bundle $\mathcal{O}_{X}$. Let $\delta$ be a
polynomial of degree at most $\dim X-1=n-1$ and positive leading
coefficient.

Let $(E,\varphi)$ be a $\delta$-unstable rank $2$ tensor. Let $m_{0}$ be an integer as in Theorem
\ref{GIT-deltark2} (i.e. such that the $\delta$-stability and the GIT stability coincide) and also such that $E$
is $m_{0}$ regular (choosing a larger integer, if necessary). Choose an integer $m\geq m_{0}$ and let $V$ be a
vector space of dimension $P_{E}(m)=h^{0}(E(m))$.

Given a filtration of vector subspaces $0\subset V_{1}\subset \cdots \subset V_{t+1}= V$ and positive numbers
 $n_{1},\cdots,n_{t}>0$, i.e., given a weighted filtration, we define the following function
$$\mu(V_{\bullet},n_{\bullet})=\frac{\sum_{i=1}^{t}  n_{i} (r\dim V_{i}-r_{i}\dim V)+
\frac{a_{2}}{a_{1}}\sum_{i=1}^{t} n_{i} \big(s\dim
V_{i}-\epsilon_{i}(\overline\Phi)\dim V \big)}
{\sqrt{\sum_{i=1}^{t+1} {\dim V^{i}_{}} \Gamma_{i}^{2}}}\; ,$$
which is a \emph{Kempf function} for this problem (c.f. Definition
\ref{kempffunction}), where the numerator of the function
coincides with the numerical function in Proposition
\ref{GITstabrk2} and the denominator is a length $||\Gamma ||$ in
the space of $1$-parameter subgroups (c.f. Definition
\ref{length}).

Let
\begin{equation}
\label{vkempf-filtrk2}
0\subset V_{1}\subset \cdots \subset V_{t+1}= V
\end{equation} be the \emph{Kempf filtration} of $V$
(c.f. Theorem \ref{kempftheoremsheaves}), and let
\begin{equation}
\label{ekempf-filtrk2}
0\subseteq (E^{m}_{1},\varphi|_{E_{1}^{m}})\subseteq
(E^{m}_{2},\varphi|_{E_{2}^{m}})\subseteq\cdots
(E^{m}_{t},\varphi|_{E_{t}^{m}})\subseteq
(E^{m}_{t+1},\varphi|_{E_{t+1}^{m}})\subseteq (E,\varphi)
\end{equation}
be the \emph{$m$-Kempf filtration} of the rank $2$ tensor $(E,\varphi)$, where $E_{i}^{m}\subset E$ is the subsheaf
generated by $V_{i}$ under the evaluation map.

We will prove the following
\begin{thm}
\label{kempfstabilizesrk2} There exists an integer $m'\gg 0$ such that the $m$-Kempf filtration of the $\rk 2$
tensor $(E,\varphi)$ is independent of $m$, for $m\geq m'$.
\end{thm}

\subsection{The $m$-Kempf filtration stabilizes with $m$}

Let us define the graph associated to the $m$-Kempf filtration of $(E,\varphi)$.

\begin{dfn}
\label{graphrk2} Let $m\geq m_{0}$. Given $0\subset V_{1}\subset \cdots \subset V_{t+1}= V$ a filtration of vector spaces
of $V$, let
$$v_{m,i}=m^{n+1}\cdot \frac{1}{\dim V^{i}\dim V}\big[r^{i}\dim V-r\dim V^{i}+\dfrac{a_{2}}{a_{1}}(\epsilon^{i}(\overline\Phi)\dim V-s\dim V^{i})\big]\; ,$$
$$b_{m}^{i}=\dfrac{1}{m^{n}}\dim V^{i}>0\; ,$$
$$w_{m}^{i}=-b_{m}^{i}\cdot v_{m,i}=m\cdot \frac{1}{\dim V}\big[r\dim V^{i}-r^{i}\dim V+\dfrac{a_{2}}{a_{1}}(s\dim V^{i}-\epsilon^{i}(\overline\Phi)\dim V) \big]\; .$$
Also let
$$b_{m,i}=b_{m}^{1}+\ldots +b_{m}^{i}=\dfrac{1}{m}\dim V_{i}\; ,$$
$$w_{m,i}=w_{m}^{1}+\ldots +w_{m}^{i}=m\cdot \frac{1}{\dim V}\big[r\dim V_{i}-r_{i}\dim V+\dfrac{a_{2}}{a_{1}}(s\dim V_{i}-\epsilon_{i}(\overline\Phi)\dim V)\big]\; .$$
We call the graph defined by points $(b_{m,i},w_{m,i})$ the \emph{graph associated to the filtration}
$V_{\bullet}\subset V$.
\end{dfn}

Now we prove a crucial Lemma which will let us prove Theorem
\ref{kempfstabilizesrk2} using the same method than in previous
sections.

\begin{lem}
\label{independenceofweightsrk2} The symbols
$\epsilon_i(\overline\Phi)=\epsilon_i(\overline\Phi,V_{\bullet},n_{\bullet})$ do not depend on the weights
$n_{\bullet}$. Therefore, the graph associated to the filtration only depends on the data $V_{\bullet}\subset
V$, not the weights $n_{\bullet}$.
\end{lem}
\begin{pr}
Note that $\rk E_{1}\geq 1$ because it is generated by, at least, a
non zero global section. Suppose that $\rk E^{m}_{1}=\rk
E^{m}_{2}=\ldots=\rk E^{m}_{k}=1$ and $\rk E^{m}_{k+1}=\ldots=\rk
E^{m}_{t}=\rk E=2$. Then, for example, $E^{m}_{1}$ coincide with
$E^{m}_{2}$ on an open set and, generically, the behavior with
respect to $\varphi$ is the same, i.e.
$$\overline\Phi|_{V_{1}\otimes\cdots\otimes
V_{1}}=0\Leftrightarrow \varphi|_{E^{m}_{1}\otimes\cdots\otimes
E^{m}_{1}}=0\Leftrightarrow \varphi|_{E^{m}_{2}\otimes
E^{m}_{1}\cdots \otimes E^{m}_{1}}=0\; .$$ Therefore, the values
$\epsilon_i(\overline\Phi,V_{\bullet},n_{\bullet})$ only depend on
the filters $E^{m}_{i}$ but not on the specific values of the
$\Gamma_{i}$. In fact, they will only depend on $\Gamma_{1}$ and
$\Gamma_{k+1}$, because they are the minimal ones among the
filters of the same rank (c.f. (\ref{rightmu}) and
(\ref{rightmuV})). In this case we will just write
$\epsilon_i(\overline\Phi,V_{\bullet})$, or
$\epsilon_i(\overline\Phi)$, when the filtration is clear from the
context.
\end{pr}

Next, we can identify the Kempf function in Theorem \ref{kempftheoremsheaves}
$$\mu(V_{\bullet},n_{\bullet})=\frac{\sum_{i=1}^{t}  n_{i} (r\dim V_{i}-r_{i}\dim V)+\frac{a_{2}}{a_{1}}
\sum_{i=1}^{t} n_{i} \big(s\dim
V_{i}-\epsilon_{i}(\overline\Phi)\dim V
\big)}{\sqrt{\sum_{i=1}^{t+1} {\dim V^{i}} \Gamma_{i}^{2}}}=$$
$$=\frac {\sum_{i=1}^{t+1} \frac{\Gamma_i}{\dim V} ( r^i \dim V - r\dim V^i)+\frac{a_2}{a_1}
 \sum_{i=1}^{t+1} \frac{\Gamma_{i}}{\dim V} \big(\epsilon^i(\overline\Phi)\dim V -s\dim V^i\big)}
 {\sqrt{\sum_{i=1}^{t+1} {\dim V^{i}} \Gamma_{i}^{2}}}\; ,$$
where $n_{i}=\frac{\Gamma_{i}-\Gamma_{i-1}}{\dim V}$, with the function in Theorem \ref{maxconvexenvelope} (c.f.
Proposition \ref{identification}). Precisely, we use Lemma \ref{independenceofweightsrk2} to assure that the
data of the filters $V_{\bullet}\subset V$, and the data of the weights $n_{\bullet}$ are independent, so we can maximize the Kempf function with
respect to each of them, independently, as in Theorem \ref{maxconvexenvelope}.

\begin{prop}
\label{identificationrk2} For every integer $m$, the following equality holds
$$\mu(V_{\bullet},n_{\bullet})=m^{(-\frac{n}{2}-1)}\cdot \mu_{v_{m}}(\Gamma)$$
between the Kempf function on Theorem \ref{kempftheoremsheaves} and the function in Theorem
\ref{maxconvexenvelope}.
\end{prop}
\begin{pr}
By Lemma \ref{independenceofweightsrk2}, we can fix a vector $v_{m}$ and look for the maximum of the function
$\mu_{v_{m}}$ among the corresponding convex cone.
\end{pr}

In the following, we will omit the subindex $m$ for the numbers $v_{m,i}$, $b_{m,i}$, $w_{m,i}$ in the
definition of the graph associated to the filtration of vector spaces, where it is clear from the context.
Recall Remark \ref{growth} to understand the meaning of the factors in $m$ in Definition \ref{graphrk2}.

Now we use Lemmas \ref{lemmaA} and \ref{lemmaB} to give the
analogous to Propositions \ref{boundedness} and \ref{task} in this
case.

Let us define
\begin{equation}
\label{C_constant_rk2} C=\max\{r|\mu _{\max }(E)|+\frac{d}{r}+r|B|+|A|+s\delta_{n-1}(n-1)!+1\;,\;1\},
\end{equation}
a positive constant, where $\delta_{n-1}$ is the leading coefficient of the polynomial $\delta(m)$, of degree $\leq n-1$ (if $\deg \delta< n-1$, set $\delta_{n-1}=0$).

\begin{prop}
\label{boundednessrk2} Given a sufficiently large $m$, each filter in the $m$-Kempf filtration of the $\rk 2$ tensor $(E,\varphi)$ has slope
$\mu(E^{m}_{i})\geq \dfrac{d}{r}-C$.
\end{prop}

\begin{pr}
Choose an $m_{1}$ such that for $m\geq m_{1}$
$$[\mu_{max}(E)+gm+B]_{+}=\mu_{max}(E)+gm+B$$
and
$$[\frac{d}{r}-C+gm+B]_{+}=\frac{d}{r}-C+gm+B\; .$$
Let $m_{2}$ be such that $P_{E}(m)-s\delta(m)>0$ for $m\geq
m_{2}$. Now consider $m\geq \max\{m_{0},m_{1},m_{2}\}$ and let
$$0\subseteq (E^{m}_{1},\varphi|_{E_{1}^{m}})\subseteq
(E^{m}_{2},\varphi|_{E_{2}^{m}})\subseteq\cdots
(E^{m}_{t},\varphi|_{E_{t}^{m}})\subseteq
(E^{m}_{t+1},\varphi|_{E_{t+1}^{m}})\subseteq (E,\varphi)$$
be the $m$-Kempf filtration.

Suppose that we have a filter $E^{m}_{i}\subseteq E$, of rank $r_{i}$ and degree $d_{i}$, such that $\mu
(E^{m}_{i})<\frac{d}{r}-C$. Again, the subsheaf $E_{i}^{m}(m)\subset E(m)$ satisfies the estimate in Lemma \ref{Simpson},
$$h^{0}(E^{m}_{i}(m))\leq \frac{1}{g^{n-1}n!}\big ((r_{i}-1)([\mu_{max}(E^{m}_{i})+gm+B]_{+})^{n}+([\mu_{min}(E^{m}_{i})+gm+B]_{+})^{n}\big)\; .$$
Hence, analogously to Proposition \ref{boundedness},
$$h^{0}(E^{m}_{i}(m))\leq \frac{1}{g^{n-1}n!}\big ((r_{i}-1)(\mu_{max}(E)+gm+B)^{n}+(\frac{d}{r}-C+gm+B)^{n}\big)=G(m)\; ,$$
where
$$G(m)=\frac{1}{g^{n-1}n!}\big
[r_{i}g^{n}m^{n}+ng^{n-1}\big((r_{i}-1)\mu_{max}(E)+\frac{d}{r}-C+r_{i}B\big)m^{n-1}+\cdots
\big ]\; .$$

By Definition \ref{graphrk2}, to the $m$-Kempf filtration we associate a graph with heights, for each $j$
$$w_{j}=w^{1}+\ldots +w^{j}=m\cdot \frac{1}{\dim V}\big[r\dim V_{j}-r_{j}\dim V+\dfrac{a_{2}}{a_{1}}(s\dim V_{j}-\epsilon_{j}
(\overline\Phi)\dim V )\big]\; .$$
We will get a contradiction by showing that $w_{i}<0$ (c.f. Proposition \ref{boundedness}).

Since $E^{m}_{i}(m)$ is generated by $V_{i}$ under the evaluation map, it is $\dim V_{i}\leq
H^{0}(E^{m}_{i}(m))$, hence
$$w_{i}=\frac{m}{\dim V}\big[r\dim V_{i}-r_{i}\dim V+\dfrac{a_{2}}{a_{1}}(s\dim V_{i}-\epsilon_{i}(\overline\Phi)\dim V )\big]\leq$$
$$\frac{m}{P_{E}(m)}\big[rh^{0}(E^{m}_{i}(m))-r_{i}P_{E}(m)+\dfrac{r\delta(m)}{P_{E}(m)-s\delta(m)}(sh^{0}(E^{m}_{i}(m))-
\epsilon_{i}(\overline\Phi)P_{E}(m))\big]\leq$$
$$\frac{m}{P_{E}(m)}\big[rG(m)-r_{i}P_{E}(m)+\dfrac{r\delta(m)}{P_{E}(m)-s\delta(m)}
(sG(m)-\epsilon_{i}(\overline\Phi)P_{E}(m))\big]=$$
$$m\cdot \frac{\big[(P_{E}(m)-s\delta(m))(rG(m)-r_{i}P_{E}(m))+
(r\delta(m))(sG(m)-\epsilon_{i}(\overline\Phi)P_{E}(m))\big]}{P_{E}(m)(P_{E}(m)-s\delta(m))}\;
.$$

Then, $w_{i}<0$ is equivalent to
$$\Psi(m)=(P_{E}(m)-s\delta(m))(rG(m)-r_{i}P_{E}(m))+(r\delta(m))(sG(m)-\epsilon_{i}(\overline\Phi)P_{E}(m))<0\; ,$$
and $\Psi(m)=\xi_{2n}m^{2n}+\xi_{2n-1}m^{2n-1}+\cdots
+\xi_{1}m+\xi_{0}$ is a $(2n)^{th}$-order polynomial, whose higher order coefficient is
$$\xi_{2n}=(P_{E}(m)-s\delta(m))_{n}(rG(m)-r_{i}P_{E}(m))_{n}+(r\delta(m))_{n}(sG(m)-
\epsilon_{i}(\overline\Phi)P_{E}(m))_{n}=$$
$$(P_{E}(m)-s\delta(m))_{n}(r\frac{r_{i}g}{n!}-r_{i}\frac{rg}{n!})+0=0\; .$$
The $(2n-1)^{th}$-order coefficient is
$$\xi_{2n-1}=(P_{E}(m)-s\delta(m))_{n}(rG(m)-r_{i}P_{E}(m))_{n-1}+(r\delta(m))_{n-1}(sG(m)-
\epsilon_{i}(\overline\Phi)P_{E}(m))_{n}=$$
$$\frac{rg}{n!}(rG_{n-1}-r_{i}\frac{A}{(n-1)!})+r\delta_{n-1}(s\frac{r_{i}g}{n!}-
\epsilon_{i}(\overline\Phi)\frac{rg}{n!})$$
where $G_{n-1}$ is the $(n-1)^{th}$-coefficient of the polynomial $G(m)$,
$$G_{n-1}=\frac{1}{g^{n-1}n!}ng^{n-1}((r_{i}-1)\mu_{max}(E)+\frac{d}{r}-C+r_{i}B)=$$
$$\frac{1}{(n-1)!}((r_{i}-1)\mu_{max}(E)+\frac{d}{r}-C+r_{i}B)\leq$$
$$\frac{1}{(n-1)!}((r_{i}-1)|\mu_{max}(E)|+\frac{d}{r}-C+r_{i}|B|)\leq$$
$$\frac{1}{(n-1)!}(r|\mu_{max}(E)|+\frac{d}{r}-C+r|B|)<\frac{-|A|}{(n-1)!}-s\delta_{n-1}\; ,$$
last inequality coming from the definition of $C$ in (\ref{C_constant_rk2}). Then
$$\xi_{2n-1}<\frac{rg}{n!}\big(r(\frac{-|A|}{(n-1)!}-s\delta_{n-1})-r_{i}\frac{A}{(n-1)!}\big )+r\delta_{n-1}\big(\frac{r_{i}g}{n!}-
\epsilon_{i}(\overline\Phi)\frac{rg}{n!}\big)=$$
$$\frac{rg}{n!}\big [ \big (\frac{-r|A|-r_{i}A}{(n-1)!}\big )
-rs\delta_{n-1}+\delta_{n-1}(r_{i}-\epsilon_{i}(\overline\Phi)r)\big
]=$$
$$\frac{rg}{n!}\big [ \big (\frac{-r|A|-r_{i}A}{(n-1)!}\big )
+\delta_{n-1}(-rs +r_{i}s-\epsilon_{i}(\overline\Phi)r)\big ]<$$
$$\frac{rg}{n!}\delta_{n-1}(-rs +r_{i}s-\epsilon_{i}(\overline\Phi)r)\; ,$$
because $-r|A|-r_{i}A<0$. Last expression is either zero if
$r_{i}=\rk E=2$ (because in that case it is $\epsilon_{i}(\overline\Phi)=\epsilon_{t+1}(\overline\Phi)=s$), or
negative if $r_{i}=1$. Hence, $\xi_{2n-1}<0$.

Therefore $\Psi(m)=\xi_{2n-1}m^{2n-1}+\cdots +\xi_{1}m+\xi_{0}$
with $\xi_{2n-1}<0$, so there exists an integer $m_{3}$ such that for $m\geq
\{m_{0},m_{1},m_{2},m_{3}\}$ we have $\Psi(m)<0$ and
$w_{i}<0$, then the contradiction.
\end{pr}

Similarly to Proposition \ref{regular}, we prove

\begin{prop}
\label{regularrk2} There exists an integer $m_{4}$ such that for
$m\geq m_{4}$ the sheaves $E^{m}_{i}$ and
$E^{m,i}=E^m_i/E^m_{i-1}$ are $m_{4}$-regular. In particular their
higher cohomology groups, after twisting with
$\mathcal{O}_{X}(m_{4})$, vanish and they are generated by global
sections.
\end{prop}

\begin{prop}
\label{taskrk2}
Let $m\geq m_{4}$. For each filter $E_{i}^{m}$ in the $m$-Kempf filtration of the $\rk 2$ tensor $(E,\varphi)$, we
have $\dim V_{i}=h^{0}(E_{i}^{m}(m))$, therefore $V_{i}\cong H^{0}(E_{i}^{m}(m))$.
\end{prop}

\begin{pr}
Let $V_{\bullet}\subseteq V$ be the Kempf filtration of $V$ (cf. Theorem \ref{kempftheoremsheaves}) and let
$(E_{\bullet}^{m},\varphi|_{E_{\bullet}^{m}})\subseteq (E,\varphi)$ be the $m$-Kempf filtration of $(E,\varphi)$.
Analogously to Proposition \ref{task} we can construct two filtrations

\begin{equation}
 \label{filtrationVrk2}
    \begin{array}{ccccccccccccc}
    0 & \subset & \cdots & \subset & V_{i} & \subset & V_{i+1} & \subset & V_{i+2} & \subset & \cdots & \subset & V \\
    & & & & \cap & & || & & ||& & & &   \\
    & & & & H^{0}(E_{i}^{m}(m)) & \subset & H^{0}(E_{i+1}^{m}(m)) & \subset & H^{0}(E_{i+2}^{m}(m)) &  &  & &
    \end{array}
\end{equation}

and

\begin{equation}
\label{filtrationV'rk2}
    \begin{array}{ccccccccccccc}
    0 & \subset & \cdots & \subset & V_{i} & \subset & H^{0}(E_{i}^{m}(m)) & \subset & V_{i+1} & \subset & \cdots & \subset &  V\\
    & & & & || & & || & & || & & & &   \\
    & & & & V'_{i} & &  V'_{i+1} & & V'_{i+2} & & & &
    \end{array}
\end{equation}
to be in situation of Lemma
\ref{lemmaB}, where $W=H^{0}(E_{i}^{m}(m))$, filtration $V_{\bullet}$ is \eqref{filtrationVrk2} and filtration
$V'_{\bullet}$ is \eqref{filtrationV'rk2}.

Now, the graph associated to filtration $V_{\bullet}$ is given, by Definition \ref{graphrk2}, by the points
$$(b_{i},w_{i})=(\dfrac{\dim V_{i}}{m^{n}},\frac{m}{\dim V}\big(r\dim V_{i}-r_{i}\dim V+\dfrac{a_{2}}{a_{1}}(s\dim V_{i}-\epsilon_{i}
(\overline\Phi,V_{\bullet})\dim V))\big)\; ,$$
the slopes $-v_{i}$ of the graph given by

$$-v_{i}=\frac{w^{i}}{b^{i}}=\frac{w_{i}-w_{i-1}}{b_{i}-b_{i-1}}=$$
$$\frac{m^{n+1}}{\dim V}\big(r-r^{i}\frac{\dim V}{\dim V^{i}}+\frac{a_{2}}{a_{1}}(s-\epsilon^{i}(\overline\Phi,V_{\bullet})\frac{\dim V}
{\dim V^{i}})\big)\leq$$
$$\frac{m^{n+1}}{\dim V}\big(r+s\frac{a_{2}}{a_{1}}\big):=R$$
and equality holds if and only if $r^{i}=0$ (note that $r^{i}=0$ implies $\epsilon^{i}(\overline\Phi,V_{\bullet})=0$).

The new point which appears in graph of the filtration $V'_{\bullet}$ is
$$Q=\big(\dfrac{h^{0}(E_{i}^{m}(m))}{m^{n}},\frac{m}{\dim V}(rh^{0}(E_{i}^{m}(m))-r_{i}\dim V+\dfrac{a_{2}}{a_{1}}
(sh^{0}(E_{i}^{m}(m))- \epsilon_{i}(\overline\Phi,V_{\bullet})\dim V))\big)\; ,$$ where we write
$\epsilon_{i}(\overline\Phi,V_{\bullet})$ instead of $\epsilon_{i}(\overline\Phi,V'_{\bullet})$, by the same
argument used in proof of Proposition \ref{taskpairs} (c.f. \eqref{epsilon}).

The slope of the segment between $(b_{i},w_{i})$ and $Q$ is, similarly,
$$-v'_{i}=\dfrac{m^{n+1}}{\dim V}(r+s\dfrac{a_{2}}{a_{1}})=R\; .$$

By Lemma \ref{lemmaA}, the graph is convex, so $v_{1}<v_{2}<\ldots<v_{t+1}$. Besides, $r^{1}=r_{1}>0$, then
$-R<v_{1}$, because $E$ is torsion free, hence also the subsheaf $E_{1}^{m}$, and a rank $0$ torsion free sheaf
is the zero sheaf. On the other hand, by Lemma \ref{lemmaB}, $v'_{i}\geq v_{i}$. Hence,
$$-R<v_{1}<v_{2}<\ldots<v_{i}\leq v'_{i}=-R\; ,$$
which is a contradiction.

Therefore, $\dim V_{i}=h^{0}(E_{i}^{m}(m))$, for every filter in the $m$-Kempf filtration.
\end{pr}

\begin{cor}
\label{rankrk2} Let $m\geq m_{4}$. For every filter $E_{i}^{m}$ in
the $m$-Kempf filtration of the $\rk 2$ tensor $(E,\varphi)$, it
is $r^{i}>0$. Therefore, the $m$-Kempf filtration consists on a
rank $1$ subsheaf, $0\subset (L^{m},\varphi|_{L^{m}})\subset (E,\varphi)$.
\end{cor}
\begin{pr}
We have seen that $r^{i}=0$ is equivalent to $-v_{i}=R$. Then the result follows from Proposition \ref{taskrk2}
because it is $r^{1}=r_{1}>0$ and $-R<v_{1}<v_{2}<\ldots<v_{t+1}$.
\end{pr}

For any $m\geq m_{4}$, by Corollary \ref{rankrk2} there is only
one filter $(L^{m},\varphi|_{L^{m}})$ in the $m$-Kempf filtration
and, by Proposition \ref{regularrk2}, $L^{m}$ is $m_{4}$-regular.
Hence, $L^{m}(m_{4})$ is generated by the subspace
$H^{0}(L^{m}(m_{4}))\subset H^{0}(E(m_{4}))$ by the evaluation map
(c.f. Lemma \ref{mregularity}). Note that the dimension of the
vector space $H^{0}(E(m_{4}))$ does not depend on $m$.

The \emph{$m$-type} of the $m$-Kempf filtration $0\subset (L^{m},\varphi|_{L^{m}})\subset (E,\varphi)$
is the Hilbert polynomial $P_{L^{m}}$ (c.f. Definition \ref{mtype}). The set of possible $m$-types
$$\mathcal{P}=\big\{P_{L^{m}}\big\}$$
is finite, for all integers $m\geq m_{4}$ (c.f. Proposition
\ref{Pisfinite}).

Rewrite the graph associated to the $m$-Kempf filtration (c.f. Definition \ref{graphrk2})
$$v_{m,i}=\frac{m^{n+1}}{\dim V^{i}\dim V}\big[r^{i}\dim V-r\dim V^{i}+\dfrac{a_{2}}{a_{1}}(\epsilon^{i}
(\overline\Phi)\dim V-s\dim V^{i})\big]\; ,$$
$$b_{m}^{i}=\frac{1}{m^{n}}\cdot \dim V^{i}\; ,$$
as
$$v_{m,i}=\frac{m^{n+1}}{P_{m}^{i}(m)P(m)}\big[r^{i}P(m)-rP_{m}^{i}(m)+\dfrac{r\delta(m)}{P(m)-s\delta(m)}(\epsilon^{i}
(\overline\Phi)P(m)-sP_{m}^{i}(m))\big]\; ,$$
$$b_{m}^{i}=\frac{1}{m^{n}}\cdot P_{m}^{i}(m)\; ,$$
by Propositions \ref{regularrk2} and \ref{taskrk2}.

Note that, by Corollary \ref{rankrk2}, the graph has only two
slopes given by
$$v_{m,1}=\frac{m^{n+1}}{P_{L^{m}}(m)P(m)}\big[P(m)-2P_{L^{m}}(m)+\dfrac{2\delta(m)}{P(m)-s\delta(m)}(\epsilon_{L^{m}}
P(m)-sP_{L^{m}}(m))\big]\; ,$$
$$v_{m,2}=\frac{m^{n+1}}{P_{E/L^{m}}(m)P(m)}\big[P(m)-2P_{E/L^{m}}(m)+\dfrac{2\delta(m)}{P(m)-s\delta(m)}((s-\epsilon_{L^{m}})
P(m)-sP_{E/L^{m}}(m))\big]\; ,$$ where $\epsilon(L^{m})$ is the number of times that the subsheaf $L^{m}$
appears on the minimal multi-index (c.f. (\ref{rightmuV}) in section \ref{exampletensors}).

The set
$$
\mathcal{A}=\{\Theta_{m}:m\geq m_{4}\}
$$
is finite (c.f. Proposition \ref{Pisfinite}), where
$$
\Theta_{m}(l)=(\mu_{v_{m}(l)}(\Gamma_{v_{m}(l)}))^{2}=||v_{m}(l)||^{2}
$$
(c.f. (\ref{maxvalue})). Let $K$ be the maximal function in $\mathcal{A}$ as in Lemma \ref{uniquefunction}) for which
$\exists\; m_{5}$ such that for all $m\geq m_{5}$ it is $\Theta_{m}=K$.

\begin{prop}
\label{eventuallyrk2} Let $l_{1}$ and $l_{2}$ be integers with
$l_{1}\geq l_{2} \geq m_{5}$. Then the $l_{1}$-Kempf filtration of
$E$ is equal to the $l_{2}$-Kempf filtration of $E$.
\end{prop}

\begin{pr}
C.f. Proposition \ref{eventually}.
\end{pr}

Therefore, Theorem \ref{kempfstabilizesrk2} follows from Proposition \ref{eventuallyrk2}. Hence, eventually, the
Kempf filtration of the $\rk 2$ tensor $(E,\varphi)$ does not depend on the integer $m$.

\begin{dfn}
If $m\geq m_{5}$, the $m$-Kempf filtration of the $\rk 2$ tensor $(E,\varphi)$
$$0\subset (L,\varphi|_{L}) \subset (E,\varphi)$$ is called the \textbf{Kempf filtration} or the \textbf{Kempf subsheaf} of $(E,\varphi)$.
\end{dfn}

\subsection{Harder-Narasimhan filtration for $\rk 2$ tensors}

Kempf theorem (c.f. Theorem \ref{kempftheoremsheaves}) says that, given an integer $m$ and $V\simeq
H^{0}(E(m))$, there exists a unique weighted filtration of vector spaces $V_{\bullet}\subseteq V$ which gives
maximum for the Kempf function
$$\mu(V_{\bullet},n_{\bullet})=\frac {\sum_{i=1}^{t+1} \frac{\Gamma_i}{\dim V} ( r^i \dim V - r\dim V^i)+\frac{a_2}{a_1}
 \sum_{i=1}^{t+1} \frac{\Gamma_{i}}{\dim V} \big(\epsilon^i(\overline\Phi)\dim V -s\dim V^i\big)}
 {\sqrt{\sum_{i=1}^{t+1} {\dim V^{i}} \Gamma_{i}^{2}}}\; .$$ This filtration induces a unique rank $1$ subsheaf $L\subset E$ called
the Kempf subsheaf of the $\rk 2$ tensor $(E,\varphi)$. By
Proposition \ref{eventuallyrk2}, the subsheaf $L$ does not depend
on $m$, for $m\geq m_{5}$.

The Kempf function is a function on $m$ (c.f. Proposition \ref{identificationrk2}). Consider the
function
$$K(m)=m^{\frac{n}{2}+1}\cdot \mu(V_{\bullet},m_{\bullet})=\mu_{v_{m}}(\Gamma)$$
and, making the substitutions for $m$ sufficiently large
$$\dim V_{1}=\dim V^{1}=h^{0}(L(m))=P_{L}(m)\; ,$$
$$\dim V^{2}=\dim V-\dim V_{1}=h^{0}(E/L(m))=P_{E/L}(m)$$
we get
$$K(m)=m^{\frac{n}{2}+1}\cdot \frac {\sum_{i=1}^{2} \frac{\gamma_{i}}{r}[(r^i P - rP^i)+\frac{r\delta}{P-s\delta}
 (\epsilon^i P -sP^i)]}
 {\sqrt{\sum_{i=1}^{2} P^{i}\frac{P^{2}}{r^{2}} \gamma_{i}^{2}}}\; ,$$
where we put $P=P_{E}(m)$, $P^{1}=P_{L}(m)$, $P^{2}=P_{E/L}(m)$, $\epsilon^{1}=\epsilon(L)$, $\epsilon^{2}=s-\epsilon(L)$.
Note that
$\epsilon^{i}=\epsilon^i(\overline\Phi)=\epsilon^{i}(\varphi)$ and
recall the relations
$$\gamma_{i}=\frac{r}{P}\Gamma_{i}\; ,$$
$$\frac{a_{2}}{a_{1}}=\frac{r\delta}{P-s\delta}\; .$$
Also recall
$$\frac{\gamma_{i+1}-\gamma_{i}}{r}=n_{i}\; ,$$
$$\sum r^{i}\gamma_{i}=\gamma_{1}+\gamma_{2}=0\; ,$$ which gives in our case $\gamma_{1}=-n_{1}$, $\gamma_{2}=n_{1}$. Substituting we get
$$K(m)=m^{\frac{n}{2}+1}\cdot \frac{1}{P-s\delta}\frac {-n_{1}[2(\delta \epsilon^{1}-P^{1})+(P-\delta s)]+n_{1}[2(\delta\epsilon^{2}-P^{2})+(P-\delta s)]}
 {\sqrt{P^{1}n_{1}^{2}+P^{2}n_{1}^{2}}}=$$
$$m^{\frac{n}{2}+1}\cdot \frac{r}{\sqrt{P}(P-s\delta)}[2P_{L}-P_{E}+\delta(s-2\epsilon(L))]\; .$$
Note that the unique weight $n_{1}$ does not appear in the function later from the substitutions, as it was
expected from a one-step filtration. Also note that the denominator of the function $K$ is positive (c.f.
choice of $m_{2}$ in proof of Proposition \ref{boundednessrk2}). Hence, we can state the following theorem.

\begin{thm}
\label{finalfunction} Given a $\delta$-unstable $\rk 2$ tensor $
(E,\varphi:\overbrace{E\otimes\cdots \otimes E}^{\text{s
times}}\too \mathcal{O}_{X}) $, there exists a unique line
subsheaf $L \subset E$ which gives maximum for the polynomial
function
$$K(m)=2P_{L}(m)-P_{E}(m)+\delta(m)(s-2\epsilon(L))\; .$$
\end{thm}

If $X$ is a one dimensional complex projective variety, i.e. a
smooth projective complex curve, we can simplify the function $K$.
Recall that, by Riemann-Roch, the Hilbert polynomial of a sheaf
$E$ of rank $r$ and degree $d$ over a curve of genus $g$ is
$$P_{E}(m)=rm+d+r(1-g)\; ,$$
and the polynomial $\delta(m)$ becomes a positive constant $\tau$.
In this case, a coherent torsion free sheaf of rank $2$ is a
vector bundle of rank $2$ over $X$, and the Kempf subsheaf will be
a line subbundle.

\begin{thm}
\label{finalfunctioncurves} Given a $\tau$-unstable $\rk 2$ tensor
$ (E,\varphi:\overbrace{E\otimes\cdots \otimes E}^{\text{s
times}}\too \mathcal{O}_{X}) $ over a smooth projective complex
curve, there exists a unique line subbundle $L \subset E$ which
maximizes the quantity
$$2\deg{L}-\deg{E}+\tau(s-2\epsilon(L))\; .$$
\end{thm}

Note that, if the tensor is unstable, such quantity will be positive, and the graph corresponding to the filtration will be \textit{a cusp} which is a convex graph.

If we define the \emph{corrected Hilbert polynomials} of $(E,\varphi)$ and $(L,\varphi|_{L})$ (c.f. Definition \ref{corrected}) as
$$\overline{P}_{E}=P_{E}-\delta s\; ,$$
$$\overline{P}_{L}=P_{L}-\delta\epsilon(L)\; ,$$
we recover the notion of stability for $\rk 2$ tensors (c.f.
Definition \ref{stabilityrk2}).  A $\rk 2$ tensor $(E,\varphi)$ is
\emph{$\delta$-unstable} if there exists a line subsheaf $L\subset E$
such that
$$\frac{\overline{P}_{L}}{\rk L}>\frac{\overline{P}_{E}}{\rk E}\Leftrightarrow \overline{P}_{L}>
\frac{\overline{P}_{E}}{2}\; .$$

Hence, this procedure allows us define a notion of a Harder-Narasimhan filtration for $\delta$-unstable $\rk 2$ tensors.

\begin{dfn}
\label{HNrk2tensors} If $(E,\varphi)$ is a $\delta$-unstable $\rk 2$ tensor, there exists a unique line subsheaf
maximizing
$$2\cdot \overline{P}_{L}-\overline{P}_{E}>0\; .$$
We call
$$0\subset (L,\varphi|_{L})\subset (E,\varphi)$$ the \emph{Harder-Narasimhan filtration} of $(E,\varphi)$, and we
call $L$ the \emph{Harder-Narasimhan subsheaf} of $(E,\varphi)$.
\end{dfn}

\begin{rem}
We do not know, in principle, how to define a quotient tensor $(E/L,\overline{\varphi}|_{E/L})$, because we do not know, a priori, how to define
$\overline{\varphi}|_{E/L}$. This is why we cannot talk about quotient tensors, as in Definition \ref{subpairs}.

Given the exact sequence of sheaves, $0\rightarrow L\rightarrow E\rightarrow E/L\rightarrow 0$, we define the
corrected Hilbert polynomial of the quotient as $\overline{P}_{E/L}=\overline{P}_{E}-\overline{P}_{L}$, and we have, trivially, the additivity of the corrected polynomials on exact sequences of sheaves. This way we can consider that
Definition \ref{HNrk2tensors} contains the analogous to conditions of Definition \ref{HNdef} for $\rk 2$ tensors. Indeed,
$$2\cdot \overline{P}_{L}-\overline{P}_{E}>0\Leftrightarrow \overline{P}_{L}>\overline{P}_{E/L}\; ,$$
and the semistability of $(L,\varphi|_{L})$ and $(E/L,\overline{\varphi}|_{E/L})$ (whichever definition of $\overline{\varphi}|_{E/L}$ we impose), would follow trivially from the
fact of they are rank $1$ tensors.

Therefore, Definition \ref{HNrk2tensors} gives a notion of Harder-Narasimhan filtration for these objects.
\end{rem}

\subsection{Stable coverings of a projective curve}
In this section we use the previous notions for $\rk 2$
tensors over curves where the morphism is symmetric, and the Definition \ref{HNrk2tensors} of the
Harder-Narasimhan subsheaf, to define stable coverings of a
projective curve and, for the unstable ones, a maximally
destabilizing object, in geometrical terms.

In the following, we shall consider rank $2$ tensors $(E,\varphi)$ where $E$ is a $\rk 2$ vector bundle over a smooth
projective complex curve $X$, and $$\varphi:\overbrace{E\otimes\cdots \otimes E}^{\text{s times}}\too
\mathcal{O}_{X}$$ is a symmetric non degenerate morphism. We call it a \emph{symmetric non degenerate rank $2$ tensor}. Let $\tau$ be a positive real number. Let
$\mathbb{P}(E)$ be the projective space bundle of the vector bundle $E$, which is a ruled algebraic surface
(c.f. \cite[Section V.2]{Ha}).

The morphism $\varphi$ is, fiberwise, a symmetric multilinear map
$$\varphi_{x}:\overbrace{V\otimes\cdots \otimes V}^{\text{s times}}\too
\mathbb{C}\; ,$$ where $V\simeq \mathbb{C}^{2}$. Then, $\varphi_{x}$
factors through $\Sym^{s}(V)$, isomorphic to the
$(s+1)$-dimensional vector space of homogeneous polynomials of
degree $s$ in two variables. Hence, fiberwise, $\varphi$ can be
represented by a polynomial
\begin{equation}
\label{tensorfiber} \varphi_{x}\equiv
\sum_{i=0}^{s}a_{i}(x)X_{0}^{i}X_{1}^{s-i}
\end{equation}  which
vanishes on $s$ points in
$\mathbb{P}(V)\simeq\mathbb{P}_{\mathbb{C}}^{1}$. Therefore, as
$\varphi$ varies on $X$, it defines a degree $s$ covering
$$\mathbb{P}(E)\supset X'\rightarrow
X\; .$$

Suppose that $(E,\varphi)$ is a $\tau$-unstable $\rk 2$ tensor. Then, by
Theorem \ref{finalfunctioncurves}, there exists a line subbundle
$L\subset E$, the \emph{Harder-Narasimhan subbundle}, giving maximum for the quantity
\begin{equation}
\label{rk2stabexpression}
2\deg(L)-\deg (E)+\tau(s-2\epsilon(L))\; .
\end{equation}

The subbundle $L$ can be seen as a section of $\mathbb{P}(E)$,
each fiber $L_{x}$ corresponding to a point $P=\{L_{x}\}\in
\mathbb{P}_{\mathbb{C}}^{1}$. Recall from Definition
\ref{stabilitytensorsrk2} that $\epsilon (L)=k$ if
$\varphi|_{L^{\otimes (k+1)}\otimes E^{\otimes (s-k-1)}}= 0$ and
$\varphi|_{L^{\otimes k}\otimes E^{\otimes (s-k)}}\neq 0$. Note
that here we use the symmetry of the morphism $\varphi$.
Therefore, $\epsilon(L)=k$ means that, generically, $P=\{L_{x}\}$
is a zero of multiplicity $s-k$ and, by definition of the covering
$X'\rightarrow X$, $s-\epsilon(L)$ is exactly the number of
branches of $X'$ which generically do coincide with the section
defined by $L$, counted with multiplicity.

Recall Examples \ref{exG} and \ref{exG2} in Section \ref{modulis}.
There, a homogeneous polynomial of degree $N$, $P=\underset{i}\sum
a_{i}X_{0}^{i}X_{1}^{N-i}$, was unstable if it contained a linear
factor of degree greater that $\frac{N}{2}$. Now, observe that the
restriction of a rank $2$ tensor to a point $x\in X$ in
(\ref{tensorfiber}), passing to the projectivization
$\mathbb{P}(E)$ hence fibers are isomorphic to
$\mathbb{P}_{\mathbb{C}}^{1}$, is precisely one of the homogeneous
polynomials in Examples \ref{exG} and \ref{exG2}. Fiberwise, the
morphism $\varphi$ defines a set of $s$ points in
$\mathbb{P}_{\mathbb{C}}^{1}$. See that, from the point of view of
Examples \ref{exG} and \ref{exG2}, letting $s=N$, the set of
points is unstable if there exists a point with multiplicity
greater that $\frac{s}{2}$.

Then, as $s-\epsilon(L)$ is the multiplicity of the point defined
by the line $L_{x}$ (the fiber of the Harder-Narasimhan subbundle
over $x$), in the set of $s$ points defined by the morphism
$\varphi$, following the previous argument, this point $\{L_{x}\}$
will destabilize the set if
$$s-\epsilon(L)>\frac{s}{2}\Leftrightarrow s-2\epsilon (L)>0\; ,$$
which is the second summand in (\ref{rk2stabexpression}). Hence,
the positivity of $s-2\epsilon(L)$ is equivalent to find a line
subbundle $L$ defining a point in the fiber
$\mathbb{P}_{\mathbb{C}}^{1}$, which coincides with one of the
zeroes of $\varphi$ in the fiber, and such that it has
multiplicity greater that $\frac{s}{2}$.

To conclude, we can say that the expression
(\ref{rk2stabexpression}) consists of two summands weighted by the
parameter $\tau$. First one, $2\deg(L)-\deg(E)$, is measuring the
stability of the vector bundle $E$. Second one, $s-2\epsilon(L)$,
is measuring the stability of the morphism or, with the previous
observations, the generic stability of the set of points defined in
$\mathbb{P}_{\mathbb{C}}^{1}$, fiberwise, as in Examples \ref{exG}
and \ref{exG2}, when varying along the covering. Therefore, an object destabilizing a rank $2$
tensor is an object which contradicts these two stabilities,
weighted by $\tau$, and the Harder-Narasimhan subbundle is the
unique one which maximally does, for a $\tau$-unstable tensor.

The sets of points in each fiber defined by $\varphi$ give a
covering of degree $s$,
$$\mathbb{P}(E)\supset X'\rightarrow X\; .$$
In the following, we rewrite the stability of the sets of points, fiberwise, as stability for the covering, using intersection theory for ruled surfaces.

\begin{prop}\cite[Proposition V.2.8]{Ha}
\label{normalized}
Given a ruled surface $\mathbb{P}(E)$, there
exists $E'\simeq E\otimes N$, with $N$ line bundle, such that
$H^{0}(E')\neq 0$ but for all line bundles $N'$ with negative
degree we have $H^{0}(E'\otimes N')=0$. Therefore,
$\mathbb{P}(E)=\mathbb{P}(E')$ and the integer $e=-\deg E'$ is an
invariant of the ruled surface. Furthermore, in this case, there
exists a section $\sigma_{0}:X\rightarrow \mathbb{P}(E')$ with
image $C_{0}$, such that $\mathcal{L}(C_{0})\simeq
\mathcal{O}_{X}(1)$.
\end{prop}

For a ruled surface $\mathbb{P}(E')$ we say that $E'$ is
\emph{normalized} if it satisfies the conditions of the
Proposition \ref{normalized}.

Let $\mathbb{P}(E')$ be a ruled surface with $E'$ normalized. Let
$\sigma:X\rightarrow \mathbb{P}(E)$ be a section, and let $D=\im
\sigma$ be a divisor on $\mathbb{P}(E)$. It can be proved that
$\deg(L)=-e-C_{0}\cdot D$, with these conventions (c.f.
\cite[Proposition V.2.9]{Ha}). Let us define, by analogy,
$\epsilon(\sigma)=\epsilon(D)$ as the number of branches of $X'$
which generically do coincide with $D$, the section defined by
$\sigma$, counted with multiplicity.

\begin{dfn}
Let $(E,\varphi:\overbrace{E\otimes\cdots \otimes E}^{\text{s times}}\too \mathcal{O}_{X}))$ be a symmetric non degenerate rank $2$ tensor over $X$. Let $f:X'\rightarrow X$
be the covering defined by $(E,\varphi)$, $X'\subset \mathbb{P}(E)$. Let $\tau$ be a positive number. We say
that $f$ is \emph{$\tau$-unstable} if there exists a section $\sigma:X\rightarrow \mathbb{P}(E)$ with image $D$,
i.e. there exists a line subbundle $L\subset E$, such that the following holds
$$-2C_{0}\cdot D-e+\tau(s-2\epsilon(D))>0$$
\end{dfn}

\begin{prop}
Let $\tau$ be a positive number. A symmetric non degenerate $\rk 2$ tensor $(E,\varphi)$ is $\tau$-unstable if and only if the associated
covering $f:X'\rightarrow X$ is $\tau$-unstable.
\end{prop}
\begin{pr}
It is only needed to check that we can assume
$X'\subset\mathbb{P}(E')$ with $E'$ normalized (c.f. Proposition
\ref{normalized}), in the definition of stability of $f$. Let $N$
be a line bundle over $X$. If we change $E$ by $E'=E\otimes N$,
then we have the line subbundle $L\otimes N\subset E'$ (by
exactness of the tensor product with locally free sheaves), and
$$\deg(E')=\deg(E\otimes N)=\deg(E)+2\deg(N)\; ,$$
$$\deg(L\otimes N)=\deg(L)+\deg(N)\; ,$$
so the quantity $2\deg(L)-\deg(E)$ is invariant by tensoring $E$
with a line bundle.

For the invariance of the rest of the formula, also note that we
can trivially extend the definition of the morphism $\varphi$,
$$\varphi':(E')^{\otimes s}=E^{\otimes
s}\otimes N^{\otimes s}\rightarrow \mathcal{O}_{X}$$ and, then, it
is $\epsilon'(L\otimes N)=\epsilon(L)$.
\end{pr}

\begin{thm}
If $f:X'\rightarrow X$ is a degree $s$ covering coming from a symmetric non degenerate rank $2$ tensor
$(E,\varphi)$ which is $\tau$-unstable, then there exists a unique section
$\sigma:X\rightarrow \mathbb{P}(E)$ with image $D$, giving maximum
for
$$-2C_{0}\cdot D-e+\tau(s-2\epsilon(D))>0\, .$$
We call $\sigma$ the \emph{Harder-Narasimhan section} of the covering.
\end{thm}

\section{Rank $3$ tensors and beyond}
\label{rk3beyond}

This final section of chapter \ref{chaptersheaves} contains
some observations about the rank $3$ tensors case, which is the
first one we cannot apply directly the techniques used in the
previous sections. The crucial point will be the impossibility of
rewriting the Kempf function (c.f. Definition \ref{kempffunction})
in this case as a geometrical function as in Proposition
\ref{identification} because, as we will see, the argument
$\Gamma$ in that geometrical function (which represents the
weights $n_{\bullet}$ in Definition \ref{stabilityfortensors})
depends on the vector $v$ (which represents the filters
$E_{\bullet}$ in Definition \ref{stabilityfortensors}) which does
not allow us to apply results of subsection \ref{convexcones}.

\subsection{Independence between multi-indexes and weights}

In the previous sections we have been able to carry out the
program designed for torsion free coherent sheaves in different
cases of tensors: holomorphic pairs in section \ref{kempfpairs}
and $\rk 2$ tensors in section \ref{kempfrk2}. The proof of the
correspondence between the $1$-parameter subgroup of Kempf and the
Harder-Narasimhan filtration in that case is based on proving
properties related with the convexity for an arbitrary filtration
of subobjects, to show that the candidates to be the Kempf
filtration are very particular, so this filtration is unique. For
holomorphic pairs we previously know about the Harder-Narasimhan
filtration so, by uniqueness, it coincides with the Kempf
filtration. For rank $2$ tensors, we define the Harder-Narasimhan filtration as the unique Kempf filtration. 

We know that the Kempf filtration gives maximum for some function,
the Kempf function, which depends on the data of the filters
(ranks, Hilbert polynomials, etc.) and the data of the weights or the
exponents of the $1$-parameter subgroup. And the key point is to
rewrite the Kempf function as a geometrical function in an
Euclidean space (c.f. definition of $\mu_v$ in (\ref{eq:mu})).
There we prove that, if we think of the filtration (referring just to the flag $V_{\bullet}\subset V$ without the weights $n_{\bullet}$) as a graph, the weights giving the maximum for
the function are given by the convex envelope of the graph (c.f.
Theorem \ref{maxconvexenvelope}).

Thanks to the independence between the vector $v$ giving the graph
of the $m$-Kempf filtration $E^{m}_{\bullet}\subset E$ and the weights
$\Gamma_{i}$, we are able to rewrite the Kempf function as a
geometrical function, then it is possible to interpret this
geometrical function as depending on two values, $v$ and $\Gamma$,
independently. For holomorphic pairs this is easily proved by
checking that
$\epsilon_{i}(\overline{\Phi},V_{\bullet})=\epsilon(E_{i})$ is
independent of $\Gamma$, because the symbol $\epsilon_{i}$ only
depends on the vanishing of the restriction of the morphism
$\varphi$ to $E_{i}$ (c.f. Lemma \ref{independenceofweightspairs}).
For $\rk 2$ tensors, this is proven in Lemma
\ref{independenceofweightsrk2}. Nevertheless, for tensors in general,
this is not possible as we are going to show.

Recall Definition \ref{deftensor} and expression (\ref{rightmu}) in the stability condition for tensors,
\begin{equation}
\label{rightmurk3} \mu(\varphi,E_{\bullet},n_{\bullet})=\min_{I\in
\mathcal{I}} \{\gamma_{r_{i_1}}+\cdots+\gamma_{r_{i_s}}:
\,\varphi|_{(E_{i_1}\otimes\cdots\otimes E_{i_s})^{\oplus c}}\neq
0 \}\; ,
\end{equation}
where $\mathcal{I}=\{1,...,t+1\}^{\times s}$ is the set of all multi-indexes $I=(i_{i},...,i_{s})$ and
$(E_{\bullet},n_{\bullet})$ is a weighted filtration of $E$. Recall that the previous quantity was expressed in another
way in (\ref{rightmu2}),
$$\mu(\varphi,E_\bullet,n_{\bullet})=\sum_{i=1}^{t}n_{i}(sr_{i}-\epsilon_{i}(E_{\bullet})r)\; .$$
Also recall that the data of the multi-index giving minimum in (\ref{rightmurk3}) is equivalent to the data of
the $\epsilon_{i}(E_{\bullet})$ in the second expression.

For holomorphic pairs (i.e. $s=1$ in Definition \ref{deftensor}),
Lemma \ref{independenceofweightspairs} guarantees that the
multi-index does not depend on the weights $\gamma_{r_{i_{j}}}$ of the $1$-parameter
subgroup. The multi-index is $i$ if $i=\min\{
i:\varphi_{E_{i}}\neq 0\}$, i.e., $\varphi_{E_{i-1}}=0$ and
$\varphi_{E_{i}}\neq 0$. Lemma \ref{independenceofweightsrk2} does
the analogous for $\rk 2$ tensors, being the multi-index (in the symmetric case)
$(i_{1},\ldots,i_{s})=(1,\ldots,1,2,\ldots,2)$ where the number of
$1$'s is $\epsilon (L)=k$ if $\varphi|_{L^{\otimes (k+1)}\otimes
E^{\otimes (s-k-1)}}= 0$ and $\varphi|_{L^{\otimes k}\otimes
E^{\otimes (s-k)}}\neq 0$. However, beyond this cases, we cannot
assure the independence of the multi-index with the weights.

Suppose the case $s=2$, $r=3$ in Definition \ref{deftensor}, and
suppose, for simplicity, that the morphism $\varphi$ is symmetric.
The multi-index in this case will be $(i_{1},i_{2})$, where $0\leq
i_{j}\leq 3$ because recall that, in (\ref{rightmurk3}), the
vanishing of the morphism is checked generically. Then, the
multi-index is checking whether $\varphi|_{E_{i_{1}}\otimes
E_{i_{2}}}=0$ or not, hence by the symmetry of $\varphi$, the only
multi-indexes which can appear are $(1,1)$, $(1,2)$, $(1,3)$,
$(2,2)$, $(2,3)$. Therefore, for a given filtration $0\subset
L\subset F\subset E$ with weights $\gamma_{1},
\gamma_{2},\gamma_{3}$, the following six situations can
occur:
\begin{enumerate}
 \item $\varphi|_{L\otimes L}\neq 0$
\item $\varphi|_{F\otimes L}\neq 0$ and $\varphi|_{L\otimes L}=0$
\item $\varphi|_{E\otimes L}\neq 0$ and $\varphi|_{F\otimes F}=0$
\item $\varphi|_{F\otimes F}\neq 0$ and $\varphi|_{E\otimes L}=0$
\item $\varphi|_{E\otimes F}\neq 0$ and $\varphi|_{E\otimes L}=0$
\item $\varphi|_{F\otimes L}=0\; , \varphi|_{F\otimes F}\neq 0\; , \varphi|_{E\otimes L}\neq 0$
\end{enumerate}
Cases $1-5$ give a fixed multi-index, $(1,1)$, $(1,2)$,
$(1,3)$, $(2,2)$ and $(2,3)$ respectively. However, in case number $6$, the
multi-index will be $(1,3)$ if $\gamma_{1}+\gamma_{3}\leq
2\gamma_{2}$ or $(2,2)$  if $\gamma_{1}+\gamma_{3}\geq
2\gamma_{2}$. Hence, this is the simplest case where the
multi-index actually depends on the weights $\gamma_{i}$. Therefore, setting $s=2$ and $r=3$ in Definition \ref{deftensor} we get the first case for which these
features can occur.

If this happens, we are not able to rewrite the Kempf function as a geometrical function (c.f. (\ref{eq:mu}))
and prove an analogous to Proposition \ref{identification} to apply the argument of the convex envelope to look
for the vector $\gamma$ giving maximum in Theorem \ref{maxconvexenvelope}. This is the reason why the general
method described in this thesis breaks down and does not apply beyond rank $2$ tensors.

\subsection{Considerations for $\rk 3$ tensors}

Let us consider \emph{symmetric rank $3$ tensors of two arguments} over a smooth projective curve, i.e. the morphism
$\varphi: E\otimes E\rightarrow \mathcal{O}_{X}$ being symmetric. This case is obtained from section
\ref{exampletensors} by setting $s=2$ $c=1$, $b=0$, $R=\Spec \mathbb{C}$ and $\mathcal{D}=\mathcal{O}_{X}$, the
structure sheaf over $X\times R\simeq X$, in Definition \ref{deftensor}.

Let $\tau$ be a positive constant and consider a tensor $(E,\varphi)$ which is $\tau$-unstable (c.f. Definition
\ref{slopestabilitytensors}). Let $m_{0}$ be an integer as in Theorem \ref{GIT-delta} (such that
$\delta$-stability and GIT stability coincide) and such that $E$ is $m_{0}$-regular (picking a larger integer,
if necessary). Let $m\geq m_{0}$ and let $V\simeq H^{0}(E(m))$. In this case, the \emph{Kempf function} (c.f.
Definition \ref{kempffunction}) is
$$\mu(V_{\bullet},n_{\bullet})=\frac{\sum_{i=1}^{t}  n_{i} (r\dim V_{i}-r_{i}\dim V)+\frac{a_{2}}{a_{1}}
\sum_{i=1}^{t} n_{i} \big(2\dim V_{i}-\epsilon_{i}(\overline\Phi)\dim V \big)}{\sqrt{\sum_{i=1}^{t+1} {\dim
V^{i}} \Gamma_{i}^{2}}}\; ,$$ where
$$\frac{a_{1}}{a_{2}}=\frac{r\tau}{P_{E}(m)-s\tau}=\frac{3\tau}{P_{E}(m)-2\tau}\; .$$ Let
$$0\subset V_{1}\subset \cdots \subset V_{t+1}= V$$
be the \emph{Kempf filtration of $V$} (c.f. Theorem \ref{kempftheoremsheaves}) and let
$$
0\subseteq (E^{m}_{1},\varphi|_{E_{1}^{m}})\subseteq
(E^{m}_{2},\varphi|_{E_{2}^{m}})\subseteq\cdots
(E^{m}_{t},\varphi|_{E_{t}^{m}})\subseteq
(E^{m}_{t+1},\varphi|_{E_{t+1}^{m}})\subseteq (E,\varphi)\; .
$$
be the \emph{$m$-Kempf filtration of $(E,\varphi)$}, by evaluating the $V_{i}$. Suppose that we are able to
prove properties satisfied by the filters of the $m$-Kempf filtration (i.e. analogous to Propositions
\ref{boundedness} and \ref{task}) to rewrite the Kempf function as
$$K=\frac {\sum_{i=1}^{t+1} \Gamma_{i}[(r^{i}d - rd^i)+\frac{r\tau}{P-2\tau}
 (\epsilon^i P -2P^i)]} {\sqrt{\sum_{i=1}^{t+1} P^{i} \Gamma_{i}^{2}}}=$$
$$\frac {\sum_{i=1}^{t+1} \Gamma_{i}[(r^{i}d - rd^{i})+\tau
 (r\epsilon^{i}-2r^{i})]} {\sqrt{\sum_{i=1}^{t+1} r^{i} \Gamma_{i}^{2}}}=
\frac {\sum_{i=1}^{t} n_{i}[(rd_{i} -r_{i}d)+\tau
 (2r_{i}-\epsilon_{i}r)]} {\sqrt{\sum_{i=1}^{t+1} r^{i} \Gamma_{i}^{2}}}$$
(c.f. Proposition \ref{finalfunctionpairs}). Recall that we are considering the case $s=2$ and $\dim X=1$ in
Definition \ref{deftensor}.

Observe that, in order to achieve a maximum of the function, it is enough to consider saturated filtrations.
Indeed, note that $n_{i}>0$, $\deg E_{i}\leq \deg\overline{E_{i}}$ and $\rk E_{i}= \rk \overline{E_{i}}$, hence
the value of the function is greater on saturated filtrations. Also, by similar reasons, it is enough to
consider filtrations with increasing ranks (c.f. Remark \ref{saturatedfiltrations}). Therefore, in order to look
for the Kempf filtration, i.e., the filtration which maximizes the previous function, we can restrict our attention 
to filtrations of the form
\begin{equation}
\label{filtrk3} 0\subset (L,\varphi|_{L})\subset
(F,\varphi|_{F})\subset (E,\varphi)\; ,
\end{equation}
for which the Kempf function is
$$\frac {\sum_{i=1}^{3} \Gamma_{i}[(r^{i}d - rd^{i})+\tau
 (r\epsilon^{i}-r^{i}s)]} {\sqrt{\sum_{i=1}^{3} r^{i} \Gamma_{i}^{2}}}\; .$$
Let $(i_{1},i_{2})$ be the multi-index in the definition of
(\ref{rightmurk3}). Consider a filtration as in (\ref{filtrk3})
which is $\tau$-destabilizing, i.e. contradicting Definition
\ref{stabilityfortensors}. We want to check if this filtration is
the Kempf filtration, and for that we would like to ask ourselves
for the best $1$-parameter subgroup giving maximum for the Kempf
function. Note that, in the definition of stability for $\rk 3$
tensors, it is not enough to consider one-step filtrations, i.e.
subobjects, hence asking for the weights of the filtration is a
meaningful question.

The crucial fact is that the coefficients of the Kempf function
(understood as the function in Theorem \ref{maxconvexenvelope}, a
function on the exponents $\Gamma_{1}$, $\Gamma_{2}$ and
$\Gamma_{3}$) vary with the multi-index $(i_{1},i_{2})$. The
filtration (\ref{filtrk3}) will give a multi-index in
(\ref{rightmurk3}). If the multi-index we obtain falls into one of the cases $1-5$
in the list of the previous subsection, then the multi-index does
not depend on the weights. However, if
$$\varphi|_{F\otimes L}=0\; , \varphi|_{F\otimes F}\neq 0\; , \varphi|_{E\otimes L}\neq 0\; ,$$
we are in case $6$, and the multi-index will be $(1,3)$ if $\Gamma_{1}+\Gamma_{3}\leq 2\Gamma_{2}$, or $(2,2)$
otherwise. In this case, the vector $v$ (the vector of the graph associated to the filtration) will depend on
the multi-index, so we can have two possible vectors associated to the filtration (\ref{filtrk3}). Call the
vectors $x=(x_{1},x_{2},x_{3})$ and $y=(y_{1},y_{2},y_{3})$, then the coordinates of these two vectors are

$$x_{1}=r^{1}d-d^{1}r+\tau(r\epsilon^{(1,3)}_{1}-2r^{1})=v_{1}+\tau$$
$$x_{2}=r^{2}d-d^{2}r+\tau(r\epsilon^{(1,3)}_{2}-2r^{2})=v_{2}-2\tau$$
$$x_{3}=r^{3}d-d^{3}r+\tau(r\epsilon^{(1,3)}_{3}-2r^{3})=v_{3}+\tau\; ,$$

$$y_{1}=r^{1}d-d^{1}r+\tau(r\epsilon^{(2,2)}_{1}-2r^{1})=v_{1}-2\tau$$
$$y_{2}=r^{2}d-d^{2}r+\tau(r\epsilon^{(2,2)}_{2}-2r^{2})=v_{2}+4\tau$$
$$y_{3}=r^{3}d-d^{3}r+\tau(r\epsilon^{(2,2)}_{3}-2r^{3})=v_{3}-2\tau\; .$$
Note that $r^{1}=r^{2}=r^{3}=1$, $r=3$, and note that we denote with the upper index the different symbols
$\epsilon_{i}$ for each multi-index. Also we call $v_{i}=r^{i}d-d^{i}r$, for each $i$. 

Note that, in both cases, the following holds
$$\sum_{i=1}^{3}x_{i}=\sum_{i=1}^{3}y_{i}=\sum_{i=1}^{3}v_{i}=0\; .$$

Suppose that $(\Gamma_{1},\Gamma_{2},\Gamma_{3})$ is the vector giving maximum in the Kempf function. And
suppose that it verifies $\Gamma_{1}+\Gamma_{3}\leq 2\Gamma_{2}$, hence the multi-index is $(1,3)$, and the
vector of the graph associated to the filtration is $x$. Taking into account that
$\Gamma_{1}+\Gamma_{2}+\Gamma_{3}=0$, because $\Gamma$ is a $1$-parameter subgroup of $SL(N)$, the Kempf function will
be
$$K=\frac{\Gamma_{1}x_{1}+\Gamma_{2}x_{2}+\Gamma_{3}x_{3}}{\sqrt{\Gamma_{1}^{2}+\Gamma_{2}^{2}+\Gamma_{3}^{2}}}
=\frac{\Gamma_{1}v_{1}+\Gamma_{2}v_{2}+\Gamma_{3}v_{3}+\tau(\Gamma_{1}+\Gamma_{3}-2\Gamma_{2})}{\sqrt{\Gamma_{1}^{2}+\Gamma_{2}^{2}+\Gamma_{3}^{2}}}=$$
$$\frac{\Gamma_{1}(v_{1}-v_{2}+3\tau)+\Gamma_{3}(v_{3}-v_{2}+3\tau)}{\sqrt{2(\Gamma_{1}^{2}+\Gamma_{1}\Gamma_{2}+\Gamma_{3}^{2})}}\; ,$$
which is a function on two arguments. To maximize the function with respect to $\Gamma_{1}$ and $\Gamma_{3}$ we
set the gradient of $K$ equal to $0$ which gives
$$\Gamma_{1}=\Gamma_{3}(\frac{v_{1}+\tau}{v_{3}+\tau})\; .$$

If we suppose that the vector giving maximum in the Kempf function verifies $\Gamma_{1}+\Gamma_{3}\geq
2\Gamma_{2}$, we obtain the Kempf function
$$K=\frac{\Gamma_{1}y_{1}+\Gamma_{2}y_{2}+\Gamma_{3}y_{3}}{\sqrt{\Gamma_{1}^{2}+\Gamma_{2}^{2}+\Gamma_{3}^{2}}}\; ,$$
which we maximize with respect to $\Gamma_{1}$ and $\Gamma_{3}$ as well, obtaining
$$\Gamma_{1}=\Gamma_{3}(\frac{v_{1}-2\tau}{v_{3}-2\tau})\; .$$
Observe that, in both cases, the vector $\Gamma=(\Gamma_{1},\Gamma_{2},\Gamma_{3})$ which maximize the Kempf
function is exactly given by the vectors $x$ and $y$, because the $\Gamma_{1}$ and $\Gamma_{3}$ obtained are
precisely multiples of their coordinates.

Note that, $\Gamma_{1}+\Gamma_{3}\leq 2\Gamma_{2}$ implies $\Gamma_{1}\leq -\Gamma_{3}$, so to be congruent, in
the first case, it has to be $v_{1}+v_{3}+2\tau\leq 0$ (here we use that $v_{3}>0$, which has to hold by
convexity in Lemma \ref{lemmaA}). Similarly, in the second case, it has to be $v_{1}+v_{3}-4\tau\geq 0$
(whenever $v_{3}-2\tau\geq 0$). Observe that, as $\tau>0$, both conditions cannot hold at the same time, then
necessarily the multi-index is one of two, either $(1,3)$ or $(2,2)$.

Finally, consider the following example. Let $(E,\varphi)$ be a rank $3$ tensor over $X=\mathbb{P}_{\mathbb{C}}^{1}$,
where $E=\mathcal{O}_{\mathbb{P}_{\mathbb{C}}^{1}}(2)\oplus \mathcal{O}_{\mathbb{P}_{\mathbb{C}}^{1}}\oplus
\mathcal{O}_{\mathbb{P}_{\mathbb{C}}^{1}}(-1)$. Consider the filtration $0\subset (L,\varphi|_{L})\subset
(F,\varphi|_{F})\subset (E,\varphi)$ where
$L=\mathcal{O}_{\mathbb{P}_{\mathbb{C}}^{1}}(2)$, $F=\mathcal{O}_{\mathbb{P}_{\mathbb{C}}^{1}}(2)\oplus
\mathcal{O}_{\mathbb{P}_{\mathbb{C}}^{1}}$, and suppose that the matrix of the morphism $\varphi$, adapted to
the filtration $0\subset L\subset F\subset E$, is
$$\left(
    \begin{array}{ccc}
      0 & 0 & X  \\
      0 & X & X  \\
      X & X & X\\
    \end{array}
  \right)$$
where $X$ represents a non zero element. Let $\tau=\frac{1}{3}$ and observe that $(E,\varphi)$ is
$\tau$-unstable, because we can find weights $\Gamma_{1}$, $\Gamma_{2}$ and $\Gamma_{3}$ such that the filtration
 $0\subset L\subset F\subset E$ contradicts Definition \ref{stabilityfortensors}. We consider such filtration and look for the best weights in order to maximize the Kempf function.
Because of how $\varphi$ looks like for this filtration, i.e. not knowing if the multi-index is $(1,3)$ or
$(2,2)$, we have to apply the previous analysis.

See that 
$$v_{1}=\rk L\cdot \deg E-\deg L\cdot \rk E=-5$$
$$v_{2}=\rk F/L\cdot \deg E-\deg F/L\cdot \rk E=1$$
$$v_{3}=\rk E/F\cdot \deg E-\deg E/F\cdot \rk E=4\; ,$$
hence we can only be in the first case, where the
multi-index is given by $(1,3)$. Substituting, we get that the best $1$-parameter subgroup is given by the
vector,
$$\Gamma=(\frac{-5+\tau}{4+\tau},\frac{1-2\tau}{4+\tau},1)=(\frac{-14}{13},\frac{1}{13},1)\; .$$
Note that we set $\Gamma_{3}=1$, and recall that the Kempf function
is invariant by rescaling the $\Gamma_{i}$ (c.f. subsection
\ref{convexcones}).

To know if the filtration $0\subset (L,\varphi|_{L})\subset
(F,\varphi|_{F})\subset (E,\varphi)$ is the Kempf filtration, i.e. to know if it is the filtration
giving maximum for the Kempf function, we would have to check all possible destabilizing filtrations and the
values the Kempf function achieves for them, to choose the greatest value which has to correspond to the Kempf
filtration, the candidate to be defined as the Harder-Narasimhan filtration.

In view of this, we can define a class of tensors for which, the
ambiguity of two or more multi-indexes which can give the minimum
in (\ref{rightmu}) or (\ref{rightmuV}), cannot appear.

\begin{dfn}
\label{determinedmultiindex}
Let $\delta$ be a polynomial with positive leading coefficient of degree at most
$n-1$. Let $(E,\varphi)$ be a $\delta$-unstable tensor over an $n$-dimensional projective variety. Suppose that
there does not exist a destabilizing weighted filtration $(E_{\bullet},n_{\bullet})$ (i.e. contradicting the expression of Definition \ref{stabilityfortensors}) 
such that for it, the symbols
$\epsilon_{i}(\varphi)$ do depend on the weights $\Gamma_{i}$ (because, for example, the expression of the morphism $\varphi$ adapted to the filtration is 
particularly easy). We call these tensors, \emph{determined
multi-index} tensors.
\end{dfn}

It is clear that we can develop the same techniques of this chapter to show that the $m$-Kempf filtration
stabilizes with the integer $m$, and to construct a Harder-Narasimhan filtration for determined multi-index
tensors as in Definition \ref{determinedmultiindex}.

\chapter[Correspondence for representations of quivers]{Correspondence for representations of quivers}
\chaptermark{Corresp. for rep. of quivers}
\label{chapterquivers}

\section{Representations of quivers on vector spaces}
\label{qvectorspaces}

Let $Q$ be a finite quiver, given by a finite set of vertices and arrows between them, and a representation of
$Q$ on finite dimensional $k$-vector spaces, where $k$ is an algebraically closed field of arbitrary
characteristic. There exists a notion of stability for such representations given by King (c.f. \cite{Ki}) and,
more generally by Reineke (c.f. \cite{Re}) (both particular cases of the abstract notion of stability for an
abelian category that we can find in \cite{Ru}), and a notion of the existence of a unique Harder-Narasimhan
filtration with respect to that stability condition.

We consider the construction of a moduli space for these objects by King (c.f. \cite{Ki}) and associate to an
unstable representation an unstable point, in the sense of Geometric Invariant Theory, in a parameter space
where a group acts. Then, the $1$-parameter subgroup given by Kempf (c.f. Theorem \ref{kempftheorem0}), which is
maximally destabilizing in the GIT sense, gives a filtration of subrepresentations and we prove that it
coincides with the Harder-Narasimhan filtration for that representation.

The proof follows the argument given in chapter \ref{chaptersheaves} to establish the correspondence between the
$1$-parameter subgroup of Kempf and the Harder-Narasimhan filtration for the different cases studied there.
However, for representations of quivers on the category of vector spaces, the proof is much simpler, as there is
no need of proving that there exists an integer $m$, sufficiently large, such that the $m$-Kempf filtration does
not depend on $m$ (c.f. subsection \ref{kempfstabilizes}), because in this construction of the moduli space
there is no such integer $m$ involved.


The definition of stability for a representation of a quiver (c.f. Definition \ref{Qstability}) contains two sets of parameters,
the coefficients of the linear functions $\Theta$ and $\sigma$. In
\cite{Ke}, the $1$-parameter subgroup is taken to maximize certain
function which depends on the choice of a linearization of the action of
the group we are taking the quotient by,
and a \textit{length} in the set of $1$-parameter subgroups
(c.f. Definition \ref{length}). In the case of sheaves the group
is $SL(N)$, which is simple, so any such length is unique
up to multiplication by a scalar, whereas for finite dimensional
representations of quivers we quotient by a product of general
linear groups, so we have to choose a scalar for each factor in
the choice of a length. Hence, we set the positive coefficients of $\sigma$ precisely as these scalars
and consider a particular linearization depending on $\sigma$ and $\Theta$,
in order to relate the Harder-Narasimhan filtration of a representation with the $1$-parameter subgroup given by Kempf in \cite{Ke} (c.f. Theorem \ref{kempfHNquivers}).

\subsection{Harder-Narasimhan filtration for representations of quivers}
A \emph{finite quiver} $Q$ is given by a finite set of vertices
$Q_{0}$ and a finite set of arrows $Q_{1}$. The arrows will be
denoted by $(\alpha: v_{i}\rightarrow v_{j})\in Q_{1}$. We denote
by $\mathbb{Z}Q_{0}$ the free abelian group generated by $Q_{0}$.

The following figures show different examples of finite quivers:

$$\xymatrix{
\bullet
\ar@(ur,dr)}\;\;\;\;\;\;\;\;\;\;\;\;\;\;\;\;\;\;\;\;\;\;\;\;\;\;\;\;\;\;\;\xymatrix{
\bullet \ar[r] &
\bullet}\;\;\;\;\;\;\;\;\;\;\;\;\;\;\;\;\;\;\;\;\;\;\;\;\;\;\;\;\;\;\;\xymatrix{
\bullet \ar @/^/[r] \ar @/_/[r] & \bullet}$$
$$ $$
$$\xymatrix{
 & \bullet\ar[dl] \ar[dr] & \\
 \bullet \ar[rr] & & \bullet}\;\;\;\;\;\;\;\;\;\;\;\;\;\;\;\;\;\;\;\;\;\;\;\;\;\;\;\;\;\;\;\;\;\;\;\;\xymatrix{
\bullet \ar[r] \ar[dr] & \bullet\ar[dr] & \bullet\ar@/^/[d]\\
& \bullet\ar[r] \ar[u]& \bullet \ar[u]}$$

Fix $k$, an algebraically closed field of arbitrary
characteristic. Let $\mod kQ$ be the category of finite
dimensional representations of $Q$ over $k$. Such category is an
abelian category and its objects are given by tuples
$$M=((M_{v})_{v\in Q_{0}},(M_{\alpha}:M_{v_{i}}\rightarrow M_{v_{j}})_{\alpha:v_{i}\rightarrow v_{j}})$$
of finite dimensional $k$-vector spaces and $k$-linear maps
between them. The \emph{dimension vector} of a representation is given by
$\underline{\dim}M=\sum_{v\in Q_{0}}\dim_{k}M_{v}\cdot
v\in\mathbb{N}Q_{0}$.

For example, in the previous figures, a representation of the
first quiver on the top left on finite dimensional vector spaces is an endomorphism of a vector space,
and a representation of the one in the top center is a
homomorphism between two vector spaces.

Let $\Theta$ be a set of numbers $\Theta_{v}$ for each $v\in
Q_{0}$ and define a linear function
$\Theta:\mathbb{Z}Q_{0}\rightarrow \mathbb{Z}$, by
$$\Theta(M):=\Theta(\underline{\dim}M)=\sum_{v\in Q_{0}}\Theta_{v}\dim_{k}M_{v}\; .$$

Let $\sigma$ be a set of strictly positive numbers $\sigma_{v}$
for each $v\in Q_{0}$, and define a (strictly positive) linear
function $\sigma:\mathbb{Z}Q_{0}\rightarrow \mathbb{Z}$, by
$$\sigma(M):=\sigma(\underline{\dim}M)=\sum_{v\in Q_{0}}\sigma_{v}\dim_{k}M_{v}\; .$$
We call $\sigma(M)$ \emph{the total dimension of $M$}. we will refer to $\Theta$ and $\sigma$ indistinctly
meaning the sets of numbers $\Theta_{v}$ and $\sigma_{v}$ or the linear functions.

For a non-zero representation $M$ of $Q$ over $k$, define its
\emph{slope} by
$$\mu_{(\Theta,\sigma)}(M):=\frac{\Theta(M)}{\sigma(M)}\; .$$

\begin{dfn}
\label{Qstability}
A representation $M$ of $Q$ over $k$ is \emph{$(\Theta,\sigma)$-semistable} if for all non-zero
subrepresentations $M'$ of $M$, we have
$$\mu_{(\Theta,\sigma)}(M')\leq\mu_{(\Theta,\sigma)}(M)\; .$$
If the inequality is strict for every non-zero subrepresentation, we say that $M$ is \emph{$(\Theta,\sigma)$-stable}. If $M$ is not $(\Theta,\sigma)$-semistable 
we say that it is \emph{$(\Theta,\sigma)$-unstable}.
\end{dfn}

\begin{lem}
\label{changesthetasigma} If we multiply the linear function
$\Theta$ by a non-negative integer, or if we add an integer
multiple of the strictly positive linear function $\sigma$ to
$\Theta$, the semistable (resp. stable) representations remain
semistable (resp. stable).
\end{lem}
\begin{pr}
Let $\Theta'=a\cdot \Theta+b\cdot \sigma,\; a,b\in\mathbb{Z},\; a>0,$ be another linear function and note that
$$\frac{\Theta'(M')}{\sigma(M')}\leq \frac{\Theta'(M)}{\sigma(M)}\Leftrightarrow \frac{a\cdot \Theta(M')+b\cdot \sigma(M)}{\sigma(M')}\leq
\frac{a\cdot \Theta(M)+b\cdot \sigma(M)}{\sigma(M)}$$
$$\Leftrightarrow \frac{\Theta(M')}{\sigma(M')}\leq \frac{\Theta(M)}{\sigma(M)}\; .$$\end{pr}

\begin{rem}
\label{alastair}
In \cite{Ki}, the stability condition (c.f
\cite[Definition 1.1]{Ki}) is formulated by not considering representations with different dimension vectors.
This leads to the construction of a moduli space and
$S$-filtrations (or Jordan-H\"{o}lder filtrations) but not to define a Harder-Narasimhan filtration,
for which is needed a slope condition as in Definition
\ref{Qstability}.

This slope stability condition, the $(\Theta,\sigma)$-stability
(c.f. Definition \ref{Qstability}), can be turned out into a
stability condition as in \cite{Ki}, by clearing denominators
$$\theta(M')=\Theta(M)\sigma(M') - \sigma(M)\Theta(M')\; ,$$
where $\theta$ is the function in \cite[Definition 1.1]{Ki} (observe that
$\theta(M)=0$), $\Theta$ and $\sigma$ are as in Definition
\ref{Qstability}, and $M'\subset M$ is a subrepresentation.

We will apply this in Proposition \ref{GITstab-stabquivers}, to
relate $(\Theta,\sigma)$-stability with GIT stability.
\end{rem}

\begin{rem}
The definition of stability which appears in \cite{Re} considers
$\sigma_{v}=1$ for each $v\in Q_{0}$, although we consider a
strictly positive linear function $\sigma$ in general. The
notation of $\sigma$ agrees with \cite{AC}, \cite{ACGP},
\cite{Sch}, while $\Theta$ agrees with \cite{Re} but in the other
references it is substituted by different notations closer to
classical moduli problems where the stability notion depends on
parameters ($\tau$-stability or $\rho$-stability).
\end{rem}

\begin{lem}\cite[Definition 1]{Ru},
\cite[Lemma 4.1]{Re} \label{stabilitycategory} Let $0\rightarrow
X\rightarrow Y\rightarrow Z\rightarrow 0$ be a short exact
sequence of non-zero representations of $Q$ over $k$. Then
$\mu_{(\Theta,\sigma)}(X)< \mu_{(\Theta,\sigma)}(Y)$ if and only
if $\mu_{(\Theta,\sigma)}(X)<\mu_{(\Theta,\sigma)}(Z)$ if and only
if  $\mu_{(\Theta,\sigma)}(Y)<\mu_{(\Theta,\sigma)}(Z)$.
\end{lem}
\begin{pr}
 Note that $\sigma(Y)=\sigma(X)+\sigma(Z)$ and, therefore
$$\mu_{(\Theta,\sigma)}(Y)=\frac{\Theta(Y)}{\sigma(Y)}=\frac{\Theta(X)+\Theta(Z)}{\sigma(X)+\sigma(Z)}\; ,$$
from which the statement follows.
\end{pr}

\begin{thm}\cite[Theorem 2]{Ru},
\cite[Lemma 4.7]{Re}
\label{HNquivers}
Given linear functions $\Theta$ and $\sigma$, (being $\sigma$ strictly positive), every representation $M$ of $Q$ over $k$ has a unique filtration
$$0\subset M_{1}\subset M_{2}\subset \ldots \subset M_{t}\subset M_{t+1}=M$$
verifying the following properties, where $M^{i}:=M_{i}/M_{i-1}$,
\begin{enumerate}
\item $\mu_{(\Theta,\sigma)}(M^{1})>\mu_{(\Theta,\sigma)}(M^{2})>\ldots >\mu_{(\Theta,\sigma)}(M^{t})>\mu_{(\Theta,\sigma)}(M^{t+1})$
\item The quotients $M^{i}$ are $(\Theta,\sigma)$-semistable
\end{enumerate}
This filtration is called the \emph{Harder-Narasimhan filtration} of $M$ (with respect to $\Theta$ and $\sigma$).
\end{thm}
\begin{pr}
The proof follows the usual argument to show the existence and uniqueness of the Harder-Narasimhan filtration.

Using Lemma \ref{stabilitycategory} we can prove the existence of a unique subrepresentation $M_{1}$, whose
slope is maximal among all the subrepresentations of $M$, and of maximal total dimension $\sigma(M_{1})$ among
those of maximal slope (c.f. \cite[Proposition 1.9]{Ru}, \cite[Lemma 4.4]{Re}). Then, proceed by recursion on
the quotient $M/M_{1}$.
\end{pr}

\subsection{Moduli space of representations of quivers}
Fix $k$ an algebraically closed field of arbitrary characteristic.
Fix a dimension vector $d\in\mathbb{Z}Q_{0}$ and fix $k$-vector
spaces $M_{v}$ of dimension $d_{v}$ for all $v\in Q_{0}$. Fix
linear functions $\Theta, \sigma:\mathbb{Z}Q_{0}\rightarrow
\mathbb{Z}$, being $\sigma$ strictly positive. We recall the
construction by King (c.f. \cite{Ki}) of a moduli space for
representations of $Q$ over $k$ with dimension vector $d$.

Consider the affine $k$-space
$$\mathcal{R}_{d}(Q)=\bigoplus_{\alpha:v_{i}\rightarrow v_{j}}\Hom_{k}(M_{v_{i}},M_{v_{j}})\; ,$$
whose points parametrize representations of $Q$ on the $k$-vector
spaces $M_{v}$. The reductive linear algebraic group
$$G_{d}=\prod_{v\in Q_{0}}GL(M_{v})$$
acts on $\mathcal{R}_{d}(Q)$ by
$$(g_{v_{i}})_{v_{i}}\cdot (M_{\alpha})_{\alpha}=(g_{v_{j}}M_{\alpha}g_{v_{i}}^{-1})_{\alpha:v_{i}\rightarrow
v_{j}}\; ,$$ and the $G_{d}$-orbits of $M$ in $\mathcal{R}_{d}(Q)$
correspond bijectively to the isomorphism classes $[M]$ of
$k$-representations of $Q$ with dimension vector $d$. We will use
Geometric Invariant Theory to take the quotient of
$\mathcal{R}_{d}(Q)$ by $G_{d}$ and construct a moduli space of
representations of the quiver $Q$ on the $k$-vector spaces
$M_{v}$.

The action of $G_{d}$ on the affine space $\mathcal{R}_{d}(Q)$ can
be lifted by a character $\chi$ to the (necessarily trivial) line
bundle $L$ required by the Geometric Invariant Theory. Note that
the subgroup of the diagonal scalar matrices in $G_{d}$,
$$\Delta=\{(t1,\ldots,t1):t\in k^{\ast}\}\; ,$$ acts
trivially on $\mathcal{R}_{d}(Q)$. Then, we have to choose $\chi$
in such a way that $\Delta$ acts trivially on the fiber, in other
words, $\chi(\Delta)=1$.

Then, using the linear functions $\Theta$ and $\sigma$, consider
the character
$$\chi_{(\Theta,\sigma)}((g_{v})_{v\in Q_{0}}):=\prod_{v\in Q_{0}}\det(g_{v})^{(\Theta(d)\sigma_{v}-\sigma(d)\Theta_{v})}$$
of $G_{d}$, and note that $\chi_{(\Theta,\sigma)}(\Delta)=1$,
because $\sum_{v\in
Q_{0}}(\Theta(d)\sigma_{v}-\sigma(d)\Theta_{v})\cdot d_{v}=0$.

Given a linearization of an action by a character $\chi$, we say that $f$ is a \emph{relative invariant of weight $\chi^{n}$} if $f(g\cdot x)=\chi^{n}(g)\cdot f(x)\;\forall x$. 

\begin{dfn}\cite[Definition 2.1]{Ki}
A point $x\in\mathcal{R}_{d}(Q)$ is \emph{$\chi_{(\Theta,\sigma)}$-semistable} if there is
a relative invariant of weight $\chi_{(\Theta,\sigma)}^{n}$, $f\in
k[\mathcal{R}_{d}(Q)]^{G_{d},\chi^{n}_{(\Theta,\sigma)}}$ with
$n\geq 1$, such that $f(x)\neq 0$.
\end{dfn}

The algebraic quotient will be given by
$$\mathcal{R}_{d}(Q)/\!\!/(G_{d},\chi_{(\Theta,\sigma)})=\Proj\big(\bigoplus_{n\geq
0}k[\mathcal{R}_{d}(Q)]^{G_{d},\chi^{n}_{(\Theta,\sigma)}}\big)\;
.$$

\begin{prop}
\label{GITstab-stabquivers} A point $x_{M}\in \mathcal{R}_{d}(Q)$
corresponding to a representation $M\in \mod kQ$ is
$\chi_{(\Theta,\sigma)}$-semistable (resp.
$\chi_{(\Theta,\sigma)}$-stable) for the action of $G_{d}$ if and
only if $M$ is $(\Theta,\sigma)$-semistable (resp.
$(\Theta,\sigma)$-stable).
\end{prop}
\begin{pr}
It follows easily from \cite[Proposition 3.1]{Ki} and the
observation in Remark \ref{alastair}. In \cite{Ki}, given a linear
function $\theta$, a representation $M$ is $\theta$-semistable if
$\theta(M)=0$ and, for every subrepresentation $M'\subset M$, we
have $\theta(M')\geq 0$ (c.f. \cite[Definition 1.1]{Ki}). Then,
\cite[Proposition 3.1]{Ki} relates the $\theta$-stability with the
$\chi_{\theta}$-stability, where the character is
$$\chi_{\theta}((g_{v})_{v}):=\prod_{v\in Q_{0}}\det(g_{v})^{\theta_{v}}\; .$$
Hence, the $\chi_{(\Theta,\sigma)}$-stability with the character
given by
$$\chi_{(\Theta,\sigma)}((g_{v})_{v}):=\prod_{v\in Q_{0}}\det(g_{v})^{(\Theta(d)\sigma_{v}-\sigma(d)\Theta_{v})}\; ,$$
is equivalent to the $(\Theta,\sigma)$-stability in Definition
\ref{Qstability} because, given a subrepresentation $M'\subset M$,
the expression
$$\sum_{v\in
Q_{0}}(\Theta(M)\sigma_{v}-\sigma(M)\Theta_{v})\cdot \dim M'_{v}=\Theta(M)\sigma(M')-\sigma(M)\Theta(M')\geq 0$$
is equivalent to
$$\frac{\Theta(M')}{\sigma(M')}\leq \frac{\Theta(M)}{\sigma(M)}\; .$$\end{pr}

Now denote by $\mathcal{R}^{(\Theta,\sigma)-ss}_{d}(Q)$ the set of
$\chi_{(\Theta,\sigma)}$-semistable points.

\begin{thm}\cite[Proposition 4.3]{Ki},
\cite[Corollary 3.7]{Re} The moduli space
$\mathfrak{M}^{(\Theta,\sigma)}_{d}(Q)=\mathcal{R}^{(\Theta,\sigma)-ss}_{d}(Q)/\!\!/G_{d}$
is a projective variety which parametrizes $S$-equivalence classes
of $(\Theta,\sigma)$-semistable representations of $Q$ of
dimension vector $d$.
\end{thm}

By the Hilbert-Mumford criterion we can characterize
$\chi_{(\Theta,\sigma)}$-semistable points by its behavior under
the action of $1$-parameter subgroups. A $1$-parameter subgroup of
$G_{d}=\prod_{v\in Q_{0}}GL(M_{v})$ is a non-trivial homomorphism
$\Gamma:k^{\ast}\rightarrow G_{d}$. There exist bases of the
vector spaces $M_{v}$ such that $\Gamma$ takes the diagonal form

$$\left(
  \begin{array}{ccc}
    t^{\Gamma_{v_{1},1}} &  & \\
     & \ddots & \\
     &  & t^{\Gamma_{v_{1},t_{1}+1}} \\
  \end{array}
\right)\times \cdots \times
\left(
  \begin{array}{ccc}
    t^{\Gamma_{v_{s},1}} &  & \\
     & \ddots & \\
     &  & t^{\Gamma_{v_{s},t_{s}+1}} \\
  \end{array}
\right)$$
where $v_{1},\ldots ,v_{s}\in Q_{0}$ are the vertices of the quiver.

Let $x\in\mathcal{R}_{d}(Q)$ and suppose that $\lim_{t\rightarrow
0}\Gamma\cdot x$ exists and is equal to $x_{0}$. Then $x_{0}$ is a
fixed point for the action of $\Gamma$, and $\Gamma$ acts on the
fiber of the trivial line bundle over $x_{0}$ as multiplication by
$t^{a}$. Define the following numerical function,
$$\mu_{\chi_{(\Theta, \sigma)}}(x,\Gamma)=-a\; .$$
The next proposition establishes a variant of the Hilbert-Mumford criterion given in Theorem
\ref{HMcrit}.

\begin{prop}\cite[Proposition 2.5]{Ki}
\label{mrwquivers} A point $x_{M}\in \mathcal{R}_{d}(Q)$ corresponding to a representation $M$ is
$\chi_{(\Theta,\sigma)}$-semistable if and only if every
$1$-parameter subgroup $\Gamma$ of $G_{d}$, for which
$\lim_{t\rightarrow 0} \Gamma(t)\cdot x_{M}$ exists, satisfies
$\mu_{\chi_{(\Theta, \sigma)}}(x_{M},\Gamma)\leq 0$.
\end{prop}

\begin{rem}
Note that in Proposition \ref{mrwquivers} we change the sign of the numerical function $\mu_{\chi_{(\Theta,
\sigma)}}(x_{M},\Gamma)$ with respect to \cite{Ki} (as we did when changing the character in the proof of
Proposition \ref{GITstab-stabquivers}), in congruence with \cite{Ke} and the numerical function in Theorem
\ref{HMcrit}.
\end{rem}

The action of a $1$-parameter subgroup $\Gamma$ of $G_{d}$
provides a decomposition of each vector space $M_{v}$, associated
to each vertex $v\in Q_{0}$, in weight spaces
$$M_{v}=\bigoplus_{n\in \mathbb{Z}}M_{v}^{n}\; ,$$
where $\Gamma(t)$ acts on the weight space $M_{v}^{n}$ as
multiplication by $t^{n}$. Every $1$-parameter subgroup, for
which $\lim_{t\rightarrow 0} \Gamma(t)\cdot x$ exists, determines
a weighted filtration $M_{\bullet}\subset M$ of subrepresentations
(c.f. \cite{Ki}),
$$0\subset M_{1}\subset M_{2}\subset \ldots \subset M_{t}\subset M_{t+1}=M\; ,$$
where $M_{i}$ is the subrepresentation with vector spaces
$M_{v,i}:=\bigoplus _{n\leq i}M_{v}^{n}$ for each vertex $v\in
Q_{0}$, and the weight corresponding to each quotient
$M^{i}=M_{i}/M_{i-1}$ is $\Gamma_{i}$. Note that two $1$-parameter
subgroups giving the same filtration are conjugated by an element
of the parabolic subgroup of $G_{d}$ defined by the filtration.
Therefore, the numerical function $\mu_{\chi_{(\Theta,
\sigma)}}(x_{M},\Gamma)$, has a simple expression in terms of the
filtration $M_{\bullet}\subset M$ (c.f. calculation in \cite{Ki}):
\begin{equation}
\label{pairing} \mu_{\chi_{(\Theta, \sigma)}}(x_{M},\Gamma)=\sum_{v\in
Q_{0}}\big[\big(\Theta(M)\sigma_{v}-\sigma(M)\Theta_{v}\big)\cdot
\sum_{i=1}^{t_{v}+1}\Gamma_{v,i}\dim M_{v}^{i}\big]\; .
\end{equation}
Let $d_{i}$, $d^{i}$ be the dimension vectors of the
subrepresentation $M_{i}$ and the quotient $M^{i}=M_{i}/M_{i-1}$,
respectively. The action of $\Gamma$ on the point corresponding to
a representation $M$ has different weights for each vertex $v\in
Q_{0}$, but collect all different weights $\Gamma_{i}$ corresponding to any vertex and form
the vector
$$\Gamma=(\Gamma_{1},\Gamma_{2},\ldots,\Gamma_{t},\Gamma_{t+1})$$ verifying $\Gamma_{1}<\Gamma_{2}<\ldots<\Gamma_{t}<\Gamma_{t+1}$.
Hence, (\ref{pairing}) turns out to be
\begin{equation}
\label{pairing2} \mu_{\chi_{(\Theta,
\sigma)}}(x_{M},\Gamma)=\sum_{i=1}^{t+1}\Gamma_{i}\cdot
[\Theta(M)\cdot \sigma(M^{i})-\sigma(M)\cdot \Theta (M^{i})]\;,
\end{equation}
and Proposition \ref{mrwquivers} can be rewritten in terms of filtrations of $M$.

\begin{prop}
\label{mrw2} A point $x_{M}\in \mathcal{R}_{d}(Q)$ corresponding
to a representation $M$ of $Q$ over $k$, is
$\chi_{(\Theta,\sigma)}$-semistable if and only if every
$1$-parameter subgroup $\Gamma$ of $G_{d}$, defining a filtration
of subrepresentations of $M$ $$0\subset M_{1}\subset M_{2}\subset
\ldots \subset M_{t}\subset M_{t+1}=M\; ,$$ satisfies that
$$\mu_{\chi_{(\Theta,
\sigma)}}(x_{M},\Gamma)=\sum_{i=1}^{t+1}\Gamma_{i}\cdot
[\Theta(M)\cdot \sigma(M^{i})-\sigma(M)\cdot \Theta (M^{i})]\leq
0\; .$$
\end{prop}

\subsection{Kempf theorem}

Given a weighted filtration of $M$,
$$0\subset M_{1}\subset M_{2}\subset \ldots \subset M_{t}\subset M_{t+1}=M\; ,$$
and $\Gamma_{1}<\Gamma_{2}<\ldots<\Gamma_{t}<\Gamma_{t+1}$, define the following function which is a \emph{Kempf function} (c.f. Definition \ref{kempffunction}) for this problem,
\begin{equation}
\label{kempffunctionquivers}
K(M_{\bullet},\Gamma)=\frac{\sum_{i=1}^{t+1}\Gamma_{i}\cdot [\Theta(M)\cdot
\sigma(M^{i})-\sigma(M)\cdot \Theta (M^{i})]}{\sqrt{\sum_{i=1}^{t+1}\sigma(M^{i})\cdot \Gamma_{i}^{2}}}
\end{equation}
It is a function whose numerator is equal to the numerical
function $\mu_{\chi_{(\Theta, \sigma)}}(x_{M},\Gamma)$ and the denominator is a \textit{length} of the $1$-parameter
subgroup $\Gamma$. Given a reductive linear algebraic group $G$, recall the notion of \emph{length} in $\Gamma(G)$,
the set of all $1$-parameter subgroups (c.f. Definition \ref{length}).

If $G$ is simple, in characteristic zero all choices of length will be multiples of the Killing form in the
Lie algebra $\mathfrak{g}$ (note that in this case $\Gamma(G)\subset \mathfrak{g}$). For an almost simple
group in arbitrary characteristic (a group $G$ whose center $Z$ is finite and $G/Z$ is simple, e.g. $SL(N)$ in positive characteristic), all lengths differ also by a scalar.

However, in this case, the group is a product of general linear groups, which is not simple. Then, there are several simple factors in the group and we can take a 
different multiple of the Killing form for each factor. Hence, observe that in the Kempf function (\ref{kempffunctionquivers}), the denominator of the expression is a
function verifying the properties of the definition of a length (c.f. Definition \ref{length}). The different multiples we take for each factor appear to be related to the choice
of the strictly positive linear function $\sigma$. 

Therefore, we can rewrite Theorem
\ref{kempftheoremsheaves} in our case as follows:

\begin{thm}
\label{kempftheoremquivers} Given a $\chi_{(\Theta,\sigma)}$-unstable point
$x_{M}\in \mathcal{R}_{d}(Q)$ corresponding to a representation $M$, there exists a unique weighted
filtration, i.e. $0\subset M_{1}\subset \cdots \subset M_{t+1}= M$
and real numbers
$\Gamma_{1}<\Gamma_{2}<\ldots<\Gamma_{t}<\Gamma_{t+1}$, called the
\emph{Kempf filtration of M}, such that the Kempf function $K(M_{\bullet},\Gamma)$
achieves the maximum among all filtrations and weights verifying
$\Gamma_{1}<\Gamma_{2}<\ldots<\Gamma_{t}<\Gamma_{t+1}$.
\end{thm}

Note that the length we are considering depends on the choice of
$\sigma$ and the Kempf function depends both on the length and the
linearization of the group action, hence depends both on $\Theta$
and $\sigma$. In order to relate the Kempf filtration of $M$ with
the Harder-Narasimhan filtration, which also depends on $\Theta$
and $\sigma$, we have set the parameters conveniently in the expression
of the stability condition (c.f. Proposition
\ref{GITstab-stabquivers}).

\subsection{Kempf filtration is Harder-Narasimhan filtration}

Finally, we close the section by relating the Kempf filtration in Theorem \ref{kempftheoremquivers} and the
Harder-Narasimhan filtration in Theorem \ref{HNquivers}. We study the geometrical properties of the Kempf
filtration by associating to it a graph which encodes the two properties satisfied by the Harder-Narasimhan
filtration. We will use the results of subsection \ref{convexcones}.

Let $\Theta:\mathbb{Z}Q_{0}\rightarrow \mathbb{Z}$ be a linear function and let
$\sigma:\mathbb{Z}Q_{0}\rightarrow \mathbb{Z}$ be a strictly positive linear function. Let $M$ be a
representation of $Q$ over an algebraically closed field $k$ of arbitrary characteristic, which is
$(\Theta,\sigma)$-unstable. Consider the $\chi_{(\Theta,\sigma)}$-unstable point $x_{M}\in \mathcal{R}_{d}(Q)$
associated to $M$, by Proposition \ref{GITstab-stabquivers}. Let $0\subset M_{1}\subset \cdots \subset M_{t+1}=
M$ and $\Gamma_{1}<\Gamma_{2}<\ldots<\Gamma_{t}<\Gamma_{t+1}$ be the Kempf filtration of $M$, by Theorem
\ref{kempftheoremquivers}.

Let $M^{i}=M_{i}/M_{i-1}$ be the quotients of the filtration.
Consider the inner product in $\mathbb{R}^{t+1}$ given by the
matrix $$
 \left(
 \begin{array}{ccc}
 \sigma(M^1) & & 0 \\
  & \ddots & \\
 0 & & \sigma(M^{t+1})\\
 \end{array}
 \right)
 $$
where $\sigma(M^{i})>0$.
\begin{dfn}
\label{graphquivers} Given a filtration $0\subset M_{1}\subset \cdots
\subset M_{t+1}= M$ of subrepresentations of $M$, define
$v=(v_{1},...,v_{t+1})$, where
$$v_{i}=\Theta(M)-\frac{\sigma(M)}{\sigma(M^{i})}\Theta(M^{i})\; ,$$
the \emph{graph} or the \emph{vector associated to the filtration}.
\end{dfn}

Now we can identify the Kempf function (c.f. (\ref{kempffunctionquivers})) with the function in Theorem
\ref{maxconvexenvelope},
$$K(M_{\bullet},\Gamma)=\frac{\sum_{i=1}^{t+1}\Gamma_{i}\cdot [\Theta(M)\sigma(M^{i})-\sigma(M)\Theta(M^{i})]}{\sqrt{\sum_{i=1}^{t+1}\sigma(M^{i})\cdot
\Gamma_{i}^{2}}}=$$
$$=\frac{\sum_{i=1}^{t+1}\sigma(M^{i})\Gamma_{i}\cdot [\Theta(M)-\frac{\sigma(M)}{\sigma(M^{i})}\Theta (M^{i})]}{\sqrt{\sum_{i=1}^{t+1}
\sigma(M^{i})\cdot
\Gamma_{i}^{2}}}=\frac{(\Gamma,v)}{\|\Gamma\|}=\mu_{v}(\Gamma)\;
.$$ Note that $\sum_{i=1}^{t+1}b^{i}v_{i}=0$.

\begin{thm}
\label{kempfHNquivers} The Kempf filtration of $M$ is the
Harder-Narasimhan filtration of $M$.
\end{thm}
\begin{pr}
The vector $v$ associated to the Kempf filtration of $M$ in Definition \ref{graphquivers} verifies properties in
Lemmas \ref{lemmaA} and \ref{lemmaB}, which are precisely properties
$1$ and $2$ in Theorem \ref{HNquivers}, respectively. Lemma \ref{lemmaA} implies that $v_{i}<v_{i+1}$, for each $i$, hence
$$\Theta(M)-\frac{\sigma(M)}{\sigma(M^{i})}\Theta(M^{i})<\Theta(M)-\frac{\sigma(M)}{\sigma(M^{i+1})}\Theta(M^{i+1})\Leftrightarrow 
\frac{\Theta(M^{i})}{\sigma(M^{i})}>\frac{\Theta(M^{i+1})}{\sigma(M^{i+1})}\; ,$$ and Lemma \ref{lemmaB} implies the $(\Theta,\sigma)$-semistability 
of each quotient $M^{i}=M_{i}/M_{i-1}$. By uniqueness
of the Harder-Narasimhan filtration of $M$, both filtrations do
coincide.
\end{pr}

\section{Representations of quivers on coherent sheaves}
\label{sectionqsheaves}

In this final section we will show the correspondence between the $1$-parameter subgroup of Kempf in Theorem
\ref{kempftheorem0} and the filtration of Harder-Narasimhan in Theorem \ref{HNunique} through the language of
representations of quivers. A coherent sheaf will be a representation of a one vertex quiver on the category of
coherent sheaves. It will have a $H$-Kronecker associated, and we will associate to it a representation of
another quiver on vector spaces, to use the results of section \ref{qvectorspaces}.

We first present the $Q$-sheaves which are representations of a
quiver on the category of coherent sheaves, and its relation with
the Kronecker modules.

\subsection{Quiver sheaves and Kronecker modules}

Let $Q$ be a quiver and let $X$ be a projective variety. A
\emph{$Q$-sheaf} over $X$ is a representation $E$ of $Q$ in the category
of coherent sheaves over $X$, given by the data of a coherent
sheaf $E_{v}$ for all $v\in Q_{0}$ and a morphism of sheaves
$\phi_{\alpha}: E_{v_{i}}\rightarrow E_{v_{j}}$ for all $\alpha\in
Q_{1}$. Let $E$ be a $Q$-sheaf over $X$. Let $P$ be a set of
polynomials $P_{v}\in \mathbb{Q}[m]$, indexed by the vertices
$v\in Q_{0}$. A $Q$-sheaf $E$ has Hilbert polynomial vector $P$ if
$P_{v}$ is the Hilbert polynomial of each sheaf $E_{v}$.

Let $\kappa$ be a set of polynomials $\kappa_{v}\in \mathbb{Q}[m]$, for $v\in Q_{0}$, such that
$\kappa_{v}(m)>0$ for $m\gg 0$ and $\deg \kappa_{v}=t$ for all $v\in Q_{0}$, for a fixed integer $t\geq 0$
independent of $v$. Let us call $\sigma, \tau$ the sets of rational numbers $\sigma_{v}>0$, $\tau_{v}$,  indexed
by the vertices $v\in Q_{0}$, such that
$$\kappa_{v}(m)=\sigma_{v}m^{t}+\tau_{v}m^{t-1}+\ldots \; .$$
The \emph{$\kappa$-Hilbert polynomial} of a $Q$-sheaf $E$ is the
polynomial $P_{\kappa}(E)\in \mathbb{Q}[m]$ given by
$$P_{\kappa}(E,m):=\sum_{v\in Q_{0}}\kappa_{v}(m)P_{E_{v}}(m)\; ,$$ where each $P_{E_{v}}(m)$ is the Hilbert polynomial of the coherent sheaf $E_{v}$.

A $Q$-sheaf $E$ is called \emph{pure of dimension $e$}, if the sheaf $E_{v}$ is pure of dimension $e$
(independent of $v$), for all $v\in Q_{0}$. A \emph{$Q$-subsheaf} of a $Q$-sheaf $E$ is given by a subsheaf
$E'_{v}\subset E_{v}$ for each vertex such that the restrictions of the morphisms are compatible.

\begin{dfn}
\label{kstability} A $Q$-sheaf $E$ over $X$ is \emph{Gieseker
$\kappa$-semistable} if it is pure (of any dimension $e$) and
$$\frac{P_{\kappa}(E',m)}{\sum_{v\in _{0}}\sigma_{v} \rk
E'_{v}}\leq \frac{P_{\kappa}(E,m)}{\sum_{v\in _{0}}\sigma_{v} \rk E_{v}}\;\; \text{for  } m\gg 0\; ,$$ for each
non-zero $Q$-subsheaf $E'\subset E$, and \emph{Gieseker $\kappa$-stable} if, furthermore, the inequality is
strict for all proper $Q$-subsheaves $E'\subset E$. If $E$ is not Gieseker $\kappa$-semistable we say that it is \emph{Gieseker $\kappa$-unstable}. 
\end{dfn}

By $\kappa$-semistable, $\kappa$-stable and $\kappa$-unstable we mean, in the following, Gieseker $\kappa$-semistable, Gieseker
$\kappa$-stable and Gieseker $\kappa$-unstable. We can rewrite Definition \ref{kstability} as it appears in \cite{AC}.

\begin{lem}\cite[Lemma 7]{AC}
\label{kstability2} A $Q$-sheaf $E$ over $X$ is $\kappa$-semistable if and only if for all non-zero $E'\subset
E$,
$$\frac{P_{\kappa}(E',m)}{P_{\kappa}(E',l)}\leq \frac{P_{\kappa}(E,m)}{P_{\kappa}(E,l)}\;\; \text{  for  } l\gg m\gg 0\; ,$$
for each non-zero $Q$-subsheaf $E'\subset E$, and Gieseker $\kappa$-stable if, furthermore, the inequality is
strict for all proper $E'\subset E$.
\end{lem}

\'{A}lvarez-C\'{o}nsul and King in \cite{ACK} give a functorial construction of the moduli space of coherent
sheaves over a projective variety by associating to a sheaf a \emph{Kronecker module} which is a representation of
a particular quiver on vector spaces.

Let $l>m$ be integers and consider the sheaf
$T=\mathcal{O}(-l)\oplus \mathcal{O}(-m)$ together with a finite
dimensional $k$-algebra
$$\left(
                          \begin{array}{cc}
                            k & H \\
                            0 & k \\
                          \end{array}
                        \right)$$
of operators on $T$, where $A=L\oplus H\subset \End_{X}(T)$, $L=k\cdot e_{0}\oplus k\cdot e_{1}$ is the
semisimple algebra generated by the two projection operators onto the summands of $T$, and $H
=H^{0}(\mathcal{O}(l-m))= \Hom(\mathcal{O}_{X}(m),\mathcal{O}_{X}(l))$, acting on $T$ in the off-diagonal way.

We can give a right $A$-module structure on $M$ by giving a right $L$-module structure and a right $L$-module
map $M\otimes_{L} H\rightarrow V$. The first one is equivalent to a direct sum decomposition $M=V\oplus W$, being
$V=M\cdot e_{0}$ and $W=M\cdot e_{1}$, and the second one given by the map
$$\alpha:V\otimes
H\rightarrow W\; .$$ This structure given on $V$ is called a
\emph{$H$-Kronecker module}. We can also say that $A$ is the path
algebra of the quiver with two vertices and, after choosing a
basis for $H$, a number of $\dim H$ arrows between them. A
representation of this quiver is also a \emph{$H$-Kronecker
module}, equivalent to the previous definition by the standard
equivalence between representations of quivers and modules for
their path algebras.

Given a sheaf $E$, $\Hom_{X}(T,E)$ can be given a structure of
$H$-Kronecker module. Indeed, it has a natural right module
structure over $A\subset \Hom_{X}(T,T)$, given by composition of
maps, and we have the decomposition
$\Hom_{X}(T,E)=H^{0}(E(m))\oplus H^{0}(E(l))$ together with the
multiplication map $\alpha_{E}: H^{0}(E(m))\otimes H\rightarrow
H^{0}(E(l))$.

Given an $A$-module $M=V\oplus W$, an \emph{$A$-submodule $M'$} is given by
$V'\subset V$ and $W'\subset W$ such that $\alpha(V'\otimes
H)\subset W'$.

\begin{dfn}\cite[Definition 2.3]{ACK}
\label{QAstability} 
An $A$ module $M=V\oplus W$ is \emph{semistable} if
$$\frac{\dim V'}{\dim W'}\leq \frac{\dim V}{\dim W}$$
for every submodule $M'=V'\oplus W'\subset M$. If the previous
inequality is strict for every submodule, we say that $M$ is
\emph{stable}. If $M$ is not semistable, we say that it is \emph{unstable}.
\end{dfn}

We associate to a sheaf $E$ a Kronecker module in this way. An observation in \cite{ACK} points out that the GIT
semistability of the orbit of $E$ is equivalent to the natural semistability of the Kronecker module associated
(c.f. \cite[Remark 2.4]{ACK}). In the following theorem we also relate the stability of the $A$-module $M$ with
the stability of a representation of the quiver $\tilde{Q}$
$$\xymatrix{
\bullet \ar[r] &
\bullet}$$
on vector spaces, to use the results of section \ref{qvectorspaces}.

\begin{thm}
\label{equivstability}
Let $E$ be a coherent sheaf over $X$, pure
of dimension $e$, with Hilbert polynomial $P$. There exists $l\gg
m\gg 0$, such that the following are equivalent:
\begin{enumerate}
\item $E$ is semistable as in Definition \ref{giesekerstab}.
\item  $E$ is $m$-regular and the $A$-module $M=H^{0}(E(m))\oplus H^{0}(E(l))$
is semistable as in Definition \ref{QAstability}.
\item The representation $M$ of the two vertex quiver $\tilde{Q}=\{v_{0},v_{1}\}$ and one arrow between them
on $k$-vector spaces, where $M_{v_{0}}=H^{0}(E(m))$, $M_{v_{1}}=H^{0}(E(l))$, is $(\Theta,\sigma)$-semistable as
in Definition \ref{Qstability}, where the linear functions $\Theta$ and $\sigma$ are defined as $\Theta(M)=\dim
M_{v_{0}}$, $\sigma(M)=\dim M_{v_{0}}+\dim M_{v_{1}}$.
\item The point $x_{M}\in \mathcal{R}^{(\Theta,\sigma)}_{d}(\tilde{Q})$ is $\chi_{(\Theta,\sigma)}$-semistable, where $d$ is
the dimension vector of the representation $M$, $d_{v_{i}}=\dim
M_{v_{i}}$.
\end{enumerate}
\end{thm}
\begin{pr}
The equivalence between 1 and 2 follows from \cite[Theorem
5.10]{ACK}. For the equivalence between 2 and 3 note that,
defining the linear functions $\Theta$ and $\sigma$ as in the
statement $3$, it is
$$\frac{\dim M'_{v_{0}}}{\dim M'_{v_{1}}}
\leq \frac{\dim M_{v_{0}}}{\dim M_{v_{1}}}\Leftrightarrow$$
$$\frac{\dim M'_{v_{0}}}{\dim M'_{v_{0}}+\dim M'_{v_{1}}} \leq \frac{\dim
M_{v_{0}}}{\dim M_{v_{0}}+\dim M_{v_{1}}}\Leftrightarrow
\frac{\Theta(M')}{\sigma(M')}\leq\frac{\Theta(M)}{\sigma(M)}\; .$$
The equivalence between 3 and 4 follows from Proposition
\ref{GITstab-stabquivers}.
\end{pr}

\subsection{Kempf filtration for $Q$-sheaves}

Here we use Theorem \ref{equivstability} to show the correspondence between the Kempf Theorem and the
Harder-Narasimhan filtration for coherent sheaves (c.f. Theorem \ref{HNqsheaves}), passing through stability for
Kronecker modules and stability for representations of quivers on vector spaces. In this way, all ideas involved
in this thesis, all different notions of stability, all correspondences of maximal unstability, appear together.
Theorem \ref{HNqsheaves} gives another proof of Theorem \ref{kempfisHN}.

Let $Q=\{v\}$ be a one vertex quiver without arrows. Then, a $Q$-sheaf $E$ is a coherent sheaf $E$ and note that the stability
condition in Definition \ref{kstability} is independent of the choice of a polynomial $\kappa$. Let $E$ be an
unstable $Q$-sheaf with Hilbert polynomial $P$, i.e. an unstable coherent sheaf. Choose $l\gg m\gg 0$ such that
Theorem \ref{equivstability} holds, and such that $E$ is $m$-regular. By Theorem \ref{equivstability}, the
representation $M=H^{0}(E(m))\oplus H^{0}(E(l))$ of $\tilde{Q}=\{v_{0},v_{1}\}$ in $k$-vector spaces is
$(\Theta,\sigma)$-unstable, and the corresponding point $x_{M}\in \mathcal{R}^{(\Theta,\sigma)}_{d}(\tilde{Q})$
is $\chi_{(\Theta,\sigma)}$-unstable. By Theorem \ref{kempftheoremquivers}, let $0\subset M_{1}\subset \cdots
\subset M_{t+1}= M$ and $\Gamma_{1}<\Gamma_{2}<\ldots<\Gamma_{t}<\Gamma_{t+1}$ be the Kempf filtration of $M$
(depending on $m$ and $l$) which, by Theorem \ref{kempfHNquivers}, is the Harder-Narasimhan filtration of $M$,
defined in Theorem \ref{HNquivers}. Recall that we denote $M^{i}=M_{i}/M_{i-1}$ for each $i$.

The Kempf function in this case is
$$K(M_{\bullet},\Gamma)=\frac{\sum_{i=1}^{t+1}\Gamma_{i}
[\Theta(M)\sigma(M^{i})-\sigma(M)\Theta(M^{i})]}{\sqrt{\sum_{i=1}^{t+1}\Gamma_{i}^{2}\sigma(M^{i})}}
=\frac{(\Gamma,v)}{\|\Gamma\|}\; ,$$ (c.f. Theorem \ref{kempffunctionquivers}) where the coordinates of the vector
$v=(v_{1},...,v_{t+1})$ are given by $v_{i}=\Theta(M)-\frac{\sigma(M)}{\sigma(M^{i})}\Theta(M^{i})$, and the
scalar product in $\mathbb{R}^{t+1}$ is given by the diagonal matrix with elements $\sigma(M^{i})$ (c.f.
Definition \ref{graphquivers}).

\begin{dfn}\cite[Definition 5.3]{ACK}
\label{tight} Let $M'=V'_{0}\oplus V'_{1}$ and $M''=V''_{0}\oplus
V''_{1}$ be submodules of an $A$-module $M$. We say that $M'$ is
\emph{subordinate} to $M''$ if
$$V'_{0}\subset V''_{0}\text{         and         }V''_{1}\subset V'_{1}\; .$$
We say that $M'$ is \emph{tight} if it is subordinate to no submodule other than itself.
\end{dfn}

\begin{prop}
\label{alltight} The submodules $M_{i}$ appearing on the Kempf
filtration of $M$ are tight submodules of $M$.
\end{prop}
\begin{pr}
Given $m$ and $l$, the Kempf filtration of $M$ is, by Theorem \ref{kempfHNquivers}, the Harder-Narasimhan
filtration of $M$. Note that, whenever a module $M'$ is subordinate to $M''$, the slopes verify
$\mu_{(\Theta,\sigma)}(M'')\geq \mu_{(\Theta,\sigma)}(M')$. In the construction of the Harder-Narasimhan
filtration in Theorem \ref{HNquivers} we look for the unique representation $M_{1}$ of maximal slope
and of maximal total dimension $\sigma(M_{1})$ among those of maximal slope and, then, proceed by recursion. Hence, by that construction, all
the submodules have to be tight.
\end{pr}

Fix $m,l$ and consider the Kempf filtration of $M$. Using Proposition \ref{alltight} and \cite[Lemma 5.5]{ACK},
all submodules appearing in the Kempf filtration of $M$ are of the form
$M_{i}=\Hom_{X}(T,E_{i})=H^{0}(E_{i}(m))\oplus H^{0}(E_{i}(l))$ for some subsheaves $E_{i}\subset E$, where
$E_{i}(m)$ is globally generated for each $i$. Then, define the filtration
\begin{equation}
\label{kfiltqsheaves} 0\subset E_{1}\subset \cdots \subset E_{t+1}= E\; ,
\end{equation} which depends on $m$ and $l$. We call
it the \emph{$m$-Kempf filtration of the $Q$-sheaf $E$}.

Now we give the analogous to Proposition \ref{boundedness} for
this case. Fix the positive constant
\begin{equation}
\label{C_constantqsheaves} C=\max\{r|\mu _{\max }(E)|+\frac{d}{r}+r|B|+|A|+1\;,\;1\}.
\end{equation}

\begin{prop}
\label{boundednessqsheaves} Given integers $m$, $l$, let $E_{\bullet}\subset E$ be the $m$-Kempf filtration of the 
$Q$-sheaf $E$ as in (\ref{kfiltqsheaves}). There exists integers $m_{2}$, $l_{2}$ such that for $m\geq
m_{2},l\geq l_{2}$, each subsheaf $E_{i}\subset E$ in the $m$-Kempf filtration has slope $\mu(E_{i})\geq
\dfrac{d}{r}-C$.
\end{prop}

\begin{pr}
The proof follows similarly to Proposition \ref{boundedness}. Choose an integer $m_{1}\geq m_{0}$ such that for
every $m\geq m_{1}$, if we have a filter $E^{m}_{i}\subseteq E$ verifying $\mu (E^{m}_{i})<\frac{d}{r}-C$ (hence
it satisfies the estimate in Lemma \ref{Simpson}), it is
$$h^{0}(E^{m}_{i}(m))\leq \frac{1}{g^{n-1}n!}\big ((r_{i}-1)(\mu_{max}(E)+gm+B)^{n}+(\frac{d}{r}-C+gm+B)^{n}\big)=G(m)\; ,$$
where
$$G(m)=\frac{1}{g^{n-1}n!}\big
[r_{i}g^{n}m^{n}+ng^{n-1}\big((r_{i}-1)\mu_{max}(E)+\frac{d}{r}-C+r_{i}B\big)m^{n-1}+\cdots \big ]\; .$$

Recall that, by Definition \ref{graphquivers}, to any such filtration we associate a graph where $w^{i}=-b^{i}\cdot
v_{i}=-\sigma(M^{i})\cdot (\Theta(M)-\frac{\sigma(M)}{\sigma(M^{i})}\Theta(M^{i}))$. Then, the heights of the
graph, for each $i$, are
$$w_{i}=w^{1}+\ldots +w^{i}=\Theta(M_{i})\sigma(M)-\Theta(M)\sigma(M_{i})=$$
$$\dim M_{v_{0},i}(\dim M_{v_{0}}+\dim M_{v_{1}})-\dim M_{v_{0}}(\dim M_{v_{1},i}+\dim
M_{v_{1},i})=$$
$$\dim M_{v_{0},i}\dim M_{v_{1}}-\dim M_{v_{0}}\dim M_{v_{1},i}\; .$$
Again, to reach a contradiction, it is enough to show that $w_{i}<0$ because, in that case, we get $w_{t+1}<0$.
But it is
$$w_{t+1}=\dim M_{v_{0},t+1}\dim M_{v_{1}}-\dim M_{v_{0}}\dim M_{v_{1},t+1}=0\; ,$$
because $M_{v_{0},t+1}=M_{v_{0}}$ and $M_{v_{1},t+1}=M_{v_{1}}$, then the contradiction.

Using Proposition \ref{tight}, \cite[Lemma 5.5]{ACK}, and the $m$-regularity of $E$, we get
$$w_{i}=h^{0}(E_{i}(m))P_{E}(l)-P_{E}(m)h^{0}(E_{i}(l))\; .$$
Given $m$ and $l$, and using \cite[Lemma 5.4 b)]{ACK}, the negativity of the numerical expression given by $w_{i}$ for each $l$ is equivalent
to the negativity of the polynomial expression
$$h^{0}(E_{i}(m))P_{E}-P_{E}(m)P_{E_{i}}\; .$$

Let us show that $w_{i}(m,l)<0$, for sufficiently larges $m$ and $l$. By the previous calculations
$$w_{i}(m,l)=h^{0}(E_{i}(m))P_{E}(l)-P_{E}(m)P_{E_{i}}(l)
\leq$$
$$G(m)P_{E}(l)-P_{E}(m)P_{E_{i}}(l))=:\Psi(m,l)\; .
$$
where $\Psi(m,l)$ can be seen as an $n^{th}$-order polynomial on
$l$, whose coefficients are polynomials in $m$,
$$\Psi(m,l)=\psi_{n}(m)l^{n}+\psi_{n-1}(m)l^{n-1}+\cdots
+\psi_{1}(m)l+\psi_{0}(m)\: .$$ Hence, it is sufficient to show
that $\psi_{n}(m)=rG(m)-r_{i}P_{E}(m)<0$ for sufficiently large
$m$.

Note that $\psi_{n}(m)=\xi_{n}m^{n}+\xi_{n-1}m^{n-1}+\cdots
+\xi_{1}m+\xi_{0}$ is an $n^{th}$-order polynomial. The coefficient in
order $n^{th}$ vanishes,
$$\xi_{n}=(rG(m)-r_{i}P(m))_{n}=r\frac{r_{i}g}{n!}-r_{i}\frac{rg}{n!}=0\; .$$
Let us calculate the $(n-1)^{th}$-coefficient:
$$\xi_{n-1}=(rG(m)-r_{i}P(m))_{n-1}=(rG_{n-1}-r_{i}\frac{A}{(n-1)!})$$
where $G_{n-1}$ is the $(n-1)^{th}$-coefficient of the polynomial $G(m)$,
$$G_{n-1}=\frac{1}{g^{n-1}n!}ng^{n-1}((r_{i}-1)\mu_{max}(E)+\frac{d}{r}-C+r_{i}B)=$$
$$\frac{1}{(n-1)!}((r_{i}-1)\mu_{max}(E)+\frac{d}{r}-C+r_{i}B)\leq$$
$$\frac{1}{(n-1)!}((r_{i}-1)|\mu_{max}(E)|+\frac{d}{r}-C+r_{i}|B|)\leq$$
$$\frac{1}{(n-1)!}(r|\mu_{max}(E)|+\frac{d}{r}-C+r|B|)<\frac{-|A|}{(n-1)!}\; ,$$
last inequality coming from the definition of $C$ in \eqref{C_constantqsheaves}. Then
$$\xi_{n-1}<r\big(\frac{-|A|}{(n-1)!}\big)-r_{i}\frac{A}{(n-1)!}=
\frac{-r|A|-r_{i}A}{(n-1)!}<0\; ,$$ because $-r|A|-r_{i}A<0$.

Therefore $\psi_{n}(m)=\xi_{n-1}m^{n-1}+\cdots +\xi_{1}m+\xi_{0}$ with $\xi_{n-1}<0$, so there exists $m_{2}\geq
m_{1}$ such that for every $m\geq m_{2}$ we have $\psi_{n}(m)<0$ and $w_{i}(m,l)<0$, for $l\gg 0$, hence the
contradiction.
\end{pr}

\begin{prop}
\label{regularqsheaf} There exists an integer $m_{3}$ such that
for $m\geq m_{3}$ the sheaves $E^{m}_{i}$ and
$E^{m,i}=E^m_i/E^m_{i-1}$ are $m_{3}$-regular. In particular their
higher cohomology groups, after twisting with
$\mathcal{O}_{X}(m_{3})$, vanish and they are generated by global
sections.
\end{prop}
\begin{pr}
The argument follows analogously to Proposition \ref{regular}.
\end{pr}

By Proposition \ref{regularqsheaf}, for any $m\geq m_{3}$, all the filters $E_{i}$ of the $m$-Kempf filtration
of $E$ are $m_{3}$-regular and hence, the $m$-Kempf filtration of sheaves
$$0\subset E_{1} \subset E_{2} \subset \cdots \subset E_{t_{m}} \subset E_{t+1}=E\; ,$$
is obtained from the filtration of vector subspaces
$$0\subset H^{0}(E_{1}(m_{3})) \subset H^{0}(E_{2}(m_{3})) \subset \cdots \subset H^{0}(E_{t}(m_{3})) \subset H^{0}(E_{t+1}(m_{3}))=H^{0}(E(m_{3}))$$
by the evaluation map, of a unique vector space $H^{0}(E(m_{3}))$, whose dimension is independent of $m$.

Let $m\geq m_{3}$ and let
$$(P_{E_{1}},\ldots,P_{E_{t+1}})$$ be 
the \emph{$m$-type} of the $m$-Kempf filtration of $E$ (c.f. Definition \ref{mtype}) and let
$$\mathcal{P}=\big\{(P_{E_{1}},\ldots,P_{E_{t+1}})\big\}$$ be the finite set of possible vectors for $m\geq m_{3}$ (c.f. Proposition \ref{Pisfinite}).

By Definition \ref{graphquivers} we associate a graph to the
$m$-Kempf filtration, which can be rewritten, by Proposition
\ref{regularqsheaf}, as
$$v_{m,i}(l)=\Theta(M)-\frac{\sigma(M)}{\sigma(M^{i})}\cdot \Theta(M^{i})=\dim M_{v_{0}}-\frac{\dim M_{v_{0}}+\dim M_{v_{1}}}{\dim M_{v_{0}}^{i}+\dim M_{v_{1}}^{i}}\cdot \dim
M_{v_{0}}^{i}=$$
$$P_{E}(m)-\frac{P_{E}(m)+P_{E}(l)}{P_{E^{i}}(m)+P_{E^{i}}(l)}\cdot P_{E^{i}}(m)\; ,$$
and
$$b_{m}^{i}(l)=\dim M_{v_{0}}^{i}+\dim
M_{v_{1}}^{i}=P_{E^{i}}(m)+P_{E^{i}}(l)\; .$$ We use notations $v_{m,i}(l)$ and $b_{m}^{i}(l)$ because, given
the $m$-Kempf filtration, its $m$-type is fixed, for $m\geq m_{3}$, hence the coordinates of the graph can be
seen as rational functions in $l$, whose coefficients are fixed functions in $m$.

Now we follow the argument in subsection \ref{sectionkempfstabilizes} with the particularities of section
\ref{kempfhiggs}. Define the functional in $\mathcal{P}$,
$$\Phi_{m}(l)=(\mu_{v_{m}}(\Gamma_{v_{m}(l)}))^{2}=\|v_{m}(l)\|^{2}\; ,$$
which is a rational function on $l$ (c.f. (\ref{maxvalue})). By finiteness of $\mathcal{P}$ there is a finite
list of such possible functions
$$
\mathcal{A}=\{\Phi_{m}:m\geq m_{3}\}
$$
and we can choose $K$ among them, such that there exist integers $m_{4}$ and $l_{4}$ with $\Phi_{m}(l)=K(l)$,
for all $m\geq m_{4}$ and $l\geq l_{4}$ (c.f. Lemma \ref{uniquefunction}).

\begin{prop}
\label{eventuallyqsheaf} Let $a_{1},a_{2}$ be integers with $a_{1}\geq a_{2} \geq m_{4}$. The $a_{1}$-Kempf
filtration of $E$ is equal to the $a_{2}$-Kempf filtration of $E$.
\end{prop}
\begin{pr}
C.f. proof of Proposition \ref{eventually}.
\end{pr}

\begin{dfn}
If $m\geq m_{4}$, the $m$-Kempf filtration of $E$ is called \emph{the Kempf filtration of $E$},
$$0\subset E_{1} \subset E_{2} \subset \cdots \subset E_{t} \subset E_{t+1}=E\; .$$
Note that it does not depend on $m$ by Proposition
\ref{eventuallyqsheaf}.
\end{dfn}

\begin{thm}
\label{HNqsheaves} Given a one vertex quiver $Q$, every $Q$-sheaf $E$ over $X$ pure of dimension $e$, (i.e. a
coherent sheaf pure of dimension $e$) has a unique filtration
$$0\subset E_{1}\subset E_{2}\subset \ldots \subset E_{t}\subset E_{t+1}=E$$
verifying the following properties, where $E^{i}:=E_{i}/E_{i-1}$,
\begin{itemize}
\item $\frac{P_{E^{1}}(m)}{\rk E^{1}}>\frac{P_{E^{2}}(m)}{\rk E^{2}}>\ldots
>\frac{P_{E^{t}}(m)}{\rk E^{t}}>\frac{P_{E^{t+1}}(m)}{\rk E^{t+1}}$
\item The quotients $E^{i}$ are semistable
\end{itemize}
This filtration is the \emph{Harder-Narasimhan filtration} of $E$ defined in Theorem \ref{HNunique}.
\end{thm}
\begin{pr}
Let $\tilde{Q}=\{v_{0},v_{1}\}$. By Theorem \ref{equivstability},
choosing $m\geq m_{4}$, we associate to an unstable $Q$-sheaf $E$
a point in the parameter space, $x_{M}\in
\mathcal{R}^{(\Theta,\sigma)}_{d}(\tilde{Q})$ which is
$\chi_{(\Theta,\sigma)}$-unstable. By uniqueness of Theorem
\ref{kempftheoremquivers}, there exists a unique filtration of $M$
$$0\subset M_{1}\subset \cdots \subset M_{t+1}= M$$ verifying the two conditions of Theorem \ref{HNquivers}, this is
$\mu(M^{1})>\mu(M^{2})>\ldots
>\mu(M^{t})>\mu(M^{t+1})$ and that the quotients $M^{i}$ are
$(\Theta,\sigma)$-semistable. Consider the Kempf filtration
$$0\subset E_{1}\subset E_{2}\subset \ldots \subset E_{t}\subset
E_{t+1}=E\; ,$$ which does not depend on $m$, by Proposition
\ref{eventuallyqsheaf}.

To the Kempf filtration we associate a graph $v$, by Definition
\ref{graphquivers}, and by Lemma \ref{lemmaA}, the coordinates
$v_{i}$ are in increasing order, hence
$$v_{i}<v_{i+1}\Leftrightarrow\Theta(M)-\frac{\sigma
(M)}{\sigma(M^{i})}\cdot \Theta(M^{i})<\Theta(M)-\frac{\sigma (M)}{\sigma(M^{i+1})}\cdot \Theta(M^{i+1})$$
$$\Leftrightarrow
\dim M_{v_{0}}-\frac{\dim M_{v_{0}}+\dim M_{v_{1}}}{\dim M_{v_{0}}^{i}+\dim M_{v_{1}}^{i}}\cdot \dim
M_{v_{0}}^{i}<\dim M_{v_{0}}-\frac{\dim M_{v_{0}}+\dim M_{v_{1}}}{\dim M_{v_{0}}^{i}+\dim M_{v_{1}}^{i+1}}\cdot
\dim M_{v_{0}}^{i+1}$$
$$\Leftrightarrow \frac{\dim M^{i}_{v_{0}}}{\dim M^{i}_{v_{1}}}>
\frac{\dim M^{i+1}_{v_{0}}}{\dim M^{i+1}_{v_{1}}}\; .$$ Using
Theorem \ref{equivstability}, the last expression is equivalent to
$$\frac{P_{E^{i}}(m)}{P_{E^{i}}(l)}> \frac{P_{E^{i+1}}(m)}{P_{E^{i+1}}(l)}
\Leftrightarrow \frac{P_{E^{i}}(m)}{\rk E^{i}}>
\frac{P_{E^{i+1}}(m)}{\rk E^{i+1}}$$ where the last equivalence
follows from Lemma \ref{kstability2}. Using Lemma \ref{lemmaB} as
in Proposition \ref{blocksemistability}, we can see that the Kempf
filtration of $E$ verifies the second property of the
Harder-Narasimhan filtration as well.
\end{pr}

\end{document}